\documentclass{article}
\usepackage{geometry}
\usepackage{pdflscape}

\usepackage[numbers]{natbib}

\usepackage[]{amssymb}
\usepackage[cmex10]{amsmath}
\usepackage{multirow}
\usepackage{booktabs}
\usepackage{algorithm,algorithmic}
\usepackage{tikz}
\usepackage{ieeetrantools}

\DeclareMathOperator{\rank}{rank}
\DeclareMathOperator{\diag}{diag}

\begin{document}

\title{Real-Value Power-Voltage Formulations of, and Bounds for, Three-Wire Unbalanced Optimal Power Flow}                      




\author{Frederik Geth\footnote{Dr. Geth is with GridQube, Springfield (Brisbane), Queensland, Australia. \texttt{frederik.geth [at] gridqube.com}} \,and\, Hakan Ergun\footnote{Dr. Ergun is with KU Leuven, Kasteelpark Arenberg 10, Heverlee, Belgium and EnergyVille, Thor Park 8310, Belgium. \texttt{hakan.ergun [at] kuleuven.be} }}




\maketitle

\begin{abstract}
Unbalanced \emph{optimal} power flow refers to a class of optimization problems subject to the steady state physics of three-phase power grids with nonnegligible phase unbalance. 
Significant progress on this problem has been made on the mathematical modelling side of unbalanced OPF, however there is a lack of information on implementation aspects as well as data sets for benchmarking. One of the key problems is the lack of definitions of current and voltage bounds across different classes of representations of the power flow equations. 
Therefore, this tutorial-style paper summarizes the structural features of the unbalanced (optimal) power problem for three-phase systems. 
The resulting nonlinear complex-value matrix formulations are presented for both the bus injection and branch flow formulation frameworks, which typically cannot be implemented as-is in optimization toolboxes. 
Therefore, this paper  also  derives the equivalent real-value  formulations, and discusses challenges related to the implementation in optimization modeling toolboxes. 
The derived formulations can be re-used easily for  continuous and discrete optimization problems in distribution networks for a variety of operational and planning problems. 
Finally, bounds are derived for all variables involved, to further the development of benchmarks for unbalanced optimal power flow, where consensus on bound semantics is a pressing need. We believe benchmarks remain a cornerstone for the development and validation of scalable and reproducible optimization models and tools. 
The soundness of the derivations is confirmed through numerical experiments, validated w.r.t. \textsc{OpenDSS} for IEEE test feeders with $3\times3$ impedance matrices. 
\end{abstract}





%
%
%
%



\newcommand{\hermitiantranspose}{\text{H}}
\newcommand{\transpose}{\text{T}}
\newcommand{\complexconj}{*}

\newcommand{\Arrow}[1]{%
\parbox{#1}{\tikz{\draw[->](0,0)--(#1,0);}}
}
\newcommand{\Arrowright}[1]{%
\parbox{#1}{\tikz{\draw[<-](0,0)--(#1,0);}}
}

\newcommand{\indexScenario}{s}
\newcommand{\indexZone}{z}
\newcommand{\indexStakeholder}{b}
\newcommand{\indexUnitTech}{{x}_{\unitss}}
\newcommand{\indexPNTech}{{x}_{\pnss}}
\newcommand{\indexGITech}{{x}_{\giss}}
\newcommand{\indexGETech}{{x}_{\gess}}
\newcommand{\indexUnitRating}{{g}_{\unitss}}
\newcommand{\indexPNRating}{{g}_{\pnss}}
\newcommand{\indexGIRating}{{g}_{\giss}}
\newcommand{\indexGERating}{{g}_{\gess}}
\newcommand{\indexNodeRating}{{g}_{\nodess}}
\newcommand{\indexGridNode}{{i}}
\newcommand{\indexGridNodeTwo}{{j}}
\newcommand{\indexGridNodeThree}{{k}}
\newcommand{\indexGridLines}{{l}}
\newcommand{\indexGridLinesTwo}{{m}}
\newcommand{\indexGridLoad}{{h}}
\newcommand{\indexPhases}{p}
\newcommand{\indexPhasesTwo}{q}
\newcommand{\indexTimestep}{{k}}
\newcommand{\indexDay}{{d}}
\newcommand{\indexWeek}{{w}}
\newcommand{\indexUnit}{{u}}
\newcommand{\indexShunt}{{h}}
\newcommand{\indexPowerNode}{{e}}
\newcommand{\indexAdd}{{a}}
\newcommand{\indexRemove}{{r}}
\newcommand{\indexCirculation}{{\indexAdd\indexRemove}}
\newcommand{\indexFlows}{{f}}
\newcommand{\indexConductor}{n}
\newcommand{\indexdim}{n}
\newcommand{\indexdimtwo}{m}
\newcommand{\indexIteration}{k}
\newcommand{\indexContingency}{c}


\newcommand{\symPower}{P}
\newcommand{\symPowerToEnergy}{\gamma}
\newcommand{\symReactivePower}{Q}
\newcommand{\symApparentPower}{S}
\newcommand{\symVoltage}{U}
\newcommand{\symVoltageSOCP}{W}
\newcommand{\symVoltageAngle}{\theta}
\newcommand{\symPhaseDifference}{\varphi}
\newcommand{\symAnnualCost}{K}
\newcommand{\symPenalty}{Y}
\newcommand{\symPenaltyWeight}{N}
\newcommand{\symMass}{m}
\newcommand{\symSpecificMass}{\dot{\symMass}}
\newcommand{\symVolume}{V}
\newcommand{\symSpecificVolume}{v}
\newcommand{\symEnergy}{E}
\newcommand{\symEnergyFlow}{\dot{E}}
\newcommand{\symCurrent}{I}
\newcommand{\symCurrentSOCP}{L}
\newcommand{\symCurrentSquared}{L}
\newcommand{\symLineCurrent}{J}
\newcommand{\symObjectiveWeight}{w}
\newcommand{\symProbability}{\lambda}
\newcommand{\symBinary}{\alpha}
\newcommand{\symDecision}{{d}}
\newcommand{\symPhaseA}{{a}}
\newcommand{\symPhaseB}{{b}}
\newcommand{\symPhaseC}{{c}}
\newcommand{\symPhaseP}{{p}}
\newcommand{\symPhaseN}{{n}}
\newcommand{\symPhaseG}{{g}}
\newcommand{\symHarmonic}{{h}}
\newcommand{\symLossFactor}{{\rho}}
\newcommand{\symAvailability}{{b}}
\newcommand{\symDirectionality}{{d}}
\newcommand{\symPrice}{p}
\newcommand{\symInvestment}{C}
\newcommand{\symSpecificInvestment}{c}
\newcommand{\symLifetime}{\tau}
\newcommand{\symResistance}{r}
\newcommand{\symReactance}{x}
\newcommand{\symImpedance}{z}
\newcommand{\symAdmittance}{y}
\newcommand{\symSusceptance}{b}
\newcommand{\symConductance}{g}
\newcommand{\symLength}{l}
\newcommand{\symTime}{t}
\newcommand{\symRatio}{T}
\newcommand{\symLocation}{l}
\newcommand{\symEnthalpy}{H}
\newcommand{\symConnectivity}{C}

\newcommand{\symVolumeFlow}{\bar{\symVolume}}
\newcommand{\symMassFlow}{\bar{\symMass}}
\newcommand{\symDensity}{{\rho}}
\newcommand{\symGravitationalConstant}{\textcolor{\paramcolor}{g}}
\newcommand{\symElevation}{{h}}
\newcommand{\symElevationDifference}{\Delta{\symElevation}}
\newcommand{\symSetting}{\sigma}

\newcommand{\symRampRate}{\dot{\symPower}}
\newcommand{\symSpecRampRate}{{\zeta}}

\newcommand{\symPressure}[0]{p}
\newcommand{\symPressureSquared}[0]{\beta}
\newcommand{\symTemperature}[0]{   T }
\newcommand{\symHeatFlow}[0]{   Q }
\newcommand{\symViscosity}[0]{   \mu }
\newcommand{\symHeatResistance}[0]{   R }
\newcommand{\symFluidPipeDiameter}[0]{  \textcolor{\paramcolor}{D }    }
\newcommand{\symFluidPipeArea}[0]{  A     }
\newcommand{\symFluidPipeFriction}[0]{  \textcolor{\paramcolor}{f }    }
\newcommand{\symHeatCapacity}[0]{ c   }

\newcommand{\symRadiality}[0]{ \beta   }
\newcommand{\symSlack}{\epsilon}



\newcommand{\basess}{\text{{base}}}
\newcommand{\seriesss}{\text{{s}}}
\newcommand{\shuntss}{\text{{sh}}}
\newcommand{\invss}{\text{{inv}}}
\newcommand{\unbss}{\text{{unb}}}
\newcommand{\totss}{\text{{tot}}}
\newcommand{\batss}{\text{{bat}}}
\newcommand{\gridss}{\text{{grid}}}
\newcommand{\usabless}{\text{{usable}}}
\newcommand{\giss}{\text{{gi}}}
\newcommand{\pnss}{\text{{pn}}}
\newcommand{\lossss}{\text{{loss}}}
\newcommand{\cyclesss}{\text{{cycles}}}
\newcommand{\standbyss}{\text{{standby}}}
\newcommand{\operationss}{\text{{op}}}
\newcommand{\Pss}{\text{\sc{\symPower}}}
\newcommand{\Vss}{\text{\sc{\symVoltage}}}
\newcommand{\Iss}{\text{\sc{\symCurrent}}}
\newcommand{\Sss}{\text{\sc{\symApparentPower}}}
\newcommand{\Qss}{\text{\sc{\symReactivePower}}}
\newcommand{\Zss}{\text{\sc{\symImpedance}}}
\newcommand{\Ess}{\text{\sc{\symEnergy}}}
\newcommand{\nomss}{\text{nom}}
\newcommand{\ratedss}{\text{rated}}
\newcommand{\absss}{\text{abs}}
\newcommand{\maxss}{\text{max}}
\newcommand{\minss}{\text{min}}
\newcommand{\effss}{\text{eff}}
\newcommand{\chargess}{\text{c}}
\newcommand{\dischargess}{\text{d}}
\newcommand{\exoaddss}{\text{+}}
\newcommand{\exoremss}{\text{-}}
\newcommand{\orthss}{\perp}
\newcommand{\deprss}{\text{depr}}
\newcommand{\PFconvexss}{\text{PFconvex}}
\newcommand{\GIconvexss}{\text{GIconvex}}
\newcommand{\iterationss}{{(q)}}
\newcommand{\refss}{\text{ref}}
\newcommand{\addss}{\indexAdd,\text{add}}
\newcommand{\addtextss}{\text{add}}
\newcommand{\remss}{\indexRemove,\text{rem}}
\newcommand{\remtextss}{\text{rem}}
\newcommand{\circss}{\indexAdd\indexRemove,\text{circ}}
\newcommand{\flowss}{\indexFlows,\text{flow}}
\newcommand{\switchss}{\text{switch}}
\newcommand{\rotatedss}{\text{rot}}
\newcommand{\circumss}{\text{circum}}
\newcommand{\sbss}{\text{SB}}
\newcommand{\nonsbss}{\text{nonSB}}
\newcommand{\gess}{\text{ge}}
\newcommand{\unitss}{\text{unit}}
\newcommand{\nodess}{\text{node}}
\newcommand{\linss}{\text{lin}}
\newcommand{\realss}{\text{re}}
\newcommand{\imagss}{\text{im}}
\newcommand{\hagenposeuilless}[0]{  \text{ HP }}
\newcommand{\darcyweisbachss}[0]{ \text{  DC }}
\newcommand{\bentallinearization}[0]{ \text{  BenTal }}

\newcommand{\fortescuess}{\text{012}}

\newcommand{\plusss}{\text{+}}
\newcommand{\subsss}{\text{-}}
\newcommand{\combustionss}{\text{comb}}
\newcommand{\formationss}{\text{form}}



\newcommand{\PowerBase}{\symApparentPower_{\basess}}
\newcommand{\CurrentBase}{\symCurrent_{\basess}}
\newcommand{\VoltageBase}{\symVoltage_{\basess}^{\text{line}}}
\newcommand{\VoltageBasePhase}{\symVoltage_{\basess}^{\text{phase}}}
\newcommand{\TimeBase}{\symTime_{\basess}}
\newcommand{\EnergyBase}{\symEnergy_{\basess}}
\newcommand{\ImpedanceBase}{\symImpedance_{\basess}}
\newcommand{\AdmittanceBase}{\symAdmittance_{\basess}}
\newcommand{\CurrencyBase}{\euro_{\basess}}

\newcommand{\Knul}{\textcolor{\paramcolor}{\symAnnualCost_{0}}}
\newcommand{\Ynul}{\textcolor{\paramcolor}{\symPenalty_{0}}}



\newcommand{\yearunit}{a}
\newcommand{\VA}{VA}
\newcommand{\var}{var}
\newcommand{\kWh}{kWh}
\newcommand{\MWh}{MWh}
\newcommand{\si}[1]{#1}
\newcommand{\Si}[2]{\unit{#1}{#2}}
\newcommand{\percent}{\%}
\newcommand{\num}{}
\newcommand{\invtan}{\mathrm{atan2}}
\newcommand{\onen}[1]{\textcolor{\paramcolor}{\mathbf{1}_{#1}}}
\newcommand{\onenon}{\textcolor{\paramcolor}{\mathbf{1}}}

\newcommand{\nk}{n_{\setTimesteps}}
\newcommand{\nunits}{n_\setUnits}
\newcommand{\nn}{n_{\omegaGridNodes}}
\newcommand{\nc}{n_{\omegaGridLines}}
\newcommand{\nslackbuses}{n_{\omegaSlackbuses}}
\newcommand{\nnonslackbuses}{n_{\omegaNonSlackbuses}}

\newcommand{\setTimesteps}{\mathcal{K}}
\newcommand{\setTimestepsprev}{\mathcal{K}_0}
\newcommand{\omegaPNTech}{\mathcal{X}_{\pnss}}
\newcommand{\omegaGITech}{\mathcal{X}_{\giss}}
\newcommand{\omegaGETech}{\mathcal{X}_{\gess}}
\newcommand{\omegaUnitTech}{\mathcal{X}_{\unitss}}
\newcommand{\setPNRating}{\mathcal{G}_{\pnss}}
\newcommand{\setGIRating}{\mathcal{G}_{\giss}}
\newcommand{\setGERating}{\mathcal{G}_{\gess}}
\newcommand{\setUnitRating}{\mathcal{G}_{\unitss}}  
\newcommand{\setNodeRating}{\mathcal{G}_{\nodess}}  %
\newcommand{\omegaTransformers}{\textcolor{\paramcolor}{\mathcal{F}}}
\newcommand{\omegaGridLines}{\textcolor{\paramcolor}{\mathcal{J}}}
\newcommand{\omegaGridNodes}{\textcolor{\paramcolor}{\mathcal{I}}}
\newcommand{\omegaGridNodesDomain}[1]{\omegaGridNodes_{#1}}
\newcommand{\omegaSlackbuses}{\omegaGridNodes_{\indexDomains, \sbss}}
\newcommand{\omegaNonSlackbuses}{\omegaGridNodes_{\indexDomains,\nonsbss}}
\newcommand{\setReferenceNodes}{\omegaGridNodes_{\text{ref}}}
\newcommand{\setReferenceNodesPrime}{\omegaGridNodes'_{\text{ref}}}
\newcommand{\setReferenceNodesRankOne}{\omegaGridNodes^{\text{R-1}}_{\text{ref}}}
\newcommand{\setReferenceNodesRankOnePrime}{\omegaGridNodes'^{\text{R-1}}_{\text{ref}}}
\newcommand{\setReferenceNodesBalanced}{\omegaGridNodes^{\text{bal}}_{\text{ref}}}
\newcommand{\setReferenceNodesFixed}{\omegaGridNodes^{\text{fix}}_{\text{ref}}}

\newcommand{\setHarmonics}{\mathcal{H}}
\newcommand{\setHarmonicsTwo}{\mathcal{H}_1}

\newcommand{\setPowerNodes}{\mathcal{E}}
\newcommand{\setAdd}{\mathcal{A}}
\newcommand{\setRemove}{\mathcal{R}}
\newcommand{\setCirculation}{\mathcal{C}}
\newcommand{\setFlows}{\mathcal{F}}
\newcommand{\setContingency}{\mathcal{C}}
\newcommand{\setTopology}{\textcolor{\paramcolor}{\mathcal{T}^{\Arrow{.1cm}}}}
\newcommand{\setTopologyBoth}{\textcolor{\paramcolor}{\mathcal{T}}}
\newcommand{\setTopologyReverse}{\textcolor{\paramcolor}{\mathcal{T}^{\Arrowright{.1cm}}}}
\newcommand{\setUnits}{\textcolor{\paramcolor}{\mathcal{U}}}
\newcommand{\setUnitTopology}{\textcolor{\paramcolor}{\mathcal{T^{\text{units}}}}}
\newcommand{\setShuntTopology}{\textcolor{\paramcolor}{\mathcal{T^{\text{shunts}}}}}

\newcommand{\setBinary}[0]{\textcolor{\setscolor}{ \left\{ 0,1 \right\} }}
\newcommand{\setComplex}[0]{\textcolor{\setscolor}{ \mathbb{C}   }}
\newcommand{\setReal}[0]{\textcolor{\setscolor}{ \mathbb{R}   }}
\newcommand{\setIntegers}[0]{\textcolor{\setscolor}{ \mathbb{Z}   }}
\newcommand{\setNaturalNumbers}[0]{\textcolor{\setscolor}{ \mathbb{N}   }}
\newcommand{\setNaturalNumbersPos}[0]{\textcolor{\setscolor}{ \mathbb{N}^{+}   }}
\newcommand{\setRealMatrix}[2]{\textcolor{\setscolor}{ \mathbb{R}^{#1 \times #2}   }}
\newcommand{\setComplexMatrix}[2]{\textcolor{\setscolor}{ \mathbb{C}^{#1\times #2}   }}
\newcommand{\setHermitianMatrix}[1]{\textcolor{\setscolor}{ \mathbb{H}^{#1}   }}
\newcommand{\setRealPos}[0]{\textcolor{\setscolor}{  \setReal_{+}    }}

\newcommand{\setPhases}{\textcolor{\paramcolor}{\mathcal{P}}}
\newcommand{\setPhasesTwo}{{\Phi}}
\newcommand{\setFortescue}{\textcolor{\paramcolor}{\mathcal{F}}}
\newcommand{\setCircuit}{\textcolor{\paramcolor}{\mathcal{C}}}
\newcommand{\setConductors}{\textcolor{\paramcolor}{\mathcal{N}}}
\newcommand{\setBuspairs}{\textcolor{\paramcolor}{\mathcal{B}}}
\newcommand{\setZIP}{\mathcal{Z}}

\newcommand{\setPSDConen}[1]{\textcolor{\setscolor}{  \mathbb{S}_{+}^{#1}    }}
\newcommand{\setNNOrthant}[0]{\textcolor{\setscolor}{  \setReal_{+}^{\indexdim}    }}
\newcommand{\setSOCone}[0]{\textcolor{\setscolor}{  \mathbb{Q}_{}^{\indexdim}    }}
\newcommand{\setPSDCone}[0]{\textcolor{\setscolor}{  \mathbb{S}_{+}^{\indexdim}    }}
\newcommand{\setPDCone}[0]{\textcolor{\setscolor}{  \mathbb{S}_{++}^{\indexdim}    }}
\newcommand{\setScenario}{\mathcal{S}}
\newcommand{\setZone}{\mathcal{Z}}
\newcommand{\setStakeholder}{\mathcal{B}}

\newcommand{\ineuro}[1]{\EUR{\ensuremath{#1}}}
\newcommand{\eurounit}{\euro}
\newcommand{\ineuroperyear}[1]{\ensuremath{#1}\,\eurounit \si{\per\yearunit}}

\newcommand{\variables}{\textcolor{black}{real-valued variables}}
\newcommand{\parameters}{\textcolor{\paramcolor}{real-valued parameters}}
\newcommand{\binaryparameters}{\textcolor{\decisioncolor}{binary parameters}}
\newcommand{\bounds}{\textcolor{\boundscolor}{generic bounds}}
\newcommand{\sizingvariables}{\textcolor{\sizingcolor}{sizing parameters}}
\newcommand{\sizingbounds}{\textcolor{\ratingboundscolor}{sizing bounds}}
\newcommand{\sizingparameters}{\textcolor{\investmentcolor}{sizing parameters}}
\newcommand{\binaryvariables}{\textcolor{\binarycolor}{binary variables}}
\newcommand{\timeparameters}{\textcolor{\timecolor}{time parameters}}
\newcommand{\priceparameters}{\textcolor{\pricecolor}{price parameters}}
\newcommand{\complexvariables}{\textcolor{\complexcolor}{complex variables}}
\newcommand{\complexparameters}{\textcolor{\complexparamcolor}{complex parameters}}

\newcommand{\paramcolor}{red}
\newcommand{\sizingcolor}{red}
\newcommand{\investmentcolor}{teal}
\newcommand{\decisioncolor}{brown}
\newcommand{\boundscolor}{violet}
\newcommand{\binarycolor}{darkgray}
\newcommand{\timecolor}{cyan}
\newcommand{\pricecolor}{brown}
\newcommand{\ratingboundscolor}{olive}
\newcommand{\setscolor}{darkgray}
\newcommand{\complexcolor}{blue}  
\newcommand{\complexparamcolor}{brown}  


\newcommand{\costheta}{{c}_{\indexGridNode \indexGridNodeTwo}   }
\newcommand{\sintheta}{{s}_{\indexGridNode \indexGridNodeTwo}   }
\newcommand{\usquared}{{u}_{\indexGridNode \indexGridNode}   }
\newcommand{\umult}{{u}_{\indexGridNode \indexGridNodeTwo}   }
\newcommand{\uucos}{{r}_{\indexGridNode \indexGridNodeTwo}   }
\newcommand{\uusin}{{t}_{\indexGridNode \indexGridNodeTwo}   }

\newcommand{\phasesi}{\textcolor{\paramcolor}{{\setPhasesTwo}_{\indexGridNode }   }}
\newcommand{\phasesj}{\textcolor{\paramcolor}{{\setPhasesTwo}_{\indexGridNodeTwo }   }}
\newcommand{\phasesij}{\textcolor{\paramcolor}{{\setPhasesTwo}_{\indexGridNode \indexGridNodeTwo}   }}

\newcommand{\fortescue}{\textcolor{\complexparamcolor}{\mathbf{F}}}
\newcommand{\alpharot}{\textcolor{\complexparamcolor}{\alpha}}

\newcommand{\Ts}{\textcolor{\timecolor}{T_{\text{s}} }}
\newcommand{\thor}{\textcolor{\timecolor}{t_{\text{h}}}}
\newcommand{\etac}{\textcolor{\paramcolor}{\eta_{\chargess} }}
\newcommand{\etad}{\textcolor{\paramcolor}{\eta_{\dischargess} }}
\newcommand{\etap}{\textcolor{\paramcolor}{\eta_{\exoaddss} }}
\newcommand{\etam}{\textcolor{\paramcolor}{\eta_{\exoremss} }}

\newcommand{\etaadd}{\textcolor{\paramcolor}{\eta_{\addss} }}
\newcommand{\etarem}{\textcolor{\paramcolor}{\eta_{\remss} }}


\newcommand{\Pck}[0]{\symPower_{\chargess} }
\newcommand{\Pdk}[0]{\symPower_{\dischargess} }
\newcommand{\Ppk}[0]{\symPower_{\exoaddss} }
\newcommand{\Pmk}[0]{\symPower_{\exoremss} }

\newcommand{\Paddk}[0]{\symPower_{\addss} }
\newcommand{\Premk}[0]{\symPower_{\remss} }
\newcommand{\Pflowk}[0]{\symPower_{\flowss} }
\newcommand{\Pflowkprev}[0]{\symPower_{\flowss} prev}
\newcommand{\Pflowrampmax}[0]{\textcolor{\paramcolor}{\symRampRate_{\flowss}^{\maxss}}}
\newcommand{\Pflowrampspecmax}[0]{\textcolor{\paramcolor}{\symSpecRampRate_{\flowss}}}

\newcommand{\Pclossk}[0]{\symPower_{\chargess, \lossss} }
\newcommand{\Pdlossk}[0]{\symPower_{\dischargess, \lossss} }
\newcommand{\Pplossk}[0]{\symPower_{\exoaddss, \lossss} }
\newcommand{\Pmlossk}[0]{\symPower_{\exoremss, \lossss} }
\newcommand{\Pgilossk}[0]{\symPower_{\exoremss, \lossss} }
\newcommand{\PPNlossk}[0]{\symPower_{\indexPowerNode, \lossss} }
\newcommand{\PPNexolossk}[0]{\symPower_{\indexPowerNode, \text{exo},\lossss} }
\newcommand{\PPNeleclossk}[0]{\symPower_{\indexPowerNode, \text{elec}, \lossss} }
\newcommand{\Punitlossk}[0]{\symPower_{\indexUnit, \lossss} }

\newcommand{\Paddlossk}[0]{\symPower_{\addss, \lossss} }
\newcommand{\Premlossk}[0]{\symPower_{\remss, \lossss} }

\newcommand{\Kswitch}[0]{\symAnnualCost^{\switchss}_{}}
\newcommand{\Kswitchl}[0]{\symAnnualCost^{\switchss}_{\indexGridLines}}

\newcommand{\Kgridloss}[0]{\symAnnualCost^{\gridss\lossss}_{}}
\newcommand{\Kgridlossl}[0]{\symAnnualCost^{\gridss\lossss}_{\indexGridLines}}

\newcommand{\Kc}[0]{\symAnnualCost^{\chargess}_{}}
\newcommand{\Kd}[0]{\symAnnualCost^{\dischargess}_{}}
\newcommand{\Kp}[0]{\symAnnualCost^{\exoaddss}_{}}
\newcommand{\Km}[0]{\symAnnualCost^{\exoremss}_{}}
\newcommand{\KE}[0]{\symAnnualCost^{\Ess}_{}}
\newcommand{\Kgi}[0]{\symAnnualCost^{\giss}_{}}
\newcommand{\KPN}[0]{\symAnnualCost_{\indexPowerNode}}
\newcommand{\Kunit}[0]{\symAnnualCost_{\indexUnit}}
\newcommand{\Kunits}[0]{\symAnnualCost_{\setUnits}}

\newcommand{\Kflow}[0]{\symAnnualCost^{\flowss}_{}}

\newcommand{\Kopex}[0]{\symAnnualCost^{\operationss}_{}}

\newcommand{\Kcoper}[0]{\symAnnualCost^{\chargess,\operationss}}
\newcommand{\Kdoper}[0]{\symAnnualCost^{\dischargess,\operationss}}
\newcommand{\Kctotoper}[0]{\symAnnualCost^{\chargess,\totss,\operationss}}
\newcommand{\Kdtotoper}[0]{\symAnnualCost^{\dischargess,\totss,\operationss}}
\newcommand{\Kpoper}[0]{\symAnnualCost^{\exoaddss,\operationss}}
\newcommand{\Kmoper}[0]{\symAnnualCost^{\exoremss,\operationss}}
\newcommand{\KPNoper}[0]{\symAnnualCost_{\indexPowerNode}^{\operationss}}
\newcommand{\Kunitoper}[0]{\symAnnualCost_{\indexUnit}^{\operationss}}
\newcommand{\Kunitsoper}[0]{\symAnnualCost_{\setUnits}^{\operationss}}

\newcommand{\Kflowoper}[0]{\symAnnualCost^{\flowss,\operationss}}

\newcommand{\Ec}[0]{\symEnergy_{\chargess}}
\newcommand{\Ed}[0]{\symEnergy_{\dischargess}}
\newcommand{\Ep}[0]{\symEnergy_{\exoaddss}}
\newcommand{\Em}[0]{\symEnergy_{\exoremss}}
\newcommand{\Egi}[0]{\symEnergy_{\giss}}
\newcommand{\EE}[0]{n_{\cyclesss}}

\newcommand{\Eflow}[0]{\symEnergy_{\flowss}}

\newcommand{\Ecflow}[0]{\symEnergyFlow_{\chargess}}
\newcommand{\Edflow}[0]{\symEnergyFlow_{\dischargess}}
\newcommand{\Epflow}[0]{\symEnergyFlow_{\exoaddss}}
\newcommand{\Emflow}[0]{\symEnergyFlow_{\exoremss}}
\newcommand{\Egiflow}[0]{\symEnergyFlow_{\giss}}
\newcommand{\EEflow}[0]{\dot{n}_{\cyclesss}}

\newcommand{\Eflowflow}[0]{\symEnergyFlow_{\flowss}}

\newcommand{\Ecmax}[0]{\textcolor{\paramcolor}{\symEnergy_{\chargess}^{\maxss}}}
\newcommand{\Edmax}[0]{\textcolor{\paramcolor}{\symEnergy_{\dischargess}^{\maxss}}}
\newcommand{\Epmax}[0]{\textcolor{\paramcolor}{\symEnergy_{\exoaddss}^{\maxss}}}
\newcommand{\Emmax}[0]{\textcolor{\paramcolor}{\symEnergy_{\exoremss}^{\maxss}}}
\newcommand{\Egimax}[0]{\textcolor{\paramcolor}{\symEnergy_{\giss}^{\maxss}}}
\newcommand{\EEmax}[0]{\textcolor{\paramcolor}{n_{\cyclesss}^{\maxss}}}

\newcommand{\Eflowmax}[0]{\textcolor{\paramcolor}{\symEnergy_{\flowss}^{\maxss}}}

\newcommand{\Ecmaxflow}[0]{\textcolor{\paramcolor}{\dot{\symEnergy}_{\chargess}^{\maxss}}}
\newcommand{\Edmaxflow}[0]{\textcolor{\paramcolor}{\dot{\symEnergy}_{\dischargess}^{\maxss}}}
\newcommand{\Epmaxflow}[0]{\textcolor{\paramcolor}{\dot{\symEnergy}_{\exoaddss}^{\maxss}}}
\newcommand{\Emmaxflow}[0]{\textcolor{\paramcolor}{\dot{\symEnergy}_{\exoremss}^{\maxss}}}
\newcommand{\Egimaxflow}[0]{\textcolor{\paramcolor}{\dot{\symEnergy}_{\giss}^{\maxss}}}
\newcommand{\EEmaxflow}[0]{\textcolor{\paramcolor}{\dot{n}_{\cyclesss}^{\maxss}}}

\newcommand{\Eflowmaxflow}[0]{\textcolor{\paramcolor}{\dot{\symEnergy}_{\flowss}^{\maxss}}}

\newcommand{\Kcdepr}[0]{\symAnnualCost^{\chargess,\deprss}}
\newcommand{\Kddepr}[0]{\symAnnualCost^{\dischargess,\deprss}}
\newcommand{\Kpdepr}[0]{\symAnnualCost^{\exoaddss,\deprss}}
\newcommand{\Kmdepr}[0]{\symAnnualCost^{\exoremss,\deprss}}
\newcommand{\KEdepr}[0]{\symAnnualCost^{\Ess,\deprss}}
\newcommand{\Kgidepr}[0]{\symAnnualCost^{\giss,\deprss}}
\newcommand{\KPNdepr}[0]{\symAnnualCost_{\indexPowerNode}^{\deprss}}
\newcommand{\Kunitdepr}[0]{\symAnnualCost_{\indexUnit}^{\deprss}}
\newcommand{\Kunitsdepr}[0]{\symAnnualCost_{\setUnits}^{\deprss}}

\newcommand{\Kflowdepr}[0]{\symAnnualCost^{\flowss,\deprss}}

\newcommand{\Cc}[0]{\textcolor{\investmentcolor}{\symInvestment_{\chargess}}}
\newcommand{\Cd}[0]{\textcolor{\investmentcolor}{\symInvestment_{\dischargess}}}
\newcommand{\Cp}[0]{\textcolor{\investmentcolor}{\symInvestment_{\exoaddss}}}
\newcommand{\Cm}[0]{\textcolor{\investmentcolor}{\symInvestment_{\exoremss}}}
\newcommand{\Cgi}[0]{\textcolor{\investmentcolor}{\symInvestment_{\Sss}}}
\newcommand{\CE}[0]{\textcolor{\investmentcolor}{\symInvestment_{\Ess}}}
\newcommand{\Cunit}[0]{\textcolor{\investmentcolor}{\symInvestment_{\indexUnit}}}
\newcommand{\CPN}[0]{\textcolor{\investmentcolor}{\symInvestment_{\indexPowerNode}}}
\newcommand{\Cunitmax}[0]{\textcolor{\boundscolor}{\symInvestment^{\maxss}_{\indexUnit}}}

\newcommand{\Cflow}[0]{\textcolor{\investmentcolor}{\symInvestment_{\flowss}}}

\newcommand{\cpc}[0]{\textcolor{\pricecolor}{\symSpecificInvestment_{\chargess}}}
\newcommand{\cpd}[0]{\textcolor{\pricecolor}{\symSpecificInvestment_{\dischargess}}}
\newcommand{\cpp}[0]{\textcolor{\pricecolor}{\symSpecificInvestment_{\exoaddss}}}
\newcommand{\cpm}[0]{\textcolor{\pricecolor}{\symSpecificInvestment_{\exoremss}}}
\newcommand{\cpgi}[0]{\textcolor{\pricecolor}{\symSpecificInvestment_{\Sss}}}
\newcommand{\cpE}[0]{\textcolor{\pricecolor}{\symSpecificInvestment_{\Ess}}}

\newcommand{\cpflow}[0]{\textcolor{\pricecolor}{\symSpecificInvestment_{\flowss}}}

\newcommand{\Vc}[0]{\textcolor{\investmentcolor}{\symVolume_{\chargess}}}
\newcommand{\Vd}[0]{\textcolor{\investmentcolor}{\symVolume_{\dischargess}}}
\newcommand{\Vp}[0]{\textcolor{\investmentcolor}{\symVolume_{\exoaddss}}}
\newcommand{\Vm}[0]{\textcolor{\investmentcolor}{\symVolume_{\exoremss}}}
\newcommand{\Vgi}[0]{\textcolor{\investmentcolor}{\symVolume_{\Sss}}}
\newcommand{\VE}[0]{\textcolor{\investmentcolor}{\symVolume_{\Ess}}}
\newcommand{\Vunit}[0]{\textcolor{\investmentcolor}{\symVolume_{\indexUnit}}}
\newcommand{\Vunitmax}[0]{\textcolor{\boundscolor}{\symVolume^{\maxss}_{\indexUnit}}}
\newcommand{\VPN}[0]{\textcolor{\investmentcolor}{\symVolume_{\indexPowerNode}}}

\newcommand{\Vflow}[0]{\textcolor{\investmentcolor}{\symVolume_{\flowss}}}

\newcommand{\vpc}[0]{\textcolor{\pricecolor}{\symSpecificVolume_{\chargess}}}
\newcommand{\vpd}[0]{\textcolor{\pricecolor}{\symSpecificVolume_{\dischargess}}}
\newcommand{\vpp}[0]{\textcolor{\pricecolor}{\symSpecificVolume_{\exoaddss}}}
\newcommand{\vpm}[0]{\textcolor{\pricecolor}{\symSpecificVolume_{\exoremss}}}
\newcommand{\vpgi}[0]{\textcolor{\pricecolor}{\symSpecificVolume_{\Sss}}}
\newcommand{\vpE}[0]{\textcolor{\pricecolor}{\symSpecificVolume_{\Ess}}}

\newcommand{\vpflow}[0]{\textcolor{\pricecolor}{\symSpecificVolume_{\flowss}}}

\newcommand{\Mc}[0]{\textcolor{\investmentcolor}{\symMass_{\chargess}}}
\newcommand{\Md}[0]{\textcolor{\investmentcolor}{\symMass_{\dischargess}}}
\newcommand{\Mp}[0]{\textcolor{\investmentcolor}{\symMass_{\exoaddss}}}
\newcommand{\Mm}[0]{\textcolor{\investmentcolor}{\symMass_{\exoremss}}}
\newcommand{\Mgi}[0]{\textcolor{\investmentcolor}{\symMass_{\Sss}}}
\newcommand{\ME}[0]{\textcolor{\investmentcolor}{\symMass_{\Ess}}}
\newcommand{\Munit}[0]{\textcolor{\investmentcolor}{\symMass_{\indexUnit}}}
\newcommand{\Munitmax}[0]{\textcolor{\boundscolor}{\symMass^{\maxss}_{\indexUnit}}}
\newcommand{\MPN}[0]{\textcolor{\investmentcolor}{\symMass_{\indexPowerNode}}}

\newcommand{\Mflow}[0]{\textcolor{\investmentcolor}{\symMass_{\flowss}}}

\newcommand{\mpc}[0]{\textcolor{\pricecolor}{\symSpecificMass_{\chargess}}}
\newcommand{\mpd}[0]{\textcolor{\pricecolor}{\symSpecificMass_{\dischargess}}}
\newcommand{\mpp}[0]{\textcolor{\pricecolor}{\symSpecificMass_{\exoaddss}}}
\newcommand{\mpm}[0]{\textcolor{\pricecolor}{\symSpecificMass_{\exoremss}}}
\newcommand{\mpgi}[0]{\textcolor{\pricecolor}{\symSpecificMass_{\Sss}}}
\newcommand{\mpE}[0]{\textcolor{\pricecolor}{\symSpecificMass_{\Ess}}}

\newcommand{\mpflow}[0]{\textcolor{\pricecolor}{\symSpecificMass_{\flowss}}}

\newcommand{\tauc}[0]{\textcolor{black}{\symLifetime_{\chargess}}}
\newcommand{\taud}[0]{\textcolor{black}{\symLifetime_{\dischargess}}}
\newcommand{\taup}[0]{\textcolor{black}{\symLifetime_{\exoaddss}}}
\newcommand{\taum}[0]{\textcolor{black}{\symLifetime_{\exoremss}}}
\newcommand{\taugi}[0]{\textcolor{black}{\symLifetime_{\Sss}}}
\newcommand{\tauE}[0]{\textcolor{black}{\symLifetime{_\Ess}}}

\newcommand{\tauflow}[0]{\textcolor{black}{\symLifetime_{\flowss}}}

\newcommand{\taucmax}[0]{\textcolor{\paramcolor}{\symLifetime^{\maxss}_{\chargess}}}
\newcommand{\taudmax}[0]{\textcolor{\paramcolor}{\symLifetime{^{\maxss}_\dischargess}}}
\newcommand{\taupmax}[0]{\textcolor{\paramcolor}{\symLifetime^{\maxss}_{\exoaddss}}}
\newcommand{\taummax}[0]{\textcolor{\paramcolor}{\symLifetime^{\maxss}_{\exoremss}}}
\newcommand{\taugimax}[0]{\textcolor{\paramcolor}{\symLifetime^{\maxss}_{\Sss}}}
\newcommand{\tauEmax}[0]{\textcolor{\paramcolor}{\symLifetime^{\maxss}_{\Ess}}}

\newcommand{\tauflowmax}[0]{\textcolor{\paramcolor}{\symLifetime^{\maxss}_{\flowss}}}

\newcommand{\ppc}[0]{\textcolor{\pricecolor}{\symPrice_{\chargess} }}
\newcommand{\ppd}[0]{\textcolor{\pricecolor}{\symPrice_{\dischargess} }}
\newcommand{\ppp}[0]{\textcolor{\pricecolor}{\symPrice_{\exoaddss} }}
\newcommand{\ppm}[0]{\textcolor{\pricecolor}{\symPrice_{\exoremss} }}

\newcommand{\ppflow}[0]{\textcolor{\pricecolor}{\symPrice_{\flowss} }}

\newcommand{\ppswitch}[0]{\textcolor{\pricecolor}{\symPrice_{\switchss} }}
\newcommand{\ppswitchl}[0]{\textcolor{\pricecolor}{\symPrice_{\indexGridLines,\switchss} }}

\newcommand{\ppPgridlossl}[0]{\textcolor{\pricecolor}{\symPrice_{\indexGridLines\gridss\lossss}^{\Pss}    }}
\newcommand{\ppQgridlossl}[0]{\textcolor{\pricecolor}{\symPrice_{\indexGridLines\gridss\lossss}^{\Qss} }}
\newcommand{\ppQabsgridlossl}[0]{\textcolor{\pricecolor}{\symPrice_{\indexGridLines\gridss\lossss}^{|\Qss|} }}

\newcommand{\Pckref}[0]{\textcolor{\paramcolor}{\symPower^{\refss}_{\chargess} }}
\newcommand{\Pdkref}[0]{\textcolor{\paramcolor}{\symPower^{\refss}_{\dischargess} }}
\newcommand{\Ppkref}[0]{\textcolor{\paramcolor}{\symPower^{\refss}_{\exoaddss} }}
\newcommand{\Pmkref}[0]{\textcolor{\paramcolor}{\symPower^{\refss}_{\exoremss} }}

\newcommand{\Pckmax}[0]{\textcolor{\boundscolor}{\symPower^{\maxss}_{\chargess} }}
\newcommand{\Pdkmax}[0]{\textcolor{\boundscolor}{\symPower^{\maxss}_{\dischargess} }}
\newcommand{\Ppkmax}[0]{\textcolor{\boundscolor}{\symPower^{\maxss}_{\exoaddss} }}
\newcommand{\Pmkmax}[0]{\textcolor{\boundscolor}{\symPower^{\maxss}_{\exoremss} }}

\newcommand{\Pckmin}[0]{\textcolor{\boundscolor}{\symPower^{\minss}_{\chargess} }}
\newcommand{\Pdkmin}[0]{\textcolor{\boundscolor}{\symPower^{\minss}_{\dischargess} }}
\newcommand{\Ppkmin}[0]{\textcolor{\boundscolor}{\symPower^{\minss}_{\exoaddss} }}
\newcommand{\Pmkmin}[0]{\textcolor{\boundscolor}{\symPower^{\minss}_{\exoremss} }}
\newcommand{\Pflowkref}[0]{\textcolor{\paramcolor}{\symPower^{\refss}_{\flowss} }}

\newcommand{\Paddkmin}[0]{\textcolor{\boundscolor}{\symPower^{\minss}_{\addss} }}
\newcommand{\Paddkmax}[0]{\textcolor{\boundscolor}{\symPower^{\maxss}_{\addss} }}
\newcommand{\Premkmin}[0]{\textcolor{\boundscolor}{\symPower^{\minss}_{\remss} }}
\newcommand{\Premkmax}[0]{\textcolor{\boundscolor}{\symPower^{\maxss}_{\remss} }}
\newcommand{\Pflowkmin}[0]{\textcolor{\boundscolor}{\symPower^{\minss}_{\flowss} }}
\newcommand{\Pflowkmax}[0]{\textcolor{\boundscolor}{\symPower^{\maxss}_{\flowss} }}

\newcommand{\Pcratedmax}[0]{\textcolor{\ratingboundscolor}{\symPower'^{\maxss}_{\chargess}}}
\newcommand{\Pdratedmax}[0]{\textcolor{\ratingboundscolor}{\symPower'^{\maxss}_{\dischargess}}}
\newcommand{\Ppratedmax}[0]{\textcolor{\ratingboundscolor}{\symPower'^{\maxss}_{\exoaddss}}}
\newcommand{\Pmratedmax}[0]{\textcolor{\ratingboundscolor}{\symPower'^{\maxss}_{\exoremss}}}

\newcommand{\Pcratedmin}[0]{\textcolor{\ratingboundscolor}{\symPower'^{\minss}_{\chargess}}}
\newcommand{\Pdratedmin}[0]{\textcolor{\ratingboundscolor}{\symPower'^{\minss}_{\dischargess}}}
\newcommand{\Ppratedmin}[0]{\textcolor{\ratingboundscolor}{\symPower'^{\minss}_{\exoaddss}}}
\newcommand{\Pmratedmin}[0]{\textcolor{\ratingboundscolor}{\symPower'^{\minss}_{\exoremss}}}

\newcommand{\Paddratedmin}[0]{\textcolor{\ratingboundscolor}{\symPower'^{\minss}_{\addss}}}
\newcommand{\Paddratedmax}[0]{\textcolor{\ratingboundscolor}{\symPower'^{\maxss}_{\addss}}}
\newcommand{\Premratedmin}[0]{\textcolor{\ratingboundscolor}{\symPower'^{\minss}_{\remss}}}
\newcommand{\Premratedmax}[0]{\textcolor{\ratingboundscolor}{\symPower'^{\maxss}_{\remss}}}

\newcommand{\Pflowrated}[0]{\textcolor{\sizingcolor}{\symPower'^{}_{\flowss}}}
\newcommand{\Pflowratedmin}[0]{\textcolor{\ratingboundscolor}{\symPower'^{\minss}_{\flowss}}}
\newcommand{\Pflowratedmax}[0]{\textcolor{\ratingboundscolor}{\symPower'^{\maxss}_{\flowss}}}

\newcommand{\Eratedmaxmax}[0]{\textcolor{\ratingboundscolor}{\symEnergy'^{\maxss\maxss}}}
\newcommand{\Eratedmaxmin}[0]{\textcolor{\ratingboundscolor}{\symEnergy'^{\maxss\minss}}}
\newcommand{\Eratedmax}[0]{\textcolor{\sizingcolor}{\symEnergy'^{\maxss}}}
\newcommand{\Eratedminmax}[0]{\textcolor{\ratingboundscolor}{\symEnergy'^{\minss\maxss}}}
\newcommand{\Eratedminmin}[0]{\textcolor{\ratingboundscolor}{\symEnergy'^{\minss\minss}}}
\newcommand{\Eratedmin}[0]{\textcolor{\sizingcolor}{\symEnergy'^{\minss}}}

\newcommand{\Pcrated}[0]{\textcolor{\sizingcolor}{\symPower'^{}_{\chargess}}}
\newcommand{\Pdrated}[0]{\textcolor{\sizingcolor}{\symPower'^{}_{\dischargess}}}
\newcommand{\Pprated}[0]{\textcolor{\sizingcolor}{\symPower'^{}_{\exoaddss}}}
\newcommand{\Pmrated}[0]{\textcolor{\sizingcolor}{\symPower'^{}_{\exoremss}}}

\newcommand{\Paddrated}[0]{\textcolor{\sizingcolor}{\symPower'^{}_{\addss}}}
\newcommand{\Premrated}[0]{\textcolor{\sizingcolor}{\symPower'^{}_{\remss}}}

\newcommand{\Pcorthk}[0]{\symPower^{\orthss}_{\chargess} }
\newcommand{\Pdorthk}[0]{\symPower^{\orthss}_{\dischargess} }
\newcommand{\Pporthk}[0]{\symPower^{\orthss}_{\exoaddss} }
\newcommand{\Pmorthk}[0]{\symPower^{\orthss}_{\exoremss} }
\newcommand{\Porthk}[0]{\symPower^{\orthss}_{} }

\newcommand{\Paddorthk}[0]{\symPower^{\orthss}_{\addss} }
\newcommand{\Premorthk}[0]{\symPower^{\orthss}_{\remss} }
\newcommand{\Pfloworthk}[0]{\symPower^{\orthss}_{\flowss} }
\newcommand{\Paddorthtotk}[0]{\symPower^{\orthss}_{\addtextss} }
\newcommand{\Premorthtotk}[0]{\symPower^{\orthss}_{\remtextss} }

\newcommand{\Pcorthkmax}[0]{\textcolor{\boundscolor}{\symPower^{\orthss, \maxss}_{\chargess} }}
\newcommand{\Pdorthkmax}[0]{\textcolor{\boundscolor}{\symPower^{\orthss, \maxss}_{\dischargess} }}
\newcommand{\Pporthkmax}[0]{\textcolor{\boundscolor}{\symPower^{\orthss, \maxss}_{\exoaddss} }}
\newcommand{\Pmorthkmax}[0]{\textcolor{\boundscolor}{\symPower^{\orthss, \maxss}_{\exoremss} }}

\newcommand{\Paddorthkmax}[0]{\textcolor{\boundscolor}{\symPower^{\orthss, \maxss}_{\addss} }}
\newcommand{\Premorthkmax}[0]{\textcolor{\boundscolor}{\symPower^{\orthss, \maxss}_{\remss} }}
\newcommand{\Pfloworthkmax}[0]{\textcolor{\boundscolor}{\symPower^{\orthss, \maxss}_{\flowss} }}

\newcommand{\Pcdk}[0]{\symPower_{\chargess\dischargess} }
\newcommand{\Pcmk}[0]{\symPower_{\chargess\exoremss} }
\newcommand{\Ppmk}[0]{\symPower_{\exoaddss\exoremss} }
\newcommand{\Ppdk}[0]{\symPower_{\exoaddss\dischargess} }

\newcommand{\Pcirck}[0]{\symPower_{\circss} }
\newcommand{\Pcircmax}[0]{\symPower^{\maxss}_{\circss}}

\newcommand{\Pcdmax}[0]{\symPower^{\maxss}_{\chargess\dischargess}}
\newcommand{\Pcmmax}[0]{\symPower^{\maxss}_{\chargess\exoremss}}
\newcommand{\Ppmmax}[0]{\symPower^{\maxss}_{\exoaddss\exoremss}}
\newcommand{\Ppdmax}[0]{\symPower^{\maxss}_{\exoaddss\dischargess}}

\newcommand{\Pcpk}[0]{\symPower_{\chargess\exoaddss} }
\newcommand{\Pdmk}[0]{\symPower_{\dischargess\exoremss} }
\newcommand{\gammacp}[0]{\textcolor{\paramcolor}{\symPowerToEnergy_{\chargess\exoaddss} }}
\newcommand{\gammadm}[0]{\textcolor{\paramcolor}{\symPowerToEnergy_{\dischargess\exoremss} }}

\newcommand{\gammaadd}[0]{\textcolor{\paramcolor}{\symPowerToEnergy_{\addss} }}
\newcommand{\gammarem}[0]{\textcolor{\paramcolor}{\symPowerToEnergy_{\remss} }}

\newcommand{\dcd}[0]{\textcolor{\decisioncolor}{\symDecision_{\chargess\dischargess} }}
\newcommand{\dcm}[0]{\textcolor{\decisioncolor}{\symDecision_{\chargess\exoremss} }}
\newcommand{\dpm}[0]{\textcolor{\decisioncolor}{\symDecision_{\exoaddss\exoremss} }}
\newcommand{\dpd}[0]{\textcolor{\decisioncolor}{\symDecision_{\exoaddss\dischargess} }}

\newcommand{\dcirc}[0]{\textcolor{\decisioncolor}{\symDecision_{\circss} }}

\newcommand{\dcorth}[0]{\textcolor{\decisioncolor}{\symDecision^{\orthss}_{\chargess} }}
\newcommand{\ddorth}[0]{\textcolor{\decisioncolor}{\symDecision^{\orthss}_{\dischargess} }}
\newcommand{\dporth}[0]{\textcolor{\decisioncolor}{\symDecision^{\orthss}_{\exoaddss} }}
\newcommand{\dmorth}[0]{\textcolor{\decisioncolor}{\symDecision^{\orthss}_{\exoremss} }}

\newcommand{\daddorth}[0]{\textcolor{\decisioncolor}{\symDecision^{\orthss}_{\addss} }}
\newcommand{\dremorth}[0]{\textcolor{\decisioncolor}{\symDecision^{\orthss}_{\remss} }}
\newcommand{\dfloworth}[0]{\textcolor{\decisioncolor}{\symDecision^{\orthss}_{\flowss} }}

\newcommand{\dclevels}[0]{\textcolor{\decisioncolor}{\symDecision^{}_{\chargess}  }}
\newcommand{\ddlevels}[0]{\textcolor{\decisioncolor}{\symDecision^{}_{\dischargess}  }}
\newcommand{\dplevels}[0]{\textcolor{\decisioncolor}{\symDecision^{}_{\exoaddss}  }}
\newcommand{\dmlevels}[0]{\textcolor{\decisioncolor}{\symDecision^{}_{\exoremss}  }}
\newcommand{\dflowlevels}[0]{\textcolor{\decisioncolor}{\symDecision^{}_{\flowss}  }}

\newcommand{\dcmaxlevels}[0]{\textcolor{\paramcolor}{\symDecision^{\maxss}_{\chargess}  }}
\newcommand{\ddmaxlevels}[0]{\textcolor{\paramcolor}{\symDecision^{\maxss}_{\dischargess}  }}
\newcommand{\dpmaxlevels}[0]{\textcolor{\paramcolor}{\symDecision^{\maxss}_{\exoaddss}  }}
\newcommand{\dmmaxlevels}[0]{\textcolor{\paramcolor}{\symDecision^{\maxss}_{\exoremss}  }}
\newcommand{\dflowmaxlevels}[0]{\textcolor{\paramcolor}{\symDecision^{\maxss}_{\flowss}  }}


\newcommand{\dE}[0]{\textcolor{\decisioncolor}{\symDecision_{\Ess}}}

\newcommand{\Emax}[0]{\textcolor{\sizingcolor}{\symEnergy'^{\maxss}}}
\newcommand{\Emin}[0]{\textcolor{\sizingcolor}{\symEnergy'^{\minss}}}
\newcommand{\Eeff}[0]{\textcolor{\sizingcolor}{\symEnergy'^{\usabless}}}
\newcommand{\Ekmax}[0]{\textcolor{\boundscolor}{\symEnergy^{\maxss} }}
\newcommand{\Ekmin}[0]{\textcolor{\boundscolor}{\symEnergy^{\minss} }}
\newcommand{\Ekeff}[0]{\textcolor{\boundscolor}{\symEnergy^{\usabless} }}
\newcommand{\Ek}[0]{\symEnergy }
\newcommand{\Ekprev}[0]{\symEnergy prev}
\newcommand{\iEk}[0]{\textcolor{\binarycolor}{\symBinary_{\Ess} }}
\newcommand{\iEkprev}[0]{\textcolor{\binarycolor}{\symBinary_{\Ess} prev}}

\newcommand{\iijk}[0]{\textcolor{\binarycolor}{\symBinary_{\indexGridNode \indexGridNodeTwo} }}
\newcommand{\ijik}[0]{\textcolor{\binarycolor}{\symBinary_{\indexGridNodeTwo \indexGridNode} }}
\newcommand{\ilk}[0]{\textcolor{\binarycolor}{\symBinary_{\indexGridLines} }}
\newcommand{\ilkprev}[0]{\textcolor{\binarycolor}{\symBinary_{\indexGridLines} prev}}
\newcommand{\ilkmin}[0]{\textcolor{\boundscolor}{\symBinary^{\minss}_{\indexGridLines} }}
\newcommand{\ilkmax}[0]{\textcolor{\boundscolor}{\symBinary^{\maxss}_{\indexGridLines} }}

\newcommand{\betaijk}[0]{\textcolor{\binarycolor}{\symRadiality_{\indexGridNode \indexGridNodeTwo} }}
\newcommand{\betajik}[0]{\textcolor{\binarycolor}{\symRadiality_{\indexGridNodeTwo \indexGridNode} }}


\newcommand{\Sgik}[0]{\textcolor{\complexcolor}{\symApparentPower_{\giss} }}
\newcommand{\Sauxk}[0]{\symApparentPower_{\text{aux}} }
\newcommand{\SauxkA}[0]{\symApparentPower_{\text{aux},\symPhaseA} }
\newcommand{\SauxkB}[0]{\symApparentPower_{\text{aux},\symPhaseB} }
\newcommand{\SauxkC}[0]{\symApparentPower_{\text{aux},\symPhaseC} }
\newcommand{\Sgirated}[0]{\textcolor{\sizingcolor}{\symApparentPower_{\giss}^{'}}}
\newcommand{\lgik}[0]{\textcolor{\paramcolor}{\symLocation_{\giss} }}
\newcommand{\Pgik}[0]{\symPower_{\giss} }
\newcommand{\Qgik}[0]{\symReactivePower_{\giss} }
\newcommand{\Qgikref}[0]{\textcolor{\boundscolor}{\symReactivePower^{\refss}_{\giss} }}
\newcommand{\Qgikmax}[0]{\textcolor{\boundscolor}{\symReactivePower^{\maxss}_{\giss} }}
\newcommand{\Qgikmin}[0]{\textcolor{\boundscolor}{\symReactivePower^{\minss}_{\giss} }}
\newcommand{\Pctotk}[0]{\symPower_{\chargess,\totss} }
\newcommand{\Pdtotk}[0]{\symPower_{\dischargess,\totss} }
\newcommand{\Pstandbyk}[0]{\symPower_{\standbyss} }
\newcommand{\PSlossk}[0]{\symPower_{\Sss,\lossss} }
\newcommand{\Punblossk}[0]{\symPower_{\unbss,\lossss} }
\newcommand{\etagi}{\textcolor{\paramcolor}{\eta_\giss  }}

\newcommand{\Pstandbykmin}[0]{\textcolor{\boundscolor}{\symPower_{\standbyss}^{\minss}  }}
\newcommand{\Pstandbykmax}[0]{\textcolor{\boundscolor}{\symPower_{\standbyss}^{\maxss}  }}

\newcommand{\aPZk}[0]{\textcolor{\paramcolor}{a_{  }^{ \Zss\Pss }  }}
\newcommand{\aPIk}[0]{\textcolor{\paramcolor}{a_{  }^{ \Iss\Pss }  }}
\newcommand{\aPPk}[0]{\textcolor{\paramcolor}{a_{  }^{ \Pss\Pss }  }}
\newcommand{\aQZk}[0]{\textcolor{\paramcolor}{a_{  }^{ \Zss\Qss }  }}
\newcommand{\aQIk}[0]{\textcolor{\paramcolor}{a_{  }^{ \Iss\Qss }  }}
\newcommand{\aQPk}[0]{\textcolor{\paramcolor}{a_{  }^{ \Pss\Qss }  }}

\newcommand{\aPZ}[0]{\textcolor{\paramcolor}{a_{ \indexZIP }^{ \Zss \Pss} }}
\newcommand{\aPI}[0]{\textcolor{\paramcolor}{a_{ \indexZIP }^{ \Iss \Pss} }}
\newcommand{\aPP}[0]{\textcolor{\paramcolor}{a_{\indexZIP  }^{ \Pss\Pss } }}
\newcommand{\aQZ}[0]{\textcolor{\paramcolor}{a_{  \indexZIP}^{ \Zss\Qss } }}
\newcommand{\aQI}[0]{\textcolor{\paramcolor}{a_{ \indexZIP }^{ \Iss \Qss} }}
\newcommand{\aQP}[0]{\textcolor{\paramcolor}{a_{ \indexZIP }^{ \Pss\Qss } }}

\newcommand{\PgikA}[0]{\symPower_{\giss,\symPhaseA} }
\newcommand{\QgikA}[0]{\symReactivePower_{\giss,\symPhaseA} }
\newcommand{\SgikA}[0]{\textcolor{\complexcolor}{\symApparentPower_{\giss,\symPhaseA} }}
\newcommand{\PgikB}[0]{\symPower_{\giss,\symPhaseB} }
\newcommand{\QgikB}[0]{\symReactivePower_{\giss,\symPhaseB} }
\newcommand{\SgikB}[0]{\textcolor{\complexcolor}{\symApparentPower_{\giss,\symPhaseB} }}
\newcommand{\PgikC}[0]{\symPower_{\giss,\symPhaseC} }
\newcommand{\QgikC}[0]{\symReactivePower_{\giss,\symPhaseC} }
\newcommand{\SgikC}[0]{\textcolor{\complexcolor}{\symApparentPower_{\giss,\symPhaseC} }}
\newcommand{\SgiratedA}[0]{\textcolor{\sizingcolor}{\symApparentPower_{\giss,\symPhaseA}^{'}}}
\newcommand{\SgiratedB}[0]{\textcolor{\sizingcolor}{\symApparentPower_{\giss,\symPhaseB}^{'}}}
\newcommand{\SgiratedC}[0]{\textcolor{\sizingcolor}{\symApparentPower_{\giss,\symPhaseC}^{'}}}

\newcommand{\Sgiratedmin}[0]{\textcolor{\ratingboundscolor}{\symApparentPower_{\giss}^{\minss}}}
\newcommand{\SgiratedminA}[0]{\textcolor{\ratingboundscolor}{\symApparentPower_{\giss,\symPhaseA}^{\minss}}}
\newcommand{\SgiratedminB}[0]{\textcolor{\ratingboundscolor}{\symApparentPower_{\giss,\symPhaseB}^{\minss}}}
\newcommand{\SgiratedminC}[0]{\textcolor{\ratingboundscolor}{\symApparentPower_{\giss,\symPhaseC}^{\minss}}}
\newcommand{\Sgiratedmax}[0]{\textcolor{\ratingboundscolor}{\symApparentPower_{\giss}^{\maxss}}}
\newcommand{\SgiratedmaxA}[0]{\textcolor{\ratingboundscolor}{\symApparentPower_{\giss,\symPhaseA}^{\maxss}}}
\newcommand{\SgiratedmaxB}[0]{\textcolor{\ratingboundscolor}{\symApparentPower_{\giss,\symPhaseB}^{\maxss}}}
\newcommand{\SgiratedmaxC}[0]{\textcolor{\ratingboundscolor}{\symApparentPower_{\giss,\symPhaseC}^{\maxss}}}

\newcommand{\Pgikunb}[0]{\symPower_{\giss,\unbss} }
\newcommand{\Qgikunb}[0]{\symReactivePower_{\giss,\unbss} }

\newcommand{\rhoPgikunb}[0]{\textcolor{\paramcolor}{\symLossFactor_{\giss,\Pss, \unbss} }}
\newcommand{\rhoQgikunb}[0]{\textcolor{\paramcolor}{\symLossFactor_{\giss,\Qss,\unbss} }}

\newcommand{\bgik}[0]{\textcolor{\paramcolor}{\symAvailability_{\giss} }}
\newcommand{\bgikA}[0]{\textcolor{\paramcolor}{\symAvailability_{\giss, \symPhaseA} }}
\newcommand{\bgikB}[0]{\textcolor{\paramcolor}{\symAvailability_{\giss, \symPhaseB} }}
\newcommand{\bgikC}[0]{\textcolor{\paramcolor}{\symAvailability_{\giss, \symPhaseC} }}

\newcommand{\dPtogik}[0]{\textcolor{\paramcolor}{\symDirectionality_{\giss}^{\Pss\plusss}   }}
\newcommand{\dPtogikA}[0]{\textcolor{\paramcolor}{\symDirectionality_{\giss, \symPhaseA}^{\Pss\plusss} }}
\newcommand{\dPtogikB}[0]{\textcolor{\paramcolor}{\symDirectionality_{\giss, \symPhaseB}^{\Pss\plusss} }}
\newcommand{\dPtogikC}[0]{\textcolor{\paramcolor}{\symDirectionality_{\giss, \symPhaseC}^{\Pss\plusss} }}

\newcommand{\dQtogik}[0]{\textcolor{\paramcolor}{\symDirectionality_{\giss}^{\Qss\plusss}  }}
\newcommand{\dQtogikA}[0]{\textcolor{\paramcolor}{\symDirectionality_{\giss, \symPhaseA}^{\Qss\plusss}  }}
\newcommand{\dQtogikB}[0]{\textcolor{\paramcolor}{\symDirectionality_{\giss, \symPhaseB}^{\Qss\plusss}  }}
\newcommand{\dQtogikC}[0]{\textcolor{\paramcolor}{\symDirectionality_{\giss, \symPhaseC}^{\Qss\plusss}  }}

\newcommand{\dPfromgik}[0]{\textcolor{\paramcolor}{\symDirectionality_{\giss}^{\Pss\subsss} }}
\newcommand{\dPfromgikA}[0]{\textcolor{\paramcolor}{\symDirectionality_{\giss, \symPhaseA}^{\Pss\subsss} }}
\newcommand{\dPfromgikB}[0]{\textcolor{\paramcolor}{\symDirectionality_{\giss, \symPhaseB}^{\Pss\subsss} }}
\newcommand{\dPfromgikC}[0]{\textcolor{\paramcolor}{\symDirectionality_{\giss, \symPhaseC}^{\Pss\subsss} }}

\newcommand{\dQfromgik}[0]{\textcolor{\paramcolor}{\symDirectionality_{\giss}^{\Qss\subsss}  }}
\newcommand{\dQfromgikA}[0]{\textcolor{\paramcolor}{\symDirectionality_{\giss, \symPhaseA}^{\Qss\subsss}  }}
\newcommand{\dQfromgikB}[0]{\textcolor{\paramcolor}{\symDirectionality_{\giss, \symPhaseB}^{\Qss\subsss}  }}
\newcommand{\dQfromgikC}[0]{\textcolor{\paramcolor}{\symDirectionality_{\giss, \symPhaseC}^{\Qss\subsss}  }}

\newcommand{\rhoDODmax}[0]{\textcolor{\paramcolor}{\symLossFactor_{\Ess}^{\maxss}}}
\newcommand{\rhoDODmin}[0]{\textcolor{\paramcolor}{\symLossFactor_{\Ess}^{\minss}}}
\newcommand{\rhoDODeff}[0]{\textcolor{\paramcolor}{\symLossFactor_{\Ess}^{\usabless}}}

\newcommand{\Eijk}[0]{\textcolor{\complexcolor}{\symEnergy_{\indexGridNode \indexGridNodeTwo}^{}  }}
\newcommand{\Eijkre}[0]{\textcolor{black}{\symEnergy_{\indexGridNode \indexGridNodeTwo}^{\realss}  }}

\newcommand{\Sij}[0]{\textcolor{\complexcolor}{\symApparentPower_{\indexGridNode \indexGridNodeTwo}^{} }}

\newcommand{\Pzref}[0]{\textcolor{\paramcolor}{\symPower_{\indexZIP }^{\refss} }}
\newcommand{\Pinode}[0]{\textcolor{black}{\symPower_{\indexGridNode}^{} }}
\newcommand{\Pz}[0]{\textcolor{black}{\symPower_{\indexZIP}^{} }}

\newcommand{\Pij}[0]{\textcolor{black}{\symPower_{\indexGridNode \indexGridNodeTwo}^{} }}
\newcommand{\Pji}[0]{\textcolor{black}{\symPower_{\indexGridNodeTwo \indexGridNode }^{} }}
\newcommand{\Pijloss}[0]{\textcolor{black}{\symPower_{ \indexGridLines}^{\lossss} }}
\newcommand{\Pijmax}[0]{\textcolor{\paramcolor}{\symPower_{\indexGridNode \indexGridNodeTwo}^{\maxss} }}
\newcommand{\Pijacc}[0]{\textcolor{black}{\symPower_{\indexGridNode \indexGridNodeTwo}^{\star} }}
\newcommand{\Pijrated}[0]{\textcolor{\sizingcolor}{\symPower_{\indexGridNode \indexGridNodeTwo}^{\ratedss} }}
\newcommand{\Pjirated}[0]{\textcolor{\sizingcolor}{\symPower_{\indexGridNodeTwo\indexGridNode }^{\ratedss} }}

\newcommand{\Qij}[0]{\textcolor{black}{\symReactivePower_{\indexGridNode \indexGridNodeTwo}^{} }}
\newcommand{\Qji}[0]{\textcolor{black}{\symReactivePower_{\indexGridNodeTwo \indexGridNode }^{} }}
\newcommand{\Qijrated}[0]{\textcolor{\sizingcolor}{\symReactivePower_{\indexGridNode \indexGridNodeTwo}^{\ratedss} }}
\newcommand{\Qjirated}[0]{\textcolor{\sizingcolor}{\symReactivePower_{\indexGridNodeTwo\indexGridNode }^{\ratedss} }}
\newcommand{\Qijmax}[0]{\textcolor{\paramcolor}{\symReactivePower_{\indexGridNode \indexGridNodeTwo}^{\maxss} }}
\newcommand{\Qijmin}[0]{\textcolor{\paramcolor}{\symReactivePower_{\indexGridNode \indexGridNodeTwo}^{\minss} }}
\newcommand{\Qjimax}[0]{\textcolor{\paramcolor}{\symReactivePower_{\indexGridNodeTwo \indexGridNode }^{\maxss} }}
\newcommand{\Qjimin}[0]{\textcolor{\paramcolor}{\symReactivePower_{\indexGridNodeTwo\indexGridNode }^{\minss} }}
\newcommand{\Qijloss}[0]{\textcolor{black}{\symReactivePower_{ \indexGridLines}^{\lossss} }}

\newcommand{\Sjiks}[0]{\textcolor{\complexcolor}{\symApparentPower_{\indexGridNodeTwo \indexGridNode , \seriesss}^{}  }}
\newcommand{\Sijks}[0]{\textcolor{\complexcolor}{\symApparentPower_{\indexGridNode \indexGridNodeTwo, \seriesss}^{}  }}
\newcommand{\Pijks}[0]{\textcolor{black}{\symPower_{\indexGridNode \indexGridNodeTwo, \seriesss}^{}  }}
\newcommand{\Qijks}[0]{\textcolor{black}{\symReactivePower_{\indexGridNode \indexGridNodeTwo, \seriesss}^{}  }}
\newcommand{\Pjiks}[0]{\textcolor{black}{\symPower_{\indexGridNodeTwo \indexGridNode , \seriesss}^{}  }}
\newcommand{\Qjiks}[0]{\textcolor{black}{\symReactivePower_{\indexGridNodeTwo \indexGridNode , \seriesss}^{}  }}

\newcommand{\Sijlossk}[0]{\textcolor{\complexcolor}{\symApparentPower_{ \indexGridLines}^{\lossss}  }}
\newcommand{\Slossk}[0]{\textcolor{black}{\symApparentPower_{}^{\lossss}  }}

\newcommand{\Sijlosssk}[0]{\textcolor{\complexcolor}{\symApparentPower_{ \indexGridLines, \seriesss}^{\lossss}  }}

\newcommand{\Plossk}[0]{\textcolor{black}{\symPower_{}^{\lossss}  }}
\newcommand{\Pijk}[0]{\textcolor{black}{\symPower_{\indexGridNode \indexGridNodeTwo}^{}  }}
\newcommand{\Pijlossk}[0]{\textcolor{black}{\symPower_{ \indexGridLines}^{\lossss}  }}
\newcommand{\Pijlosssk}[0]{\textcolor{black}{\symPower_{ \indexGridLines, \seriesss}^{\lossss}  }}
\newcommand{\Pijlossshk}[0]{\textcolor{black}{\symPower_{ \indexGridLines, \shuntss}^{\lossss}  }}
\newcommand{\Pilossshk}[0]{\textcolor{black}{\symPower_{\indexGridNode \indexGridNodeTwo, \shuntss }^{\lossss}  }}
\newcommand{\Pjlossshk}[0]{\textcolor{black}{\symPower_{ \indexGridNodeTwo\indexGridNode, \shuntss}^{\lossss}  }}
\newcommand{\Pjik}[0]{\textcolor{black}{\symPower_{\indexGridNodeTwo \indexGridNode}^{}    }}
\newcommand{\Pinodek}[0]{\textcolor{black}{\symPower_{ \indexGridNode}^{}    }}
\newcommand{\Qijk}[0]{\textcolor{black}{\symReactivePower_{\indexGridNode \indexGridNodeTwo}^{}  }}
\newcommand{\Qijkdelta}[0]{\textcolor{black}{\symReactivePower_{\indexGridNode \indexGridNodeTwo}^{\Delta}  }}
\newcommand{\Qijlossk}[0]{\textcolor{black}{\symReactivePower_{ \indexGridLines}^{\lossss}  }}
\newcommand{\Qijlosssk}[0]{\textcolor{black}{\symReactivePower_{ \indexGridLines,\seriesss}^{\lossss}  }}
\newcommand{\Qijlossshk}[0]{\textcolor{black}{\symReactivePower_{ \indexGridLines,\shuntss}^{\lossss}  }}
\newcommand{\Qilossshk}[0]{\textcolor{black}{\symReactivePower_{\indexGridNode\indexGridNodeTwo,\shuntss }^{\lossss}  }}
\newcommand{\Qjlossshk}[0]{\textcolor{black}{\symReactivePower_{ \indexGridNodeTwo\indexGridNode,\shuntss}^{\lossss}  }}
\newcommand{\Qjik}[0]{\textcolor{black}{\symReactivePower_{\indexGridNodeTwo \indexGridNode}^{}    }}
\newcommand{\Qinodek}[0]{\textcolor{black}{\symReactivePower_{ \indexGridNode}^{}    }}
\newcommand{\Sinodek}[0]{\textcolor{\complexcolor}{\symApparentPower_{ \indexGridNode}^{}    }}

\newcommand{\Wijk}[0]{\textcolor{\complexcolor}{W_{\indexGridNode \indexGridNodeTwo}^{}  }} 
\newcommand{\WijkSDP}[0]{\textcolor{\complexcolor}{\mathbf{W}_{}^{}  }} 
\newcommand{\Wiik}[0]{\textcolor{black}{W_{\indexGridNode \indexGridNode}^{}  }} 
\newcommand{\Wjjk}[0]{\textcolor{black}{W_{\indexGridNodeTwo \indexGridNodeTwo}^{}  }} 
\newcommand{\Wjik}[0]{\textcolor{\complexcolor}{W_{\indexGridNodeTwo \indexGridNode}^{}  }} 

\newcommand{\Rijk}[0]{\textcolor{black}{R_{\indexGridNode \indexGridNodeTwo}^{}  }} 
\newcommand{\Tijk}[0]{\textcolor{black}{T_{\indexGridNode \indexGridNodeTwo}^{}  }}  
\newcommand{\Rjik}[0]{\textcolor{black}{R_{\indexGridNodeTwo\indexGridNode }^{}  }} 
\newcommand{\Tjik}[0]{\textcolor{black}{T_{\indexGridNodeTwo\indexGridNode }^{}  }} 
\newcommand{\Rlk}[0]{\textcolor{black}{R_{\indexGridLines }^{}  }}    
\newcommand{\rijk}[0]{\textcolor{black}{r_{\indexGridNode \indexGridNodeTwo}^{}  }}   

\newcommand{\SikA}[0]{\textcolor{\complexcolor}{\symApparentPower_{\indexGridNode , \symPhaseA}^{}  }}
\newcommand{\SikB}[0]{\textcolor{\complexcolor}{\symApparentPower_{\indexGridNode , \symPhaseB}^{}  }}
\newcommand{\SikC}[0]{\textcolor{\complexcolor}{\symApparentPower_{\indexGridNode , \symPhaseC}^{}  }}

\newcommand{\IhshuntSDPrated}[0]{\textcolor{\paramcolor}{\mathbf{\symCurrent}_{ \indexShunt}^{\ratedss}  }}

\newcommand{\PijkA}[0]{\textcolor{black}{\symPower_{\indexGridNode \indexGridNodeTwo, \symPhaseA}^{}  }}
\newcommand{\PijkB}[0]{\textcolor{black}{\symPower_{\indexGridNode \indexGridNodeTwo, \symPhaseB}^{}  }}
\newcommand{\PijkC}[0]{\textcolor{black}{\symPower_{\indexGridNode \indexGridNodeTwo, \symPhaseC}^{}  }}
\newcommand{\QijkA}[0]{\textcolor{black}{\symReactivePower_{\indexGridNode \indexGridNodeTwo, \symPhaseA}^{}  }}
\newcommand{\QijkB}[0]{\textcolor{black}{\symReactivePower_{\indexGridNode \indexGridNodeTwo, \symPhaseB}^{}  }}
\newcommand{\QijkC}[0]{\textcolor{black}{\symReactivePower_{\indexGridNode \indexGridNodeTwo, \symPhaseC}^{}  }}

\newcommand{\PjikA}[0]{\textcolor{black}{\symPower_{ \indexGridNodeTwo\indexGridNode, \symPhaseA}^{}  }}
\newcommand{\PjikB}[0]{\textcolor{black}{\symPower_{ \indexGridNodeTwo\indexGridNode, \symPhaseB}^{}  }}
\newcommand{\PjikC}[0]{\textcolor{black}{\symPower_{ \indexGridNodeTwo\indexGridNode, \symPhaseC}^{}  }}
\newcommand{\QjikA}[0]{\textcolor{black}{\symReactivePower_{ \indexGridNodeTwo\indexGridNode, \symPhaseA}^{}  }}
\newcommand{\QjikB}[0]{\textcolor{black}{\symReactivePower_{ \indexGridNodeTwo\indexGridNode, \symPhaseB}^{}  }}
\newcommand{\QjikC}[0]{\textcolor{black}{\symReactivePower_{ \indexGridNodeTwo\indexGridNode, \symPhaseC}^{}  }}

\newcommand{\PijlosskA}[0]{\textcolor{black}{\symPower_{\indexGridNode \indexGridNodeTwo,\symPhaseA}^{\lossss}  }}
\newcommand{\PijlosskB}[0]{\textcolor{black}{\symPower_{\indexGridNode \indexGridNodeTwo,\symPhaseB}^{\lossss}  }}
\newcommand{\PijlosskC}[0]{\textcolor{black}{\symPower_{\indexGridNode \indexGridNodeTwo,\symPhaseC}^{\lossss}  }}
\newcommand{\QijlosskA}[0]{\textcolor{black}{\symReactivePower_{\indexGridNode \indexGridNodeTwo,\symPhaseA}^{\lossss}  }}
\newcommand{\QijlosskB}[0]{\textcolor{black}{\symReactivePower_{\indexGridNode \indexGridNodeTwo,\symPhaseB}^{\lossss}  }}
\newcommand{\QijlosskC}[0]{\textcolor{black}{\symReactivePower_{\indexGridNode \indexGridNodeTwo,\symPhaseC}^{\lossss}  }}

\newcommand{\PijlosskminA}[0]{\textcolor{\boundscolor}{\symPower_{\indexGridNode \indexGridNodeTwo,\symPhaseA}^{\lossss,\minss}  }}
\newcommand{\PijlosskminB}[0]{\textcolor{\boundscolor}{\symPower_{\indexGridNode \indexGridNodeTwo,\symPhaseB}^{\lossss,\minss}  }}
\newcommand{\PijlosskminC}[0]{\textcolor{\boundscolor}{\symPower_{\indexGridNode \indexGridNodeTwo,\symPhaseC}^{\lossss,\minss}  }}
\newcommand{\QijlosskminA}[0]{\textcolor{\boundscolor}{\symReactivePower_{\indexGridNode \indexGridNodeTwo,\symPhaseA}^{\lossss,\minss}  }}
\newcommand{\QijlosskminB}[0]{\textcolor{\boundscolor}{\symReactivePower_{\indexGridNode \indexGridNodeTwo,\symPhaseB}^{\lossss,\minss}  }}
\newcommand{\QijlosskminC}[0]{\textcolor{\boundscolor}{\symReactivePower_{\indexGridNode \indexGridNodeTwo,\symPhaseC}^{\lossss,\minss}  }}

\newcommand{\PijlosskmaxA}[0]{\textcolor{\boundscolor}{\symPower_{\indexGridNode \indexGridNodeTwo,\symPhaseA}^{\lossss,\maxss}  }}
\newcommand{\PijlosskmaxB}[0]{\textcolor{\boundscolor}{\symPower_{\indexGridNode \indexGridNodeTwo,\symPhaseB}^{\lossss,\maxss}  }}
\newcommand{\PijlosskmaxC}[0]{\textcolor{\boundscolor}{\symPower_{\indexGridNode \indexGridNodeTwo,\symPhaseC}^{\lossss,\maxss}  }}
\newcommand{\QijlosskmaxA}[0]{\textcolor{\boundscolor}{\symReactivePower_{\indexGridNode \indexGridNodeTwo,\symPhaseA}^{\lossss,\maxss}  }}
\newcommand{\QijlosskmaxB}[0]{\textcolor{\boundscolor}{\symReactivePower_{\indexGridNode \indexGridNodeTwo,\symPhaseB}^{\lossss,\maxss}  }}
\newcommand{\QijlosskmaxC}[0]{\textcolor{\boundscolor}{\symReactivePower_{\indexGridNode \indexGridNodeTwo,\symPhaseC}^{\lossss,\maxss}  }}

\newcommand{\Iijkabs}[0]{\textcolor{black}{|\symCurrent_{\indexGridLines\indexGridNode \indexGridNodeTwo}^{}  |}}
\newcommand{\Sijkabs}[0]{\textcolor{black}{|\symApparentPower_{\indexGridLines\indexGridNode \indexGridNodeTwo}^{}  |}}
\newcommand{\Sijk}[0]{\textcolor{\complexcolor}{\symApparentPower_{\indexGridLines\indexGridNode \indexGridNodeTwo}^{}  }}

\newcommand{\SijkSDPrated}[0]{\textcolor{\paramcolor}{\mathbf{\symApparentPower}_{\indexGridLines\indexGridNode \indexGridNodeTwo}^{\ratedss}  }}
\newcommand{\IijkSDPrated}[0]{\textcolor{\paramcolor}{\mathbf{\symCurrent}_{\indexGridLines\indexGridNode \indexGridNodeTwo}^{\ratedss}  }}
\newcommand{\IijksSDPrated}[0]{\textcolor{\paramcolor}{\mathbf{\symCurrent}_{\indexGridLines }^{\seriesss, \ratedss}  }}

\newcommand{\SijkSDPseq}[0]{\textcolor{\complexcolor}{\mathbf{\symApparentPower}_{\indexGridLines\indexGridNode \indexGridNodeTwo}^{\fortescuess}  }}
\newcommand{\SjikSDPseq}[0]{\textcolor{\complexcolor}{\mathbf{\symApparentPower}_{\indexGridLines\indexGridNodeTwo \indexGridNode }^{\fortescuess} }}
\newcommand{\SijkSDPH}[0]{\textcolor{\complexcolor}{(\mathbf{\symApparentPower}_{\indexGridLines\indexGridNode \indexGridNodeTwo})^{\hermitiantranspose}  }}
\newcommand{\SijkSDP}[0]{\textcolor{\complexcolor}{\mathbf{\symApparentPower}_{\indexGridLines\indexGridNode \indexGridNodeTwo}^{}  }}
\newcommand{\SjikSDP}[0]{\textcolor{\complexcolor}{\mathbf{\symApparentPower}_{\indexGridLines\indexGridNodeTwo \indexGridNode }^{}  }}
\newcommand{\SjkkSDP}[0]{\textcolor{\complexcolor}{\mathbf{\symApparentPower}_{ \indexGridLines\indexGridNodeTwo \indexGridNodeThree}^{}  }}
\newcommand{\PjkkSDP}[0]{\textcolor{black}{\mathbf{\symPower}_{ \indexGridLines\indexGridNodeTwo \indexGridNodeThree}^{}  }}
\newcommand{\QjkkSDP}[0]{\textcolor{black}{\mathbf{\symReactivePower}_{ \indexGridLines\indexGridNodeTwo \indexGridNodeThree}^{}  }}

\newcommand{\SijkSDPAdot}[0]{\textcolor{\complexcolor}{\mathbf{\symApparentPower}_{\indexGridLines\indexGridNode \indexGridNodeTwo}^{\symPhaseA \! \cdot}  }}
\newcommand{\SijkSDPBdot}[0]{\textcolor{\complexcolor}{\mathbf{\symApparentPower}_{\indexGridLines\indexGridNode \indexGridNodeTwo}^{\symPhaseB \! \cdot}  }}
\newcommand{\SijkSDPCdot}[0]{\textcolor{\complexcolor}{\mathbf{\symApparentPower}_{\indexGridLines\indexGridNode \indexGridNodeTwo}^{\symPhaseC \! \cdot}  }}
\newcommand{\SijkSDPNdot}[0]{\textcolor{\complexcolor}{\mathbf{\symApparentPower}_{\indexGridLines\indexGridNode \indexGridNodeTwo}^{\symPhaseN \! \cdot}  }}
\newcommand{\SijkSDPGdot}[0]{\textcolor{\complexcolor}{\mathbf{\symApparentPower}_{\indexGridLines\indexGridNode \indexGridNodeTwo}^{\symPhaseG \! \cdot}  }}

\newcommand{\SjikSDPAdot}[0]{\textcolor{\complexcolor}{\mathbf{\symApparentPower}_{\indexGridLines\indexGridNodeTwo\indexGridNode }^{\symPhaseA \! \cdot}  }}
\newcommand{\SjikSDPBdot}[0]{\textcolor{\complexcolor}{\mathbf{\symApparentPower}_{\indexGridLines\indexGridNodeTwo\indexGridNode }^{\symPhaseB \! \cdot}  }}
\newcommand{\SjikSDPCdot}[0]{\textcolor{\complexcolor}{\mathbf{\symApparentPower}_{\indexGridLines\indexGridNodeTwo\indexGridNode }^{\symPhaseC \! \cdot}  }}
\newcommand{\SjikSDPNdot}[0]{\textcolor{\complexcolor}{\mathbf{\symApparentPower}_{\indexGridLines\indexGridNodeTwo\indexGridNode }^{\symPhaseN \! \cdot}  }}

\newcommand{\PijkSDPAdot}[0]{\textcolor{black}{\mathbf{\symPower}_{\indexGridLines\indexGridNode \indexGridNodeTwo}^{\symPhaseA \! \cdot}  }}
\newcommand{\PijkSDPBdot}[0]{\textcolor{black}{\mathbf{\symPower}_{\indexGridLines\indexGridNode \indexGridNodeTwo}^{\symPhaseB \! \cdot}  }}
\newcommand{\PijkSDPCdot}[0]{\textcolor{black}{\mathbf{\symPower}_{\indexGridLines\indexGridNode \indexGridNodeTwo}^{\symPhaseC \! \cdot}  }}
\newcommand{\PijkSDPNdot}[0]{\textcolor{black}{\mathbf{\symPower}_{\indexGridLines\indexGridNode \indexGridNodeTwo}^{\symPhaseN \! \cdot}  }}

\newcommand{\QijkSDPAdot}[0]{\textcolor{black}{\mathbf{\symReactivePower}_{\indexGridLines\indexGridNode \indexGridNodeTwo}^{\symPhaseA \! \cdot}  }}
\newcommand{\QijkSDPBdot}[0]{\textcolor{black}{\mathbf{\symReactivePower}_{\indexGridLines\indexGridNode \indexGridNodeTwo}^{\symPhaseB \! \cdot}  }}
\newcommand{\QijkSDPCdot}[0]{\textcolor{black}{\mathbf{\symReactivePower}_{\indexGridLines\indexGridNode \indexGridNodeTwo}^{\symPhaseC \! \cdot}  }}
\newcommand{\QijkSDPNdot}[0]{\textcolor{black}{\mathbf{\symReactivePower}_{\indexGridLines\indexGridNode \indexGridNodeTwo}^{\symPhaseN \! \cdot}  }}

\newcommand{\PijksSDPAdot}[0]{\textcolor{black}{\mathbf{\symPower}_{\indexGridLines\indexGridNode \indexGridNodeTwo}^{\seriesss,\symPhaseA \! \cdot}  }}
\newcommand{\PijksSDPBdot}[0]{\textcolor{black}{\mathbf{\symPower}_{\indexGridLines\indexGridNode \indexGridNodeTwo}^{\seriesss,\symPhaseB \! \cdot}  }}
\newcommand{\PijksSDPCdot}[0]{\textcolor{black}{\mathbf{\symPower}_{\indexGridLines\indexGridNode \indexGridNodeTwo}^{\seriesss,\symPhaseC \! \cdot}  }}
\newcommand{\PijksSDPNdot}[0]{\textcolor{black}{\mathbf{\symPower}_{\indexGridLines\indexGridNode \indexGridNodeTwo}^{\seriesss,\symPhaseN \! \cdot}  }}

\newcommand{\QijksSDPAdot}[0]{\textcolor{black}{\mathbf{\symReactivePower}_{\indexGridLines\indexGridNode \indexGridNodeTwo}^{\seriesss,\symPhaseA \! \cdot}  }}
\newcommand{\QijksSDPBdot}[0]{\textcolor{black}{\mathbf{\symReactivePower}_{\indexGridLines\indexGridNode \indexGridNodeTwo}^{\seriesss,\symPhaseB \! \cdot}  }}
\newcommand{\QijksSDPCdot}[0]{\textcolor{black}{\mathbf{\symReactivePower}_{\indexGridLines\indexGridNode \indexGridNodeTwo}^{\seriesss,\symPhaseC \! \cdot}  }}
\newcommand{\QijksSDPNdot}[0]{\textcolor{black}{\mathbf{\symReactivePower}_{\indexGridLines\indexGridNode \indexGridNodeTwo}^{\seriesss,\symPhaseN \! \cdot}  }}

\newcommand{\SijkSDPdotA}[0]{\textcolor{\complexcolor}{\mathbf{\symApparentPower}_{\indexGridLines\indexGridNode \indexGridNodeTwo}^{\cdot\!\symPhaseA}  }}
\newcommand{\SijkSDPdotB}[0]{\textcolor{\complexcolor}{\mathbf{\symApparentPower}_{\indexGridLines\indexGridNode \indexGridNodeTwo}^{\cdot\!\symPhaseB }  }}
\newcommand{\SijkSDPdotC}[0]{\textcolor{\complexcolor}{\mathbf{\symApparentPower}_{\indexGridLines\indexGridNode \indexGridNodeTwo}^{\cdot\!\symPhaseC }  }}
\newcommand{\SijkSDPdotN}[0]{\textcolor{\complexcolor}{\mathbf{\symApparentPower}_{\indexGridLines\indexGridNode \indexGridNodeTwo}^{\cdot \!\symPhaseN }  }}

\newcommand{\PijkSDPdotA}[0]{\textcolor{black}{\mathbf{\symPower}_{\indexGridLines\indexGridNode \indexGridNodeTwo}^{\cdot\!\symPhaseA}  }}
\newcommand{\PijkSDPdotB}[0]{\textcolor{black}{\mathbf{\symPower}_{\indexGridLines\indexGridNode \indexGridNodeTwo}^{\cdot\!\symPhaseB }  }}
\newcommand{\PijkSDPdotC}[0]{\textcolor{black}{\mathbf{\symPower}_{\indexGridLines\indexGridNode \indexGridNodeTwo}^{\cdot\!\symPhaseC }  }}
\newcommand{\PijkSDPdotN}[0]{\textcolor{black}{\mathbf{\symPower}_{\indexGridLines\indexGridNode \indexGridNodeTwo}^{\cdot \!\symPhaseN }  }}

\newcommand{\QijkSDPdotA}[0]{\textcolor{black}{\mathbf{\symReactivePower}_{\indexGridLines\indexGridNode \indexGridNodeTwo}^{\cdot\!\symPhaseA}  }}
\newcommand{\QijkSDPdotB}[0]{\textcolor{black}{\mathbf{\symReactivePower}_{\indexGridLines\indexGridNode \indexGridNodeTwo}^{\cdot\!\symPhaseB }  }}
\newcommand{\QijkSDPdotC}[0]{\textcolor{black}{\mathbf{\symReactivePower}_{\indexGridLines\indexGridNode \indexGridNodeTwo}^{\cdot\!\symPhaseC }  }}
\newcommand{\QijkSDPdotN}[0]{\textcolor{black}{\mathbf{\symReactivePower}_{\indexGridLines\indexGridNode \indexGridNodeTwo}^{\cdot \!\symPhaseN }  }}

\newcommand{\PijksSDPdotA}[0]{\textcolor{black}{\mathbf{\symPower}_{\indexGridLines\indexGridNode \indexGridNodeTwo}^{\seriesss,\cdot\!\symPhaseA}  }}
\newcommand{\PijksSDPdotB}[0]{\textcolor{black}{\mathbf{\symPower}_{\indexGridLines\indexGridNode \indexGridNodeTwo}^{\seriesss,\cdot\!\symPhaseB }  }}
\newcommand{\PijksSDPdotC}[0]{\textcolor{black}{\mathbf{\symPower}_{\indexGridLines\indexGridNode \indexGridNodeTwo}^{\seriesss,\cdot\!\symPhaseC }  }}
\newcommand{\PijksSDPdotN}[0]{\textcolor{black}{\mathbf{\symPower}_{\indexGridLines\indexGridNode \indexGridNodeTwo}^{\seriesss,\cdot \!\symPhaseN }  }}

\newcommand{\QijksSDPdotA}[0]{\textcolor{black}{\mathbf{\symReactivePower}_{\indexGridLines\indexGridNode \indexGridNodeTwo}^{\seriesss,\cdot\!\symPhaseA}  }}
\newcommand{\QijksSDPdotB}[0]{\textcolor{black}{\mathbf{\symReactivePower}_{\indexGridLines\indexGridNode \indexGridNodeTwo}^{\seriesss,\cdot\!\symPhaseB }  }}
\newcommand{\QijksSDPdotC}[0]{\textcolor{black}{\mathbf{\symReactivePower}_{\indexGridLines\indexGridNode \indexGridNodeTwo}^{\seriesss,\cdot\!\symPhaseC }  }}
\newcommand{\QijksSDPdotN}[0]{\textcolor{black}{\mathbf{\symReactivePower}_{\indexGridLines\indexGridNode \indexGridNodeTwo}^{\seriesss,\cdot \!\symPhaseN }  }}

\newcommand{\SijskSDPAdot}[0]{\textcolor{\complexcolor}{\mathbf{\symApparentPower}_{\indexGridLines\indexGridNode \indexGridNodeTwo}^{\seriesss, \symPhaseA \! \cdot}  }}
\newcommand{\SijskSDPBdot}[0]{\textcolor{\complexcolor}{\mathbf{\symApparentPower}_{\indexGridLines\indexGridNode \indexGridNodeTwo}^{\seriesss, \symPhaseB \! \cdot}  }}
\newcommand{\SijskSDPCdot}[0]{\textcolor{\complexcolor}{\mathbf{\symApparentPower}_{\indexGridLines\indexGridNode \indexGridNodeTwo}^{\seriesss, \symPhaseC \! \cdot}  }}
\newcommand{\SijskSDPNdot}[0]{\textcolor{\complexcolor}{\mathbf{\symApparentPower}_{\indexGridLines\indexGridNode \indexGridNodeTwo}^{\seriesss, \symPhaseN \! \cdot}  }}
\newcommand{\SijskSDPGdot}[0]{\textcolor{\complexcolor}{\mathbf{\symApparentPower}_{\indexGridLines\indexGridNode \indexGridNodeTwo}^{\seriesss, \symPhaseG \! \cdot}  }}

\newcommand{\SijskSDPdotA}[0]{\textcolor{\complexcolor}{\mathbf{\symApparentPower}_{\indexGridLines\indexGridNode \indexGridNodeTwo}^{\seriesss, \cdot\!\symPhaseA}  }}
\newcommand{\SijskSDPdotB}[0]{\textcolor{\complexcolor}{\mathbf{\symApparentPower}_{\indexGridLines\indexGridNode \indexGridNodeTwo}^{\seriesss, \cdot\!\symPhaseB }  }}
\newcommand{\SijskSDPdotC}[0]{\textcolor{\complexcolor}{\mathbf{\symApparentPower}_{\indexGridLines\indexGridNode \indexGridNodeTwo}^{\seriesss, \cdot\!\symPhaseC }  }}
\newcommand{\SijskSDPdotN}[0]{\textcolor{\complexcolor}{\mathbf{\symApparentPower}_{\indexGridLines\indexGridNode \indexGridNodeTwo}^{\seriesss, \cdot \!\symPhaseN }  }}

\newcommand{\SijksSDPmax}[0]{\textcolor{\paramcolor}{\mathbf{\symApparentPower}_{\indexGridLines \indexGridNode \indexGridNodeTwo}^{\seriesss, \maxss}  }}

\newcommand{\SijksSDP}[0]{\textcolor{\complexcolor}{\mathbf{\symApparentPower}_{\indexGridLines \indexGridNode \indexGridNodeTwo}^{\seriesss}  }}
\newcommand{\SjiksSDP}[0]{\textcolor{\complexcolor}{\mathbf{\symApparentPower}_{\indexGridLines \indexGridNodeTwo\indexGridNode }^{\seriesss}  }}

\newcommand{\SijksSDPH}[0]{\textcolor{\complexcolor}{\mathbf{(\symApparentPower}_{\indexGridLines \indexGridNode \indexGridNodeTwo}^{\seriesss})^{\hermitiantranspose}  }}

\newcommand{\SijksSDPdiag}[0]{\textcolor{\complexcolor}{\mathbf{\symApparentPower}_{\indexGridLines\indexGridNode \indexGridNodeTwo}^{\seriesss, \text{diag}}  }}
\newcommand{\PijksSDPdiag}[0]{\textcolor{black}{\mathbf{\symPower}_{\indexGridLines\indexGridNode \indexGridNodeTwo}^{\seriesss, \text{diag}}  }}
\newcommand{\QijksSDPdiag}[0]{\textcolor{black}{\mathbf{\symReactivePower}_{\indexGridLines\indexGridNode \indexGridNodeTwo}^{\seriesss, \text{diag}}  }}

\newcommand{\PijkSDP}[0]{\textcolor{black}{\mathbf{\symPower}_{\indexGridLines\indexGridNode \indexGridNodeTwo}^{}  }}
\newcommand{\QijkSDP}[0]{\textcolor{black}{\mathbf{\symReactivePower}_{\indexGridLines\indexGridNode \indexGridNodeTwo}^{}  }}
\newcommand{\PijkSDPH}[0]{\textcolor{black}{\mathbf{\symPower}_{\indexGridLines\indexGridNode \indexGridNodeTwo}^{\hermitiantranspose}  }}
\newcommand{\QijkSDPH}[0]{\textcolor{black}{\mathbf{\symReactivePower}_{\indexGridLines\indexGridNode \indexGridNodeTwo}^{\hermitiantranspose}  }}
\newcommand{\PjikSDP}[0]{\textcolor{black}{\mathbf{\symPower}_{\indexGridLines\indexGridNodeTwo\indexGridNode }^{}  }}
\newcommand{\QjikSDP}[0]{\textcolor{black}{\mathbf{\symReactivePower}_{\indexGridLines\indexGridNodeTwo\indexGridNode }^{}  }}

\newcommand{\PijksSDPH}[0]{\textcolor{black}{\mathbf{\symPower}_{\indexGridLines\indexGridNode \indexGridNodeTwo, \seriesss}^{\hermitiantranspose}  }}
\newcommand{\QijksSDPH}[0]{\textcolor{black}{\mathbf{\symReactivePower}_{\indexGridLines\indexGridNode \indexGridNodeTwo, \seriesss}^{\hermitiantranspose}  }}

\newcommand{\PijksSDP}[0]{\textcolor{black}{\mathbf{\symPower}_{\indexGridLines\indexGridNode \indexGridNodeTwo}^{\seriesss}  }}
\newcommand{\QijksSDP}[0]{\textcolor{black}{\mathbf{\symReactivePower}_{\indexGridLines\indexGridNode \indexGridNodeTwo}^{\seriesss}  }}
\newcommand{\PjiksSDP}[0]{\textcolor{black}{\mathbf{\symPower}_{\indexGridLines\indexGridNodeTwo\indexGridNode }^{\seriesss}  }}
\newcommand{\QjiksSDP}[0]{\textcolor{black}{\mathbf{\symReactivePower}_{\indexGridLines\indexGridNodeTwo\indexGridNode }^{\seriesss}  }}

\newcommand{\SijklossSDP}[0]{\textcolor{\complexcolor}{\mathbf{\symApparentPower}_{\indexGridNode \indexGridNodeTwo}^{\lossss}  }}
\newcommand{\SjiklossSDP}[0]{\textcolor{\complexcolor}{\mathbf{\symApparentPower}_{\indexGridNodeTwo\indexGridNode }^{\lossss}  }}
\newcommand{\SlklossSDP}[0]{\textcolor{\complexcolor}{\mathbf{\symApparentPower}_{\indexGridLines}^{\lossss}  }}
\newcommand{\PijklossSDP}[0]{\textcolor{black}{\mathbf{\symPower}_{\indexGridNode \indexGridNodeTwo}^{\lossss}  }}
\newcommand{\QijklossSDP}[0]{\textcolor{black}{\mathbf{\symReactivePower}_{\indexGridNode \indexGridNodeTwo}^{\lossss}  }}

\newcommand{\PijkslossSDP}[0]{\textcolor{black}{\mathbf{\symPower}_{\indexGridNode \indexGridNodeTwo, \seriesss}^{\lossss}  }}
\newcommand{\QijkslossSDP}[0]{\textcolor{black}{\mathbf{\symReactivePower}_{\indexGridNode \indexGridNodeTwo, \seriesss}^{\lossss}  }}

\newcommand{\SijkshlossSDP}[0]{\textcolor{\complexcolor}{\mathbf{\symApparentPower}_{\indexGridLines\indexGridNode \indexGridNodeTwo}^{\lossss, \shuntss}  }}
\newcommand{\SjikshlossSDP}[0]{\textcolor{\complexcolor}{\mathbf{\symApparentPower}_{\indexGridLines\indexGridNodeTwo\indexGridNode}^{\lossss, \shuntss}  }}
\newcommand{\SlkshlossSDP}[0]{\textcolor{\complexcolor}{\mathbf{\symApparentPower}_{\indexGridLines}^{\lossss, \shuntss}  }}

\newcommand{\SijkslossSDP}[0]{\textcolor{\complexcolor}{\mathbf{\symApparentPower}_{\indexGridLines }^{\lossss,\seriesss}  }}
\newcommand{\SjikslossSDP}[0]{\textcolor{\complexcolor}{\mathbf{\symApparentPower}_{\indexGridLines}^{\lossss,\seriesss}  }}
\newcommand{\SlkslossSDP}[0]{\textcolor{\complexcolor}{\mathbf{\symApparentPower}_{\indexGridLines}^{\lossss,\seriesss}  }}

\newcommand{\Sinode}[0]{\textcolor{black}{\symApparentPower_{\indexGridNode }^{}  }}
\newcommand{\Sjnode}[0]{\textcolor{black}{\symApparentPower_{ \indexGridNodeTwo}^{}  }}
\newcommand{\SjnodeSDP}[0]{\textcolor{\complexcolor}{\mathbf{\symApparentPower}_{ \indexGridNodeTwo}^{}  }}
\newcommand{\QjnodeSDP}[0]{\textcolor{black}{\mathbf{\symReactivePower}_{ \indexGridNodeTwo}^{}  }}
\newcommand{\PjnodeSDP}[0]{\textcolor{black}{\mathbf{\symPower}_{ \indexGridNodeTwo}^{}  }}

\newcommand{\SinodeSDP}[0]{\textcolor{\complexcolor}{\mathbf{\symApparentPower}_{ \indexGridNode}^{}  }}
\newcommand{\QinodeSDP}[0]{\textcolor{black}{\mathbf{\symReactivePower}_{ \indexGridNode}^{}  }}
\newcommand{\PinodeSDP}[0]{\textcolor{black}{\mathbf{\symPower}_{ \indexGridNode}^{}  }}

\newcommand{\IinodeSDP}[0]{\textcolor{\complexcolor}{\mathbf{\symCurrent}_{ \indexGridNode}^{}  }}

\newcommand{\IuunitSDPrated}[0]{\textcolor{\paramcolor}{\mathbf{\symCurrent}_{ \indexUnit}^{\ratedss}  }}
\newcommand{\IuunitSDPdeltarated}[0]{\textcolor{\paramcolor}{\mathbf{\symCurrent}_{ \indexUnit}^{\Delta\ratedss}  }}

\newcommand{\IuunitSDPref}[0]{\textcolor{\complexparamcolor}{\mathbf{\symCurrent}_{ \indexUnit}^{\refss}  }}

\newcommand{\IuunitSDP}[0]{\textcolor{\complexcolor}{\mathbf{\symCurrent}_{ \indexUnit}^{}  }}
\newcommand{\IuunitSDPreal}[0]{\textcolor{black}{\mathbf{\symCurrent}_{ \indexUnit}^{\text{re}}  }}
\newcommand{\IuunitSDPimag}[0]{\textcolor{black}{\mathbf{\symCurrent}_{ \indexUnit}^{\text{im}}  }}

\newcommand{\IuunitSDPdelta}[0]{\textcolor{\complexcolor}{\mathbf{\symCurrent}_{ \indexUnit}^{\Delta}  }}
\newcommand{\IuunitSDPdeltareal}[0]{\textcolor{black}{\mathbf{\symCurrent}_{ \indexUnit}^{\Delta, \realss}  }}
\newcommand{\IuunitSDPdeltaimag}[0]{\textcolor{black}{\mathbf{\symCurrent}_{ \indexUnit}^{\Delta, \imagss}  }}

\newcommand{\UuunitSDP}[0]{\textcolor{\complexcolor}{\mathbf{\symVoltage}_{ \indexUnit}^{}  }}

\newcommand{\IbSDP}[0]{\textcolor{\complexcolor}{\mathbf{\symCurrent}_{ \indexShunt}^{}  }}
\newcommand{\IbSDPreal}[0]{\textcolor{black}{\mathbf{\symCurrent}_{ \indexShunt}^{\text{re}}  }}
\newcommand{\IbSDPimag}[0]{\textcolor{black}{\mathbf{\symCurrent}_{ \indexShunt}^{\text{im}}  }}

\newcommand{\SbSDP}[0]{\textcolor{\complexcolor}{\mathbf{\symApparentPower}_{ \indexShunt} }}
\newcommand{\PbSDP}[0]{\textcolor{black}{\mathbf{\symPower}_{ \indexShunt} }}
\newcommand{\QbSDP}[0]{\textcolor{black}{\mathbf{\symReactivePower}_{ \indexShunt} }}

\newcommand{\SuunitSDPint}[0]{\textcolor{\complexcolor}{\mathbf{\symApparentPower}_{ \indexUnit}^{'}  }}

\newcommand{\SuunitSDPref}[0]{\textcolor{\complexparamcolor}{\mathbf{\symApparentPower}_{ \indexUnit}^{\refss}  }}
\newcommand{\SuunitSDPrefdelta}[0]{\textcolor{\complexparamcolor}{\mathbf{\symApparentPower}_{ \indexUnit}^{\Delta\refss}  }}
\newcommand{\PuunitSDPref}[0]{\textcolor{\paramcolor}{\mathbf{\symPower}_{ \indexUnit}^{\refss}  }}
\newcommand{\QuunitSDPref}[0]{\textcolor{\paramcolor}{\mathbf{\symReactivePower}_{ \indexUnit}^{\refss}  }}

\newcommand{\PuunitSDPmin}[0]{\textcolor{\paramcolor}{\mathbf{\symPower}_{ \indexUnit}^{\minss}  }}
\newcommand{\QuunitSDPmin}[0]{\textcolor{\paramcolor}{\mathbf{\symReactivePower}_{ \indexUnit}^{\minss}  }}
\newcommand{\PuunitSDPmax}[0]{\textcolor{\paramcolor}{\mathbf{\symPower}_{ \indexUnit}^{\maxss}  }}
\newcommand{\QuunitSDPmax}[0]{\textcolor{\paramcolor}{\mathbf{\symReactivePower}_{ \indexUnit}^{\maxss}  }}
\newcommand{\SuunitSDPmax}[0]{\textcolor{\paramcolor}{\mathbf{\symApparentPower}_{ \indexUnit}^{\maxss}  }}

\newcommand{\PuunitSDPmindelta}[0]{\textcolor{\paramcolor}{\mathbf{\symPower}_{ \indexUnit}^{\Delta\minss}  }}
\newcommand{\QuunitSDPmindelta}[0]{\textcolor{\paramcolor}{\mathbf{\symReactivePower}_{ \indexUnit}^{\Delta\minss}  }}
\newcommand{\PuunitSDPmaxdelta}[0]{\textcolor{\paramcolor}{\mathbf{\symPower}_{ \indexUnit}^{\Delta\maxss}  }}
\newcommand{\QuunitSDPmaxdelta}[0]{\textcolor{\paramcolor}{\mathbf{\symReactivePower}_{ \indexUnit}^{\Delta\maxss}  }}

\newcommand{\SuunitSDP}[0]{\textcolor{\complexcolor}{\mathbf{\symApparentPower}_{ \indexUnit}^{}  }}
\newcommand{\PuunitSDP}[0]{\textcolor{black}{\mathbf{\symPower}_{ \indexUnit}^{}  }}
\newcommand{\QuunitSDP}[0]{\textcolor{black}{\mathbf{\symReactivePower}_{ \indexUnit}^{}  }}

\newcommand{\SuunitSDPdelta}[0]{\textcolor{\complexcolor}{\mathbf{\symApparentPower}_{ \indexUnit}^{\Delta}  }}
\newcommand{\PuunitSDPdelta}[0]{\textcolor{black}{\mathbf{\symPower}_{ \indexUnit}^{\Delta}  }}
\newcommand{\QuunitSDPdelta}[0]{\textcolor{black}{\mathbf{\symReactivePower}_{ \indexUnit}^{\Delta}  }}

\newcommand{\Tuunitdelta}[0]{\textcolor{\paramcolor}{\mathbf{T}^{\Delta}  }}
\newcommand{\Nwye}[0]{\textcolor{\paramcolor}{\mathbf{N}^{\text{pn}}  }}

\newcommand{\XljiSDPdelta}[0]{\textcolor{\complexcolor}{\mathbf{X}_{ lji}  }}
\newcommand{\XlijSDPdelta}[0]{\textcolor{\complexcolor}{\mathbf{X}_{ lij}  }}

\newcommand{\XuunitSDPdelta}[0]{\textcolor{\complexcolor}{\mathbf{X}_{ \indexUnit}^{\Delta}  }}
\newcommand{\XuunitSDPdeltareal}[0]{\textcolor{black}{\mathbf{X}_{ \indexUnit}^{\Delta\realss}  }}
\newcommand{\XuunitSDPdeltaimag}[0]{\textcolor{black}{\mathbf{X}_{ \indexUnit}^{\Delta\imagss}  }}

\newcommand{\IuunitAN}[0]{\textcolor{\complexcolor}{{\symCurrent}_{ \indexUnit, \symPhaseA}^{}  }}
\newcommand{\IuunitBN}[0]{\textcolor{\complexcolor}{{\symCurrent}_{ \indexUnit, \symPhaseB}^{}  }}
\newcommand{\IuunitCN}[0]{\textcolor{\complexcolor}{{\symCurrent}_{ \indexUnit, \symPhaseC}^{}  }}
\newcommand{\IuunitN}[0]{\textcolor{\complexcolor}{{\symCurrent}_{ \indexUnit, \symPhaseN}^{}  }}
\newcommand{\IuunitG}[0]{\textcolor{\complexcolor}{{\symCurrent}_{ \indexUnit, \symPhaseG}^{}  }}

\newcommand{\IuunitArated}[0]{\textcolor{\paramcolor}{{\symCurrent}_{ \indexUnit, \symPhaseA}^{\ratedss}  }}
\newcommand{\IuunitBrated}[0]{\textcolor{\paramcolor}{{\symCurrent}_{ \indexUnit, \symPhaseB}^{\ratedss}  }}
\newcommand{\IuunitCrated}[0]{\textcolor{\paramcolor}{{\symCurrent}_{ \indexUnit, \symPhaseC}^{\ratedss}  }}
\newcommand{\IuunitNrated}[0]{\textcolor{\paramcolor}{{\symCurrent}_{ \indexUnit, \symPhaseN}^{\ratedss}  }}

\newcommand{\IuunitABrated}[0]{\textcolor{\paramcolor}{{\symCurrent}_{ \indexUnit, \symPhaseA\symPhaseB}^{\ratedss}  }}
\newcommand{\IuunitBCrated}[0]{\textcolor{\paramcolor}{{\symCurrent}_{ \indexUnit, \symPhaseB\symPhaseC}^{\ratedss}  }}
\newcommand{\IuunitCArated}[0]{\textcolor{\paramcolor}{{\symCurrent}_{ \indexUnit, \symPhaseC\symPhaseA}^{\ratedss}  }}

\newcommand{\IbA}[0]{\textcolor{\complexcolor}{{\symCurrent}_{ \indexShunt, \symPhaseA}^{}  }}
\newcommand{\IbB}[0]{\textcolor{\complexcolor}{{\symCurrent}_{ \indexShunt, \symPhaseB}^{}  }}
\newcommand{\IbC}[0]{\textcolor{\complexcolor}{{\symCurrent}_{ \indexShunt, \symPhaseC}^{}  }}
\newcommand{\IbN}[0]{\textcolor{\complexcolor}{{\symCurrent}_{ \indexShunt, \symPhaseN}^{}  }}

\newcommand{\IuunitNreal}[0]{\textcolor{black}{{\symCurrent}_{ \indexUnit, \symPhaseN}^{\text{re}}  }}
\newcommand{\IuunitNimag}[0]{\textcolor{black}{{\symCurrent}_{ \indexUnit, \symPhaseN}^{\text{im}}  }}

\newcommand{\IuunitAB}[0]{\textcolor{\complexcolor}{{\symCurrent}_{ \indexUnit, \symPhaseA\symPhaseB}^{\Delta}  }}
\newcommand{\IuunitBC}[0]{\textcolor{\complexcolor}{{\symCurrent}_{ \indexUnit, \symPhaseB\symPhaseC}^{\Delta}  }}
\newcommand{\IuunitCA}[0]{\textcolor{\complexcolor}{{\symCurrent}_{ \indexUnit, \symPhaseC\symPhaseA}^{\Delta}  }}

\newcommand{\IuunitBA}[0]{\textcolor{\complexcolor}{{\symCurrent}_{ \indexUnit,\symPhaseB \symPhaseA}^{\Delta}  }}
\newcommand{\IuunitCB}[0]{\textcolor{\complexcolor}{{\symCurrent}_{ \indexUnit, \symPhaseC\symPhaseB}^{\Delta}  }}
\newcommand{\IuunitAC}[0]{\textcolor{\complexcolor}{{\symCurrent}_{ \indexUnit, \symPhaseA\symPhaseC}^{\Delta}  }}

\newcommand{\SuunitAAdelta}[0]{\textcolor{\complexcolor}{{\symApparentPower}_{ \indexUnit, \symPhaseA\symPhaseA}^{\Delta}  }}
\newcommand{\SuunitBBdelta}[0]{\textcolor{\complexcolor}{{\symApparentPower}_{ \indexUnit, \symPhaseB\symPhaseB}^{\Delta}  }}
\newcommand{\SuunitCCdelta}[0]{\textcolor{\complexcolor}{{\symApparentPower}_{ \indexUnit, \symPhaseC\symPhaseC}^{\Delta}  }}

\newcommand{\PuunitAAdelta}[0]{\textcolor{black}{{\symPower}_{ \indexUnit, \symPhaseA\symPhaseA}^{\Delta}  }}
\newcommand{\PuunitBBdelta}[0]{\textcolor{black}{{\symPower}_{ \indexUnit, \symPhaseB\symPhaseB}^{\Delta}  }}
\newcommand{\PuunitCCdelta}[0]{\textcolor{black}{{\symPower}_{ \indexUnit, \symPhaseC\symPhaseC}^{\Delta}  }}

\newcommand{\QuunitAAdelta}[0]{\textcolor{black}{{\symReactivePower}_{ \indexUnit, \symPhaseA\symPhaseA}^{\Delta}  }}
\newcommand{\QuunitBBdelta}[0]{\textcolor{black}{{\symReactivePower}_{ \indexUnit, \symPhaseB\symPhaseB}^{\Delta}  }}
\newcommand{\QuunitCCdelta}[0]{\textcolor{black}{{\symReactivePower}_{ \indexUnit, \symPhaseC\symPhaseC}^{\Delta}  }}

\newcommand{\SuunitABdelta}[0]{\textcolor{\complexcolor}{{\symApparentPower}_{ \indexUnit, \symPhaseA\symPhaseB}^{\Delta}  }}
\newcommand{\SuunitBCdelta}[0]{\textcolor{\complexcolor}{{\symApparentPower}_{ \indexUnit, \symPhaseB\symPhaseC}^{\Delta}  }}
\newcommand{\SuunitCAdelta}[0]{\textcolor{\complexcolor}{{\symApparentPower}_{ \indexUnit, \symPhaseC\symPhaseA}^{\Delta}  }}

\newcommand{\SuunitBAdelta}[0]{\textcolor{\complexcolor}{{\symApparentPower}_{ \indexUnit,\symPhaseB \symPhaseA}^{\Delta}  }}
\newcommand{\SuunitCBdelta}[0]{\textcolor{\complexcolor}{{\symApparentPower}_{ \indexUnit, \symPhaseC\symPhaseB}^{\Delta}  }}
\newcommand{\SuunitACdelta}[0]{\textcolor{\complexcolor}{{\symApparentPower}_{ \indexUnit, \symPhaseA\symPhaseC}^{\Delta}  }}

\newcommand{\SuunitABCdelta}[0]{\textcolor{\complexcolor}{{\symApparentPower}_{ \indexUnit, \symPhaseA\symPhaseB\symPhaseC}^{\Delta}  }}
\newcommand{\SuunitBCAdelta}[0]{\textcolor{\complexcolor}{{\symApparentPower}_{ \indexUnit, \symPhaseB\symPhaseC\symPhaseA}^{\Delta}  }}
\newcommand{\SuunitCABdelta}[0]{\textcolor{\complexcolor}{{\symApparentPower}_{ \indexUnit, \symPhaseC\symPhaseA\symPhaseB}^{\Delta}  }}

\newcommand{\PuunitABdeltamax}[0]{\textcolor{\paramcolor}{{\symPower}_{ \indexUnit, \symPhaseA\symPhaseB}^{\Delta, \maxss}  }}
\newcommand{\PuunitBCdeltamax}[0]{\textcolor{\paramcolor}{{\symPower}_{ \indexUnit, \symPhaseB\symPhaseC}^{\Delta, \maxss}  }}
\newcommand{\PuunitCAdeltamax}[0]{\textcolor{\paramcolor}{{\symPower}_{ \indexUnit, \symPhaseC\symPhaseA}^{\Delta, \maxss}  }}

\newcommand{\QuunitABdeltamax}[0]{\textcolor{\paramcolor}{{\symReactivePower}_{ \indexUnit, \symPhaseA\symPhaseB}^{\Delta, \maxss}  }}
\newcommand{\QuunitBCdeltamax}[0]{\textcolor{\paramcolor}{{\symReactivePower}_{ \indexUnit, \symPhaseB\symPhaseC}^{\Delta, \maxss}  }}
\newcommand{\QuunitCAdeltamax}[0]{\textcolor{\paramcolor}{{\symReactivePower}_{ \indexUnit, \symPhaseC\symPhaseA}^{\Delta, \maxss}  }}

\newcommand{\PuunitABdeltamin}[0]{\textcolor{\paramcolor}{{\symPower}_{ \indexUnit, \symPhaseA\symPhaseB}^{\Delta, \minss}  }}
\newcommand{\PuunitBCdeltamin}[0]{\textcolor{\paramcolor}{{\symPower}_{ \indexUnit, \symPhaseB\symPhaseC}^{\Delta, \minss}  }}
\newcommand{\PuunitCAdeltamin}[0]{\textcolor{\paramcolor}{{\symPower}_{ \indexUnit, \symPhaseC\symPhaseA}^{\Delta, \minss}  }}

\newcommand{\QuunitABdeltamin}[0]{\textcolor{\paramcolor}{{\symReactivePower}_{ \indexUnit, \symPhaseA\symPhaseB}^{\Delta, \minss}  }}
\newcommand{\QuunitBCdeltamin}[0]{\textcolor{\paramcolor}{{\symReactivePower}_{ \indexUnit, \symPhaseB\symPhaseC}^{\Delta, \minss}  }}
\newcommand{\QuunitCAdeltamin}[0]{\textcolor{\paramcolor}{{\symReactivePower}_{ \indexUnit, \symPhaseC\symPhaseA}^{\Delta, \minss}  }}

\newcommand{\PuunitAAmax}[0]{\textcolor{\paramcolor}{{\symPower}_{ \indexUnit, \symPhaseA\symPhaseA}^{\maxss}  }}
\newcommand{\PuunitBBmax}[0]{\textcolor{\paramcolor}{{\symPower}_{ \indexUnit, \symPhaseB\symPhaseB}^{ \maxss}  }}
\newcommand{\PuunitCCmax}[0]{\textcolor{\paramcolor}{{\symPower}_{ \indexUnit, \symPhaseC\symPhaseC}^{ \maxss}  }}

\newcommand{\QuunitAAmax}[0]{\textcolor{\paramcolor}{{\symReactivePower}_{ \indexUnit, \symPhaseA\symPhaseA}^{\maxss}  }}
\newcommand{\QuunitBBmax}[0]{\textcolor{\paramcolor}{{\symReactivePower}_{ \indexUnit, \symPhaseB\symPhaseB}^{ \maxss}  }}
\newcommand{\QuunitCCmax}[0]{\textcolor{\paramcolor}{{\symReactivePower}_{ \indexUnit, \symPhaseC\symPhaseC}^{ \maxss}  }}

\newcommand{\PuunitAAmin}[0]{\textcolor{\paramcolor}{{\symPower}_{ \indexUnit, \symPhaseA\symPhaseA}^{\minss}  }}
\newcommand{\PuunitBBmin}[0]{\textcolor{\paramcolor}{{\symPower}_{ \indexUnit, \symPhaseB\symPhaseB}^{ \minss}  }}
\newcommand{\PuunitCCmin}[0]{\textcolor{\paramcolor}{{\symPower}_{ \indexUnit, \symPhaseC\symPhaseC}^{ \minss}  }}

\newcommand{\QuunitAAmin}[0]{\textcolor{\paramcolor}{{\symReactivePower}_{ \indexUnit, \symPhaseA\symPhaseA}^{\minss}  }}
\newcommand{\QuunitBBmin}[0]{\textcolor{\paramcolor}{{\symReactivePower}_{ \indexUnit, \symPhaseB\symPhaseB}^{ \minss}  }}
\newcommand{\QuunitCCmin}[0]{\textcolor{\paramcolor}{{\symReactivePower}_{ \indexUnit, \symPhaseC\symPhaseC}^{ \minss}  }}

\newcommand{\PuunitABdelta}[0]{\textcolor{black}{{\symPower}_{ \indexUnit, \symPhaseA\symPhaseB}^{\Delta}  }}
\newcommand{\PuunitBCdelta}[0]{\textcolor{black}{{\symPower}_{ \indexUnit, \symPhaseB\symPhaseC}^{\Delta}  }}
\newcommand{\PuunitCAdelta}[0]{\textcolor{black}{{\symPower}_{ \indexUnit, \symPhaseC\symPhaseA}^{\Delta}  }}

\newcommand{\PuunitBAdelta}[0]{\textcolor{black}{{\symPower}_{ \indexUnit,\symPhaseB \symPhaseA}^{\Delta}  }}
\newcommand{\PuunitCBdelta}[0]{\textcolor{black}{{\symPower}_{ \indexUnit, \symPhaseC\symPhaseB}^{\Delta}  }}
\newcommand{\PuunitACdelta}[0]{\textcolor{black}{{\symPower}_{ \indexUnit, \symPhaseA\symPhaseC}^{\Delta}  }}

\newcommand{\QuunitABdelta}[0]{\textcolor{black}{{\symReactivePower}_{ \indexUnit, \symPhaseA\symPhaseB}^{\Delta}  }}
\newcommand{\QuunitBCdelta}[0]{\textcolor{black}{{\symReactivePower}_{ \indexUnit, \symPhaseB\symPhaseC}^{\Delta}  }}
\newcommand{\QuunitCAdelta}[0]{\textcolor{black}{{\symReactivePower}_{ \indexUnit, \symPhaseC\symPhaseA}^{\Delta}  }}

\newcommand{\QuunitBAdelta}[0]{\textcolor{black}{{\symReactivePower}_{ \indexUnit,\symPhaseB \symPhaseA}^{\Delta}  }}
\newcommand{\QuunitCBdelta}[0]{\textcolor{black}{{\symReactivePower}_{ \indexUnit, \symPhaseC\symPhaseB}^{\Delta}  }}
\newcommand{\QuunitACdelta}[0]{\textcolor{black}{{\symReactivePower}_{ \indexUnit, \symPhaseA\symPhaseC}^{\Delta}  }}

\newcommand{\PuunitAN}[0]{\textcolor{black}{{\symPower}_{ \indexUnit, \symPhaseA}^{}  }}
\newcommand{\PuunitBN}[0]{\textcolor{black}{{\symPower}_{ \indexUnit, \symPhaseB}^{}  }}
\newcommand{\PuunitCN}[0]{\textcolor{black}{{\symPower}_{ \indexUnit, \symPhaseC}^{}  }}

\newcommand{\SuunitrefAA}[0]{\textcolor{\complexparamcolor}{{\symApparentPower}_{ \indexUnit, \symPhaseA\symPhaseA}^{\refss}  }}
\newcommand{\SuunitrefBB}[0]{\textcolor{\complexparamcolor}{{\symApparentPower}_{ \indexUnit, \symPhaseB\symPhaseB}^{\refss}  }}
\newcommand{\SuunitrefCC}[0]{\textcolor{\complexparamcolor}{{\symApparentPower}_{ \indexUnit, \symPhaseC\symPhaseC}^{\refss}  }}

\newcommand{\SuunitrefAB}[0]{\textcolor{\complexparamcolor}{{\symApparentPower}_{ \indexUnit, \symPhaseA\symPhaseB}^{\refss}  }}
\newcommand{\SuunitrefBC}[0]{\textcolor{\complexparamcolor}{{\symApparentPower}_{ \indexUnit, \symPhaseB\symPhaseC}^{\refss}  }}
\newcommand{\SuunitrefCA}[0]{\textcolor{\complexparamcolor}{{\symApparentPower}_{ \indexUnit, \symPhaseC\symPhaseA}^{\refss}  }}

\newcommand{\PuunitrefAB}[0]{\textcolor{\paramcolor}{{\symPower}_{ \indexUnit, \symPhaseA\symPhaseB}^{\refss}  }}
\newcommand{\PuunitrefBC}[0]{\textcolor{\paramcolor}{{\symPower}_{ \indexUnit, \symPhaseB\symPhaseC}^{\refss}  }}
\newcommand{\PuunitrefCA}[0]{\textcolor{\paramcolor}{{\symPower}_{ \indexUnit, \symPhaseC\symPhaseA}^{\refss}  }}

\newcommand{\QuunitrefAB}[0]{\textcolor{\paramcolor}{{\symReactivePower}_{ \indexUnit, \symPhaseA\symPhaseB}^{\refss}  }}
\newcommand{\QuunitrefBC}[0]{\textcolor{\paramcolor}{{\symReactivePower}_{ \indexUnit, \symPhaseB\symPhaseC}^{\refss}  }}
\newcommand{\QuunitrefCA}[0]{\textcolor{\paramcolor}{{\symReactivePower}_{ \indexUnit, \symPhaseC\symPhaseA}^{\refss}  }}

\newcommand{\PuunitrefAA}[0]{\textcolor{\paramcolor}{{\symPower}_{ \indexUnit, \symPhaseA\symPhaseA}^{\refss}  }}
\newcommand{\PuunitrefBB}[0]{\textcolor{\paramcolor}{{\symPower}_{ \indexUnit, \symPhaseB\symPhaseB}^{\refss}  }}
\newcommand{\PuunitrefCC}[0]{\textcolor{\paramcolor}{{\symPower}_{ \indexUnit, \symPhaseC\symPhaseC}^{\refss}  }}

\newcommand{\QuunitrefAA}[0]{\textcolor{\paramcolor}{{\symReactivePower}_{ \indexUnit, \symPhaseA\symPhaseA}^{\refss}  }}
\newcommand{\QuunitrefBB}[0]{\textcolor{\paramcolor}{{\symReactivePower}_{ \indexUnit, \symPhaseB\symPhaseB}^{\refss}  }}
\newcommand{\QuunitrefCC}[0]{\textcolor{\paramcolor}{{\symReactivePower}_{ \indexUnit, \symPhaseC\symPhaseC}^{\refss}  }}

\newcommand{\SuunitAdot}[0]{\textcolor{\complexcolor}{\mathbf{\symApparentPower}_{ \indexUnit}^{\symPhaseA\!\cdot}  }}
\newcommand{\SuunitBdot}[0]{\textcolor{\complexcolor}{\mathbf{\symApparentPower}_{ \indexUnit}^{\symPhaseB\!\cdot}  }}
\newcommand{\SuunitCdot}[0]{\textcolor{\complexcolor}{\mathbf{\symApparentPower}_{ \indexUnit}^{\symPhaseC\!\cdot}  }}
\newcommand{\SuunitNdot}[0]{\textcolor{\complexcolor}{\mathbf{\symApparentPower}_{ \indexUnit}^{\symPhaseN\!\cdot}  }}
\newcommand{\SuunitGdot}[0]{\textcolor{\complexcolor}{\mathbf{\symApparentPower}_{ \indexUnit}^{\symPhaseG\!\cdot}  }}

\newcommand{\SuunitdotN}[0]{\textcolor{\complexcolor}{\mathbf{\symApparentPower}_{ \indexUnit}^{\cdot\! \symPhaseN}  }}

\newcommand{\PuunitNdot}[0]{\textcolor{black}{\mathbf{\symPower}_{ \indexUnit}^{\symPhaseN\!\cdot}  }}
\newcommand{\QuunitNdot}[0]{\textcolor{black}{\mathbf{\symReactivePower}_{ \indexUnit}^{\symPhaseN\!\cdot}  }}

\newcommand{\XuunitAA}[0]{\textcolor{\complexcolor}{{X}_{ \indexUnit, \symPhaseA\symPhaseA}^{}  }}
\newcommand{\XuunitAB}[0]{\textcolor{\complexcolor}{{X}_{ \indexUnit, \symPhaseA\symPhaseB}^{}  }}
\newcommand{\XuunitAC}[0]{\textcolor{\complexcolor}{{X}_{ \indexUnit, \symPhaseA\symPhaseC}^{}  }}

\newcommand{\XuunitBA}[0]{\textcolor{\complexcolor}{{X}_{ \indexUnit, \symPhaseB\symPhaseA}^{}  }}
\newcommand{\XuunitBB}[0]{\textcolor{\complexcolor}{{X}_{ \indexUnit, \symPhaseB\symPhaseB}^{}  }}
\newcommand{\XuunitBC}[0]{\textcolor{\complexcolor}{{X}_{ \indexUnit, \symPhaseB\symPhaseC}^{}  }}

\newcommand{\XuunitCA}[0]{\textcolor{\complexcolor}{{X}_{ \indexUnit, \symPhaseC\symPhaseA}^{}  }}
\newcommand{\XuunitCB}[0]{\textcolor{\complexcolor}{{X}_{ \indexUnit, \symPhaseC\symPhaseB}^{}  }}
\newcommand{\XuunitCC}[0]{\textcolor{\complexcolor}{{X}_{ \indexUnit, \symPhaseC\symPhaseC}^{}  }}

\newcommand{\XuunitAAreal}[0]{\textcolor{black}{{X}_{ \indexUnit, \symPhaseA\symPhaseA}^{\realss}  }}
\newcommand{\XuunitABreal}[0]{\textcolor{black}{{X}_{ \indexUnit, \symPhaseA\symPhaseB}^{\realss}  }}
\newcommand{\XuunitACreal}[0]{\textcolor{black}{{X}_{ \indexUnit, \symPhaseA\symPhaseC}^{\realss}  }}

\newcommand{\XuunitBAreal}[0]{\textcolor{black}{{X}_{ \indexUnit, \symPhaseB\symPhaseA}^{\realss}  }}
\newcommand{\XuunitBBreal}[0]{\textcolor{black}{{X}_{ \indexUnit, \symPhaseB\symPhaseB}^{\realss}  }}
\newcommand{\XuunitBCreal}[0]{\textcolor{black}{{X}_{ \indexUnit, \symPhaseB\symPhaseC}^{\realss}  }}

\newcommand{\XuunitCAreal}[0]{\textcolor{black}{{X}_{ \indexUnit, \symPhaseC\symPhaseA}^{\realss}  }}
\newcommand{\XuunitCBreal}[0]{\textcolor{black}{{X}_{ \indexUnit, \symPhaseC\symPhaseB}^{\realss}  }}
\newcommand{\XuunitCCreal}[0]{\textcolor{black}{{X}_{ \indexUnit, \symPhaseC\symPhaseC}^{\realss}  }}

\newcommand{\XuunitAAimag}[0]{\textcolor{black}{{X}_{ \indexUnit, \symPhaseA\symPhaseA}^{\imagss}  }}
\newcommand{\XuunitABimag}[0]{\textcolor{black}{{X}_{ \indexUnit, \symPhaseA\symPhaseB}^{\imagss}  }}
\newcommand{\XuunitACimag}[0]{\textcolor{black}{{X}_{ \indexUnit, \symPhaseA\symPhaseC}^{\imagss}  }}

\newcommand{\XuunitBAimag}[0]{\textcolor{black}{{X}_{ \indexUnit, \symPhaseB\symPhaseA}^{\imagss}  }}
\newcommand{\XuunitBBimag}[0]{\textcolor{black}{{X}_{ \indexUnit, \symPhaseB\symPhaseB}^{\imagss}  }}
\newcommand{\XuunitBCimag}[0]{\textcolor{black}{{X}_{ \indexUnit, \symPhaseB\symPhaseC}^{\imagss}  }}

\newcommand{\XuunitCAimag}[0]{\textcolor{black}{{X}_{ \indexUnit, \symPhaseC\symPhaseA}^{\imagss}  }}
\newcommand{\XuunitCBimag}[0]{\textcolor{black}{{X}_{ \indexUnit, \symPhaseC\symPhaseB}^{\imagss}  }}
\newcommand{\XuunitCCimag}[0]{\textcolor{black}{{X}_{ \indexUnit, \symPhaseC\symPhaseC}^{\imagss}  }}

\newcommand{\SuunitAA}[0]{\textcolor{\complexcolor}{{\symApparentPower}_{ \indexUnit, \symPhaseA\symPhaseA}^{}  }}
\newcommand{\SuunitAB}[0]{\textcolor{\complexcolor}{{\symApparentPower}_{ \indexUnit, \symPhaseA\symPhaseB}^{}  }}
\newcommand{\SuunitAC}[0]{\textcolor{\complexcolor}{{\symApparentPower}_{ \indexUnit, \symPhaseA\symPhaseC}^{}  }}
\newcommand{\SuunitAN}[0]{\textcolor{\complexcolor}{{\symApparentPower}_{ \indexUnit, \symPhaseA\symPhaseN}^{}  }}

\newcommand{\SuunitBA}[0]{\textcolor{\complexcolor}{{\symApparentPower}_{ \indexUnit, \symPhaseB\symPhaseA}^{}  }}
\newcommand{\SuunitBB}[0]{\textcolor{\complexcolor}{{\symApparentPower}_{ \indexUnit, \symPhaseB\symPhaseB}^{}  }}
\newcommand{\SuunitBC}[0]{\textcolor{\complexcolor}{{\symApparentPower}_{ \indexUnit, \symPhaseB\symPhaseC}^{}  }}
\newcommand{\SuunitBN}[0]{\textcolor{\complexcolor}{{\symApparentPower}_{ \indexUnit, \symPhaseB\symPhaseN}^{}  }}

\newcommand{\SuunitCA}[0]{\textcolor{\complexcolor}{{\symApparentPower}_{ \indexUnit, \symPhaseC\symPhaseA}^{}  }}
\newcommand{\SuunitCB}[0]{\textcolor{\complexcolor}{{\symApparentPower}_{ \indexUnit, \symPhaseC\symPhaseB}^{}  }}
\newcommand{\SuunitCC}[0]{\textcolor{\complexcolor}{{\symApparentPower}_{ \indexUnit, \symPhaseC\symPhaseC}^{}  }}
\newcommand{\SuunitCN}[0]{\textcolor{\complexcolor}{{\symApparentPower}_{ \indexUnit, \symPhaseC\symPhaseN}^{}  }}
\newcommand{\SuunitNN}[0]{\textcolor{\complexcolor}{{\symApparentPower}_{ \indexUnit, \symPhaseN\symPhaseN}^{}  }}
\newcommand{\SuunitGG}[0]{\textcolor{\complexcolor}{{\symApparentPower}_{ \indexUnit, \symPhaseG\symPhaseG}^{}  }}

\newcommand{\SuunitNA}[0]{\textcolor{\complexcolor}{{\symApparentPower}_{ \indexUnit, \symPhaseN\symPhaseA}^{}  }}
\newcommand{\SuunitNB}[0]{\textcolor{\complexcolor}{{\symApparentPower}_{ \indexUnit, \symPhaseN\symPhaseB}^{}  }}
\newcommand{\SuunitNC}[0]{\textcolor{\complexcolor}{{\symApparentPower}_{ \indexUnit, \symPhaseN\symPhaseC}^{}  }}

\newcommand{\PuunitNN}[0]{\textcolor{black}{{\symPower}_{ \indexUnit, \symPhaseN\symPhaseN}^{}  }}
\newcommand{\PuunitNA}[0]{\textcolor{black}{{\symPower}_{ \indexUnit, \symPhaseN\symPhaseA}^{}  }}
\newcommand{\PuunitNB}[0]{\textcolor{black}{{\symPower}_{ \indexUnit, \symPhaseN\symPhaseB}^{}  }}
\newcommand{\PuunitNC}[0]{\textcolor{black}{{\symPower}_{ \indexUnit, \symPhaseN\symPhaseC}^{}  }}

\newcommand{\PuunitAA}[0]{\textcolor{black}{{\symPower}_{ \indexUnit, \symPhaseA\symPhaseA}^{}  }}
\newcommand{\PuunitAB}[0]{\textcolor{black}{{\symPower}_{ \indexUnit, \symPhaseA\symPhaseB}^{}  }}
\newcommand{\PuunitAC}[0]{\textcolor{black}{{\symPower}_{ \indexUnit, \symPhaseA\symPhaseC}^{}  }}

\newcommand{\PuunitBA}[0]{\textcolor{black}{{\symPower}_{ \indexUnit, \symPhaseB\symPhaseA}^{}  }}
\newcommand{\PuunitBB}[0]{\textcolor{black}{{\symPower}_{ \indexUnit, \symPhaseB\symPhaseB}^{}  }}
\newcommand{\PuunitBC}[0]{\textcolor{black}{{\symPower}_{ \indexUnit, \symPhaseB\symPhaseC}^{}  }}

\newcommand{\PuunitCA}[0]{\textcolor{black}{{\symPower}_{ \indexUnit, \symPhaseC\symPhaseA}^{}  }}
\newcommand{\PuunitCB}[0]{\textcolor{black}{{\symPower}_{ \indexUnit, \symPhaseC\symPhaseB}^{}  }}
\newcommand{\PuunitCC}[0]{\textcolor{black}{{\symPower}_{ \indexUnit, \symPhaseC\symPhaseC}^{}  }}

\newcommand{\QuunitNN}[0]{\textcolor{black}{{\symReactivePower}_{ \indexUnit, \symPhaseN\symPhaseN}^{}  }}
\newcommand{\QuunitNA}[0]{\textcolor{black}{{\symReactivePower}_{ \indexUnit, \symPhaseN\symPhaseA}^{}  }}
\newcommand{\QuunitNB}[0]{\textcolor{black}{{\symReactivePower}_{ \indexUnit, \symPhaseN\symPhaseB}^{}  }}
\newcommand{\QuunitNC}[0]{\textcolor{black}{{\symReactivePower}_{ \indexUnit, \symPhaseN\symPhaseC}^{}  }}

\newcommand{\QuunitAA}[0]{\textcolor{black}{{\symReactivePower}_{ \indexUnit, \symPhaseA\symPhaseA}^{}  }}
\newcommand{\QuunitAB}[0]{\textcolor{black}{{\symReactivePower}_{ \indexUnit, \symPhaseA\symPhaseB}^{}  }}
\newcommand{\QuunitAC}[0]{\textcolor{black}{{\symReactivePower}_{ \indexUnit, \symPhaseA\symPhaseC}^{}  }}

\newcommand{\QuunitBA}[0]{\textcolor{black}{{\symReactivePower}_{ \indexUnit, \symPhaseB\symPhaseA}^{}  }}
\newcommand{\QuunitBB}[0]{\textcolor{black}{{\symReactivePower}_{ \indexUnit, \symPhaseB\symPhaseB}^{}  }}
\newcommand{\QuunitBC}[0]{\textcolor{black}{{\symReactivePower}_{ \indexUnit, \symPhaseB\symPhaseC}^{}  }}

\newcommand{\QuunitCA}[0]{\textcolor{black}{{\symReactivePower}_{ \indexUnit, \symPhaseC\symPhaseA}^{}  }}
\newcommand{\QuunitCB}[0]{\textcolor{black}{{\symReactivePower}_{ \indexUnit, \symPhaseC\symPhaseB}^{}  }}
\newcommand{\QuunitCC}[0]{\textcolor{black}{{\symReactivePower}_{ \indexUnit, \symPhaseC\symPhaseC}^{}  }}

\newcommand{\QuunitAN}[0]{\textcolor{black}{{\symReactivePower}_{ \indexUnit, \symPhaseA}^{}  }}
\newcommand{\QuunitBN}[0]{\textcolor{black}{{\symReactivePower}_{ \indexUnit, \symPhaseB}^{}  }}
\newcommand{\QuunitCN}[0]{\textcolor{black}{{\symReactivePower}_{ \indexUnit, \symPhaseC}^{}  }}

\newcommand{\SuunitrefSDP}[0]{\textcolor{\complexparamcolor}{\mathbf{\symApparentPower}_{ \indexUnit}^{\refss}  }}
\newcommand{\PuunitrefSDP}[0]{\textcolor{\paramcolor}{\mathbf{\symPower}_{ \indexUnit}^{\refss}  }}
\newcommand{\QuunitrefSDP}[0]{\textcolor{\paramcolor}{\mathbf{\symReactivePower}_{ \indexUnit}^{\refss}  }}

\newcommand{\PuunitrefAN}[0]{\textcolor{\paramcolor}{{\symPower}_{ \indexUnit, \symPhaseA}^{\refss}  }}
\newcommand{\PuunitrefBN}[0]{\textcolor{\paramcolor}{{\symPower}_{ \indexUnit, \symPhaseB}^{\refss}  }}
\newcommand{\PuunitrefCN}[0]{\textcolor{\paramcolor}{{\symPower}_{ \indexUnit, \symPhaseC}^{\refss}  }}

\newcommand{\QuunitrefAN}[0]{\textcolor{\paramcolor}{{\symReactivePower}_{ \indexUnit, \symPhaseA}^{\refss}  }}
\newcommand{\QuunitrefBN}[0]{\textcolor{\paramcolor}{{\symReactivePower}_{ \indexUnit, \symPhaseB}^{\refss}  }}
\newcommand{\QuunitrefCN}[0]{\textcolor{\paramcolor}{{\symReactivePower}_{ \indexUnit, \symPhaseC}^{\refss}  }}

\newcommand{\SbNN}[0]{\textcolor{\complexcolor}{{\symApparentPower}_{ \indexShunt, \symPhaseN\symPhaseN}^{}  }}

\newcommand{\PijkrefAA}[0]{\textcolor{\paramcolor}{\symPower_{\indexGridLines\indexGridNode \indexGridNodeTwo, \symPhaseA\symPhaseA}^{\refss}  }}
\newcommand{\PijkrefBB}[0]{\textcolor{\paramcolor}{\symPower_{\indexGridLines\indexGridNode \indexGridNodeTwo, \symPhaseB\symPhaseB}^{\refss}  }}
\newcommand{\PijkrefCC}[0]{\textcolor{\paramcolor}{\symPower_{\indexGridLines\indexGridNode \indexGridNodeTwo, \symPhaseC\symPhaseC}^{\refss}  }}

\newcommand{\Pijkpp}[0]{\textcolor{black}{\symPower_{\indexGridLines\indexGridNode \indexGridNodeTwo, \indexPhases\indexPhases}^{}  }}
\newcommand{\Qijkpp}[0]{\textcolor{black}{\symReactivePower_{\indexGridLines\indexGridNode \indexGridNodeTwo, \indexPhases\indexPhases}^{}  }}

\newcommand{\Pkpp}[0]{\textcolor{black}{\symPower_{\indexShunt, \indexPhases\indexPhases}  }}
\newcommand{\Qkpp}[0]{\textcolor{black}{\symReactivePower_{\indexShunt,  \indexPhases\indexPhases} }}

\newcommand{\PjikAA}[0]{\textcolor{black}{\symPower_{\indexGridLines\indexGridNodeTwo\indexGridNode , \symPhaseA\symPhaseA}^{}  }}
\newcommand{\PijkAAloss}[0]{\textcolor{black}{\symPower_{\indexGridLines\indexGridNode \indexGridNodeTwo, \symPhaseA\symPhaseA}^{\lossss}  }}
\newcommand{\PjikAAloss}[0]{\textcolor{black}{\symPower_{\indexGridLines\indexGridNodeTwo\indexGridNode , \symPhaseA\symPhaseA}^{\lossss}  }}

\newcommand{\SijkAA}[0]{\textcolor{\complexcolor}{\symApparentPower_{\indexGridLines\indexGridNode \indexGridNodeTwo, \symPhaseA\symPhaseA}^{}  }}
\newcommand{\SijkAB}[0]{\textcolor{\complexcolor}{\symApparentPower_{\indexGridLines\indexGridNode \indexGridNodeTwo, \symPhaseA\symPhaseB}^{}  }}
\newcommand{\SijkAC}[0]{\textcolor{\complexcolor}{\symApparentPower_{\indexGridLines\indexGridNode \indexGridNodeTwo, \symPhaseA\symPhaseC}^{}  }}
\newcommand{\SijkAN}[0]{\textcolor{\complexcolor}{\symApparentPower_{\indexGridLines\indexGridNode \indexGridNodeTwo, \symPhaseA\symPhaseN}^{}  }}
\newcommand{\SijkAG}[0]{\textcolor{\complexcolor}{\symApparentPower_{\indexGridLines\indexGridNode \indexGridNodeTwo, \symPhaseA\symPhaseG}^{}  }}
\newcommand{\SijkBA}[0]{\textcolor{\complexcolor}{\symApparentPower_{\indexGridLines\indexGridNode \indexGridNodeTwo, \symPhaseB\symPhaseA}^{}  }}
\newcommand{\SijkBB}[0]{\textcolor{\complexcolor}{\symApparentPower_{\indexGridLines\indexGridNode \indexGridNodeTwo, \symPhaseB\symPhaseB}^{}  }}
\newcommand{\SijkBC}[0]{\textcolor{\complexcolor}{\symApparentPower_{\indexGridLines\indexGridNode \indexGridNodeTwo, \symPhaseB\symPhaseC}^{}  }}
\newcommand{\SijkBN}[0]{\textcolor{\complexcolor}{\symApparentPower_{\indexGridLines\indexGridNode \indexGridNodeTwo, \symPhaseB\symPhaseN}^{}  }}
\newcommand{\SijkBG}[0]{\textcolor{\complexcolor}{\symApparentPower_{\indexGridLines\indexGridNode \indexGridNodeTwo, \symPhaseB\symPhaseG}^{}  }}
\newcommand{\SijkCA}[0]{\textcolor{\complexcolor}{\symApparentPower_{\indexGridLines\indexGridNode \indexGridNodeTwo, \symPhaseC\symPhaseA}^{}  }}
\newcommand{\SijkCB}[0]{\textcolor{\complexcolor}{\symApparentPower_{\indexGridLines\indexGridNode \indexGridNodeTwo, \symPhaseC\symPhaseB}^{}  }}
\newcommand{\SijkCC}[0]{\textcolor{\complexcolor}{\symApparentPower_{\indexGridLines\indexGridNode \indexGridNodeTwo, \symPhaseC\symPhaseC}^{}  }}
\newcommand{\SijkCN}[0]{\textcolor{\complexcolor}{\symApparentPower_{\indexGridLines\indexGridNode \indexGridNodeTwo, \symPhaseC\symPhaseN}^{}  }}
\newcommand{\SijkCG}[0]{\textcolor{\complexcolor}{\symApparentPower_{\indexGridLines\indexGridNode \indexGridNodeTwo, \symPhaseC\symPhaseG}^{}  }}

\newcommand{\SijkNA}[0]{\textcolor{\complexcolor}{\symApparentPower_{\indexGridLines\indexGridNode \indexGridNodeTwo, \symPhaseN\symPhaseA}^{}  }}
\newcommand{\SijkNB}[0]{\textcolor{\complexcolor}{\symApparentPower_{\indexGridLines\indexGridNode \indexGridNodeTwo, \symPhaseN\symPhaseB}^{}  }}
\newcommand{\SijkNC}[0]{\textcolor{\complexcolor}{\symApparentPower_{\indexGridLines\indexGridNode \indexGridNodeTwo, \symPhaseN\symPhaseC}^{}  }}
\newcommand{\SijkNN}[0]{\textcolor{\complexcolor}{\symApparentPower_{\indexGridLines\indexGridNode \indexGridNodeTwo, \symPhaseN\symPhaseN}^{}  }}
\newcommand{\SijkNG}[0]{\textcolor{\complexcolor}{\symApparentPower_{\indexGridLines\indexGridNode \indexGridNodeTwo, \symPhaseN\symPhaseG}^{}  }}
\newcommand{\PijkNG}[0]{\textcolor{black}{\symPower_{\indexGridLines\indexGridNode \indexGridNodeTwo, \symPhaseN\symPhaseG}^{}  }}
\newcommand{\QijkNG}[0]{\textcolor{black}{\symReactivePower_{\indexGridLines\indexGridNode \indexGridNodeTwo, \symPhaseN\symPhaseG}^{}  }}

\newcommand{\SijkGA}[0]{\textcolor{\complexcolor}{\symApparentPower_{\indexGridLines\indexGridNode \indexGridNodeTwo, \symPhaseG\symPhaseA}^{}  }}
\newcommand{\SijkGB}[0]{\textcolor{\complexcolor}{\symApparentPower_{\indexGridLines\indexGridNode \indexGridNodeTwo, \symPhaseG\symPhaseB}^{}  }}
\newcommand{\SijkGC}[0]{\textcolor{\complexcolor}{\symApparentPower_{\indexGridLines\indexGridNode \indexGridNodeTwo, \symPhaseG\symPhaseC}^{}  }}
\newcommand{\SijkGN}[0]{\textcolor{\complexcolor}{\symApparentPower_{\indexGridLines\indexGridNode \indexGridNodeTwo, \symPhaseG\symPhaseN}^{}  }}
\newcommand{\SijkGG}[0]{\textcolor{\complexcolor}{\symApparentPower_{\indexGridLines\indexGridNode \indexGridNodeTwo, \symPhaseG\symPhaseG}^{}  }}

\newcommand{\SjikAA}[0]{\textcolor{\complexcolor}{\symApparentPower_{\indexGridLines\indexGridNodeTwo\indexGridNode , \symPhaseA\symPhaseA}^{}  }}
\newcommand{\SjikAB}[0]{\textcolor{\complexcolor}{\symApparentPower_{\indexGridLines\indexGridNodeTwo\indexGridNode , \symPhaseA\symPhaseB}^{}  }}
\newcommand{\SjikAC}[0]{\textcolor{\complexcolor}{\symApparentPower_{\indexGridLines\indexGridNodeTwo\indexGridNode , \symPhaseA\symPhaseC}^{}  }}
\newcommand{\SjikAN}[0]{\textcolor{\complexcolor}{\symApparentPower_{\indexGridLines\indexGridNodeTwo\indexGridNode , \symPhaseA\symPhaseN}^{}  }}
\newcommand{\SjikBB}[0]{\textcolor{\complexcolor}{\symApparentPower_{\indexGridLines\indexGridNodeTwo\indexGridNode , \symPhaseB\symPhaseB}^{}  }}
\newcommand{\SjikCC}[0]{\textcolor{\complexcolor}{\symApparentPower_{\indexGridLines\indexGridNodeTwo\indexGridNode , \symPhaseC\symPhaseC}^{}  }}
\newcommand{\SjikNA}[0]{\textcolor{\complexcolor}{\symApparentPower_{\indexGridLines\indexGridNodeTwo\indexGridNode , \symPhaseN\symPhaseA}^{}  }}
\newcommand{\SjikNB}[0]{\textcolor{\complexcolor}{\symApparentPower_{\indexGridLines\indexGridNodeTwo\indexGridNode , \symPhaseN\symPhaseB}^{}  }}
\newcommand{\SjikNC}[0]{\textcolor{\complexcolor}{\symApparentPower_{\indexGridLines\indexGridNodeTwo\indexGridNode , \symPhaseN\symPhaseC}^{}  }}
\newcommand{\SjikNN}[0]{\textcolor{\complexcolor}{\symApparentPower_{\indexGridLines\indexGridNodeTwo\indexGridNode , \symPhaseN\symPhaseN}^{}  }}

\newcommand{\SikNG}[0]{\textcolor{\complexcolor}{\symApparentPower_{\indexGridNode, \symPhaseN\symPhaseG}^{}  }}
\newcommand{\PikNG}[0]{\textcolor{black}{\symPower_{\indexGridNode, \symPhaseN\symPhaseG}^{}  }}
\newcommand{\QikNG}[0]{\textcolor{black}{\symReactivePower_{\indexGridNode, \symPhaseN\symPhaseG}^{}  }}

\newcommand{\PijkAA}[0]{\textcolor{black}{\symPower_{\indexGridLines\indexGridNode \indexGridNodeTwo, \symPhaseA\symPhaseA}^{}  }}
\newcommand{\PijkAB}[0]{\textcolor{black}{\symPower_{\indexGridLines\indexGridNode \indexGridNodeTwo, \symPhaseA\symPhaseB}^{}  }}
\newcommand{\PijkAC}[0]{\textcolor{black}{\symPower_{\indexGridLines\indexGridNode \indexGridNodeTwo, \symPhaseA\symPhaseC}^{}  }}
\newcommand{\PijkBA}[0]{\textcolor{black}{\symPower_{\indexGridLines\indexGridNode \indexGridNodeTwo, \symPhaseB\symPhaseA}^{}  }}
\newcommand{\PijkBB}[0]{\textcolor{black}{\symPower_{\indexGridLines\indexGridNode \indexGridNodeTwo, \symPhaseB\symPhaseB}^{}  }}
\newcommand{\PijkBC}[0]{\textcolor{black}{\symPower_{\indexGridLines\indexGridNode \indexGridNodeTwo, \symPhaseB\symPhaseC}^{}  }}
\newcommand{\PijkCA}[0]{\textcolor{black}{\symPower_{\indexGridLines\indexGridNode \indexGridNodeTwo, \symPhaseC\symPhaseA}^{}  }}
\newcommand{\PijkCB}[0]{\textcolor{black}{\symPower_{\indexGridLines\indexGridNode \indexGridNodeTwo, \symPhaseC\symPhaseB}^{}  }}
\newcommand{\PijkCC}[0]{\textcolor{black}{\symPower_{\indexGridLines\indexGridNode \indexGridNodeTwo, \symPhaseC\symPhaseC}^{}  }}
\newcommand{\PijkNN}[0]{\textcolor{black}{\symPower_{\indexGridLines\indexGridNode \indexGridNodeTwo, \symPhaseN\symPhaseN}^{}  }}
\newcommand{\QijkAA}[0]{\textcolor{black}{\symReactivePower_{\indexGridLines\indexGridNode \indexGridNodeTwo, \symPhaseA\symPhaseA}^{}  }}
\newcommand{\QijkAB}[0]{\textcolor{black}{\symReactivePower_{\indexGridLines\indexGridNode \indexGridNodeTwo, \symPhaseA\symPhaseB}^{}  }}
\newcommand{\QijkAC}[0]{\textcolor{black}{\symReactivePower_{\indexGridLines\indexGridNode \indexGridNodeTwo, \symPhaseA\symPhaseC}^{}  }}
\newcommand{\QijkBA}[0]{\textcolor{black}{\symReactivePower_{\indexGridLines\indexGridNode \indexGridNodeTwo, \symPhaseB\symPhaseA}^{}  }}
\newcommand{\QijkBB}[0]{\textcolor{black}{\symReactivePower_{\indexGridLines\indexGridNode \indexGridNodeTwo, \symPhaseB\symPhaseB}^{}  }}
\newcommand{\QijkBC}[0]{\textcolor{black}{\symReactivePower_{\indexGridLines\indexGridNode \indexGridNodeTwo, \symPhaseB\symPhaseC}^{}  }}
\newcommand{\QijkCA}[0]{\textcolor{black}{\symReactivePower_{\indexGridLines\indexGridNode \indexGridNodeTwo, \symPhaseC\symPhaseA}^{}  }}
\newcommand{\QijkCB}[0]{\textcolor{black}{\symReactivePower_{\indexGridLines\indexGridNode \indexGridNodeTwo, \symPhaseC\symPhaseB}^{}  }}
\newcommand{\QijkCC}[0]{\textcolor{black}{\symReactivePower_{\indexGridLines\indexGridNode \indexGridNodeTwo, \symPhaseC\symPhaseC}^{}  }}
\newcommand{\QijkNN}[0]{\textcolor{black}{\symReactivePower_{\indexGridLines\indexGridNode \indexGridNodeTwo, \symPhaseN\symPhaseN}^{}  }}

\newcommand{\PjikAB}[0]{\textcolor{black}{\symPower_{\indexGridLines \indexGridNodeTwo \indexGridNode, \symPhaseA\symPhaseB}^{}  }}
\newcommand{\QjikAB}[0]{\textcolor{black}{\symReactivePower_{\indexGridLines \indexGridNodeTwo \indexGridNode, \symPhaseA\symPhaseB}^{}  }}

\newcommand{\PijksAA}[0]{\textcolor{black}{\symPower_{\indexGridLines\indexGridNode \indexGridNodeTwo, \symPhaseA\symPhaseA }^{\seriesss}  }}
\newcommand{\QijksAA}[0]{\textcolor{black}{\symReactivePower_{\indexGridLines\indexGridNode \indexGridNodeTwo, \symPhaseA\symPhaseA }^{\seriesss}  }}
\newcommand{\SijksAA}[0]{\textcolor{\complexcolor}{\symApparentPower_{\indexGridLines\indexGridNode \indexGridNodeTwo, \symPhaseA\symPhaseA }^{\seriesss}  }}
\newcommand{\PijksAB}[0]{\textcolor{black}{\symPower_{\indexGridLines\indexGridNode \indexGridNodeTwo, \symPhaseA\symPhaseB }^{\seriesss}  }}
\newcommand{\QijksAB}[0]{\textcolor{black}{\symReactivePower_{\indexGridLines\indexGridNode \indexGridNodeTwo, \symPhaseA\symPhaseB }^{\seriesss}  }}
\newcommand{\SijksAB}[0]{\textcolor{\complexcolor}{\symApparentPower_{\indexGridLines\indexGridNode \indexGridNodeTwo, \symPhaseA\symPhaseB }^{\seriesss}  }}
\newcommand{\PijksAC}[0]{\textcolor{black}{\symPower_{\indexGridLines\indexGridNode \indexGridNodeTwo, \symPhaseA\symPhaseC }^{\seriesss}  }}
\newcommand{\QijksAC}[0]{\textcolor{black}{\symReactivePower_{\indexGridLines\indexGridNode \indexGridNodeTwo, \symPhaseA\symPhaseC }^{\seriesss}  }}
\newcommand{\SijksAC}[0]{\textcolor{\complexcolor}{\symApparentPower_{\indexGridLines\indexGridNode \indexGridNodeTwo, \symPhaseA\symPhaseC }^{\seriesss}  }}
\newcommand{\PijksBA}[0]{\textcolor{black}{\symPower_{\indexGridLines\indexGridNode \indexGridNodeTwo, \symPhaseB\symPhaseA}^{\seriesss}  }}
\newcommand{\QijksBA}[0]{\textcolor{black}{\symReactivePower_{\indexGridLines\indexGridNode \indexGridNodeTwo, \symPhaseB\symPhaseA}^{\seriesss}  }}
\newcommand{\SijksBA}[0]{\textcolor{\complexcolor}{\symApparentPower_{\indexGridLines\indexGridNode \indexGridNodeTwo, \symPhaseB\symPhaseA}^{\seriesss}  }}
\newcommand{\PijksBB}[0]{\textcolor{black}{\symPower_{\indexGridLines\indexGridNode \indexGridNodeTwo, \symPhaseB\symPhaseB}^{\seriesss}  }}
\newcommand{\QijksBB}[0]{\textcolor{black}{\symReactivePower_{\indexGridLines\indexGridNode \indexGridNodeTwo, \symPhaseB\symPhaseB}^{\seriesss}  }}
\newcommand{\SijksBB}[0]{\textcolor{\complexcolor}{\symApparentPower_{\indexGridLines\indexGridNode \indexGridNodeTwo, \symPhaseB\symPhaseB}^{\seriesss}  }}
\newcommand{\PijksBC}[0]{\textcolor{black}{\symPower_{\indexGridLines\indexGridNode \indexGridNodeTwo, \symPhaseB\symPhaseC}^{\seriesss}  }}
\newcommand{\QijksBC}[0]{\textcolor{black}{\symReactivePower_{\indexGridLines\indexGridNode \indexGridNodeTwo, \symPhaseB\symPhaseC}^{\seriesss}  }}
\newcommand{\SijksBC}[0]{\textcolor{\complexcolor}{\symApparentPower_{\indexGridLines\indexGridNode \indexGridNodeTwo, \symPhaseB\symPhaseC}^{\seriesss}  }}
\newcommand{\PijksCA}[0]{\textcolor{black}{\symPower_{\indexGridLines\indexGridNode \indexGridNodeTwo, \symPhaseC\symPhaseA}^{\seriesss}  }}
\newcommand{\QijksCA}[0]{\textcolor{black}{\symReactivePower_{\indexGridLines\indexGridNode \indexGridNodeTwo, \symPhaseC\symPhaseA}^{\seriesss}  }}
\newcommand{\SijksCA}[0]{\textcolor{\complexcolor}{\symApparentPower_{\indexGridLines\indexGridNode \indexGridNodeTwo, \symPhaseC\symPhaseA}^{\seriesss}  }}
\newcommand{\PijksCB}[0]{\textcolor{black}{\symPower_{\indexGridLines\indexGridNode \indexGridNodeTwo, \symPhaseC\symPhaseB}^{\seriesss}  }}
\newcommand{\QijksCB}[0]{\textcolor{black}{\symReactivePower_{\indexGridLines\indexGridNode \indexGridNodeTwo, \symPhaseC\symPhaseB}^{\seriesss}  }}
\newcommand{\SijksCB}[0]{\textcolor{\complexcolor}{\symApparentPower_{\indexGridLines\indexGridNode \indexGridNodeTwo, \symPhaseC\symPhaseB}^{\seriesss}  }}
\newcommand{\PijksCC}[0]{\textcolor{black}{\symPower_{\indexGridLines\indexGridNode \indexGridNodeTwo, \symPhaseC\symPhaseC}^{\seriesss}  }}
\newcommand{\QijksCC}[0]{\textcolor{black}{\symReactivePower_{\indexGridLines\indexGridNode \indexGridNodeTwo, \symPhaseC\symPhaseC}^{\seriesss}  }}
\newcommand{\SijksCC}[0]{\textcolor{\complexcolor}{\symApparentPower_{\indexGridLines\indexGridNode \indexGridNodeTwo, \symPhaseC\symPhaseC}^{\seriesss}  }}

\newcommand{\PijksAN}[0]{\textcolor{black}{\symPower_{\indexGridLines\indexGridNode \indexGridNodeTwo, \symPhaseA\symPhaseN }^{\seriesss}  }}
\newcommand{\QijksAN}[0]{\textcolor{black}{\symReactivePower_{\indexGridLines\indexGridNode \indexGridNodeTwo, \symPhaseA\symPhaseN }^{\seriesss}  }}
\newcommand{\SijksAN}[0]{\textcolor{\complexcolor}{\symApparentPower_{\indexGridLines\indexGridNode \indexGridNodeTwo, \symPhaseA\symPhaseN }^{\seriesss}  }}
\newcommand{\PijksBN}[0]{\textcolor{black}{\symPower_{\indexGridLines\indexGridNode \indexGridNodeTwo, \symPhaseB\symPhaseN }^{\seriesss}  }}
\newcommand{\QijksBN}[0]{\textcolor{black}{\symReactivePower_{\indexGridLines\indexGridNode \indexGridNodeTwo, \symPhaseB\symPhaseN }^{\seriesss}  }}
\newcommand{\SijksBN}[0]{\textcolor{\complexcolor}{\symApparentPower_{\indexGridLines\indexGridNode \indexGridNodeTwo, \symPhaseB\symPhaseN }^{\seriesss}  }}
\newcommand{\PijksCN}[0]{\textcolor{black}{\symPower_{\indexGridLines\indexGridNode \indexGridNodeTwo, \symPhaseC\symPhaseN }^{\seriesss}  }}
\newcommand{\QijksCN}[0]{\textcolor{black}{\symReactivePower_{\indexGridLines\indexGridNode \indexGridNodeTwo, \symPhaseC\symPhaseN }^{\seriesss}  }}
\newcommand{\SijksCN}[0]{\textcolor{\complexcolor}{\symApparentPower_{\indexGridLines\indexGridNode \indexGridNodeTwo, \symPhaseC\symPhaseN }^{\seriesss}  }}
\newcommand{\PijksNN}[0]{\textcolor{black}{\symPower_{\indexGridLines\indexGridNode \indexGridNodeTwo, \symPhaseN\symPhaseN }^{\seriesss}  }}
\newcommand{\QijksNN}[0]{\textcolor{black}{\symReactivePower_{\indexGridLines\indexGridNode \indexGridNodeTwo, \symPhaseN\symPhaseN }^{\seriesss}  }}
\newcommand{\SijksNN}[0]{\textcolor{\complexcolor}{\symApparentPower_{\indexGridLines\indexGridNode \indexGridNodeTwo, \symPhaseN\symPhaseN }^{\seriesss}  }}
\newcommand{\PijksNA}[0]{\textcolor{black}{\symPower_{\indexGridLines\indexGridNode \indexGridNodeTwo, \symPhaseN\symPhaseA }^{\seriesss}  }}
\newcommand{\QijksNA}[0]{\textcolor{black}{\symReactivePower_{\indexGridLines\indexGridNode \indexGridNodeTwo, \symPhaseN\symPhaseA }^{\seriesss}  }}
\newcommand{\SijksNA}[0]{\textcolor{\complexcolor}{\symApparentPower_{\indexGridLines\indexGridNode \indexGridNodeTwo, \symPhaseN\symPhaseA }^{\seriesss}  }}
\newcommand{\PijksNB}[0]{\textcolor{black}{\symPower_{\indexGridLines\indexGridNode \indexGridNodeTwo, \symPhaseN\symPhaseB }^{\seriesss}  }}
\newcommand{\QijksNB}[0]{\textcolor{black}{\symReactivePower_{\indexGridLines\indexGridNode \indexGridNodeTwo, \symPhaseN\symPhaseB }^{\seriesss}  }}
\newcommand{\SijksNB}[0]{\textcolor{\complexcolor}{\symApparentPower_{\indexGridLines\indexGridNode \indexGridNodeTwo, \symPhaseN\symPhaseB }^{\seriesss}  }}
\newcommand{\PijksNC}[0]{\textcolor{black}{\symPower_{\indexGridLines\indexGridNode \indexGridNodeTwo, \symPhaseN\symPhaseC }^{\seriesss}  }}
\newcommand{\QijksNC}[0]{\textcolor{black}{\symReactivePower_{\indexGridLines\indexGridNode \indexGridNodeTwo, \symPhaseN\symPhaseC }^{\seriesss}  }}
\newcommand{\SijksNC}[0]{\textcolor{\complexcolor}{\symApparentPower_{\indexGridLines\indexGridNode \indexGridNodeTwo, \symPhaseN\symPhaseC }^{\seriesss}  }}

\newcommand{\PijklossAA}[0]{\textcolor{black}{\symPower_{\indexGridNode \indexGridNodeTwo, \symPhaseA\symPhaseA}^{\lossss}  }}
\newcommand{\PijklossAB}[0]{\textcolor{black}{\symPower_{\indexGridNode \indexGridNodeTwo, \symPhaseA\symPhaseB}^{\lossss}  }}
\newcommand{\PijklossAC}[0]{\textcolor{black}{\symPower_{\indexGridNode \indexGridNodeTwo, \symPhaseA\symPhaseC}^{\lossss}  }}
\newcommand{\PijklossBA}[0]{\textcolor{black}{\symPower_{\indexGridNode \indexGridNodeTwo, \symPhaseB\symPhaseA}^{\lossss}  }}
\newcommand{\PijklossBB}[0]{\textcolor{black}{\symPower_{\indexGridNode \indexGridNodeTwo, \symPhaseB\symPhaseB}^{\lossss}  }}
\newcommand{\PijklossBC}[0]{\textcolor{black}{\symPower_{\indexGridNode \indexGridNodeTwo, \symPhaseB\symPhaseC}^{\lossss}  }}
\newcommand{\PijklossCA}[0]{\textcolor{black}{\symPower_{\indexGridNode \indexGridNodeTwo, \symPhaseC\symPhaseA}^{\lossss}  }}
\newcommand{\PijklossCB}[0]{\textcolor{black}{\symPower_{\indexGridNode \indexGridNodeTwo, \symPhaseC\symPhaseB}^{\lossss}  }}
\newcommand{\PijklossCC}[0]{\textcolor{black}{\symPower_{\indexGridNode \indexGridNodeTwo, \symPhaseC\symPhaseC}^{\lossss}  }}
\newcommand{\QijklossAA}[0]{\textcolor{black}{\symReactivePower_{\indexGridNode \indexGridNodeTwo, \symPhaseA\symPhaseA}^{\lossss}  }}
\newcommand{\QijklossAB}[0]{\textcolor{black}{\symReactivePower_{\indexGridNode \indexGridNodeTwo, \symPhaseA\symPhaseB}^{\lossss}  }}
\newcommand{\QijklossAC}[0]{\textcolor{black}{\symReactivePower_{\indexGridNode \indexGridNodeTwo, \symPhaseA\symPhaseC}^{\lossss}  }}
\newcommand{\QijklossBA}[0]{\textcolor{black}{\symReactivePower_{\indexGridNode \indexGridNodeTwo, \symPhaseB\symPhaseA}^{\lossss}  }}
\newcommand{\QijklossBB}[0]{\textcolor{black}{\symReactivePower_{\indexGridNode \indexGridNodeTwo, \symPhaseB\symPhaseB}^{\lossss}  }}
\newcommand{\QijklossBC}[0]{\textcolor{black}{\symReactivePower_{\indexGridNode \indexGridNodeTwo, \symPhaseB\symPhaseC}^{\lossss}  }}
\newcommand{\QijklossCA}[0]{\textcolor{black}{\symReactivePower_{\indexGridNode \indexGridNodeTwo, \symPhaseC\symPhaseA}^{\lossss}  }}
\newcommand{\QijklossCB}[0]{\textcolor{black}{\symReactivePower_{\indexGridNode \indexGridNodeTwo, \symPhaseC\symPhaseB}^{\lossss}  }}
\newcommand{\QijklossCC}[0]{\textcolor{black}{\symReactivePower_{\indexGridNode \indexGridNodeTwo, \symPhaseC\symPhaseC}^{\lossss}  }}


\newcommand{\Iik}[0]{\textcolor{\complexcolor}{\symCurrent_{\indexGridNode }^{}  }}

\newcommand{\Ijik}[0]{\textcolor{\complexcolor}{\symCurrent_{\indexGridLines\indexGridNodeTwo  \indexGridNode }^{}  }}

\newcommand{\Iijkre}[0]{\textcolor{black}{\symCurrent_{\indexGridLines\indexGridNode \indexGridNodeTwo}^{\realss}  }}
\newcommand{\Iijkim}[0]{\textcolor{black}{\symCurrent_{\indexGridLines\indexGridNode \indexGridNodeTwo}^{\imagss}  }}
\newcommand{\Iijkangle}[0]{\textcolor{black}{\symCurrent_{\indexGridLines\indexGridNode \indexGridNodeTwo}^{\imagangle}  }}
\newcommand{\Iijkanglemin}[0]{\textcolor{\paramcolor}{\symCurrent_{\indexGridLines\indexGridNode \indexGridNodeTwo}^{\imagangle\minss}  }}
\newcommand{\Iijkanglemax}[0]{\textcolor{\paramcolor}{\symCurrent_{\indexGridLines\indexGridNode \indexGridNodeTwo}^{\imagangle\maxss}  }}

\newcommand{\IijkSDPseq}[0]{\textcolor{\complexcolor}{\mathbf{\symCurrent}_{\indexGridLines\indexGridNode \indexGridNodeTwo   }^{\fortescuess}  }}

\newcommand{\IijkSDPreal}[0]{\textcolor{black}{\mathbf{\symCurrent}_{\indexGridLines\indexGridNode \indexGridNodeTwo   }^{\realss}  }}
\newcommand{\IijkSDPimag}[0]{\textcolor{black}{\mathbf{\symCurrent}_{\indexGridLines\indexGridNode \indexGridNodeTwo   }^{\imagss}  }}

\newcommand{\IijkSDP}[0]{\textcolor{\complexcolor}{\mathbf{\symCurrent}_{\indexGridLines\indexGridNode \indexGridNodeTwo   }^{}  }}
\newcommand{\IijkSDPH}[0]{\textcolor{\complexcolor}{(\mathbf{\symCurrent}_{\indexGridLines\indexGridNode \indexGridNodeTwo   })^{\hermitiantranspose}  }}

\newcommand{\IijkSDPref}[0]{\textcolor{\complexparamcolor}{\mathbf{\symCurrent}_{\indexGridLines\indexGridNode \indexGridNodeTwo   }^{\refss}  }}
\newcommand{\IijkSDPrefre}[0]{\textcolor{\paramcolor}{\mathbf{\symCurrent}_{\indexGridLines\indexGridNode \indexGridNodeTwo   }^{\refss, \realss}  }}
\newcommand{\IijkSDPrefim}[0]{\textcolor{\paramcolor}{\mathbf{\symCurrent}_{\indexGridLines\indexGridNode \indexGridNodeTwo   }^{\refss, \imagss}  }}

\newcommand{\IijksSDPref}[0]{\textcolor{\complexparamcolor}{\mathbf{\symCurrent}_{\indexGridLines\indexGridNode \indexGridNodeTwo   }^{\seriesss, \refss}  }}
\newcommand{\IijksSDPrefre}[0]{\textcolor{\paramcolor}{\mathbf{\symCurrent}_{\indexGridLines\indexGridNode \indexGridNodeTwo   }^{\seriesss, \refss, \realss}  }}
\newcommand{\IijksSDPrefim}[0]{\textcolor{\paramcolor}{\mathbf{\symCurrent}_{\indexGridLines\indexGridNode \indexGridNodeTwo   }^{\seriesss, \refss, \imagss}  }}

\newcommand{\IijkshSDPref}[0]{\textcolor{\complexparamcolor}{\mathbf{\symCurrent}_{\indexGridLines\indexGridNode \indexGridNodeTwo   }^{\shuntss, \refss}  }}
\newcommand{\IijkshSDPrefre}[0]{\textcolor{\paramcolor}{\mathbf{\symCurrent}_{\indexGridLines\indexGridNode \indexGridNodeTwo   }^{\shuntss, \refss, \realss}  }}
\newcommand{\IijkshSDPrefim}[0]{\textcolor{\paramcolor}{\mathbf{\symCurrent}_{\indexGridLines\indexGridNode \indexGridNodeTwo   }^{\shuntss, \refss, \imagss}  }}

\newcommand{\IijkSDPdelta}[0]{\textcolor{\complexcolor}{\mathbf{\symCurrent}_{\indexGridLines\indexGridNode \indexGridNodeTwo   }^{\Delta}  }}
\newcommand{\IijkSDPdeltare}[0]{\textcolor{black}{\mathbf{\symCurrent}_{\indexGridLines\indexGridNode \indexGridNodeTwo   }^{\Delta, \realss}  }}
\newcommand{\IijkSDPdeltaim}[0]{\textcolor{black}{\mathbf{\symCurrent}_{\indexGridLines\indexGridNode \indexGridNodeTwo   }^{\Delta, \imagss}  }}

\newcommand{\IijksSDPdelta}[0]{\textcolor{\complexcolor}{\mathbf{\symCurrent}_{\indexGridLines\indexGridNode \indexGridNodeTwo   }^{\seriesss, \Delta}  }}
\newcommand{\IijksSDPdeltare}[0]{\textcolor{black}{\mathbf{\symCurrent}_{\indexGridLines\indexGridNode \indexGridNodeTwo   }^{\seriesss, \Delta, \realss}  }}
\newcommand{\IijksSDPdeltaim}[0]{\textcolor{black}{\mathbf{\symCurrent}_{\indexGridLines\indexGridNode \indexGridNodeTwo   }^{\seriesss, \Delta, \imagss}  }}

\newcommand{\IsqijkSDPdelta}[0]{\textcolor{\complexcolor}{\mathbf{\symCurrentSquared}_{\indexGridLines\indexGridNode \indexGridNodeTwo   }^{\Delta}  }}
\newcommand{\IsqijkSDPdeltare}[0]{\textcolor{black}{\mathbf{\symCurrentSquared}_{\indexGridLines\indexGridNode \indexGridNodeTwo   }^{\Delta, \realss}  }}
\newcommand{\IsqijkSDPdeltaim}[0]{\textcolor{black}{\mathbf{\symCurrentSquared}_{\indexGridLines\indexGridNode \indexGridNodeTwo   }^{\Delta, \imagss}  }}

\newcommand{\IsqijksSDPdelta}[0]{\textcolor{\complexcolor}{\mathbf{\symCurrentSquared}_{\indexGridLines\indexGridNode \indexGridNodeTwo   }^{\seriesss, \Delta}  }}
\newcommand{\IsqijksSDPdeltare}[0]{\textcolor{black}{\mathbf{\symCurrentSquared}_{\indexGridLines\indexGridNode \indexGridNodeTwo   }^{\seriesss, \Delta, \realss}  }}
\newcommand{\IsqijksSDPdeltaim}[0]{\textcolor{black}{\mathbf{\symCurrentSquared}_{\indexGridLines\indexGridNode \indexGridNodeTwo   }^{\seriesss, \Delta, \imagss}  }}

\newcommand{\IjikSDP}[0]{\textcolor{\complexcolor}{\mathbf{\symCurrent}_{\indexGridLines\indexGridNodeTwo\indexGridNode    }^{}  }}
\newcommand{\IjikSDPreal}[0]{\textcolor{black}{\mathbf{\symCurrent}_{\indexGridLines\indexGridNodeTwo\indexGridNode    }^{\realss}  }}
\newcommand{\IjikSDPimag}[0]{\textcolor{black}{\mathbf{\symCurrent}_{\indexGridLines\indexGridNodeTwo\indexGridNode    }^{\imagss}  }}

\newcommand{\IijkrefSDP}[0]{\textcolor{\complexparamcolor}{\mathbf{\symCurrent}_{\indexGridLines\indexGridNode \indexGridNodeTwo   }^{\refss}  }}

\newcommand{\IijksrefSDP}[0]{\textcolor{\complexparamcolor}{\mathbf{\symCurrent}_{\indexGridLines\indexGridNode \indexGridNodeTwo, \seriesss   }^{\refss}  }}

\newcommand{\IijskSDPseq}[0]{\textcolor{\complexcolor}{\mathbf{\symCurrent}_{\indexGridLines\indexGridNode \indexGridNodeTwo   }^{\seriesss, \fortescuess}  }}

\newcommand{\IijskSDPH}[0]{\textcolor{\complexcolor}{(\mathbf{\symCurrent}_{\indexGridLines\indexGridNode \indexGridNodeTwo   }^{\seriesss})^{\hermitiantranspose}  }}
\newcommand{\IijskSDP}[0]{\textcolor{\complexcolor}{\mathbf{\symCurrent}_{\indexGridLines\indexGridNode \indexGridNodeTwo   }^{\seriesss}  }}
\newcommand{\IijskSDPreal}[0]{\textcolor{black}{\mathbf{\symCurrent}_{\indexGridLines\indexGridNode \indexGridNodeTwo   }^{\seriesss, \realss}  }}
\newcommand{\IijskSDPimag}[0]{\textcolor{black}{\mathbf{\symCurrent}_{\indexGridLines\indexGridNode \indexGridNodeTwo   }^{\seriesss, \imagss}  }}

\newcommand{\IjiskSDP}[0]{\textcolor{\complexcolor}{\mathbf{\symCurrent}_{\indexGridLines\indexGridNodeTwo\indexGridNode    }^{\seriesss}  }}
\newcommand{\IjiskSDPreal}[0]{\textcolor{black}{\mathbf{\symCurrent}_{\indexGridLines\indexGridNodeTwo\indexGridNode    }^{\seriesss, \realss}  }}
\newcommand{\IjiskSDPimag}[0]{\textcolor{black}{\mathbf{\symCurrent}_{\indexGridLines\indexGridNodeTwo\indexGridNode    }^{\seriesss, \imagss}  }}

\newcommand{\IijshkSDP}[0]{\textcolor{\complexcolor}{\mathbf{\symCurrent}_{\indexGridLines\indexGridNode \indexGridNodeTwo   }^{\shuntss}  }}
\newcommand{\IjishkSDP}[0]{\textcolor{\complexcolor}{\mathbf{\symCurrent}_{ \indexGridLines\indexGridNodeTwo \indexGridNode   }^{\shuntss}  }}

\newcommand{\IijshkSDPreal}[0]{\textcolor{black}{\mathbf{\symCurrent}_{\indexGridLines\indexGridNode \indexGridNodeTwo   }^{\shuntss, \realss}  }}

\newcommand{\IijshkSDPimag}[0]{\textcolor{black}{\mathbf{\symCurrent}_{\indexGridLines\indexGridNode \indexGridNodeTwo   }^{\shuntss, \imagss}  }}

\newcommand{\Iijk}[0]{\textcolor{\complexcolor}{\symCurrent_{\indexGridLines\indexGridNode \indexGridNodeTwo}^{}  }}
\newcommand{\Iijsk}[0]{\textcolor{\complexcolor}{\symCurrent_{\indexGridLines\indexGridNode \indexGridNodeTwo,\seriesss}^{}  }}
\newcommand{\Iijshk}[0]{\textcolor{\complexcolor}{\symCurrent_{\indexGridLines\indexGridNode \indexGridNodeTwo,\shuntss}^{}  }}
\newcommand{\Iijsqk}[0]{\textcolor{black}{\symCurrent_{\indexGridLines\indexGridNode \indexGridNodeTwo}^{2}  }}
\newcommand{\Iijsqsk}[0]{\textcolor{black}{\symCurrent_{\indexGridLines\indexGridLines,\seriesss}^{2}  }}
\newcommand{\Iijsqshk}[0]{\textcolor{black}{\symCurrent_{\indexGridLines\indexGridNode \indexGridNodeTwo,\shuntss}^{2}  }}

\newcommand{\IijkH}[0]{\textcolor{\complexcolor}{\symCurrent_{\indexGridLines\indexGridNode \indexGridNodeTwo}^{\hermitiantranspose}  }}
\newcommand{\Ijisk}[0]{\textcolor{\complexcolor}{\symCurrent_{\indexGridLines\indexGridNodeTwo\indexGridNode ,\seriesss}^{}  }}
\newcommand{\Ijishk}[0]{\textcolor{\complexcolor}{\symCurrent_{\indexGridLines\indexGridNodeTwo\indexGridNode ,\shuntss}^{}  }}

\newcommand{\Iijkmax}[0]{\textcolor{\boundscolor}{\symCurrent_{\indexGridLines\indexGridNode \indexGridNodeTwo}^{\maxss}  }}
\newcommand{\Iijkmin}[0]{\textcolor{\boundscolor}{\symCurrent_{\indexGridLines\indexGridNode \indexGridNodeTwo}^{\minss}  }}
\newcommand{\Iijrated}[0]{\textcolor{\sizingcolor}{\symCurrent_{\indexGridLines\indexGridNode \indexGridNodeTwo}^{\ratedss} }}
\newcommand{\Ijirated}[0]{\textcolor{\sizingcolor}{\symCurrent_{\indexGridLines\indexGridNodeTwo  \indexGridNode }^{\ratedss} }}
\newcommand{\IijratedA}[0]{\textcolor{\sizingcolor}{\symCurrent_{\indexGridLines\indexGridNode \indexGridNodeTwo,\symPhaseA}^{\ratedss} }}
\newcommand{\IjiratedA}[0]{\textcolor{\sizingcolor}{\symCurrent_{\indexGridLines\indexGridNodeTwo  \indexGridNode,\symPhaseA}^{\ratedss} }}
\newcommand{\IijratedB}[0]{\textcolor{\sizingcolor}{\symCurrent_{\indexGridLines\indexGridNode \indexGridNodeTwo,\symPhaseB}^{\ratedss} }}
\newcommand{\IjiratedB}[0]{\textcolor{\sizingcolor}{\symCurrent_{\indexGridLines\indexGridNodeTwo  \indexGridNode,\symPhaseB}^{\ratedss} }}
\newcommand{\IijratedC}[0]{\textcolor{\sizingcolor}{\symCurrent_{\indexGridLines\indexGridNode \indexGridNodeTwo,\symPhaseC}^{\ratedss} }}
\newcommand{\IjiratedC}[0]{\textcolor{\sizingcolor}{\symCurrent_{\indexGridLines\indexGridNodeTwo  \indexGridNode,\symPhaseC}^{\ratedss} }}

\newcommand{\IijratedN}[0]{\textcolor{\sizingcolor}{\symCurrent_{\indexGridLines\indexGridNode \indexGridNodeTwo,\symPhaseN}^{\ratedss} }}

\newcommand{\IijsqskSDPreal}[0]{\textcolor{black}{\mathbf{\symCurrentSOCP}_{\indexGridLines  }^{\seriesss,\realss}  }}
\newcommand{\IijsqskSDPimag}[0]{\textcolor{black}{\mathbf{\symCurrentSOCP}_{\indexGridLines  }^{\seriesss,\imagss}  }}

\newcommand{\IijsqkSDPreal}[0]{\textcolor{black}{\mathbf{\symCurrentSOCP}_{\indexGridLines\indexGridNode \indexGridNodeTwo }^{\realss}  }}
\newcommand{\IijsqkSDPimag}[0]{\textcolor{black}{\mathbf{\symCurrentSOCP}_{\indexGridLines\indexGridNode \indexGridNodeTwo }^{\imagss}  }}

\newcommand{\IijsqkSDPAdot}[0]{\textcolor{\complexcolor}{\mathbf{\symCurrentSOCP}_{\indexGridLines\indexGridNode \indexGridNodeTwo }^{\symPhaseA\!\cdot}  }}
\newcommand{\IijsqkSDPBdot}[0]{\textcolor{\complexcolor}{\mathbf{\symCurrentSOCP}_{\indexGridLines\indexGridNode \indexGridNodeTwo }^{\symPhaseB\!\cdot}  }}
\newcommand{\IijsqkSDPCdot}[0]{\textcolor{\complexcolor}{\mathbf{\symCurrentSOCP}_{\indexGridLines\indexGridNode \indexGridNodeTwo }^{\symPhaseC\!\cdot}  }}

\newcommand{\IjisqkSDP}[0]{\textcolor{\complexcolor}{\mathbf{\symCurrentSOCP}_{\indexGridLines \indexGridNodeTwo \indexGridNode  }^{}  }}
\newcommand{\IijsqkSDP}[0]{\textcolor{\complexcolor}{\mathbf{\symCurrentSOCP}_{\indexGridLines\indexGridNode \indexGridNodeTwo }^{}  }}
\newcommand{\IijsqkmaxSDP}[0]{\textcolor{\complexparamcolor}{\mathbf{\symCurrentSOCP}_{\indexGridLines\indexGridNode \indexGridNodeTwo }^{\maxss}  }}

\newcommand{\IijkmaxSDP}[0]{\textcolor{\paramcolor}{\mathbf{\symCurrent}_{\indexGridLines\indexGridNode \indexGridNodeTwo }^{\maxss}  }}

\newcommand{\IjisqkSDPreal}[0]{\textcolor{black}{\mathbf{\symCurrentSOCP}_{\indexGridLines \indexGridNodeTwo \indexGridNode  }^{\realss}  }}

\newcommand{\IijsqshkSDP}[0]{\textcolor{\complexcolor}{\mathbf{\symCurrentSOCP}_{\indexGridLines\indexGridNode \indexGridNodeTwo }^{\shuntss}  }}
\newcommand{\IjisqshkSDP}[0]{\textcolor{\complexcolor}{\mathbf{\symCurrentSOCP}_{\indexGridLines\indexGridNodeTwo \indexGridNode  }^{\shuntss}  }}

\newcommand{\IusqkSDPAdot}[0]{\textcolor{\complexcolor}{\mathbf{\symCurrentSOCP}_{\indexUnit  }^{\symPhaseA\!\cdot}  }}
\newcommand{\IusqkSDPBdot}[0]{\textcolor{\complexcolor}{\mathbf{\symCurrentSOCP}_{\indexUnit  }^{\symPhaseB\!\cdot}  }}
\newcommand{\IusqkSDPCdot}[0]{\textcolor{\complexcolor}{\mathbf{\symCurrentSOCP}_{\indexUnit  }^{\symPhaseC\!\cdot}  }}
\newcommand{\IusqkSDPNdot}[0]{\textcolor{\complexcolor}{\mathbf{\symCurrentSOCP}_{\indexUnit  }^{\symPhaseN\!\cdot}  }}

\newcommand{\IusqkSDPdotA}[0]{\textcolor{\complexcolor}{\mathbf{\symCurrentSOCP}_{\indexUnit  }^{\!\cdot\!\symPhaseA}  }}
\newcommand{\IusqkSDPdotB}[0]{\textcolor{\complexcolor}{\mathbf{\symCurrentSOCP}_{\indexUnit  }^{\!\cdot\!\symPhaseB}  }}
\newcommand{\IusqkSDPdotC}[0]{\textcolor{\complexcolor}{\mathbf{\symCurrentSOCP}_{\indexUnit  }^{\!\cdot\!\symPhaseC}  }}
\newcommand{\IusqkSDPdotN}[0]{\textcolor{\complexcolor}{\mathbf{\symCurrentSOCP}_{\indexUnit  }^{\!\cdot\!\symPhaseN}  }}

\newcommand{\IijsqskSDPseq}[0]{\textcolor{\complexcolor}{\mathbf{\symCurrentSOCP}_{\indexGridLines  }^{\seriesss, \fortescuess}  }}
\newcommand{\IusqkSDP}[0]{\textcolor{\complexcolor}{\mathbf{\symCurrentSOCP}_{\indexUnit  }  }}
\newcommand{\IusqkSDPreal}[0]{\textcolor{black}{\mathbf{\symCurrentSOCP}_{\indexUnit  }^{\realss}  }}
\newcommand{\IusqkSDPimag}[0]{\textcolor{black}{\mathbf{\symCurrentSOCP}_{\indexUnit  }^{\imagss}  }}

\newcommand{\IusqkSDPdelta}[0]{\textcolor{\complexcolor}{\mathbf{\symCurrentSOCP}_{\indexUnit  }^{\Delta}  }}
\newcommand{\IusqkSDPdeltareal}[0]{\textcolor{black}{\mathbf{\symCurrentSOCP}_{\indexUnit  }^{\Delta\realss}  }}
\newcommand{\IusqkSDPdeltaimag}[0]{\textcolor{black}{\mathbf{\symCurrentSOCP}_{\indexUnit  }^{\Delta\imagss}  }}

\newcommand{\IijsqskSDP}[0]{\textcolor{\complexcolor}{\mathbf{\symCurrentSOCP}_{\indexGridLines  }^{\seriesss}  }}
\newcommand{\IjisqskSDP}[0]{\textcolor{\complexcolor}{\mathbf{\symCurrentSOCP}_{\indexGridLines  }^{\seriesss}  }}
\newcommand{\Iijsqsksocp}[0]{\textcolor{black}{\symCurrentSOCP_{\indexGridLines\indexGridNode \indexGridNodeTwo ,\seriesss}  }}

\newcommand{\Iijsqksocp}[0]{\textcolor{black}{\symCurrentSOCP_{\indexGridLines\indexGridNode \indexGridNodeTwo}  }}

\newcommand{\IijsqsksocpAA}[0]{\textcolor{black}{\symCurrentSOCP_{\indexGridLines\indexGridNode \indexGridNodeTwo,\seriesss, \symPhaseA\symPhaseA}  }}
\newcommand{\IijsqsksocpAB}[0]{\textcolor{\complexcolor}{\symCurrentSOCP_{\indexGridLines\indexGridNode \indexGridNodeTwo,\seriesss, \symPhaseA\symPhaseB}  }}
\newcommand{\IijsqsksocpAC}[0]{\textcolor{\complexcolor}{\symCurrentSOCP_{\indexGridLines\indexGridNode \indexGridNodeTwo,\seriesss, \symPhaseA\symPhaseC}  }}
\newcommand{\IijsqsksocpBA}[0]{\textcolor{\complexcolor}{\symCurrentSOCP_{\indexGridLines\indexGridNode \indexGridNodeTwo,\seriesss, \symPhaseB\symPhaseA}  }}
\newcommand{\IijsqsksocpBB}[0]{\textcolor{black}{\symCurrentSOCP_{\indexGridLines\indexGridNode \indexGridNodeTwo,\seriesss, \symPhaseB\symPhaseB}  }}
\newcommand{\IijsqsksocpBC}[0]{\textcolor{\complexcolor}{\symCurrentSOCP_{\indexGridLines\indexGridNode \indexGridNodeTwo,\seriesss, \symPhaseB\symPhaseC}  }}
\newcommand{\IijsqsksocpCA}[0]{\textcolor{\complexcolor}{\symCurrentSOCP_{\indexGridLines\indexGridNode \indexGridNodeTwo,\seriesss, \symPhaseC\symPhaseA}  }}
\newcommand{\IijsqsksocpCB}[0]{\textcolor{\complexcolor}{\symCurrentSOCP_{\indexGridLines\indexGridNode \indexGridNodeTwo,\seriesss, \symPhaseC\symPhaseB}  }}
\newcommand{\IijsqsksocpCC}[0]{\textcolor{black}{\symCurrentSOCP_{\indexGridLines\indexGridNode \indexGridNodeTwo,\seriesss, \symPhaseC\symPhaseC}  }}

\newcommand{\IijsqksocpAAreal}[0]{\textcolor{black}{\symCurrentSOCP_{\indexGridLines\indexGridNode \indexGridNodeTwo, \symPhaseA\symPhaseA}^{\realss}  }}
\newcommand{\IijsqksocpBBreal}[0]{\textcolor{black}{\symCurrentSOCP_{\indexGridLines\indexGridNode \indexGridNodeTwo, \symPhaseB\symPhaseB}^{\realss}  }}
\newcommand{\IijsqksocpCCreal}[0]{\textcolor{black}{\symCurrentSOCP_{\indexGridLines\indexGridNode \indexGridNodeTwo, \symPhaseC\symPhaseC}^{\realss}  }}

\newcommand{\IijsqksocpNNreal}[0]{\textcolor{black}{\symCurrentSOCP_{\indexGridLines\indexGridNode \indexGridNodeTwo, \symPhaseN\symPhaseN}^{\realss}  }}
\newcommand{\IijsqksocpGGreal}[0]{\textcolor{black}{\symCurrentSOCP_{\indexGridLines\indexGridNode \indexGridNodeTwo, \symPhaseG\symPhaseG}^{\realss}  }}
\newcommand{\kifsqksocpGGreal}[0]{\textcolor{black}{\symCurrentSOCP_{k\indexGridNode f, \symPhaseG\symPhaseG}^{\realss}  }}

\newcommand{\IijingsqksocpGG}[0]{\textcolor{\complexcolor}{\symCurrentSOCP_{\indexGridLines\indexGridNode \indexGridNodeTwo, \symPhaseN\symPhaseG}  }}

\newcommand{\kifingsqksocpGG}[0]{\textcolor{\complexcolor}{\symCurrentSOCP_{k\indexGridNode f, \symPhaseN\symPhaseG}  }}
    \newcommand{\parm}{\mathord{\color{black!33}\bullet}}%

\newcommand{\IijsqsksocpAAreal}[0]{\textcolor{black}{\symCurrentSOCP_{\indexGridLines\indexGridNode \indexGridNodeTwo, \symPhaseA\symPhaseA}^{\seriesss,\realss}  }}
\newcommand{\IijsqsksocpABreal}[0]{\textcolor{black}{\symCurrentSOCP_{\indexGridLines\indexGridNode \indexGridNodeTwo, \symPhaseA\symPhaseB}^{\seriesss,\realss}  }}
\newcommand{\IijsqsksocpACreal}[0]{\textcolor{black}{\symCurrentSOCP_{\indexGridLines\indexGridNode \indexGridNodeTwo, \symPhaseA\symPhaseC}^{\seriesss,\realss}  }}
\newcommand{\IijsqsksocpBBreal}[0]{\textcolor{black}{\symCurrentSOCP_{\indexGridLines\indexGridNode \indexGridNodeTwo, \symPhaseB\symPhaseB}^{\seriesss,\realss}  }}
\newcommand{\IijsqsksocpBCreal}[0]{\textcolor{black}{\symCurrentSOCP_{\indexGridLines\indexGridNode \indexGridNodeTwo, \symPhaseB\symPhaseC}^{\seriesss,\realss}  }}
\newcommand{\IijsqsksocpCCreal}[0]{\textcolor{black}{\symCurrentSOCP_{\indexGridLines\indexGridNode \indexGridNodeTwo, \symPhaseC\symPhaseC}^{\seriesss,\realss}  }}
\newcommand{\IijsqsksocpNNreal}[0]{\textcolor{black}{\symCurrentSOCP_{\indexGridLines\indexGridNode \indexGridNodeTwo, \symPhaseN\symPhaseN}^{\seriesss,\realss}  }}

\newcommand{\IijsqsksocpANreal}[0]{\textcolor{black}{\symCurrentSOCP_{\indexGridLines\indexGridNode \indexGridNodeTwo, \symPhaseA\symPhaseN}^{\seriesss,\realss}  }}
\newcommand{\IijsqsksocpBNreal}[0]{\textcolor{black}{\symCurrentSOCP_{\indexGridLines\indexGridNode \indexGridNodeTwo, \symPhaseB\symPhaseN}^{\seriesss,\realss}  }}
\newcommand{\IijsqsksocpCNreal}[0]{\textcolor{black}{\symCurrentSOCP_{\indexGridLines\indexGridNode \indexGridNodeTwo, \symPhaseC\symPhaseN}^{\seriesss,\realss}  }}

\newcommand{\IijsqsksocpANimag}[0]{\textcolor{black}{\symCurrentSOCP_{\indexGridLines\indexGridNode \indexGridNodeTwo, \symPhaseA\symPhaseN}^{\seriesss,\imagss}  }}
\newcommand{\IijsqsksocpBNimag}[0]{\textcolor{black}{\symCurrentSOCP_{\indexGridLines\indexGridNode \indexGridNodeTwo, \symPhaseB\symPhaseN}^{\seriesss,\imagss}  }}
\newcommand{\IijsqsksocpCNimag}[0]{\textcolor{black}{\symCurrentSOCP_{\indexGridLines\indexGridNode \indexGridNodeTwo, \symPhaseC\symPhaseN}^{\seriesss,\imagss}  }}

\newcommand{\IijsqsksocpGGreal}[0]{\textcolor{black}{\symCurrentSOCP_{\indexGridLines\indexGridNode \indexGridNodeTwo, \symPhaseG\symPhaseG}^{\seriesss,\realss}  }}
\newcommand{\IisqsksocpNGreal}[0]{\textcolor{black}{\symCurrentSOCP_{\indexGridNode , \symPhaseN\symPhaseG}  }}

\newcommand{\IjisqsksocpBBreal}[0]{\textcolor{black}{\symCurrentSOCP_{\indexGridLines\indexGridNodeTwo\indexGridNode , \symPhaseB\symPhaseB}^{\seriesss,\realss}  }}

\newcommand{\IijsqsksocpABimag}[0]{\textcolor{black}{\symCurrentSOCP_{\indexGridLines\indexGridNode \indexGridNodeTwo, \symPhaseA\symPhaseB}^{\seriesss,\imagss}  }}
\newcommand{\IijsqsksocpACimag}[0]{\textcolor{black}{\symCurrentSOCP_{\indexGridLines\indexGridNode \indexGridNodeTwo, \symPhaseA\symPhaseC}^{\seriesss,\imagss}  }}
\newcommand{\IijsqsksocpBCimag}[0]{\textcolor{black}{\symCurrentSOCP_{\indexGridLines\indexGridNode \indexGridNodeTwo, \symPhaseB\symPhaseC}^{\seriesss,\imagss}  }}

\newcommand{\VjsqksocpAAreal}[0]{\textcolor{black}{\symVoltageSOCP_{\indexGridNodeTwo, \symPhaseA\symPhaseA}^{\realss}  }}
\newcommand{\VjsqksocpBBreal}[0]{\textcolor{black}{\symVoltageSOCP_{\indexGridNodeTwo, \symPhaseB\symPhaseB}^{\realss}  }}
\newcommand{\VjsqksocpCCreal}[0]{\textcolor{black}{\symVoltageSOCP_{\indexGridNodeTwo, \symPhaseC\symPhaseC}^{\realss}  }}

\newcommand{\VisqksocpNNreal}[0]{\textcolor{black}{\symVoltageSOCP_{\indexGridNode, \symPhaseN\symPhaseN}^{\realss}  }}
\newcommand{\VisqksocpGGreal}[0]{\textcolor{black}{\symVoltageSOCP_{\indexGridNode, \symPhaseG\symPhaseG}^{\realss}  }}

\newcommand{\VisqksocpAAreal}[0]{\textcolor{black}{\symVoltageSOCP_{\indexGridNode, \symPhaseA\symPhaseA}^{\realss}  }}
\newcommand{\VisqksocpABreal}[0]{\textcolor{black}{\symVoltageSOCP_{\indexGridNode, \symPhaseA\symPhaseB}^{\realss}  }}
\newcommand{\VisqksocpACreal}[0]{\textcolor{black}{\symVoltageSOCP_{\indexGridNode, \symPhaseA\symPhaseC}^{\realss}  }}

\newcommand{\VisqksocpANreal}[0]{\textcolor{black}{\symVoltageSOCP_{\indexGridNode, \symPhaseA\symPhaseN}^{\realss}  }}
\newcommand{\VisqksocpBNreal}[0]{\textcolor{black}{\symVoltageSOCP_{\indexGridNode, \symPhaseB\symPhaseN}^{\realss}  }}
\newcommand{\VisqksocpCNreal}[0]{\textcolor{black}{\symVoltageSOCP_{\indexGridNode, \symPhaseC\symPhaseN}^{\realss}  }}
\newcommand{\VisqksocpANimag}[0]{\textcolor{black}{\symVoltageSOCP_{\indexGridNode, \symPhaseA\symPhaseN}^{\imagss}  }}
\newcommand{\VisqksocpBNimag}[0]{\textcolor{black}{\symVoltageSOCP_{\indexGridNode, \symPhaseB\symPhaseN}^{\imagss}  }}
\newcommand{\VisqksocpCNimag}[0]{\textcolor{black}{\symVoltageSOCP_{\indexGridNode, \symPhaseC\symPhaseN}^{\imagss}  }}

\newcommand{\VisqksocpBBreal}[0]{\textcolor{black}{\symVoltageSOCP_{\indexGridNode, \symPhaseB\symPhaseB}^{\realss}  }}
\newcommand{\VisqksocpBCreal}[0]{\textcolor{black}{\symVoltageSOCP_{\indexGridNode, \symPhaseB\symPhaseC}^{\realss}  }}
\newcommand{\VisqksocpCCreal}[0]{\textcolor{black}{\symVoltageSOCP_{\indexGridNode, \symPhaseC\symPhaseC}^{\realss}  }}

\newcommand{\VisqksocpABimag}[0]{\textcolor{black}{\symVoltageSOCP_{\indexGridNode, \symPhaseA\symPhaseB}^{\imagss}  }}
\newcommand{\VisqksocpACimag}[0]{\textcolor{black}{\symVoltageSOCP_{\indexGridNode, \symPhaseA\symPhaseC}^{\imagss}  }}
\newcommand{\VisqksocpBCimag}[0]{\textcolor{black}{\symVoltageSOCP_{\indexGridNode, \symPhaseB\symPhaseC}^{\imagss}  }}

\newcommand{\VisqksocpAA}[0]{\textcolor{black}{\symVoltageSOCP_{\indexGridNode, \symPhaseA\symPhaseA} }}
\newcommand{\VisqksocpAB}[0]{\textcolor{\complexcolor}{\symVoltageSOCP_{\indexGridNode, \symPhaseA\symPhaseB}}}
\newcommand{\VisqksocpAC}[0]{\textcolor{\complexcolor}{\symVoltageSOCP_{\indexGridNode, \symPhaseA\symPhaseC}  }}
\newcommand{\VisqksocpBA}[0]{\textcolor{\complexcolor}{\symVoltageSOCP_{\indexGridNode, \symPhaseB\symPhaseA} }}
\newcommand{\VisqksocpBB}[0]{\textcolor{black}{\symVoltageSOCP_{\indexGridNode, \symPhaseB\symPhaseB} }}
\newcommand{\VisqksocpBC}[0]{\textcolor{\complexcolor}{\symVoltageSOCP_{\indexGridNode, \symPhaseB\symPhaseC}}}
\newcommand{\VisqksocpCA}[0]{\textcolor{\complexcolor}{\symVoltageSOCP_{\indexGridNode, \symPhaseC\symPhaseA}  }}
\newcommand{\VisqksocpCB}[0]{\textcolor{\complexcolor}{\symVoltageSOCP_{\indexGridNode, \symPhaseC\symPhaseB}  }}
\newcommand{\VisqksocpCC}[0]{\textcolor{black}{\symVoltageSOCP_{\indexGridNode, \symPhaseC\symPhaseC}  }}

\newcommand{\VijsqksocpAAreal}[0]{\textcolor{black}{\symVoltageSOCP_{\indexGridNode\indexGridNodeTwo, \symPhaseA\symPhaseA}^{\realss}  }}
\newcommand{\VijsqksocpABreal}[0]{\textcolor{black}{\symVoltageSOCP_{\indexGridNode\indexGridNodeTwo, \symPhaseA\symPhaseB}^{\realss}  }}
\newcommand{\VijsqksocpACreal}[0]{\textcolor{black}{\symVoltageSOCP_{\indexGridNode\indexGridNodeTwo, \symPhaseA\symPhaseC}^{\realss}  }}

\newcommand{\VijsqksocpBAreal}[0]{\textcolor{black}{\symVoltageSOCP_{\indexGridNode\indexGridNodeTwo, \symPhaseB\symPhaseA}^{\realss}  }}
\newcommand{\VijsqksocpBBreal}[0]{\textcolor{black}{\symVoltageSOCP_{\indexGridNode\indexGridNodeTwo, \symPhaseB\symPhaseB}^{\realss}  }}
\newcommand{\VijsqksocpBCreal}[0]{\textcolor{black}{\symVoltageSOCP_{\indexGridNode\indexGridNodeTwo, \symPhaseB\symPhaseC}^{\realss}  }}

\newcommand{\VijsqksocpCAreal}[0]{\textcolor{black}{\symVoltageSOCP_{\indexGridNode\indexGridNodeTwo, \symPhaseC\symPhaseA}^{\realss}  }}
\newcommand{\VijsqksocpCBreal}[0]{\textcolor{black}{\symVoltageSOCP_{\indexGridNode\indexGridNodeTwo, \symPhaseC\symPhaseB}^{\realss}  }}
\newcommand{\VijsqksocpCCreal}[0]{\textcolor{black}{\symVoltageSOCP_{\indexGridNode\indexGridNodeTwo, \symPhaseC\symPhaseC}^{\realss}  }}

\newcommand{\VjsqksocpABimag}[0]{\textcolor{black}{\symVoltageSOCP_{\indexGridNodeTwo, \symPhaseA\symPhaseB}^{\imagss}  }}
\newcommand{\VjsqksocpACimag}[0]{\textcolor{black}{\symVoltageSOCP_{\indexGridNodeTwo, \symPhaseA\symPhaseC}^{\imagss}  }}
\newcommand{\VjsqksocpBCimag}[0]{\textcolor{black}{\symVoltageSOCP_{\indexGridNodeTwo, \symPhaseB\symPhaseC}^{\imagss}  }}

\newcommand{\VjsqksocpABreal}[0]{\textcolor{black}{\symVoltageSOCP_{\indexGridNodeTwo, \symPhaseA\symPhaseB}^{\realss}  }}
\newcommand{\VjsqksocpACreal}[0]{\textcolor{black}{\symVoltageSOCP_{\indexGridNodeTwo, \symPhaseA\symPhaseC}^{\realss}  }}
\newcommand{\VjsqksocpBCreal}[0]{\textcolor{black}{\symVoltageSOCP_{\indexGridNodeTwo, \symPhaseB\symPhaseC}^{\realss}  }}

\newcommand{\VijsqksocpAAimag}[0]{\textcolor{black}{\symVoltageSOCP_{\indexGridNode\indexGridNodeTwo, \symPhaseA\symPhaseA}^{\imagss}  }}
\newcommand{\VijsqksocpABimag}[0]{\textcolor{black}{\symVoltageSOCP_{\indexGridNode\indexGridNodeTwo, \symPhaseA\symPhaseB}^{\imagss}  }}
\newcommand{\VijsqksocpACimag}[0]{\textcolor{black}{\symVoltageSOCP_{\indexGridNode\indexGridNodeTwo, \symPhaseA\symPhaseC}^{\imagss}  }}

\newcommand{\VijsqksocpBAimag}[0]{\textcolor{black}{\symVoltageSOCP_{\indexGridNode\indexGridNodeTwo, \symPhaseB\symPhaseA}^{\imagss}  }}
\newcommand{\VijsqksocpBBimag}[0]{\textcolor{black}{\symVoltageSOCP_{\indexGridNode\indexGridNodeTwo, \symPhaseB\symPhaseB}^{\imagss}  }}
\newcommand{\VijsqksocpBCimag}[0]{\textcolor{black}{\symVoltageSOCP_{\indexGridNode\indexGridNodeTwo, \symPhaseB\symPhaseC}^{\imagss}  }}

\newcommand{\VijsqksocpCAimag}[0]{\textcolor{black}{\symVoltageSOCP_{\indexGridNode\indexGridNodeTwo, \symPhaseC\symPhaseA}^{\imagss}  }}
\newcommand{\VijsqksocpCBimag}[0]{\textcolor{black}{\symVoltageSOCP_{\indexGridNode\indexGridNodeTwo, \symPhaseC\symPhaseB}^{\imagss}  }}
\newcommand{\VijsqksocpCCimag}[0]{\textcolor{black}{\symVoltageSOCP_{\indexGridNode\indexGridNodeTwo, \symPhaseC\symPhaseC}^{\imagss}  }}

\newcommand{\IijsqskA}[0]{\textcolor{black}{\symCurrent_{\indexGridLines\indexGridNode \indexGridNodeTwo,\seriesss, \symPhaseA}^{2}  }}
\newcommand{\IijsqskB}[0]{\textcolor{black}{\symCurrent_{\indexGridLines\indexGridNode \indexGridNodeTwo,\seriesss, \symPhaseB}^{2}  }}
\newcommand{\IijsqskC}[0]{\textcolor{black}{\symCurrent_{\indexGridLines\indexGridNode \indexGridNodeTwo,\seriesss, \symPhaseC}^{2}  }}

\newcommand{\Sijrated}[0]{\textcolor{\sizingcolor}{\symApparentPower_{\indexGridLines\indexGridNode \indexGridNodeTwo}^{\ratedss} }}
\newcommand{\Sjirated}[0]{\textcolor{\sizingcolor}{\symApparentPower_{\indexGridLines\indexGridNodeTwo \indexGridNode }^{\ratedss} }}
\newcommand{\SijratedA}[0]{\textcolor{\sizingcolor}{\symApparentPower_{\indexGridLines\indexGridNode \indexGridNodeTwo,\symPhaseA}^{\ratedss} }}
\newcommand{\SjiratedA}[0]{\textcolor{\sizingcolor}{\symApparentPower_{\indexGridLines\indexGridNodeTwo \indexGridNode ,\symPhaseA}^{\ratedss} }}
\newcommand{\SijratedB}[0]{\textcolor{\sizingcolor}{\symApparentPower_{\indexGridLines\indexGridNode \indexGridNodeTwo,\symPhaseB}^{\ratedss} }}
\newcommand{\SjiratedB}[0]{\textcolor{\sizingcolor}{\symApparentPower_{\indexGridLines\indexGridNodeTwo \indexGridNode ,\symPhaseB}^{\ratedss} }}
\newcommand{\SijratedC}[0]{\textcolor{\sizingcolor}{\symApparentPower_{\indexGridLines\indexGridNode \indexGridNodeTwo,\symPhaseC}^{\ratedss} }}
\newcommand{\SijratedN}[0]{\textcolor{\sizingcolor}{\symApparentPower_{\indexGridLines\indexGridNode \indexGridNodeTwo,\symPhaseN}^{\ratedss} }}
\newcommand{\SjiratedC}[0]{\textcolor{\sizingcolor}{\symApparentPower_{\indexGridLines\indexGridNodeTwo \indexGridNode ,\symPhaseC}^{\ratedss} }}

\newcommand{\VikSDPmin}[0]{\textcolor{\paramcolor}{\mathbf{\symVoltage}_{\indexGridNode }^{\minss}  }}
\newcommand{\VikSDPmax}[0]{\textcolor{\paramcolor}{\mathbf{\symVoltage}_{\indexGridNode }^{\maxss}  }}
\newcommand{\VikSDPref}[0]{\textcolor{\paramcolor}{\mathbf{\symVoltage}_{\indexGridNode }^{\refss}  }}

\newcommand{\VikSDPwye}[0]{\textcolor{\complexcolor}{\mathbf{\symVoltage}_{\indexGridNode }^{\text{pn}}  }}

\newcommand{\VikSDPminwye}[0]{\textcolor{\paramcolor}{\mathbf{\symVoltage}_{\indexGridNode }^{\text{pn},\minss}  }}
\newcommand{\VikSDPmaxwye}[0]{\textcolor{\paramcolor}{\mathbf{\symVoltage}_{\indexGridNode }^{\text{pn},\maxss}  }}

\newcommand{\VikSDPmindelta}[0]{\textcolor{\paramcolor}{\mathbf{\symVoltage}_{\indexGridNode }^{\Delta\minss}  }}
\newcommand{\VikSDPmaxdelta}[0]{\textcolor{\paramcolor}{\mathbf{\symVoltage}_{\indexGridNode }^{\Delta\maxss}  }}

\newcommand{\VjkSDPmin}[0]{\textcolor{\paramcolor}{\mathbf{\symVoltage}_{\indexGridNodeTwo }^{\minss}  }}
\newcommand{\VjkSDPmax}[0]{\textcolor{\paramcolor}{\mathbf{\symVoltage}_{\indexGridNodeTwo }^{\maxss}  }}

\newcommand{\VjkSDPseq}[0]{\textcolor{\complexcolor}{\mathbf{\symVoltage}_{\indexGridNodeTwo }^{\fortescuess}  }}

\newcommand{\VikSDPacc}[0]{\textcolor{\complexcolor}{\mathbf{\bar{\symVoltage}}_{\indexGridNode }  }}
\newcommand{\VikSDPaccreal}[0]{\textcolor{black}{\mathbf{\bar{\symVoltage}}_{\indexGridNode }^{\realss}  }}
\newcommand{\VikSDPaccimag}[0]{\textcolor{black}{\mathbf{\bar{\symVoltage}}_{\indexGridNode }^{\imagss}  }}

\newcommand{\VikSDPdelta}[0]{\textcolor{\complexcolor}{\mathbf{\symVoltage}_{\indexGridNode }^{\Delta}  }}
\newcommand{\VikSDPdeltareal}[0]{\textcolor{black}{\mathbf{\symVoltage}_{\indexGridNode }^{\Delta, \realss}  }}
\newcommand{\VikSDPdeltaimag}[0]{\textcolor{black}{\mathbf{\symVoltage}_{\indexGridNode }^{\Delta, \imagss}  }}

\newcommand{\VikSDPseq}[0]{\textcolor{\complexcolor}{\mathbf{\symVoltage}_{\indexGridNode }^{\fortescuess}  }}
\newcommand{\VikSDPseqmin}[0]{\textcolor{\paramcolor}{\mathbf{\symVoltage}_{\indexGridNode }^{\fortescuess,\minss}  }}
\newcommand{\VikSDPseqmax}[0]{\textcolor{\paramcolor}{\mathbf{\symVoltage}_{\indexGridNode }^{\fortescuess,\maxss}  }}
\newcommand{\VikSDP}[0]{\textcolor{\complexcolor}{\mathbf{\symVoltage}_{\indexGridNode }^{}  }}
\newcommand{\VikSDPreal}[0]{\textcolor{black}{\mathbf{\symVoltage}_{\indexGridNode }^{\realss}  }}
\newcommand{\VikSDPimag}[0]{\textcolor{black}{\mathbf{\symVoltage}_{\indexGridNode }^{\imagss}  }}
\newcommand{\VjkSDPreal}[0]{\textcolor{black}{\mathbf{\symVoltage}_{\indexGridNodeTwo }^{\realss}  }}
\newcommand{\VjkSDPimag}[0]{\textcolor{black}{\mathbf{\symVoltage}_{\indexGridNodeTwo }^{\imagss}  }}
\newcommand{\VikSDPH}[0]{\textcolor{\complexcolor}{\mathbf{\symVoltage}_{\indexGridNode }^{\hermitiantranspose}  }}
\newcommand{\VjkSDP}[0]{\textcolor{\complexcolor}{\mathbf{\symVoltage}_{\indexGridNodeTwo }^{}  }}
\newcommand{\VzkSDP}[0]{\textcolor{\complexcolor}{\mathbf{\symVoltage}_{z }^{}  }}
\newcommand{\VjkSDPH}[0]{\textcolor{\complexcolor}{\mathbf{\symVoltage}_{\indexGridNodeTwo }^{\hermitiantranspose}  }}
\newcommand{\Vik}[0]{\textcolor{black}{\symVoltage_{\indexGridNode }^{}  }}
\newcommand{\Vjk}[0]{\textcolor{black}{\symVoltage_{\indexGridNodeTwo }^{}  }}
\newcommand{\Vikdelta}[0]{\textcolor{black}{\Delta\symVoltage_{\indexGridNode }^{}  }}
\newcommand{\Vjkdelta}[0]{\textcolor{black}{\Delta\symVoltage_{\indexGridNodeTwo }^{}  }}
\newcommand{\Viacck}[0]{\textcolor{black}{\symVoltage_{\indexGridNode }^{'}  }}
\newcommand{\Vjacck}[0]{\textcolor{black}{\symVoltage_{\indexGridNodeTwo }^{'}  }}
\newcommand{\Viacckl}[0]{\textcolor{black}{\symVoltage_{\indexGridNode }^{\star}  }}
\newcommand{\Vjacckl}[0]{\textcolor{black}{\symVoltage_{\indexGridNodeTwo }^{\star}  }}

\newcommand{\VikSDPmag}[0]{\textcolor{black}{\mathbf{\symVoltage}_{\indexGridNode }^{\text{mag}}  }}
\newcommand{\VikSDPang}[0]{\textcolor{black}{\mathbf{\Theta}_{\indexGridNode }  }}

\newcommand{\VjkSDPmag}[0]{\textcolor{black}{\mathbf{\symVoltage}_{\indexGridNodeTwo }^{\text{mag}}  }}
\newcommand{\VjkSDPang}[0]{\textcolor{black}{\mathbf{\Theta}_{\indexGridNodeTwo }  }}

\newcommand{\VukrefSDP}[0]{\textcolor{\complexparamcolor}{\mathbf{\symVoltage}_{\indexUnit }^{\refss}  }}

\newcommand{\IurefSDP}[0]{\textcolor{\complexparamcolor}{\mathbf{\symCurrent}_{\indexGridNode }^{\refss}  }}

\newcommand{\VikrefSDP}[0]{\textcolor{\complexparamcolor}{\mathbf{\symVoltage}_{\indexGridNode }^{\refss}  }}
\newcommand{\VikrefSDPreal}[0]{\textcolor{\paramcolor}{\mathbf{\symVoltage}_{\indexGridNode }^{\refss,\realss}  }}
\newcommand{\VikrefSDPimag}[0]{\textcolor{\paramcolor}{\mathbf{\symVoltage}_{\indexGridNode }^{\refss,\imagss}  }}
\newcommand{\VikdeltaSDP}[0]{\textcolor{\complexcolor}{\mathbf{\symVoltage}_{\indexGridNode }^{\Delta}  }}
\newcommand{\VikdeltaSDPreal}[0]{\textcolor{black}{\mathbf{\symVoltage}_{\indexGridNode }^{\Delta,\realss}  }}
\newcommand{\VikdeltaSDPimag}[0]{\textcolor{black}{\mathbf{\symVoltage}_{\indexGridNode }^{\Delta,\imagss}  }}

\newcommand{\VsqikdeltaSDP}[0]{\textcolor{\complexcolor}{\mathbf{\symVoltageSOCP}_{\indexGridNode }^{\Delta}  }}

\newcommand{\VikaccSDP}[0]{\textcolor{\complexcolor}{\mathbf{\symVoltage}_{\indexGridNode }^{'}  }}
\newcommand{\VjkaccSDP}[0]{\textcolor{\complexcolor}{\mathbf{\symVoltage}_{\indexGridNodeTwo }^{'}  }}
\newcommand{\Cdelta}[0]{\textcolor{\complexcolor}{\mathbf{\symConnectivity}_{\Delta }^{}  }}
\newcommand{\Czigzag}[0]{\textcolor{\complexcolor}{\mathbf{\symConnectivity}_{Z }^{}  }}
\newcommand{\Cwye}[0]{\textcolor{\complexcolor}{\mathbf{\symConnectivity}_{Y }^{}  }}

\newcommand{\CijkSDP}[0]{\textcolor{\complexcolor}{\mathbf{\symConnectivity}_{\indexGridNode\indexGridNodeTwo }^{}  }}
\newcommand{\CjikSDP}[0]{\textcolor{\complexcolor}{\mathbf{\symConnectivity}_{\indexGridNodeTwo\indexGridNode }^{}  }}

\newcommand{\Visqk}[0]{\textcolor{black}{\symVoltage_{\indexGridNode }^{2}  }}
\newcommand{\Vjsqk}[0]{\textcolor{black}{\symVoltage_{\indexGridNodeTwo }^{2}  }}

\newcommand{\Vjh}[0]{\textcolor{\complexcolor}{\mathbf{\symVoltage}_{\indexGridNodeTwo ,\symHarmonic}   }}

\newcommand{\VjAh}[0]{\textcolor{\complexcolor}{\symVoltage_{\indexGridNodeTwo,\symPhaseA ,\symHarmonic}  }}
\newcommand{\VjBh}[0]{\textcolor{\complexcolor}{\symVoltage_{\indexGridNodeTwo,\symPhaseB,\symHarmonic }  }}
\newcommand{\VjCh}[0]{\textcolor{\complexcolor}{\symVoltage_{\indexGridNodeTwo,\symPhaseC,\symHarmonic }   }}

\newcommand{\Vih}[0]{\textcolor{\complexcolor}{\mathbf{\symVoltage}_{\indexGridNode ,\symHarmonic}   }}

\newcommand{\ViAh}[0]{\textcolor{\complexcolor}{\symVoltage_{\indexGridNode,\symPhaseA ,\symHarmonic}  }}
\newcommand{\ViBh}[0]{\textcolor{\complexcolor}{\symVoltage_{\indexGridNode,\symPhaseB,\symHarmonic }  }}
\newcommand{\ViCh}[0]{\textcolor{\complexcolor}{\symVoltage_{\indexGridNode,\symPhaseC,\symHarmonic }   }}

\newcommand{\VmagiAh}[0]{\symVoltage_{\indexGridNode,\symPhaseA ,\symHarmonic}^{\text{mag}}  }
\newcommand{\VmagiBh}[0]{\symVoltage_{\indexGridNode,\symPhaseB,\symHarmonic }^{\text{mag}}  }
\newcommand{\VmagiCh}[0]{\symVoltage_{\indexGridNode,\symPhaseC,\symHarmonic }^{\text{mag}}  }

\newcommand{\VangiAh}[0]{\symVoltageAngle_{\indexGridNode,\symPhaseA ,\symHarmonic}  }
\newcommand{\VangiBh}[0]{\symVoltageAngle_{\indexGridNode,\symPhaseB,\symHarmonic }  }
\newcommand{\VangiCh}[0]{\symVoltageAngle_{\indexGridNode,\symPhaseC,\symHarmonic }  }

\newcommand{\VmagiAhmax}[0]{\textcolor{\paramcolor}{\symVoltage_{\indexGridNode,\symPhaseA ,\symHarmonic}^{\text{max}}  }}
\newcommand{\VmagiBhmax}[0]{\textcolor{\paramcolor}{\symVoltage_{\indexGridNode,\symPhaseB,\symHarmonic }^{\text{max}}  }}
\newcommand{\VmagiChmax}[0]{\textcolor{\paramcolor}{\symVoltage_{\indexGridNode,\symPhaseC,\symHarmonic }^{\text{max}}  }}

\newcommand{\VmagiAthd}[0]{\textcolor{\paramcolor}{\symVoltage_{\indexGridNode,\symPhaseA }^{\text{THD}}  }}
\newcommand{\VmagiBthd}[0]{\textcolor{\paramcolor}{\symVoltage_{\indexGridNode,\symPhaseB }^{\text{THD}}  }}
\newcommand{\VmagiCthd}[0]{\textcolor{\paramcolor}{\symVoltage_{\indexGridNode,\symPhaseC }^{\text{THD}}  }}

\newcommand{\VikAN}[0]{\textcolor{\complexcolor}{\symVoltage_{\indexGridNode,\symPhaseA\symPhaseN }^{}  }}
\newcommand{\VikBN}[0]{\textcolor{\complexcolor}{\symVoltage_{\indexGridNode,\symPhaseB \symPhaseN}^{}  }}
\newcommand{\VikCN}[0]{\textcolor{\complexcolor}{\symVoltage_{\indexGridNode,\symPhaseC\symPhaseN }^{}  }}

\newcommand{\VikAB}[0]{\textcolor{\complexcolor}{\symVoltage_{\indexGridNode,\symPhaseA\symPhaseB }^{\Delta}  }}
\newcommand{\VikBC}[0]{\textcolor{\complexcolor}{\symVoltage_{\indexGridNode,\symPhaseB\symPhaseC }^{\Delta}  }}
\newcommand{\VikCA}[0]{\textcolor{\complexcolor}{\symVoltage_{\indexGridNode,\symPhaseC\symPhaseA }^{\Delta}  }}

\newcommand{\Vikzero}[0]{\textcolor{\complexcolor}{\symVoltage_{\indexGridNode,0 }^{}  }}
\newcommand{\Vikone}[0]{\textcolor{\complexcolor}{\symVoltage_{\indexGridNode,1 }^{}  }}
\newcommand{\Viktwo}[0]{\textcolor{\complexcolor}{\symVoltage_{\indexGridNode,2 }^{}  }}

\newcommand{\VikANreal}[0]{\textcolor{black}{\symVoltage_{\indexGridNode,\symPhaseA }^{\realss}  }}
\newcommand{\VikBNreal}[0]{\textcolor{black}{\symVoltage_{\indexGridNode,\symPhaseB }^{\realss}  }}
\newcommand{\VikCNreal}[0]{\textcolor{black}{\symVoltage_{\indexGridNode,\symPhaseC }^{\realss}  }}
\newcommand{\VikNreal}[0]{\textcolor{black}{\symVoltage_{\indexGridNode,\symPhaseN }^{\realss}  }}

\newcommand{\VikANimag}[0]{\textcolor{black}{\symVoltage_{\indexGridNode,\symPhaseA }^{\imagss}  }}
\newcommand{\VikBNimag}[0]{\textcolor{black}{\symVoltage_{\indexGridNode,\symPhaseB }^{\imagss}  }}
\newcommand{\VikCNimag}[0]{\textcolor{black}{\symVoltage_{\indexGridNode,\symPhaseC }^{\imagss}  }}
\newcommand{\VikNimag}[0]{\textcolor{black}{\symVoltage_{\indexGridNode,\symPhaseN }^{\imagss}  }}

\newcommand{\VikpNreal}[0]{\textcolor{black}{\symVoltage_{\indexGridNode,\indexPhases }^{\realss}  }}
\newcommand{\VikpNimag}[0]{\textcolor{black}{\symVoltage_{\indexGridNode,\indexPhases }^{\imagss}  }}
\newcommand{\VikppNreal}[0]{\textcolor{black}{\symVoltage_{\indexGridNode,\indexPhasesTwo }^{\realss}  }}
\newcommand{\VikppNimag}[0]{\textcolor{black}{\symVoltage_{\indexGridNode,\indexPhasesTwo }^{\imagss}  }}
\newcommand{\VjkppNreal}[0]{\textcolor{black}{\symVoltage_{\indexGridNodeTwo,\indexPhasesTwo }^{\realss}  }}
\newcommand{\VjkppNimag}[0]{\textcolor{black}{\symVoltage_{\indexGridNodeTwo,\indexPhasesTwo }^{\imagss}  }}

\newcommand{\VikpNmag}[0]{\textcolor{black}{\symVoltage_{\indexGridNode,\indexPhases }^{\text{mag}}  }}
\newcommand{\VikppNmag}[0]{\textcolor{black}{\symVoltage_{\indexGridNode,\indexPhasesTwo }^{\text{mag}}  }}
\newcommand{\VjkppNmag}[0]{\textcolor{black}{\symVoltage_{\indexGridNodeTwo,\indexPhasesTwo }^{\text{mag}}  }}

\newcommand{\VikANmag}[0]{\textcolor{black}{\symVoltage_{\indexGridNode,\symPhaseA }^{\text{mag}}  }}
\newcommand{\VikBNmag}[0]{\textcolor{black}{\symVoltage_{\indexGridNode,\symPhaseB }^{\text{mag}}  }}
\newcommand{\VikCNmag}[0]{\textcolor{black}{\symVoltage_{\indexGridNode,\symPhaseC }^{\text{mag}}  }}
\newcommand{\Vikmag}[0]{\textcolor{black}{\symVoltage_{\indexGridNode }^{\text{mag}}  }}

\newcommand{\VjkANmag}[0]{\textcolor{black}{\symVoltage_{\indexGridNodeTwo,\symPhaseA }^{\text{mag}}  }}
\newcommand{\VjkBNmag}[0]{\textcolor{black}{\symVoltage_{\indexGridNodeTwo,\symPhaseB }^{\text{mag}}  }}
\newcommand{\VjkCNmag}[0]{\textcolor{black}{\symVoltage_{\indexGridNodeTwo,\symPhaseC }^{\text{mag}}  }}

\newcommand{\VjkANreal}[0]{\textcolor{black}{\symVoltage_{\indexGridNodeTwo,\symPhaseA }^{\text{re}}  }}
\newcommand{\VjkBNreal}[0]{\textcolor{black}{\symVoltage_{\indexGridNodeTwo,\symPhaseB }^{\text{re}}  }}
\newcommand{\VjkCNreal}[0]{\textcolor{black}{\symVoltage_{\indexGridNodeTwo,\symPhaseC }^{\text{re}}  }}

\newcommand{\VjkANimag}[0]{\textcolor{black}{\symVoltage_{\indexGridNodeTwo,\symPhaseA }^{\text{im}}  }}
\newcommand{\VjkBNimag}[0]{\textcolor{black}{\symVoltage_{\indexGridNodeTwo,\symPhaseB }^{\text{im}}  }}
\newcommand{\VjkCNimag}[0]{\textcolor{black}{\symVoltage_{\indexGridNodeTwo,\symPhaseC }^{\text{im}}  }}

\newcommand{\VikANangle}[0]{\textcolor{black}{\symVoltageAngle_{\indexGridNode,\symPhaseA }^{}  }}
\newcommand{\VikBNangle}[0]{\textcolor{black}{\symVoltageAngle_{\indexGridNode,\symPhaseB }^{}  }}
\newcommand{\VikCNangle}[0]{\textcolor{black}{\symVoltageAngle_{\indexGridNode,\symPhaseC }^{}  }}
\newcommand{\VikNangle}[0]{\textcolor{black}{\symVoltageAngle_{\indexGridNode,\symPhaseN }^{}  }}

\newcommand{\VjkANangle}[0]{\textcolor{black}{\symVoltageAngle_{\indexGridNodeTwo,\symPhaseA }^{}  }}
\newcommand{\VjkBNangle}[0]{\textcolor{black}{\symVoltageAngle_{\indexGridNodeTwo,\symPhaseB }^{}  }}
\newcommand{\VjkCNangle}[0]{\textcolor{black}{\symVoltageAngle_{\indexGridNodeTwo,\symPhaseC }^{}  }}
\newcommand{\VjkNangle}[0]{\textcolor{black}{\symVoltageAngle_{\indexGridNodeTwo,\symPhaseN }^{}  }}

\newcommand{\VikpNangle}[0]{\textcolor{black}{\symVoltageAngle_{\indexGridNode,\indexPhases }^{}  }}
\newcommand{\VikhNangle}[0]{\textcolor{black}{\symVoltageAngle_{\indexGridNode,\indexPhasesTwo }^{}  }}
\newcommand{\VjkpNangle}[0]{\textcolor{black}{\symVoltageAngle_{\indexGridNodeTwo,\indexPhases }^{}  }}
\newcommand{\VjkhNangle}[0]{\textcolor{black}{\symVoltageAngle_{\indexGridNodeTwo,\indexPhasesTwo }^{}  }}

\newcommand{\VikrefAN}[0]{\textcolor{\complexparamcolor}{\symVoltage_{\indexGridNode,\symPhaseA }^{\refss}  }}
\newcommand{\VikrefBN}[0]{\textcolor{\complexparamcolor}{\symVoltage_{\indexGridNode,\symPhaseB }^{\refss}  }}
\newcommand{\VikrefCN}[0]{\textcolor{\complexparamcolor}{\symVoltage_{\indexGridNode,\symPhaseC  }^{\refss}  }}

\newcommand{\VikANref}[0]{\textcolor{\paramcolor}{\symVoltage_{\indexGridNode,\symPhaseA }^{\refss}  }}
\newcommand{\VikBNref}[0]{\textcolor{\paramcolor}{\symVoltage_{\indexGridNode,\symPhaseB }^{\refss}  }}
\newcommand{\VikCNref}[0]{\textcolor{\paramcolor}{\symVoltage_{\indexGridNode,\symPhaseC }^{\refss}  }}
\newcommand{\VikNref}[0]{\textcolor{\paramcolor}{\symVoltage_{\indexGridNode,n }^{\refss}  }}

\newcommand{\thetaikANref}[0]{\textcolor{\paramcolor}{\symVoltageAngle_{\indexGridNode,\symPhaseA }^{\refss}  }}
\newcommand{\thetaikBNref}[0]{\textcolor{\paramcolor}{\symVoltageAngle_{\indexGridNode,\symPhaseB }^{\refss}  }}
\newcommand{\thetaikCNref}[0]{\textcolor{\paramcolor}{\symVoltageAngle_{\indexGridNode,\symPhaseC }^{\refss}  }}
\newcommand{\thetaikNref}[0]{\textcolor{\paramcolor}{\symVoltageAngle_{\indexGridNode,n }^{\refss}  }}

\newcommand{\VikAmin}[0]{\textcolor{\paramcolor}{\symVoltage_{\indexGridNode,\symPhaseA }^{\minss}  }}
\newcommand{\VikBmin}[0]{\textcolor{\paramcolor}{\symVoltage_{\indexGridNode,\symPhaseB }^{\minss}  }}
\newcommand{\VikCmin}[0]{\textcolor{\paramcolor}{\symVoltage_{\indexGridNode,\symPhaseC }^{\minss}  }}
\newcommand{\VikAmax}[0]{\textcolor{\paramcolor}{\symVoltage_{\indexGridNode,\symPhaseA }^{\maxss}  }}
\newcommand{\VikBmax}[0]{\textcolor{\paramcolor}{\symVoltage_{\indexGridNode,\symPhaseB }^{\maxss}  }}
\newcommand{\VikCmax}[0]{\textcolor{\paramcolor}{\symVoltage_{\indexGridNode,\symPhaseC }^{\maxss}  }}

\newcommand{\VikANmin}[0]{\textcolor{\paramcolor}{\symVoltage_{\indexGridNode,\symPhaseA\symPhaseN }^{\minss}  }}
\newcommand{\VikBNmin}[0]{\textcolor{\paramcolor}{\symVoltage_{\indexGridNode,\symPhaseB\symPhaseN }^{\minss}  }}
\newcommand{\VikCNmin}[0]{\textcolor{\paramcolor}{\symVoltage_{\indexGridNode,\symPhaseC\symPhaseN }^{\minss}  }}
\newcommand{\VikNmin}[0]{\textcolor{\paramcolor}{\symVoltage_{\indexGridNode,\symPhaseN }^{\minss}  }}
\newcommand{\VikANmax}[0]{\textcolor{\paramcolor}{\symVoltage_{\indexGridNode,\symPhaseA\symPhaseN }^{\maxss}  }}
\newcommand{\VikBNmax}[0]{\textcolor{\paramcolor}{\symVoltage_{\indexGridNode,\symPhaseB\symPhaseN }^{\maxss}  }}
\newcommand{\VikCNmax}[0]{\textcolor{\paramcolor}{\symVoltage_{\indexGridNode,\symPhaseC\symPhaseN }^{\maxss}  }}
\newcommand{\VikNmax}[0]{\textcolor{\paramcolor}{\symVoltage_{\indexGridNode,\symPhaseN }^{\maxss}  }}

\newcommand{\VikABmax}[0]{\textcolor{\paramcolor}{\symVoltage_{\indexGridNode,\symPhaseA\symPhaseB }^{\maxss}  }}
\newcommand{\VikBCmax}[0]{\textcolor{\paramcolor}{\symVoltage_{\indexGridNode,\symPhaseB\symPhaseC }^{\maxss}  }}
\newcommand{\VikCAmax}[0]{\textcolor{\paramcolor}{\symVoltage_{\indexGridNode,\symPhaseC\symPhaseA }^{\maxss}  }}
\newcommand{\VikABmin}[0]{\textcolor{\paramcolor}{\symVoltage_{\indexGridNode,\symPhaseA\symPhaseB }^{\minss}  }}
\newcommand{\VikBCmin}[0]{\textcolor{\paramcolor}{\symVoltage_{\indexGridNode,\symPhaseB\symPhaseC }^{\minss}  }}
\newcommand{\VikCAmin}[0]{\textcolor{\paramcolor}{\symVoltage_{\indexGridNode,\symPhaseC\symPhaseA }^{\minss}  }}

\newcommand{\VjkAN}[0]{\textcolor{\complexcolor}{\symVoltage_{\indexGridNodeTwo,\symPhaseA }^{}  }}
\newcommand{\VjkBN}[0]{\textcolor{\complexcolor}{\symVoltage_{\indexGridNodeTwo,\symPhaseB }^{}  }}
\newcommand{\VjkCN}[0]{\textcolor{\complexcolor}{\symVoltage_{\indexGridNodeTwo,\symPhaseC }^{}  }}
\newcommand{\VjkANmax}[0]{\textcolor{\paramcolor}{\symVoltage_{\indexGridNodeTwo,\symPhaseA }^{\maxss}  }}

\newcommand{\VjkA}[0]{\textcolor{\complexcolor}{\symVoltage_{\indexGridNodeTwo,\symPhaseA }^{}  }}
\newcommand{\VjkB}[0]{\textcolor{\complexcolor}{\symVoltage_{\indexGridNodeTwo,\symPhaseB }^{}  }}
\newcommand{\VjkC}[0]{\textcolor{\complexcolor}{\symVoltage_{\indexGridNodeTwo,\symPhaseC }^{}  }}
\newcommand{\VjkN}[0]{\textcolor{\complexcolor}{\symVoltage_{\indexGridNodeTwo,\symPhaseN }^{}  }}
\newcommand{\VjkG}[0]{\textcolor{\complexcolor}{\symVoltage_{\indexGridNodeTwo,\symPhaseG }^{}  }}

\newcommand{\VikA}[0]{\textcolor{\complexcolor}{\symVoltage_{\indexGridNode,\symPhaseA }^{}  }}
\newcommand{\VikB}[0]{\textcolor{\complexcolor}{\symVoltage_{\indexGridNode,\symPhaseB }^{}  }}
\newcommand{\VikC}[0]{\textcolor{\complexcolor}{\symVoltage_{\indexGridNode,\symPhaseC }^{}  }}
\newcommand{\VikN}[0]{\textcolor{\complexcolor}{\symVoltage_{\indexGridNode,\symPhaseN }^{}  }}
\newcommand{\VikG}[0]{\textcolor{\complexcolor}{\symVoltage_{\indexGridNode,\symPhaseG }^{}  }}

\newcommand{\Vikmin}[0]{\textcolor{\boundscolor}{\symVoltage^{\minss}_{\indexGridNode }  }}
\newcommand{\Vjkmin}[0]{\textcolor{\boundscolor}{\symVoltage^{\minss}_{\indexGridNodeTwo }  }}
\newcommand{\Vikmax}[0]{\textcolor{\boundscolor}{\symVoltage^{\maxss}_{\indexGridNode }  }}
\newcommand{\Vjkmax}[0]{\textcolor{\boundscolor}{\symVoltage^{\maxss}_{\indexGridNodeTwo }  }}
\newcommand{\bigMikmax}[0]{\textcolor{\boundscolor}{M^{}_{\indexGridNode }  }}
\newcommand{\bigMjkmax}[0]{\textcolor{\boundscolor}{M^{}_{\indexGridNodeTwo }  }}
\newcommand{\Virated}[0]{\textcolor{\sizingcolor}{\symVoltage^{\ratedss}_{\indexGridNode } }}
\newcommand{\Vjrated}[0]{\textcolor{\sizingcolor}{\symVoltage^{\ratedss}_{\indexGridNodeTwo } }}
\newcommand{\Vijrated}[0]{\textcolor{\sizingcolor}{\symVoltage^{\ratedss}_{\indexGridNode \indexGridNodeTwo} }}
\newcommand{\Vjirated}[0]{\textcolor{\sizingcolor}{\symVoltage^{\ratedss}_{\indexGridNodeTwo \indexGridNode } }}

\newcommand{\ViratedAN}[0]{\textcolor{\sizingcolor}{\symVoltage^{\ratedss}_{\indexGridNode,\symPhaseA\symPhaseN } }}
\newcommand{\ViratedBN}[0]{\textcolor{\sizingcolor}{\symVoltage^{\ratedss}_{\indexGridNode,\symPhaseB\symPhaseN } }}
\newcommand{\ViratedCN}[0]{\textcolor{\sizingcolor}{\symVoltage^{\ratedss}_{\indexGridNode,\symPhaseC\symPhaseN } }}

\newcommand{\VijratedAN}[0]{\textcolor{\sizingcolor}{\symVoltage^{\ratedss}_{\indexGridNode \indexGridNodeTwo,\symPhaseA\symPhaseN} }}
\newcommand{\VjiratedAN}[0]{\textcolor{\sizingcolor}{\symVoltage^{\ratedss}_{\indexGridNodeTwo \indexGridNode,\symPhaseA\symPhaseN } }}
\newcommand{\VijratedBN}[0]{\textcolor{\sizingcolor}{\symVoltage^{\ratedss}_{\indexGridNode \indexGridNodeTwo,\symPhaseB\symPhaseN} }}
\newcommand{\VjiratedBN}[0]{\textcolor{\sizingcolor}{\symVoltage^{\ratedss}_{\indexGridNodeTwo \indexGridNode ,\symPhaseB\symPhaseN} }}
\newcommand{\VijratedCN}[0]{\textcolor{\sizingcolor}{\symVoltage^{\ratedss}_{\indexGridNode \indexGridNodeTwo,\symPhaseC\symPhaseN} }}
\newcommand{\VjiratedCN}[0]{\textcolor{\sizingcolor}{\symVoltage^{\ratedss}_{\indexGridNodeTwo \indexGridNode,\symPhaseC\symPhaseN } }}

\newcommand{\Visqksocplin}[0]{\textcolor{black}{\symVoltage_{\indexGridNode }^{\linss}  }}

\newcommand{\VisqkSDPseq}[0]{\textcolor{\complexcolor}{\mathbf{\symVoltageSOCP}_{\indexGridNode }^{\fortescuess} }}

\newcommand{\VusqkSDP}[0]{\textcolor{\complexcolor}{\mathbf{\symVoltageSOCP}_{\indexUnit }  }}

\newcommand{\VisqkSDPdelta}[0]{\textcolor{\complexcolor}{\mathbf{\symVoltageSOCP}_{\indexGridNode }^{\Delta}  }}
\newcommand{\VisqkSDPdeltareal}[0]{\textcolor{black}{\mathbf{\symVoltageSOCP}_{\indexGridNode }^{\Delta \realss}  }}
\newcommand{\VisqkSDPdeltaimag}[0]{\textcolor{black}{\mathbf{\symVoltageSOCP}_{\indexGridNode }^{\Delta \imagss}  }}

\newcommand{\VisqkSDP}[0]{\textcolor{\complexcolor}{\mathbf{\symVoltageSOCP}_{\indexGridNode }  }}
\newcommand{\VjsqkSDP}[0]{\textcolor{\complexcolor}{\mathbf{\symVoltageSOCP}_{\indexGridNodeTwo }  }}
\newcommand{\VisqkSDPreal}[0]{\textcolor{black}{\mathbf{\symVoltageSOCP}_{\indexGridNode }^{\realss}  }}
\newcommand{\VisqkSDPimag}[0]{\textcolor{black}{\mathbf{\symVoltageSOCP}_{\indexGridNode }^{\imagss}  }}
\newcommand{\VjsqkSDPreal}[0]{\textcolor{black}{\mathbf{\symVoltageSOCP}_{\indexGridNodeTwo }^{\realss}  }}
\newcommand{\VjsqkSDPimag}[0]{\textcolor{black}{\mathbf{\symVoltageSOCP}_{\indexGridNodeTwo }^{\imagss}  }}
\newcommand{\VisqkaccSDP}[0]{\textcolor{\complexcolor}{\mathbf{\symVoltageSOCP}_{\indexGridNode }^{'}  }}

\newcommand{\VzsqkSDP}[0]{\textcolor{\complexcolor}{\mathbf{\symVoltageSOCP}_{z }  }}
\newcommand{\VsqkSDP}[0]{\textcolor{\complexcolor}{\mathbf{M}  }}
\newcommand{\VsqkSDPtworeal}[0]{\textcolor{black}{\mathbf{M}^{2\realss}  }}
\newcommand{\VsqkSDPreal}[0]{\textcolor{black}{\mathbf{M}^{\realss}  }}
\newcommand{\VsqkSDPimag}[0]{\textcolor{black}{\mathbf{M}^{\imagss}  }}
\newcommand{\MisqkSDP}[0]{\textcolor{\complexcolor}{\mathbf{M}_i  }}
\newcommand{\MisqkSDPtworeal}[0]{\textcolor{black}{\mathbf{M}^{2\realss}_{\indexGridNode}  }}

\newcommand{\VbpsqkSDP}[0]{\textcolor{\complexcolor}{\mathbf{M}_{\indexGridNode\indexGridNodeTwo}  }}
\newcommand{\VbpsqkSDPtworeal}[0]{\textcolor{black}{\mathbf{M}^{2\realss}_{\indexGridNode\indexGridNodeTwo}  }}
\newcommand{\VbpsqkSDPreal}[0]{\textcolor{black}{\mathbf{M}_{\indexGridNode\indexGridNodeTwo} ^{\realss}  }}
\newcommand{\VbpsqkSDPimag}[0]{\textcolor{black}{\mathbf{M}_{\indexGridNode\indexGridNodeTwo} ^{\imagss}  }}

\newcommand{\VlinesqkSDP}[0]{\textcolor{\complexcolor}{\mathbf{M}_{\indexGridLines\indexGridNode\indexGridNodeTwo}  }}
\newcommand{\VlinesqkSDPreal}[0]{\textcolor{black}{\mathbf{M}_{\indexGridLines\indexGridNode\indexGridNodeTwo}^{\realss}  }}
\newcommand{\VlinesqkSDPrealtwo}[0]{\textcolor{black}{\mathbf{M}_{\indexGridLines\indexGridNode\indexGridNodeTwo}^{2\realss}  }}

\newcommand{\VlinetosqkSDP}[0]{\textcolor{\complexcolor}{\mathbf{M}_{\indexGridLines\indexGridNodeTwo\indexGridNode}  }}

\newcommand{\VunitsqkSDP}[0]{\textcolor{\complexcolor}{\mathbf{M}_{\indexGridNode\indexUnit}  }}
\newcommand{\VunitsqkSDPvoltage}[0]{\textcolor{\complexcolor}{\mathbf{M}^{\symVoltage}_{\indexGridNode\indexUnit}  }}
\newcommand{\VunitsqkSDPcurrent}[0]{\textcolor{\complexcolor}{\mathbf{M}_{\indexGridNode\indexUnit}  }}
\newcommand{\VunitsqkSDPcurrenttworeal}[0]{\textcolor{black}{\mathbf{M}^{2\realss}_{\indexGridNode\indexUnit}  }}

\newcommand{\VunitsqkSDPcurrentdelta}[0]{\textcolor{\complexcolor}{\mathbf{M}^{\Delta}_{\indexGridNode\indexUnit}  }}
\newcommand{\VunitsqkSDPcurrentdeltap}[0]{\textcolor{\complexcolor}{\mathbf{M}^{\Delta'}_{\indexGridNode\indexUnit}  }}

\newcommand{\VjzsqkSDP}[0]{\textcolor{\complexcolor}{\mathbf{\symVoltageSOCP}_{jz }  }}
\newcommand{\VizsqkSDP}[0]{\textcolor{\complexcolor}{\mathbf{\symVoltageSOCP}_{iz }  }}
\newcommand{\VzjsqkSDP}[0]{\textcolor{\complexcolor}{\mathbf{\symVoltageSOCP}_{zj }  }}
\newcommand{\VzisqkSDP}[0]{\textcolor{\complexcolor}{\mathbf{\symVoltageSOCP}_{zi }  }}

\newcommand{\ViusqkSDP}[0]{\textcolor{\complexcolor}{\mathbf{\symVoltageSOCP}_{\indexGridNode \indexUnit}  }}
\newcommand{\VuisqkSDP}[0]{\textcolor{\complexcolor}{\mathbf{\symVoltageSOCP}_{ \indexUnit \indexGridNode}  }}
\newcommand{\ViusqkSDPreal}[0]{\textcolor{black}{\mathbf{\symVoltageSOCP}_{\indexGridNode \indexUnit}^{\realss}  }}
\newcommand{\VuisqkSDPreal}[0]{\textcolor{black}{\mathbf{\symVoltageSOCP}_{ \indexUnit \indexGridNode}^{\realss}  }}
\newcommand{\ViusqkSDPimag}[0]{\textcolor{black}{\mathbf{\symVoltageSOCP}_{\indexGridNode \indexUnit}^{\imagss}   }}
\newcommand{\VuisqkSDPimag}[0]{\textcolor{black}{\mathbf{\symVoltageSOCP}_{ \indexUnit \indexGridNode}^{\imagss}   }}

\newcommand{\VijsqkSDP}[0]{\textcolor{\complexcolor}{\mathbf{\symVoltageSOCP}_{\indexGridNode \indexGridNodeTwo}  }}
\newcommand{\VjisqkSDP}[0]{\textcolor{\complexcolor}{\mathbf{\symVoltageSOCP}_{\indexGridNodeTwo \indexGridNode}  }}
\newcommand{\VijsqkSDPreal}[0]{\textcolor{black}{\mathbf{\symVoltageSOCP}_{\indexGridNode \indexGridNodeTwo}^{\realss}  }}
\newcommand{\VijsqkSDPimag}[0]{\textcolor{black}{\mathbf{\symVoltageSOCP}_{\indexGridNode \indexGridNodeTwo}^{\imagss}  }}

\newcommand{\ratioij}[0]{\textcolor{\complexcolor}{\mathbf{\symRatio}_{\indexGridLines }  }}
\newcommand{\ratiosqij}[0]{\textcolor{black}{\mathbf{H}_{\indexGridLines }  }}
\newcommand{\ratioijA}[0]{\textcolor{\complexcolor}{{\symRatio}_{\indexGridLines , \symPhaseA}  }}
\newcommand{\ratioijB}[0]{\textcolor{\complexcolor}{{\symRatio}_{\indexGridLines , \symPhaseB}  }}
\newcommand{\ratioijC}[0]{\textcolor{\complexcolor}{{\symRatio}_{\indexGridLines , \symPhaseC}  }}

\newcommand{\ratioijreal}[0]{\textcolor{black}{\mathbf{\symRatio}_{\indexGridLines }^{\realss}  }}
\newcommand{\ratioijimag}[0]{\textcolor{black}{\mathbf{\symRatio}_{\indexGridLines }^{\imagss}  }}
\newcommand{\ratioijmag}[0]{\textcolor{black}{\mathbf{\symRatio}_{\indexGridLines }^{\text{mag}}  }}

\newcommand{\ratiomijA}[0]{\textcolor{black}{{\symRatio}^{\text{mag}}_{\indexGridLines , \symPhaseA}  }}
\newcommand{\ratiomijB}[0]{\textcolor{black}{{\symRatio}^{\text{mag}}_{\indexGridLines , \symPhaseB}  }}
\newcommand{\ratiomijC}[0]{\textcolor{black}{{\symRatio}^{\text{mag}}_{\indexGridLines , \symPhaseC}  }}

\newcommand{\ratiomminijA}[0]{\textcolor{\paramcolor}{{\symRatio}^{\text{min}}_{\indexGridLines , \symPhaseA}  }}
\newcommand{\ratiomminijB}[0]{\textcolor{\paramcolor}{{\symRatio}^{\text{min}}_{\indexGridLines , \symPhaseB}  }}
\newcommand{\ratiomminijC}[0]{\textcolor{\paramcolor}{{\symRatio}^{\text{min}}_{\indexGridLines , \symPhaseC}  }}

\newcommand{\ratiommaxijA}[0]{\textcolor{\paramcolor}{{\symRatio}^{\text{max}}_{\indexGridLines , \symPhaseA}  }}
\newcommand{\ratiommaxijB}[0]{\textcolor{\paramcolor}{{\symRatio}^{\text{max}}_{\indexGridLines , \symPhaseB}  }}
\newcommand{\ratiommaxijC}[0]{\textcolor{\paramcolor}{{\symRatio}^{\text{max}}_{\indexGridLines , \symPhaseC}  }}

\newcommand{\ratioaijA}[0]{\textcolor{black}{{\symRatio}^{\angle}_{\indexGridLines , \symPhaseA}  }}
\newcommand{\ratioaijB}[0]{\textcolor{black}{{\symRatio}^{\angle}_{\indexGridLines , \symPhaseB}  }}
\newcommand{\ratioaijC}[0]{\textcolor{black}{{\symRatio}^{\angle}_{\indexGridLines , \symPhaseC}  }}

\newcommand{\ratioaijAmin}[0]{\textcolor{\paramcolor}{{\symRatio}^{\angle}_{\indexGridLines , \symPhaseA}  }}
\newcommand{\ratioaijBmin}[0]{\textcolor{\paramcolor}{{\symRatio}^{\angle}_{\indexGridLines , \symPhaseB}  }}
\newcommand{\ratioaijCmin}[0]{\textcolor{\paramcolor}{{\symRatio}^{\angle}_{\indexGridLines , \symPhaseC}  }}

\newcommand{\ratioaijAmax}[0]{\textcolor{\paramcolor}{{\symRatio}^{\angle}_{\indexGridLines , \symPhaseA}  }}
\newcommand{\ratioaijBmax}[0]{\textcolor{\paramcolor}{{\symRatio}^{\angle}_{\indexGridLines , \symPhaseB}  }}
\newcommand{\ratioaijCmax}[0]{\textcolor{\paramcolor}{{\symRatio}^{\angle}_{\indexGridLines , \symPhaseC}  }}

\newcommand{\ratioaijref}[0]{\textcolor{\paramcolor}{{\symRatio}^{\angle \text{ref}}_{\indexGridLines\indexGridNode \indexGridNodeTwo}  }}

\newcommand{\ratiomaxij}[0]{\textcolor{\paramcolor}{\mathbf{\symRatio}^{\text{max}}_{\indexGridLines }  }}
\newcommand{\ratiominij}[0]{\textcolor{\paramcolor}{\mathbf{\symRatio}^{\text{min}}_{\indexGridLines }  }}

\newcommand{\ratioamaxij}[0]{\textcolor{\paramcolor}{\mathbf{\symRatio}^{\angle\text{max}}_{\indexGridLines }  }}
\newcommand{\ratioaminij}[0]{\textcolor{\paramcolor}{\mathbf{\symRatio}^{\angle\text{min}}_{\indexGridLines }  }}

\newcommand{\Visqksocp}[0]{\textcolor{black}{\symVoltageSOCP_{\indexGridNode }  }}
\newcommand{\Vjsqksocp}[0]{\textcolor{black}{\symVoltageSOCP_{\indexGridNodeTwo } }}
\newcommand{\Visqksocpacc}[0]{\textcolor{black}{\symVoltageSOCP_{\indexGridNode }^{'}  }}
\newcommand{\Vjsqksocpacc}[0]{\textcolor{black}{\symVoltageSOCP_{\indexGridNodeTwo }^{'} }}
\newcommand{\Visqksocpaccl}[0]{\textcolor{black}{\symVoltageSOCP_{\indexGridNode }^{\star}  }}
\newcommand{\Vjsqksocpaccl}[0]{\textcolor{black}{\symVoltageSOCP_{\indexGridNodeTwo }^{\star} }}
\newcommand{\Visqksocpjabr}[0]{\textcolor{black}{\symVoltageSOCP^{\star}_{\indexGridNode }  }}
\newcommand{\Vjsqksocpjabr}[0]{\textcolor{black}{\symVoltageSOCP^{\star}_{\indexGridNodeTwo } }}
\newcommand{\Visqksocpjabrl}[0]{\textcolor{black}{\symVoltageSOCP^{\star\indexGridLines}_{\indexGridNode }  }}
\newcommand{\Vjsqksocpjabrl}[0]{\textcolor{black}{\symVoltageSOCP^{\star\indexGridLines}_{\indexGridNodeTwo } }}

\newcommand{\Ajabr}[0]{\textcolor{\paramcolor}{A^{'}_{\indexGridNode \indexGridNodeTwo } }}
\newcommand{\Bjabr}[0]{\textcolor{\paramcolor}{B^{'}_{\indexGridNode \indexGridNodeTwo } }}
\newcommand{\Cjabr}[0]{\textcolor{\paramcolor}{C^{'}_{\indexGridNode \indexGridNodeTwo } }}
\newcommand{\Djabr}[0]{\textcolor{\paramcolor}{D^{'}_{\indexGridNode \indexGridNodeTwo } }}

\newcommand{\bigM}[0]{\textcolor{\paramcolor}{M^{} }}

\newcommand{\thetaikref}[0]{\textcolor{\paramcolor}{\symVoltageAngle_{\indexGridNode }^{\refss}  }}
\newcommand{\thetaikrefone}[0]{\textcolor{\paramcolor}{\symVoltageAngle_{\indexGridNode,1 }^{\refss}  }}
\newcommand{\Vikref}[0]{\textcolor{\paramcolor}{\symVoltage_{\indexGridNode }^{\refss}  }}
\newcommand{\Viref}[0]{\textcolor{\paramcolor}{\symVoltage_{\indexGridNode }^{\refss} }}

\newcommand{\Vjkref}[0]{\textcolor{\paramcolor}{\symVoltage_{\indexGridNodeTwo }^{\refss}  }}
\newcommand{\Vjref}[0]{\textcolor{\paramcolor}{\symVoltage_{\indexGridNodeTwo }^{\refss} }}

\newcommand{\thetaik}[0]{\textcolor{black}{\symVoltageAngle_{\indexGridNode }^{}  }}
\newcommand{\thetajk}[0]{\textcolor{black}{\symVoltageAngle_{\indexGridNodeTwo }^{}  }}
\newcommand{\thetaijk}[0]{\textcolor{black}{\symVoltageAngle_{\indexGridNode\indexGridNodeTwo }^{}  }}
\newcommand{\thetaijacckmin}[0]{\textcolor{\paramcolor}{\symVoltageAngle_{\indexGridNode\indexGridNodeTwo }^{'\minss}  }}
\newcommand{\thetaijacckmax}[0]{\textcolor{\paramcolor}{\symVoltageAngle_{\indexGridNode\indexGridNodeTwo }^{'\maxss}  }}
\newcommand{\thetaikmin}[0]{\textcolor{\boundscolor}{\symVoltageAngle_{\indexGridNode }^{\minss}  }}
\newcommand{\thetaikmax}[0]{\textcolor{\boundscolor}{\symVoltageAngle_{\indexGridNode }^{\maxss}  }}

\newcommand{\thetaijkmax}[0]{\textcolor{\boundscolor}{\symVoltageAngle_{\indexGridNode \indexGridNodeTwo}^{\maxss}  }}
\newcommand{\thetaijkmin}[0]{\textcolor{\boundscolor}{\symVoltageAngle_{\indexGridNode \indexGridNodeTwo}^{\minss}  }}

\newcommand{\thetaijkabs}[0]{\textcolor{\boundscolor}{\symVoltageAngle_{\indexGridNode \indexGridNodeTwo}^{\text{absmax}}  }}

\newcommand{\thetaijAAmin}[0]{\textcolor{\paramcolor}{\symVoltageAngle_{\indexGridNode\indexGridNodeTwo,\symPhaseA\symPhaseA }^{\minss}  }}
\newcommand{\thetaijAAmax}[0]{\textcolor{\paramcolor}{\symVoltageAngle_{\indexGridNode\indexGridNodeTwo,\symPhaseA\symPhaseA  }^{\maxss}  }}

\newcommand{\thetaijBBmin}[0]{\textcolor{\paramcolor}{\symVoltageAngle_{\indexGridNode\indexGridNodeTwo,\symPhaseB\symPhaseB }^{\minss}  }}
\newcommand{\thetaijBBmax}[0]{\textcolor{\paramcolor}{\symVoltageAngle_{\indexGridNode\indexGridNodeTwo,\symPhaseB\symPhaseB  }^{\maxss}  }}

\newcommand{\thetaijCCmin}[0]{\textcolor{\paramcolor}{\symVoltageAngle_{\indexGridNode\indexGridNodeTwo,\symPhaseC\symPhaseC }^{\minss}  }}
\newcommand{\thetaijCCmax}[0]{\textcolor{\paramcolor}{\symVoltageAngle_{\indexGridNode\indexGridNodeTwo,\symPhaseC\symPhaseC  }^{\maxss}  }}

\newcommand{\thetaiAAmin}[0]{\textcolor{\paramcolor}{\symVoltageAngle_{\indexGridNode,\symPhaseA\symPhaseA }^{\minss}  }}
\newcommand{\thetaiAAmax}[0]{\textcolor{\paramcolor}{\symVoltageAngle_{\indexGridNode,\symPhaseA\symPhaseA  }^{\maxss}  }}

\newcommand{\thetaiBBmin}[0]{\textcolor{\paramcolor}{\symVoltageAngle_{\indexGridNode,\symPhaseB\symPhaseB }^{\minss}  }}
\newcommand{\thetaiBBmax}[0]{\textcolor{\paramcolor}{\symVoltageAngle_{\indexGridNode,\symPhaseB\symPhaseB  }^{\maxss}  }}

\newcommand{\thetaiCCmin}[0]{\textcolor{\paramcolor}{\symVoltageAngle_{\indexGridNode,\symPhaseC\symPhaseC }^{\minss}  }}
\newcommand{\thetaiCCmax}[0]{\textcolor{\paramcolor}{\symVoltageAngle_{\indexGridNode,\symPhaseC\symPhaseC  }^{\maxss}  }}

\newcommand{\thetaijkmaxacc}[0]{\textcolor{\boundscolor}{\symVoltageAngle_{\indexGridNode \indexGridNodeTwo}^{'\maxss}  }}
\newcommand{\thetaijkminacc}[0]{\textcolor{\boundscolor}{\symVoltageAngle_{\indexGridNode \indexGridNodeTwo}^{'\minss}  }}

\newcommand{\thetaijkabsacc}[0]{\textcolor{\boundscolor}{\symVoltageAngle_{\indexGridNode \indexGridNodeTwo}^{'\text{absmax}}  }}

\newcommand{\phiijk}[0]{\textcolor{black}{\symPhaseDifference_{\indexGridNode\indexGridNodeTwo }^{}  }}
\newcommand{\phiikmax}[0]{\textcolor{\boundscolor}{\symPhaseDifference_{\indexGridNode }^{\maxss}  }}
\newcommand{\phiikmin}[0]{\textcolor{\boundscolor}{\symPhaseDifference_{\indexGridNode }^{\minss}  }}
\newcommand{\phiijkmin}[0]{\textcolor{\boundscolor}{\symPhaseDifference_{\indexGridNode   \indexGridNodeTwo }^{\minss}  }}
\newcommand{\phiijkmax}[0]{\textcolor{\boundscolor}{\symPhaseDifference_{\indexGridNode   \indexGridNodeTwo }^{\maxss}  }}

\newcommand{\phijik}[0]{\textcolor{black}{\symPhaseDifference_{\indexGridNodeTwo   \indexGridNode }^{}  }}
\newcommand{\phijikmin}[0]{\textcolor{\boundscolor}{\symPhaseDifference_{\indexGridNodeTwo   \indexGridNode    }^{\minss}  }}
\newcommand{\phijikmax}[0]{\textcolor{\boundscolor}{\symPhaseDifference_{\indexGridNodeTwo  \indexGridNode    }^{\maxss}  }}

\newcommand{\thetaiacck}[0]{\textcolor{black}{\symVoltageAngle_{\indexGridNode }^{'}  }}
\newcommand{\thetajacck}[0]{\textcolor{black}{\symVoltageAngle_{\indexGridNodeTwo }^{'}  }}
\newcommand{\thetaijacck}[0]{\textcolor{black}{\symVoltageAngle_{\indexGridNode  \indexGridNodeTwo }^{'}  }}

\newcommand{\thetaiacckl}[0]{\textcolor{black}{\symVoltageAngle_{\indexGridNode }^{\star}  }}
\newcommand{\thetajacckl}[0]{\textcolor{black}{\symVoltageAngle_{\indexGridNodeTwo }^{\star}  }}

\newcommand{\yijkshSDPseq}[0]{\textcolor{\complexparamcolor}{\mathbf{\symAdmittance}_{\indexGridLines \indexGridNode   \indexGridNodeTwo}^{\shuntss, \fortescuess}    }}

\newcommand{\yijkshSDPAdot}[0]{\textcolor{\complexparamcolor}{\mathbf{\symAdmittance}_{\indexGridLines \indexGridNode   \indexGridNodeTwo}^{\shuntss, \symPhaseA\!\cdot}    }}
\newcommand{\yijkshSDPBdot}[0]{\textcolor{\complexparamcolor}{\mathbf{\symAdmittance}_{\indexGridLines \indexGridNode   \indexGridNodeTwo}^{\shuntss, \symPhaseB\!\cdot}    }}
\newcommand{\yijkshSDPCdot}[0]{\textcolor{\complexparamcolor}{\mathbf{\symAdmittance}_{\indexGridLines \indexGridNode   \indexGridNodeTwo}^{\shuntss, \symPhaseC\!\cdot}    }}
\newcommand{\yijkshSDPNdot}[0]{\textcolor{\complexparamcolor}{\mathbf{\symAdmittance}_{\indexGridLines \indexGridNode   \indexGridNodeTwo}^{\shuntss, \symPhaseN\!\cdot}    }}

\newcommand{\yijkshSDPdotN}[0]{\textcolor{\complexparamcolor}{\mathbf{\symAdmittance}_{\indexGridLines \indexGridNode   \indexGridNodeTwo}^{\shuntss, \cdot \! \symPhaseN}    }}

\newcommand{\yijkshSDP}[0]{\textcolor{\complexparamcolor}{\mathbf{\symAdmittance}_{\indexGridLines \indexGridNode   \indexGridNodeTwo}^{\shuntss}    }}
\newcommand{\gijkshSDP}[0]{\textcolor{\paramcolor}{\mathbf{\symConductance}_{ \indexGridLines\indexGridNode   \indexGridNodeTwo}^{\shuntss}    }}
\newcommand{\bijkshSDP}[0]{\textcolor{\paramcolor}{\mathbf{\symSusceptance}_{ \indexGridLines\indexGridNode   \indexGridNodeTwo}^{\shuntss}    }}
\newcommand{\zijkshSDP}[0]{\textcolor{\complexparamcolor}{\mathbf{\symImpedance}_{\indexGridLines \indexGridNode   \indexGridNodeTwo}^{\shuntss}    }}

\newcommand{\yuunitSDP}[0]{\textcolor{\complexparamcolor}{\mathbf{\symAdmittance}_{    \indexUnit}  }}
\newcommand{\yuunitSDPDelta}[0]{\textcolor{\complexparamcolor}{\mathbf{\symAdmittance}_{    \indexUnit}^{\Delta}  }}

\newcommand{\ybSDP}[0]{\textcolor{\complexparamcolor}{\mathbf{\symAdmittance}_{    \indexShunt}  }}
\newcommand{\gbSDP}[0]{\textcolor{\paramcolor}{\mathbf{\symConductance}_{    \indexShunt}  }}
\newcommand{\bbSDP}[0]{\textcolor{\paramcolor}{\mathbf{\symSusceptance}_{    \indexShunt}  }}

\newcommand{\yjikshSDP}[0]{\textcolor{\complexparamcolor}{\mathbf{\symAdmittance}_{ \indexGridLines\indexGridNodeTwo   \indexGridNode}^{\shuntss}    }}
\newcommand{\gjikshSDP}[0]{\textcolor{\paramcolor}{\mathbf{\symConductance}_{ \indexGridLines\indexGridNodeTwo   \indexGridNode}^{\shuntss}    }}
\newcommand{\bjikshSDP}[0]{\textcolor{\paramcolor}{\mathbf{\symSusceptance}_{ \indexGridLines\indexGridNodeTwo   \indexGridNode}^{\shuntss}    }}
\newcommand{\zjikshSDP}[0]{\textcolor{\complexparamcolor}{\mathbf{\symImpedance}_{\indexGridLines \indexGridNodeTwo   \indexGridNode}^{\shuntss}    }}

\newcommand{\yikSDP}[0]{\textcolor{\complexparamcolor}{\mathbf{\symAdmittance}_{ \indexGridNode   }^{}    }}
\newcommand{\gikSDP}[0]{\textcolor{\paramcolor}{\mathbf{\symConductance}_{ \indexGridNode   }^{}    }}
\newcommand{\bikSDP}[0]{\textcolor{\paramcolor}{\mathbf{\symSusceptance}_{ \indexGridNode   }^{}    }}

\newcommand{\zijksSDPkron}[0]{\textcolor{\complexparamcolor}{\mathbf{{\symImpedance}}_{ \indexGridLines}^{\seriesss, \text{Kron}}    }}
\newcommand{\zijksSDPphase}[0]{\textcolor{\complexparamcolor}{\mathbf{{\symImpedance}}_{ \indexGridLines}^{\seriesss, \text{phase}}    }}
\newcommand{\zijksSDPprim}[0]{\textcolor{\complexparamcolor}{\mathbf{\hat{\symImpedance}}_{ \indexGridLines}^{\seriesss, \text{circ}}    }}
\newcommand{\zijksSDPvec}[0]{\textcolor{\complexparamcolor}{\mathbf{\hat{\symImpedance}}_{ \indexGridLines}^{\seriesss, \text{pn}}    }}
\newcommand{\zijksSDPvect}[0]{\textcolor{\complexparamcolor}{\mathbf{\hat{\symImpedance}}_{ \indexGridLines}^{\seriesss,  \text{np}}    }}
\newcommand{\zijksSDPhat}[0]{\textcolor{\complexparamcolor}{\mathbf{\hat{\symImpedance}}_{ \indexGridLines}^{\seriesss,  \text{pp}}    }}

\newcommand{\zijksSDPcond}[0]{\textcolor{\complexparamcolor}{\mathbf{\bar{\symImpedance}}_{ \indexGridLines}^{\seriesss, \text{cond}}    }}
\newcommand{\zijksSDPseq}[0]{\textcolor{\complexparamcolor}{\mathbf{\symImpedance}_{ \indexGridLines}^{\seriesss, \fortescuess}    }}

\newcommand{\zijksSDPseqdiag}[0]{\textcolor{\complexparamcolor}{\mathbf{\symImpedance}_{ \indexGridLines}^{\seriesss, \fortescuess, \text{diag}}    }}

\newcommand{\thetaijSDPmax}[0]{\textcolor{\paramcolor}{\mathbf{\Theta}_{\indexGridNode\indexGridNodeTwo  }^{\maxss}  }}
\newcommand{\thetaijSDPmin}[0]{\textcolor{\paramcolor}{\mathbf{\Theta}_{\indexGridNode\indexGridNodeTwo  }^{\minss}  }}

\newcommand{\thetaiSDPmax}[0]{\textcolor{\paramcolor}{\mathbf{\Theta}_{\indexGridNode  }^{\maxss}  }}
\newcommand{\thetaiSDPmin}[0]{\textcolor{\paramcolor}{\mathbf{\Theta}_{\indexGridNode  }^{\minss}  }}

\newcommand{\seqalpha}[0]{\textcolor{\complexparamcolor}{\alpha}}
\newcommand{\seq}[0]{\textcolor{\complexparamcolor}{\mathbf{A}}}
\newcommand{\seqinv}[0]{\textcolor{\complexparamcolor}{\mathbf{F}}}
\newcommand{\seqH}[0]{\textcolor{\complexparamcolor}{\mathbf{A}^{\hermitiantranspose}}}
\newcommand{\seqconj}[0]{\textcolor{\complexparamcolor}{\mathbf{A}^{*}}}
\newcommand{\seqT}[0]{\textcolor{\complexparamcolor}{\mathbf{A}^{\transpose}}}

\newcommand{\seqvec}[0]{\textcolor{\complexparamcolor}{\mathbf{a}}}
\newcommand{\seqreal}[0]{\textcolor{\paramcolor}{\mathbf{A}^{\realss}}}
\newcommand{\seqimag}[0]{\textcolor{\paramcolor}{\mathbf{A}^{\imagss}}}
\newcommand{\seqinvreal}[0]{\textcolor{\paramcolor}{\mathbf{F}^{\realss}}}
\newcommand{\seqinvimag}[0]{\textcolor{\paramcolor}{\mathbf{F}^{\imagss}}}

\newcommand{\zslh}[0]{\textcolor{\complexparamcolor}{\mathbf{\symImpedance}_{ \indexGridLines, \symHarmonic}^{\seriesss}    }}
\newcommand{\rslh}[0]{\textcolor{\paramcolor}{\mathbf{\symResistance}_{  \indexGridLines, \symHarmonic   }^{\seriesss}    }}
\newcommand{\xslh}[0]{\textcolor{\paramcolor}{\mathbf{\symReactance}_{  \indexGridLines, \symHarmonic   }^{\seriesss}    }}

\newcommand{\yshlijh}[0]{\textcolor{\complexparamcolor}{\mathbf{\symAdmittance}_{ \indexGridLines\indexGridNode\indexGridNodeTwo, \symHarmonic}^{\shuntss}    }}
\newcommand{\gshlijh}[0]{\textcolor{\paramcolor}{\mathbf{\symConductance}_{  \indexGridLines\indexGridNode\indexGridNodeTwo, \symHarmonic   }^{\shuntss}    }}
\newcommand{\bshlijh}[0]{\textcolor{\paramcolor}{\mathbf{\symSusceptance}_{  \indexGridLines\indexGridNode\indexGridNodeTwo, \symHarmonic   }^{\shuntss}    }}

\newcommand{\yshljih}[0]{\textcolor{\complexparamcolor}{\mathbf{\symAdmittance}_{ \indexGridLines\indexGridNodeTwo\indexGridNode, \symHarmonic}^{\shuntss}    }}
\newcommand{\gshljih}[0]{\textcolor{\paramcolor}{\mathbf{\symConductance}_{  \indexGridLines\indexGridNodeTwo\indexGridNode, \symHarmonic   }^{\shuntss}    }}
\newcommand{\bshljih}[0]{\textcolor{\paramcolor}{\mathbf{\symSusceptance}_{  \indexGridLines\indexGridNodeTwo\indexGridNode, \symHarmonic   }^{\shuntss}    }}

\newcommand{\zijksSDP}[0]{\textcolor{\complexparamcolor}{\mathbf{\symImpedance}_{ \indexGridLines}^{\seriesss}    }}
\newcommand{\rijksSDP}[0]{\textcolor{\paramcolor}{\mathbf{\symResistance}_{  \indexGridLines   }^{\seriesss}    }}
\newcommand{\xijksSDP}[0]{\textcolor{\paramcolor}{\mathbf{\symReactance}_{  \indexGridLines   }^{\seriesss}    }}

\newcommand{\zjiksSDP}[0]{\textcolor{\complexparamcolor}{\mathbf{\symImpedance}_{ \indexGridLines   }^{\seriesss}    }}
\newcommand{\rjiksSDP}[0]{\textcolor{\paramcolor}{\mathbf{\symResistance}_{ \indexGridLines   }^{\seriesss}    }}
\newcommand{\xjiksSDP}[0]{\textcolor{\paramcolor}{\mathbf{\symReactance}_{ \indexGridLines   }^{\seriesss}    }}

\newcommand{\zijks}[0]{\textcolor{\complexparamcolor}{\symImpedance_{ \indexGridLines,\seriesss}    }}
\newcommand{\rijks}[0]{\textcolor{\paramcolor}{\symResistance_{  \indexGridLines,\seriesss}    }}
\newcommand{\xijks}[0]{\textcolor{\paramcolor}{\symReactance_{  \indexGridLines,\seriesss}    }}
\newcommand{\yijks}[0]{\textcolor{\complexparamcolor}{\symAdmittance_{  \indexGridLines}^{\seriesss}    }}
\newcommand{\bijks}[0]{\textcolor{\paramcolor}{\symSusceptance_{  \indexGridLines}^{\seriesss}    }}
\newcommand{\gijks}[0]{\textcolor{\paramcolor}{\symConductance_{  \indexGridLines}^{\seriesss}    }}

\newcommand{\zijksSDPH}[0]{\textcolor{\complexparamcolor}{(\mathbf{\symImpedance}_{ \indexGridLines}^{\seriesss})^{\hermitiantranspose}    }}
\newcommand{\rijksSDPH}[0]{\textcolor{\paramcolor}{(\mathbf{\symResistance}_{ \indexGridLines,\seriesss})^{\transpose}    }}
\newcommand{\xijksSDPH}[0]{\textcolor{\paramcolor}{(\mathbf{\symReactance}_{ \indexGridLines,\seriesss})^{\transpose}    }}

\newcommand{\yukSDP}[0]{\textcolor{\complexparamcolor}{\mathbf{\symAdmittance}_{ \indexUnit}    }}
\newcommand{\gukSDP}[0]{\textcolor{\paramcolor}{\mathbf{\symConductance}_{ \indexUnit}    }}
\newcommand{\bukSDP}[0]{\textcolor{\paramcolor}{\mathbf{\symSusceptance}_{ \indexUnit}    }}

\newcommand{\zukSDP}[0]{\textcolor{\complexparamcolor}{\mathbf{\symImpedance}_{ \indexUnit}    }}
\newcommand{\rukSDP}[0]{\textcolor{\paramcolor}{\mathbf{\symResistance}_{ \indexUnit}    }}
\newcommand{\xukSDP}[0]{\textcolor{\paramcolor}{\mathbf{\symReactance}_{ \indexUnit}    }}

\newcommand{\yijksSDP}[0]{\textcolor{\complexparamcolor}{\mathbf{\symAdmittance}_{ \indexGridLines}^{\seriesss}    }}
\newcommand{\gijksSDP}[0]{\textcolor{\paramcolor}{\mathbf{\symConductance}_{ \indexGridLines}^{\seriesss}    }}
\newcommand{\bijksSDP}[0]{\textcolor{\paramcolor}{\mathbf{\symSusceptance}_{ \indexGridLines}^{\seriesss}    }}

\newcommand{\yjiksSDP}[0]{\textcolor{\complexparamcolor}{\mathbf{\symAdmittance}_{  \indexGridLines\indexGridNodeTwo \indexGridNode }^{\seriesss}    }}
\newcommand{\gjiksSDP}[0]{\textcolor{\paramcolor}{\mathbf{\symConductance}_{ \indexGridLines \indexGridNodeTwo \indexGridNode }^{\seriesss}    }}
\newcommand{\bjiksSDP}[0]{\textcolor{\paramcolor}{\mathbf{\symSusceptance}_{  \indexGridLines\indexGridNodeTwo \indexGridNode }^{\seriesss}    }}

\newcommand{\zijksh}[0]{\textcolor{\complexparamcolor}{\symImpedance_{ \indexGridNode \indexGridNodeTwo,\shuntss}    }}
\newcommand{\rijksh}[0]{\textcolor{\paramcolor}{\symResistance_{ \indexGridNode \indexGridNodeTwo,\shuntss}    }}
\newcommand{\xijksh}[0]{\textcolor{\paramcolor}{\symReactance_{ \indexGridNode \indexGridNodeTwo,\shuntss}    }}
\newcommand{\yijksh}[0]{\textcolor{\complexparamcolor}{\symAdmittance_{ \indexGridNode \indexGridNodeTwo,\shuntss}    }}
\newcommand{\gijksh}[0]{\textcolor{\paramcolor}{\symConductance_{ \indexGridNode \indexGridNodeTwo,\shuntss}    }}
\newcommand{\bijksh}[0]{\textcolor{\paramcolor}{\symSusceptance_{ \indexGridNode \indexGridNodeTwo,\shuntss}    }}

\newcommand{\zjiksh}[0]{\textcolor{\complexparamcolor}{\symImpedance_{ \indexGridNodeTwo \indexGridNode ,\shuntss}    }}
\newcommand{\rjiksh}[0]{\textcolor{\paramcolor}{\symResistance_{\indexGridNodeTwo \indexGridNode ,\shuntss}    }}
\newcommand{\xjiksh}[0]{\textcolor{\paramcolor}{\symReactance_{ \indexGridNodeTwo\indexGridNode ,\shuntss}    }}
\newcommand{\yjiksh}[0]{\textcolor{\complexparamcolor}{\symAdmittance_{\indexGridNodeTwo \indexGridNode ,\shuntss}    }}
\newcommand{\gjiksh}[0]{\textcolor{\paramcolor}{\symConductance_{\indexGridNodeTwo \indexGridNode ,\shuntss}    }}
\newcommand{\bjiksh}[0]{\textcolor{\paramcolor}{\symSusceptance_{ \indexGridNodeTwo \indexGridNode ,\shuntss}    }}

\newcommand{\zlksh}[0]{\textcolor{\complexparamcolor}{\symImpedance_{  \indexGridLines,\shuntss}    }}
\newcommand{\rlksh}[0]{\textcolor{\paramcolor}{\symResistance_{  \indexGridLines,\shuntss}    }}
\newcommand{\xlksh}[0]{\textcolor{\paramcolor}{\symReactance_{  \indexGridLines,\shuntss}    }}
\newcommand{\ylksh}[0]{\textcolor{\complexparamcolor}{\symAdmittance_{  \indexGridLines,\shuntss}    }}
\newcommand{\glksh}[0]{\textcolor{\paramcolor}{\symConductance_{  \indexGridLines,\shuntss}    }}
\newcommand{\blksh}[0]{\textcolor{\paramcolor}{\symSusceptance_{  \indexGridLines,\shuntss}    }}

\newcommand{\zlkshspec}[0]{\textcolor{\complexparamcolor}{\dot{\symImpedance}_{  \indexGridLines,\shuntss}    }}
\newcommand{\rlkshspec}[0]{\textcolor{\paramcolor}{\dot{\symResistance}_{  \indexGridLines,\shuntss}    }}
\newcommand{\xlkshv}[0]{\textcolor{\paramcolor}{\dot{\symReactance}_{  \indexGridLines,\shuntss}    }}
\newcommand{\ylkshspec}[0]{\textcolor{\complexparamcolor}{\dot{\symAdmittance}_{  \indexGridLines,\shuntss}    }}
\newcommand{\glkshspec}[0]{\textcolor{\paramcolor}{\dot{\symConductance}_{  \indexGridLines,\shuntss}    }}
\newcommand{\blkshspec}[0]{\textcolor{\paramcolor}{\dot{\symSusceptance}_{  \indexGridLines,\shuntss}    }}

\newcommand{\zijksspec}[0]{\textcolor{\complexparamcolor}{\dot{\symImpedance}_{  \indexGridLines,\seriesss}    }}
\newcommand{\rijksspec}[0]{\textcolor{\paramcolor}{\dot{\symResistance}_{  \indexGridLines,\seriesss}    }}
\newcommand{\xijksspec}[0]{\textcolor{\paramcolor}{\dot{\symReactance}_{  \indexGridLines,\seriesss}    }}
\newcommand{\yijkshspec}[0]{\textcolor{\complexparamcolor}{\dot{\symAdmittance}_{  \indexGridNode \indexGridNodeTwo,\shuntss}    }}
\newcommand{\bijkshspec}[0]{\textcolor{\paramcolor}{\dot{\symSusceptance}_{  \indexGridNode \indexGridNodeTwo,\shuntss}    }}
\newcommand{\gijkshspec}[0]{\textcolor{\paramcolor}{\dot{\symConductance}_{  \indexGridNode \indexGridNodeTwo,\shuntss}    }}

\newcommand{\lijk}[0]{\textcolor{\paramcolor}{{\symLength}_{  \indexGridLines}    }}
\newcommand{\vprimvsecijk}[0]{\textcolor{black}{{\symRatio}_{ \indexGridNode \indexGridNodeTwo}    }}
\newcommand{\vprimvsecijkmin}[0]{\textcolor{\boundscolor}{{\symRatio}_{ \indexGridNode \indexGridNodeTwo}^{\minss}    }}
\newcommand{\vprimvsecijkmax}[0]{\textcolor{\boundscolor}{{\symRatio}_{ \indexGridNode \indexGridNodeTwo}^{\maxss}    }}

\newcommand{\vsecvprimijk}[0]{\textcolor{black}{{\symRatio}_{\indexGridNodeTwo \indexGridNode }    }}
\newcommand{\vsecvprimijkmin}[0]{\textcolor{\boundscolor}{{\symRatio}_{ \indexGridNodeTwo\indexGridNode }^{\minss}    }}
\newcommand{\vsecvprimijkmax}[0]{\textcolor{\boundscolor}{{\symRatio}_{\indexGridNodeTwo \indexGridNode }^{\maxss}    }}

\newcommand{\YPFconvexification}[0]{\textcolor{black}{{\symPenalty}^{\PFconvexss}   }}
\newcommand{\YGIconvexification}[0]{\textcolor{black}{{\symPenalty}^{\GIconvexss}   }}
\newcommand{\YPgridloss}[0]{\textcolor{black}{{\symPenalty}^{\symPower,\lossss}   }}
\newcommand{\YSgridloss}[0]{\textcolor{black}{{\symPenalty}^{\symApparentPower,\lossss}   }}
\newcommand{\Ytotcurrent}[0]{\textcolor{black}{{\symPenalty}^{\text{\symCurrent}}   }}

\newcommand{\Yobj}[0]{\textcolor{black}{{\symPenalty}^{\text{cost}}   }}
\newcommand{\Kobj}[0]{\textcolor{black}{{\symAnnualCost}^{\text{cost}}   }}

\newcommand{\NPFconvexification}[0]{\textcolor{\paramcolor}{{\symPenaltyWeight}_{\PFconvexss}   }}
\newcommand{\NGIconvexification}[0]{\textcolor{\paramcolor}{{\symPenaltyWeight}_{\GIconvexss}   }}

\newcommand{\wPFconvexification}[0]{\textcolor{\paramcolor}{{\symObjectiveWeight}_{\PFconvexss}   }}
\newcommand{\wPFconvexificationopt}[0]{\textcolor{\paramcolor}{{\symObjectiveWeight}_{\PFconvexss}^{\text{opt}}   }}

\newcommand{\zijkszero}[0]{\textcolor{\complexparamcolor}{\symImpedance_{ \indexGridLines , 0}^{\seriesss}    }}
\newcommand{\zijksone}[0]{\textcolor{\complexparamcolor}{\symImpedance_{ \indexGridLines , 1}^{\seriesss}    }}

\newcommand{\zijkssymmzerozero}[0]{\textcolor{\complexparamcolor}{\symImpedance_{ \indexGridLines , 0 0}^{\seriesss}    }}
\newcommand{\zijkssymmzeroone}[0]{\textcolor{\complexparamcolor}{\symImpedance_{ \indexGridLines , 0 1}^{\seriesss}    }}
\newcommand{\zijkssymmzerotwo}[0]{\textcolor{\complexparamcolor}{\symImpedance_{ \indexGridLines , 0 2}^{\seriesss}    }}

\newcommand{\zijkssymmonezero}[0]{\textcolor{\complexparamcolor}{\symImpedance_{ \indexGridLines , 1 0}^{\seriesss}    }}
\newcommand{\zijkssymmoneone}[0]{\textcolor{\complexparamcolor}{\symImpedance_{ \indexGridLines , 1 1}^{\seriesss}    }}
\newcommand{\zijkssymmonetwo}[0]{\textcolor{\complexparamcolor}{\symImpedance_{ \indexGridLines , 1 2}^{\seriesss}    }}

\newcommand{\zijkssymmtwozero}[0]{\textcolor{\complexparamcolor}{\symImpedance_{ \indexGridLines , 2 0}^{\seriesss}    }}
\newcommand{\zijkssymmtwoone}[0]{\textcolor{\complexparamcolor}{\symImpedance_{ \indexGridLines , 2 1}^{\seriesss}    }}
\newcommand{\zijkssymmtwotwo}[0]{\textcolor{\complexparamcolor}{\symImpedance_{ \indexGridLines , 2 2}^{\seriesss}    }}

\newcommand{\zjNG}[0]{\textcolor{\complexparamcolor}{\symImpedance_{ \indexGridNodeTwo , \symPhaseN \symPhaseG}   }}
\newcommand{\ziNG}[0]{\textcolor{\complexparamcolor}{\symImpedance_{ \indexGridNode , \symPhaseN \symPhaseG}   }}
\newcommand{\riNG}[0]{\textcolor{\paramcolor}{\symResistance_{ \indexGridNode , \symPhaseN \symPhaseG}   }}
\newcommand{\xiNG}[0]{\textcolor{\paramcolor}{\symReactance_{ \indexGridNode , \symPhaseN \symPhaseG}   }}

\newcommand{\zijksAA}[0]{\textcolor{\complexparamcolor}{\symImpedance_{ \indexGridLines , \symPhaseA \symPhaseA}^{\seriesss}    }}
\newcommand{\zijksAB}[0]{\textcolor{\complexparamcolor}{\symImpedance_{ \indexGridLines , \symPhaseA \symPhaseB}^{\seriesss}    }}
\newcommand{\zijksAC}[0]{\textcolor{\complexparamcolor}{\symImpedance_{ \indexGridLines , \symPhaseA \symPhaseC}^{\seriesss}    }}
\newcommand{\zijksBA}[0]{\textcolor{\complexparamcolor}{\symImpedance_{ \indexGridLines , \symPhaseB \symPhaseA}^{\seriesss}    }}
\newcommand{\zijksBB}[0]{\textcolor{\complexparamcolor}{\symImpedance_{ \indexGridLines , \symPhaseB \symPhaseB}^{\seriesss}    }}
\newcommand{\zijksBC}[0]{\textcolor{\complexparamcolor}{\symImpedance_{ \indexGridLines , \symPhaseB \symPhaseC}^{\seriesss}    }}
\newcommand{\zijksCA}[0]{\textcolor{\complexparamcolor}{\symImpedance_{ \indexGridLines , \symPhaseC \symPhaseA}^{\seriesss}    }}
\newcommand{\zijksCB}[0]{\textcolor{\complexparamcolor}{\symImpedance_{ \indexGridLines , \symPhaseC \symPhaseB}^{\seriesss}    }}
\newcommand{\zijksCC}[0]{\textcolor{\complexparamcolor}{\symImpedance_{ \indexGridLines , \symPhaseC \symPhaseC}^{\seriesss}    }}

\newcommand{\zijksAN}[0]{\textcolor{\complexparamcolor}{\symImpedance_{ \indexGridLines , \symPhaseA \symPhaseN}^{\seriesss}    }}
\newcommand{\zijksBN}[0]{\textcolor{\complexparamcolor}{\symImpedance_{ \indexGridLines , \symPhaseB \symPhaseN}^{\seriesss}    }}
\newcommand{\zijksCN}[0]{\textcolor{\complexparamcolor}{\symImpedance_{ \indexGridLines , \symPhaseC \symPhaseN}^{\seriesss}    }}
\newcommand{\zijksNA}[0]{\textcolor{\complexparamcolor}{\symImpedance_{ \indexGridLines , \symPhaseN \symPhaseA }^{\seriesss}    }}
\newcommand{\zijksNB}[0]{\textcolor{\complexparamcolor}{\symImpedance_{ \indexGridLines ,\symPhaseN \symPhaseB }^{\seriesss}    }}
\newcommand{\zijksNC}[0]{\textcolor{\complexparamcolor}{\symImpedance_{ \indexGridLines , \symPhaseN\symPhaseC }^{\seriesss}    }}

\newcommand{\zijksNN}[0]{\textcolor{\complexparamcolor}{\symImpedance_{ \indexGridLines , \symPhaseN \symPhaseN}^{\seriesss}    }}

\newcommand{\zijksAG}[0]{\textcolor{\complexparamcolor}{\symImpedance_{ \indexGridLines , \symPhaseA \symPhaseG}^{\seriesss}    }}
\newcommand{\zijksBG}[0]{\textcolor{\complexparamcolor}{\symImpedance_{ \indexGridLines , \symPhaseB \symPhaseG}^{\seriesss}    }}
\newcommand{\zijksCG}[0]{\textcolor{\complexparamcolor}{\symImpedance_{ \indexGridLines , \symPhaseC \symPhaseG}^{\seriesss}    }}
\newcommand{\zijksNG}[0]{\textcolor{\complexparamcolor}{\symImpedance_{ \indexGridLines , \symPhaseN \symPhaseG}^{\seriesss}    }}
\newcommand{\zijksGA}[0]{\textcolor{\complexparamcolor}{\symImpedance_{ \indexGridLines , \symPhaseG \symPhaseA }^{\seriesss}    }}
\newcommand{\zijksGB}[0]{\textcolor{\complexparamcolor}{\symImpedance_{ \indexGridLines ,\symPhaseG \symPhaseB }^{\seriesss}    }}
\newcommand{\zijksGC}[0]{\textcolor{\complexparamcolor}{\symImpedance_{ \indexGridLines , \symPhaseG\symPhaseC }^{\seriesss}    }}
\newcommand{\zijksGN}[0]{\textcolor{\complexparamcolor}{\symImpedance_{ \indexGridLines , \symPhaseG \symPhaseN}^{\seriesss}    }}

\newcommand{\zijksGG}[0]{\textcolor{\complexparamcolor}{\symImpedance_{ \indexGridLines , \symPhaseG \symPhaseG}^{\seriesss}    }}

\newcommand{\zijksNNhat}[0]{\textcolor{\complexparamcolor}{\hat{\symImpedance}_{ \indexGridLines , \symPhaseN \symPhaseN}^{\seriesss}    }}
\newcommand{\zijksANhat}[0]{\textcolor{\complexparamcolor}{\hat{\symImpedance}_{ \indexGridLines , \symPhaseA \symPhaseN}^{\seriesss}    }}
\newcommand{\zijksBNhat}[0]{\textcolor{\complexparamcolor}{\hat{\symImpedance}_{ \indexGridLines , \symPhaseB \symPhaseN}^{\seriesss}    }}
\newcommand{\zijksCNhat}[0]{\textcolor{\complexparamcolor}{\hat{\symImpedance}_{ \indexGridLines , \symPhaseC \symPhaseN}^{\seriesss}    }}

\newcommand{\zijksAAhat}[0]{\textcolor{\complexparamcolor}{\hat{\symImpedance}_{ \indexGridLines , \symPhaseA \symPhaseA}^{\seriesss}    }}
\newcommand{\zijksBAhat}[0]{\textcolor{\complexparamcolor}{\hat{\symImpedance}_{ \indexGridLines , \symPhaseB \symPhaseA}^{\seriesss}    }}
\newcommand{\zijksCAhat}[0]{\textcolor{\complexparamcolor}{\hat{\symImpedance}_{ \indexGridLines , \symPhaseC \symPhaseA}^{\seriesss}    }}

\newcommand{\zijksABhat}[0]{\textcolor{\complexparamcolor}{\hat{\symImpedance}_{ \indexGridLines , \symPhaseA \symPhaseB}^{\seriesss}    }}
\newcommand{\zijksBBhat}[0]{\textcolor{\complexparamcolor}{\hat{\symImpedance}_{ \indexGridLines , \symPhaseB \symPhaseB}^{\seriesss}    }}
\newcommand{\zijksCBhat}[0]{\textcolor{\complexparamcolor}{\hat{\symImpedance}_{ \indexGridLines , \symPhaseC \symPhaseB}^{\seriesss}    }}

\newcommand{\zijksAChat}[0]{\textcolor{\complexparamcolor}{\hat{\symImpedance}_{ \indexGridLines , \symPhaseA \symPhaseC}^{\seriesss}    }}
\newcommand{\zijksBChat}[0]{\textcolor{\complexparamcolor}{\hat{\symImpedance}_{ \indexGridLines , \symPhaseB \symPhaseC}^{\seriesss}    }}
\newcommand{\zijksCChat}[0]{\textcolor{\complexparamcolor}{\hat{\symImpedance}_{ \indexGridLines , \symPhaseC \symPhaseC}^{\seriesss}    }}

\newcommand{\zijksNAhat}[0]{\textcolor{\complexparamcolor}{\hat{\symImpedance}_{ \indexGridLines , \symPhaseN \symPhaseA }^{\seriesss}    }}
\newcommand{\zijksNBhat}[0]{\textcolor{\complexparamcolor}{\hat{\symImpedance}_{ \indexGridLines , \symPhaseN \symPhaseB }^{\seriesss}    }}
\newcommand{\zijksNChat}[0]{\textcolor{\complexparamcolor}{\hat{\symImpedance}_{ \indexGridLines , \symPhaseN \symPhaseC }^{\seriesss}    }}

\newcommand{\yijksAA}[0]{\textcolor{\complexparamcolor}{\symAdmittance_{ \indexGridLines , \symPhaseA \symPhaseA}^{\seriesss}    }}
\newcommand{\yijksAB}[0]{\textcolor{\complexparamcolor}{\symAdmittance_{ \indexGridLines , \symPhaseA \symPhaseB}^{\seriesss}    }}
\newcommand{\yijksAC}[0]{\textcolor{\complexparamcolor}{\symAdmittance_{ \indexGridLines , \symPhaseA \symPhaseC}^{\seriesss}    }}
\newcommand{\yijksBA}[0]{\textcolor{\complexparamcolor}{\symAdmittance_{ \indexGridLines , \symPhaseB \symPhaseA}^{\seriesss}    }}
\newcommand{\yijksBB}[0]{\textcolor{\complexparamcolor}{\symAdmittance_{ \indexGridLines , \symPhaseB \symPhaseB}^{\seriesss}    }}
\newcommand{\yijksBC}[0]{\textcolor{\complexparamcolor}{\symAdmittance_{ \indexGridLines , \symPhaseB \symPhaseC}^{\seriesss}    }}
\newcommand{\yijksCA}[0]{\textcolor{\complexparamcolor}{\symAdmittance_{ \indexGridLines , \symPhaseC \symPhaseA}^{\seriesss}    }}
\newcommand{\yijksCB}[0]{\textcolor{\complexparamcolor}{\symAdmittance_{ \indexGridLines , \symPhaseC \symPhaseB}^{\seriesss}    }}
\newcommand{\yijksCC}[0]{\textcolor{\complexparamcolor}{\symAdmittance_{ \indexGridLines , \symPhaseC \symPhaseC}^{\seriesss}    }}

\newcommand{\gijksAA}[0]{\textcolor{\paramcolor}{\symConductance_{ \indexGridLines , \symPhaseA \symPhaseA}^{\seriesss}    }}
\newcommand{\gijksAB}[0]{\textcolor{\paramcolor}{\symConductance_{ \indexGridLines , \symPhaseA \symPhaseB}^{\seriesss}    }}
\newcommand{\gijksAC}[0]{\textcolor{\paramcolor}{\symConductance_{ \indexGridLines , \symPhaseA \symPhaseC}^{\seriesss}    }}
\newcommand{\gijksBA}[0]{\textcolor{\paramcolor}{\symConductance_{ \indexGridLines , \symPhaseB \symPhaseA}^{\seriesss}    }}
\newcommand{\gijksBB}[0]{\textcolor{\paramcolor}{\symConductance_{ \indexGridLines , \symPhaseB \symPhaseB}^{\seriesss}    }}
\newcommand{\gijksBC}[0]{\textcolor{\paramcolor}{\symConductance_{ \indexGridLines , \symPhaseB \symPhaseC}^{\seriesss}    }}
\newcommand{\gijksCA}[0]{\textcolor{\paramcolor}{\symConductance_{ \indexGridLines , \symPhaseC \symPhaseA}^{\seriesss}    }}
\newcommand{\gijksCB}[0]{\textcolor{\paramcolor}{\symConductance_{ \indexGridLines , \symPhaseC \symPhaseB}^{\seriesss}    }}
\newcommand{\gijksCC}[0]{\textcolor{\paramcolor}{\symConductance_{ \indexGridLines , \symPhaseC \symPhaseC}^{\seriesss}    }}

\newcommand{\bijksAA}[0]{\textcolor{\paramcolor}{\symSusceptance_{ \indexGridLines , \symPhaseA \symPhaseA}^{\seriesss}    }}
\newcommand{\bijksAB}[0]{\textcolor{\paramcolor}{\symSusceptance_{ \indexGridLines , \symPhaseA \symPhaseB}^{\seriesss}    }}
\newcommand{\bijksAC}[0]{\textcolor{\paramcolor}{\symSusceptance_{ \indexGridLines , \symPhaseA \symPhaseC}^{\seriesss}    }}
\newcommand{\bijksBA}[0]{\textcolor{\paramcolor}{\symSusceptance_{ \indexGridLines , \symPhaseB \symPhaseA}^{\seriesss}    }}
\newcommand{\bijksBB}[0]{\textcolor{\paramcolor}{\symSusceptance_{ \indexGridLines , \symPhaseB \symPhaseB}^{\seriesss}    }}
\newcommand{\bijksBC}[0]{\textcolor{\paramcolor}{\symSusceptance_{ \indexGridLines , \symPhaseB \symPhaseC}^{\seriesss}    }}
\newcommand{\bijksCA}[0]{\textcolor{\paramcolor}{\symSusceptance_{ \indexGridLines , \symPhaseC \symPhaseA}^{\seriesss}    }}
\newcommand{\bijksCB}[0]{\textcolor{\paramcolor}{\symSusceptance_{ \indexGridLines , \symPhaseC \symPhaseB}^{\seriesss}    }}
\newcommand{\bijksCC}[0]{\textcolor{\paramcolor}{\symSusceptance_{ \indexGridLines , \symPhaseC \symPhaseC}^{\seriesss}    }}

\newcommand{\Ilijh}[0]{\textcolor{\complexcolor}{\mathbf{\symCurrent}_{\indexGridLines\indexGridNode \indexGridNodeTwo,\symHarmonic}^{}  }}
\newcommand{\Iljih}[0]{\textcolor{\complexcolor}{\mathbf{\symCurrent}_{\indexGridLines\indexGridNodeTwo\indexGridNode ,\symHarmonic}^{}  }}
\newcommand{\IlijAh}[0]{\textcolor{\complexcolor}{\symCurrent_{\indexGridLines\indexGridNode \indexGridNodeTwo, \symPhaseA,\symHarmonic}^{}  }}
\newcommand{\IlijBh}[0]{\textcolor{\complexcolor}{\symCurrent_{\indexGridLines\indexGridNode \indexGridNodeTwo, \symPhaseB,\symHarmonic}^{}  }}
\newcommand{\IlijCh}[0]{\textcolor{\complexcolor}{\symCurrent_{\indexGridLines\indexGridNode \indexGridNodeTwo, \symPhaseC,\symHarmonic}^{}  }}

\newcommand{\Slijh}[0]{\textcolor{\complexcolor}{\mathbf{\symApparentPower}_{\indexGridLines\indexGridNode \indexGridNodeTwo,\symHarmonic}^{}  }}
\newcommand{\Sljih}[0]{\textcolor{\complexcolor}{\mathbf{\symApparentPower}_{\indexGridLines \indexGridNodeTwo \indexGridNode,\symHarmonic}^{}  }}

\newcommand{\Sslijh}[0]{\textcolor{\complexcolor}{\mathbf{\symApparentPower}_{\indexGridLines\indexGridNode \indexGridNodeTwo,\symHarmonic}^{\seriesss}  }}
\newcommand{\Ssljih}[0]{\textcolor{\complexcolor}{\mathbf{\symApparentPower}_{\indexGridLines \indexGridNodeTwo \indexGridNode,\symHarmonic}^{\seriesss}  }}

\newcommand{\Islijh}[0]{\textcolor{\complexcolor}{\mathbf{\symCurrent}_{\indexGridLines\indexGridNode \indexGridNodeTwo,\symHarmonic}^{\seriesss}  }}
\newcommand{\Isljih}[0]{\textcolor{\complexcolor}{\mathbf{\symCurrent}_{\indexGridLines \indexGridNodeTwo \indexGridNode,\symHarmonic}^{\seriesss}  }}
\newcommand{\IslijAh}[0]{\textcolor{\complexcolor}{\symCurrent_{\indexGridLines\indexGridNode \indexGridNodeTwo, \symPhaseA,\symHarmonic}^{\seriesss}  }}
\newcommand{\IslijBh}[0]{\textcolor{\complexcolor}{\symCurrent_{\indexGridLines\indexGridNode \indexGridNodeTwo, \symPhaseB,\symHarmonic}^{\seriesss}  }}
\newcommand{\IslijCh}[0]{\textcolor{\complexcolor}{\symCurrent_{\indexGridLines\indexGridNode \indexGridNodeTwo, \symPhaseC,\symHarmonic}^{\seriesss}  }}

\newcommand{\Ishlijh}[0]{\textcolor{\complexcolor}{\mathbf{\symCurrent}_{\indexGridLines\indexGridNode \indexGridNodeTwo,\symHarmonic}^{\shuntss}  }}
\newcommand{\Ishljih}[0]{\textcolor{\complexcolor}{\mathbf{\symCurrent}_{\indexGridLines \indexGridNodeTwo\indexGridNode ,\symHarmonic}^{\shuntss}  }}
\newcommand{\IshlijAh}[0]{\textcolor{\complexcolor}{\symCurrent_{\indexGridLines\indexGridNode \indexGridNodeTwo, \symPhaseA,\symHarmonic}^{\shuntss}  }}
\newcommand{\IshlijBh}[0]{\textcolor{\complexcolor}{\symCurrent_{\indexGridLines\indexGridNode \indexGridNodeTwo, \symPhaseB,\symHarmonic}^{\shuntss}  }}
\newcommand{\IshlijCh}[0]{\textcolor{\complexcolor}{\symCurrent_{\indexGridLines\indexGridNode \indexGridNodeTwo, \symPhaseC,\symHarmonic}^{\shuntss}  }}

\newcommand{\IijkA}[0]{\textcolor{\complexcolor}{\symCurrent_{\indexGridLines\indexGridNode \indexGridNodeTwo, \symPhaseA}^{}  }}
\newcommand{\IijkB}[0]{\textcolor{\complexcolor}{\symCurrent_{\indexGridLines\indexGridNode \indexGridNodeTwo, \symPhaseB}^{}  }}
\newcommand{\IijkC}[0]{\textcolor{\complexcolor}{\symCurrent_{\indexGridLines\indexGridNode \indexGridNodeTwo, \symPhaseC}^{}  }}
\newcommand{\IijkN}[0]{\textcolor{\complexcolor}{\symCurrent_{\indexGridLines\indexGridNode \indexGridNodeTwo, \symPhaseN}^{}  }}
\newcommand{\IijkG}[0]{\textcolor{\complexcolor}{\symCurrent_{\indexGridLines\indexGridNode \indexGridNodeTwo, \symPhaseG}^{}  }}

\newcommand{\kifkG}[0]{\textcolor{\complexcolor}{\symCurrent_{k\indexGridNode f, \symPhaseG}^{}  }}

\newcommand{\IijkNreal}[0]{\textcolor{black}{\symCurrent_{\indexGridLines\indexGridNode \indexGridNodeTwo, \symPhaseN}^{\text{re}} }}
\newcommand{\IijkNimag}[0]{\textcolor{black}{\symCurrent_{\indexGridLines\indexGridNode \indexGridNodeTwo, \symPhaseN}^{\text{im}} }}

\newcommand{\IikNG}[0]{\textcolor{\complexcolor}{\symCurrent_{\indexGridNode , \symPhaseN\symPhaseG}}}
\newcommand{\IikNGreal}[0]{\textcolor{black}{\symCurrent_{\indexGridNode , \symPhaseN\symPhaseG}^{\text{re}} }}
\newcommand{\IikNGimag}[0]{\textcolor{black}{\symCurrent_{\indexGridNode , \symPhaseN\symPhaseG}^{\text{im}}}}

\newcommand{\Iijkzero}[0]{\textcolor{\complexcolor}{\symCurrent_{\indexGridLines\indexGridNode \indexGridNodeTwo, 0}^{}  }}
\newcommand{\Iijkone}[0]{\textcolor{\complexcolor}{\symCurrent_{\indexGridLines\indexGridNode \indexGridNodeTwo, 1}^{}  }}
\newcommand{\Iijktwo}[0]{\textcolor{\complexcolor}{\symCurrent_{\indexGridLines\indexGridNode \indexGridNodeTwo, 2}^{}  }}

\newcommand{\IjikA}[0]{\textcolor{\complexcolor}{\symCurrent_{\indexGridLines\indexGridNodeTwo\indexGridNode , \symPhaseA}^{}  }}
\newcommand{\IjikB}[0]{\textcolor{\complexcolor}{\symCurrent_{\indexGridLines\indexGridNodeTwo\indexGridNode  , \symPhaseB}^{}  }}
\newcommand{\IjikC}[0]{\textcolor{\complexcolor}{\symCurrent_{\indexGridLines\indexGridNodeTwo\indexGridNode , \symPhaseC}^{}  }}
\newcommand{\IjikN}[0]{\textcolor{\complexcolor}{\symCurrent_{\indexGridLines\indexGridNodeTwo\indexGridNode , \symPhaseN}^{}  }}
\newcommand{\IjikG}[0]{\textcolor{\complexcolor}{\symCurrent_{\indexGridLines\indexGridNodeTwo\indexGridNode , \symPhaseG}^{}  }}

\newcommand{\IijkrefA}[0]{\textcolor{\complexparamcolor}{\symCurrent_{\indexGridLines\indexGridNode \indexGridNodeTwo,  \symPhaseA}^{\refss}  }}
\newcommand{\IijkrefB}[0]{\textcolor{\complexparamcolor}{\symCurrent_{\indexGridLines\indexGridNode \indexGridNodeTwo, \symPhaseB}^{\refss}  }}
\newcommand{\IijkrefC}[0]{\textcolor{\complexparamcolor}{\symCurrent_{\indexGridLines\indexGridNode \indexGridNodeTwo, \symPhaseC}^{\refss}  }}
\newcommand{\IijksrefA}[0]{\textcolor{\complexparamcolor}{\symCurrent_{\indexGridLines\indexGridNode \indexGridNodeTwo, \seriesss, \symPhaseA}^{\refss}  }}
\newcommand{\IijksrefB}[0]{\textcolor{\complexparamcolor}{\symCurrent_{\indexGridLines\indexGridNode \indexGridNodeTwo, \seriesss, \symPhaseB}^{\refss}  }}
\newcommand{\IijksrefC}[0]{\textcolor{\complexparamcolor}{\symCurrent_{\indexGridLines\indexGridNode \indexGridNodeTwo, \seriesss, \symPhaseC}^{\refss}  }}

\newcommand{\IijksA}[0]{\textcolor{\complexcolor}{\symCurrent_{\indexGridLines\indexGridNode \indexGridNodeTwo, \symPhaseA}^{ \seriesss}  }}
\newcommand{\IijksB}[0]{\textcolor{\complexcolor}{\symCurrent_{\indexGridLines\indexGridNode \indexGridNodeTwo, \symPhaseB}^{\seriesss}  }}
\newcommand{\IijksC}[0]{\textcolor{\complexcolor}{\symCurrent_{\indexGridLines\indexGridNode \indexGridNodeTwo, \symPhaseC}^{\seriesss}  }}
\newcommand{\IijksN}[0]{\textcolor{\complexcolor}{\symCurrent_{\indexGridLines\indexGridNode \indexGridNodeTwo, \symPhaseN}^{\seriesss}  }}
\newcommand{\IijksG}[0]{\textcolor{\complexcolor}{\symCurrent_{\indexGridLines\indexGridNode \indexGridNodeTwo, \symPhaseG}^{\seriesss}  }}

\newcommand{\IjiksA}[0]{\textcolor{\complexcolor}{\symCurrent_{\indexGridLines\indexGridNodeTwo\indexGridNode , \symPhaseA}^{ \seriesss}  }}
\newcommand{\IjiksB}[0]{\textcolor{\complexcolor}{\symCurrent_{\indexGridLines\indexGridNodeTwo\indexGridNode , \symPhaseB}^{\seriesss}  }}
\newcommand{\IjiksC}[0]{\textcolor{\complexcolor}{\symCurrent_{\indexGridLines\indexGridNodeTwo\indexGridNode , \symPhaseC}^{\seriesss}  }}
\newcommand{\IjiksN}[0]{\textcolor{\complexcolor}{\symCurrent_{\indexGridLines\indexGridNodeTwo\indexGridNode , \symPhaseN}^{\seriesss}  }}

\newcommand{\IijkshA}[0]{\textcolor{\complexcolor}{\symCurrent_{\indexGridLines\indexGridNode \indexGridNodeTwo, \symPhaseA}^{ \shuntss}  }}
\newcommand{\IijkshB}[0]{\textcolor{\complexcolor}{\symCurrent_{\indexGridLines\indexGridNode \indexGridNodeTwo, \symPhaseB}^{\shuntss}  }}
\newcommand{\IijkshC}[0]{\textcolor{\complexcolor}{\symCurrent_{\indexGridLines\indexGridNode \indexGridNodeTwo, \symPhaseC}^{\shuntss}  }}
\newcommand{\IijkshN}[0]{\textcolor{\complexcolor}{\symCurrent_{\indexGridLines\indexGridNode \indexGridNodeTwo, \symPhaseN}^{\shuntss}  }}
\newcommand{\IijkshG}[0]{\textcolor{\complexcolor}{\symCurrent_{\indexGridLines\indexGridNode \indexGridNodeTwo, \symPhaseG}^{\shuntss}  }}

\newcommand{\IjikshN}[0]{\textcolor{\complexcolor}{\symCurrent_{\indexGridLines \indexGridNodeTwo \indexGridNode , \symPhaseN}^{\shuntss}  }}

\newcommand{\IikA}[0]{\textcolor{\complexcolor}{\symCurrent_{\indexGridNode , \symPhaseA}^{}  }}
\newcommand{\IikB}[0]{\textcolor{\complexcolor}{\symCurrent_{\indexGridNode , \symPhaseB}^{}  }}
\newcommand{\IikC}[0]{\textcolor{\complexcolor}{\symCurrent_{\indexGridNode , \symPhaseC}^{}  }}

\newcommand{\yijkshAA}[0]{\textcolor{\complexparamcolor}{\symAdmittance_{\indexGridLines \indexGridNode \indexGridNodeTwo, \symPhaseA \symPhaseA}^{\shuntss}    }}
\newcommand{\yijkshAB}[0]{\textcolor{\complexparamcolor}{\symAdmittance_{\indexGridLines \indexGridNode \indexGridNodeTwo, \symPhaseA \symPhaseB}^{\shuntss}    }}
\newcommand{\yijkshAC}[0]{\textcolor{\complexparamcolor}{\symAdmittance_{\indexGridLines \indexGridNode \indexGridNodeTwo, \symPhaseA \symPhaseC}^{\shuntss}    }}
\newcommand{\yijkshAN}[0]{\textcolor{\complexparamcolor}{\symAdmittance_{\indexGridLines \indexGridNode \indexGridNodeTwo, \symPhaseA \symPhaseN}^{\shuntss}    }}
\newcommand{\yijkshAG}[0]{\textcolor{\complexparamcolor}{\symAdmittance_{\indexGridLines \indexGridNode \indexGridNodeTwo, \symPhaseA \symPhaseG}^{\shuntss}    }}

\newcommand{\yijkshBA}[0]{\textcolor{\complexparamcolor}{\symAdmittance_{\indexGridLines \indexGridNode \indexGridNodeTwo, \symPhaseB \symPhaseA}^{\shuntss}    }}
\newcommand{\yijkshBB}[0]{\textcolor{\complexparamcolor}{\symAdmittance_{\indexGridLines \indexGridNode \indexGridNodeTwo, \symPhaseB \symPhaseB}^{\shuntss}    }}
\newcommand{\yijkshBC}[0]{\textcolor{\complexparamcolor}{\symAdmittance_{\indexGridLines \indexGridNode \indexGridNodeTwo, \symPhaseB \symPhaseC}^{\shuntss}    }}
\newcommand{\yijkshBN}[0]{\textcolor{\complexparamcolor}{\symAdmittance_{\indexGridLines \indexGridNode \indexGridNodeTwo, \symPhaseB \symPhaseN}^{\shuntss}    }}
\newcommand{\yijkshBG}[0]{\textcolor{\complexparamcolor}{\symAdmittance_{\indexGridLines \indexGridNode \indexGridNodeTwo, \symPhaseB \symPhaseG}^{\shuntss}    }}

\newcommand{\yijkshCA}[0]{\textcolor{\complexparamcolor}{\symAdmittance_{\indexGridLines \indexGridNode \indexGridNodeTwo, \symPhaseC \symPhaseA}^{\shuntss}    }}
\newcommand{\yijkshCB}[0]{\textcolor{\complexparamcolor}{\symAdmittance_{\indexGridLines \indexGridNode \indexGridNodeTwo, \symPhaseC \symPhaseB}^{\shuntss}    }}
\newcommand{\yijkshCC}[0]{\textcolor{\complexparamcolor}{\symAdmittance_{\indexGridLines \indexGridNode \indexGridNodeTwo, \symPhaseC \symPhaseC}^{\shuntss}    }}\newcommand{\yijkshCN}[0]{\textcolor{\complexparamcolor}{\symAdmittance_{\indexGridLines \indexGridNode \indexGridNodeTwo, \symPhaseC \symPhaseN}^{\shuntss}    }}\newcommand{\yijkshCG}[0]{\textcolor{\complexparamcolor}{\symAdmittance_{\indexGridLines \indexGridNode \indexGridNodeTwo, \symPhaseC \symPhaseG}^{\shuntss}    }}

\newcommand{\yijkshNA}[0]{\textcolor{\complexparamcolor}{\symAdmittance_{\indexGridLines \indexGridNode \indexGridNodeTwo, \symPhaseN \symPhaseA}^{\shuntss}    }}
\newcommand{\yijkshNB}[0]{\textcolor{\complexparamcolor}{\symAdmittance_{\indexGridLines \indexGridNode \indexGridNodeTwo, \symPhaseN \symPhaseB}^{\shuntss}    }}
\newcommand{\yijkshNC}[0]{\textcolor{\complexparamcolor}{\symAdmittance_{\indexGridLines \indexGridNode \indexGridNodeTwo, \symPhaseN \symPhaseC}^{\shuntss}    }}\newcommand{\yijkshNN}[0]{\textcolor{\complexparamcolor}{\symAdmittance_{\indexGridLines \indexGridNode \indexGridNodeTwo, \symPhaseN \symPhaseN}^{\shuntss}    }}\newcommand{\yijkshNG}[0]{\textcolor{\complexparamcolor}{\symAdmittance_{\indexGridLines \indexGridNode \indexGridNodeTwo, \symPhaseN \symPhaseG}^{\shuntss}    }}

\newcommand{\yijkshGA}[0]{\textcolor{\complexparamcolor}{\symAdmittance_{\indexGridLines \indexGridNode \indexGridNodeTwo, \symPhaseG \symPhaseA}^{\shuntss}    }}
\newcommand{\yijkshGB}[0]{\textcolor{\complexparamcolor}{\symAdmittance_{\indexGridLines \indexGridNode \indexGridNodeTwo, \symPhaseG \symPhaseB}^{\shuntss}    }}
\newcommand{\yijkshGC}[0]{\textcolor{\complexparamcolor}{\symAdmittance_{\indexGridLines \indexGridNode \indexGridNodeTwo, \symPhaseG \symPhaseC}^{\shuntss}    }}\newcommand{\yijkshGN}[0]{\textcolor{\complexparamcolor}{\symAdmittance_{\indexGridLines \indexGridNode \indexGridNodeTwo, \symPhaseG \symPhaseN}^{\shuntss}    }}\newcommand{\yijkshGG}[0]{\textcolor{\complexparamcolor}{\symAdmittance_{\indexGridLines \indexGridNode \indexGridNodeTwo, \symPhaseG \symPhaseG}^{\shuntss}    }}

\newcommand{\gijkshAA}[0]{\textcolor{\paramcolor}{\symConductance_{ \indexGridLines\indexGridNode \indexGridNodeTwo, \symPhaseA \symPhaseA}^{\shuntss}    }}
\newcommand{\gijkshAB}[0]{\textcolor{\paramcolor}{\symConductance_{\indexGridLines \indexGridNode \indexGridNodeTwo, \symPhaseA \symPhaseB}^{\shuntss}    }}
\newcommand{\gijkshAC}[0]{\textcolor{\paramcolor}{\symConductance_{\indexGridLines \indexGridNode \indexGridNodeTwo, \symPhaseA \symPhaseC}^{\shuntss}    }}

\newcommand{\gijkshBA}[0]{\textcolor{\paramcolor}{\symConductance_{ \indexGridLines\indexGridNode \indexGridNodeTwo, \symPhaseB \symPhaseA}^{\shuntss}    }}
\newcommand{\gijkshBB}[0]{\textcolor{\paramcolor}{\symConductance_{\indexGridLines \indexGridNode \indexGridNodeTwo, \symPhaseB \symPhaseB}^{\shuntss}    }}
\newcommand{\gijkshBC}[0]{\textcolor{\paramcolor}{\symConductance_{\indexGridLines \indexGridNode \indexGridNodeTwo, \symPhaseB \symPhaseC}^{\shuntss}    }}

\newcommand{\gijkshCA}[0]{\textcolor{\paramcolor}{\symConductance_{ \indexGridLines\indexGridNode \indexGridNodeTwo, \symPhaseC \symPhaseA}^{\shuntss}    }}
\newcommand{\gijkshCB}[0]{\textcolor{\paramcolor}{\symConductance_{ \indexGridLines\indexGridNode \indexGridNodeTwo, \symPhaseC \symPhaseB}^{\shuntss}    }}
\newcommand{\gijkshCC}[0]{\textcolor{\paramcolor}{\symConductance_{\indexGridLines \indexGridNode \indexGridNodeTwo, \symPhaseC \symPhaseC}^{\shuntss}    }}

\newcommand{\bijkshAA}[0]{\textcolor{\paramcolor}{\symSusceptance_{ \indexGridLines\indexGridNode \indexGridNodeTwo, \symPhaseA \symPhaseA}^{\shuntss}    }}
\newcommand{\bijkshAB}[0]{\textcolor{\paramcolor}{\symSusceptance_{ \indexGridLines\indexGridNode \indexGridNodeTwo, \symPhaseA \symPhaseB}^{\shuntss}    }}
\newcommand{\bijkshAC}[0]{\textcolor{\paramcolor}{\symSusceptance_{ \indexGridLines\indexGridNode \indexGridNodeTwo, \symPhaseA \symPhaseC}^{\shuntss}    }}

\newcommand{\bijkshBA}[0]{\textcolor{\paramcolor}{\symSusceptance_{ \indexGridLines\indexGridNode \indexGridNodeTwo, \symPhaseB \symPhaseA}^{\shuntss}    }}
\newcommand{\bijkshBB}[0]{\textcolor{\paramcolor}{\symSusceptance_{ \indexGridLines\indexGridNode \indexGridNodeTwo, \symPhaseB \symPhaseB}^{\shuntss}    }}
\newcommand{\bijkshBC}[0]{\textcolor{\paramcolor}{\symSusceptance_{ \indexGridLines\indexGridNode \indexGridNodeTwo, \symPhaseB \symPhaseC}^{\shuntss}    }}

\newcommand{\bijkshCA}[0]{\textcolor{\paramcolor}{\symSusceptance_{ \indexGridLines\indexGridNode \indexGridNodeTwo, \symPhaseC \symPhaseA}^{\shuntss}    }}
\newcommand{\bijkshCB}[0]{\textcolor{\paramcolor}{\symSusceptance_{\indexGridLines \indexGridNode \indexGridNodeTwo, \symPhaseC \symPhaseB}^{\shuntss}    }}
\newcommand{\bijkshCC}[0]{\textcolor{\paramcolor}{\symSusceptance_{ \indexGridLines\indexGridNode \indexGridNodeTwo, \symPhaseC \symPhaseC}^{\shuntss}    }}

\newcommand{\yjikshAA}[0]{\textcolor{\complexparamcolor}{\symAdmittance_{\indexGridLines\indexGridNodeTwo \indexGridNode , \symPhaseA \symPhaseA}^{\shuntss}    }}
\newcommand{\yjikshAB}[0]{\textcolor{\complexparamcolor}{\symAdmittance_{\indexGridLines\indexGridNodeTwo \indexGridNode , \symPhaseA \symPhaseB}^{\shuntss}    }}
\newcommand{\yjikshAC}[0]{\textcolor{\complexparamcolor}{\symAdmittance_{\indexGridLines \indexGridNodeTwo\indexGridNode , \symPhaseA \symPhaseC}^{\shuntss}    }}
\newcommand{\yjikshAN}[0]{\textcolor{\complexparamcolor}{\symAdmittance_{\indexGridLines \indexGridNodeTwo\indexGridNode , \symPhaseA \symPhaseN}^{\shuntss}    }}
\newcommand{\yjikshAG}[0]{\textcolor{\complexparamcolor}{\symAdmittance_{\indexGridLines \indexGridNodeTwo\indexGridNode , \symPhaseA \symPhaseG}^{\shuntss}    }}

\newcommand{\yjikshBA}[0]{\textcolor{\complexparamcolor}{\symAdmittance_{\indexGridLines\indexGridNodeTwo \indexGridNode , \symPhaseB \symPhaseA}^{\shuntss}    }}
\newcommand{\yjikshBB}[0]{\textcolor{\complexparamcolor}{\symAdmittance_{\indexGridLines\indexGridNodeTwo \indexGridNode , \symPhaseB \symPhaseB}^{\shuntss}    }}
\newcommand{\yjikshBC}[0]{\textcolor{\complexparamcolor}{\symAdmittance_{\indexGridLines \indexGridNodeTwo\indexGridNode , \symPhaseB \symPhaseC}^{\shuntss}    }}
\newcommand{\yjikshBN}[0]{\textcolor{\complexparamcolor}{\symAdmittance_{\indexGridLines \indexGridNodeTwo\indexGridNode , \symPhaseB \symPhaseN}^{\shuntss}    }}
\newcommand{\yjikshBG}[0]{\textcolor{\complexparamcolor}{\symAdmittance_{\indexGridLines \indexGridNodeTwo\indexGridNode , \symPhaseB \symPhaseG}^{\shuntss}    }}

\newcommand{\yjikshCA}[0]{\textcolor{\complexparamcolor}{\symAdmittance_{\indexGridLines\indexGridNodeTwo \indexGridNode , \symPhaseC \symPhaseA}^{\shuntss}    }}
\newcommand{\yjikshCB}[0]{\textcolor{\complexparamcolor}{\symAdmittance_{\indexGridLines\indexGridNodeTwo \indexGridNode , \symPhaseC \symPhaseB}^{\shuntss}    }}
\newcommand{\yjikshCC}[0]{\textcolor{\complexparamcolor}{\symAdmittance_{\indexGridLines \indexGridNodeTwo\indexGridNode , \symPhaseC \symPhaseC}^{\shuntss}    }}
\newcommand{\yjikshCN}[0]{\textcolor{\complexparamcolor}{\symAdmittance_{\indexGridLines \indexGridNodeTwo\indexGridNode , \symPhaseC \symPhaseN}^{\shuntss}    }}
\newcommand{\yjikshCG}[0]{\textcolor{\complexparamcolor}{\symAdmittance_{\indexGridLines \indexGridNodeTwo\indexGridNode , \symPhaseC \symPhaseG}^{\shuntss}    }}

\newcommand{\yjikshNA}[0]{\textcolor{\complexparamcolor}{\symAdmittance_{\indexGridLines\indexGridNodeTwo \indexGridNode , \symPhaseN \symPhaseA}^{\shuntss}    }}
\newcommand{\yjikshNB}[0]{\textcolor{\complexparamcolor}{\symAdmittance_{\indexGridLines\indexGridNodeTwo \indexGridNode , \symPhaseN \symPhaseB}^{\shuntss}    }}
\newcommand{\yjikshNC}[0]{\textcolor{\complexparamcolor}{\symAdmittance_{\indexGridLines \indexGridNodeTwo\indexGridNode , \symPhaseN \symPhaseC}^{\shuntss}    }}
\newcommand{\yjikshNN}[0]{\textcolor{\complexparamcolor}{\symAdmittance_{\indexGridLines \indexGridNodeTwo\indexGridNode , \symPhaseN \symPhaseN}^{\shuntss}    }}
\newcommand{\yjikshNG}[0]{\textcolor{\complexparamcolor}{\symAdmittance_{\indexGridLines \indexGridNodeTwo\indexGridNode , \symPhaseN \symPhaseG}^{\shuntss}    }}

\newcommand{\yjikshGA}[0]{\textcolor{\complexparamcolor}{\symAdmittance_{\indexGridLines\indexGridNodeTwo \indexGridNode , \symPhaseG \symPhaseA}^{\shuntss}    }}
\newcommand{\yjikshGB}[0]{\textcolor{\complexparamcolor}{\symAdmittance_{\indexGridLines\indexGridNodeTwo \indexGridNode , \symPhaseG \symPhaseB}^{\shuntss}    }}
\newcommand{\yjikshGC}[0]{\textcolor{\complexparamcolor}{\symAdmittance_{\indexGridLines \indexGridNodeTwo\indexGridNode , \symPhaseG \symPhaseC}^{\shuntss}    }}
\newcommand{\yjikshGN}[0]{\textcolor{\complexparamcolor}{\symAdmittance_{\indexGridLines \indexGridNodeTwo\indexGridNode , \symPhaseG \symPhaseN}^{\shuntss}    }}
\newcommand{\yjikshGG}[0]{\textcolor{\complexparamcolor}{\symAdmittance_{\indexGridLines \indexGridNodeTwo\indexGridNode , \symPhaseG \symPhaseG}^{\shuntss}    }}

\newcommand{\gijkshpp}[0]{\textcolor{\paramcolor}{\symConductance_{\indexGridLines \indexGridNode \indexGridNodeTwo, \indexPhases\indexPhases}^{\shuntss}    }}
\newcommand{\gijkshph}[0]{\textcolor{\paramcolor}{\symConductance_{ \indexGridLines\indexGridNode \indexGridNodeTwo, \indexPhases\indexPhasesTwo}^{\shuntss}    }}
\newcommand{\gkshph}[0]{\textcolor{\paramcolor}{\symConductance_{ \indexShunt, \indexPhases\indexPhasesTwo}   }}
\newcommand{\bkshph}[0]{\textcolor{\paramcolor}{\symSusceptance_{ \indexShunt, \indexPhases\indexPhasesTwo}   }}

\newcommand{\gijkspp}[0]{\textcolor{\paramcolor}{\symConductance_{ \indexGridLines  , \indexPhases\indexPhases}^{\seriesss}    }}
\newcommand{\gijksph}[0]{\textcolor{\paramcolor}{\symConductance_{ \indexGridLines  , \indexPhases\indexPhasesTwo}^{\seriesss}    }}

\newcommand{\bijkshpp}[0]{\textcolor{\paramcolor}{\symSusceptance_{\indexGridLines \indexGridNode \indexGridNodeTwo, \indexPhases\indexPhases}^{\shuntss}    }}
\newcommand{\bijkspp}[0]{\textcolor{\paramcolor}{\symSusceptance_{ \indexGridLines, \indexPhases\indexPhases}^{\seriesss}    }}
\newcommand{\bijksph}[0]{\textcolor{\paramcolor}{\symSusceptance_{ \indexGridLines, \indexPhases\indexPhasesTwo}^{\seriesss}    }}

\newcommand{\bijkshph}[0]{\textcolor{\paramcolor}{\symSusceptance_{\indexGridLines \indexGridNode \indexGridNodeTwo, \indexPhases\indexPhasesTwo}^{\shuntss}    }}

\newcommand{\zbental}[0]{{z}}
\newcommand{\ybental}[0]{{y}}
\newcommand{\xbental}[0]{{x}}
\newcommand{\xonebental}[0]{\textcolor{black}{\xbental_{1} }}
\newcommand{\xtwobental}[0]{\textcolor{black}{\xbental_{2} }}
\newcommand{\xthreebental}[0]{\textcolor{black}{\xbental_{3} }}
\newcommand{\xfourbental}[0]{\textcolor{black}{\xbental_{4} }}
\newcommand{\xfivebental}[0]{\textcolor{black}{\xbental_{5} }}
\newcommand{\xsixbental}[0]{\textcolor{black}{\xbental_{6} }}
\newcommand{\xsevenbental}[0]{\textcolor{black}{\xbental_{7} }}
\newcommand{\xibental}[0]{\textcolor{black}{\xbental_{\counterbental} }}
\newcommand{\upperbental}[0]{{U}}
\newcommand{\lowerbental}[0]{{L}}
\newcommand{\xilowerbental}[0]{\textcolor{\paramcolor}{\xibental^{ \lowerbental } }}
\newcommand{\xiupperbental}[0]{\textcolor{\paramcolor}{\xibental^{ \upperbental } }}
\newcommand{\xlowerbental}[1]{\textcolor{\boundscolor}{\xbental_{#1}^{ \lowerbental } }}
\newcommand{\xupperbental}[1]{\textcolor{\boundscolor}{\xbental_{#1}^{ \upperbental } }}

\newcommand{\sigmabental}[1]{\textcolor{black}{\sigma_{(#1)}}}
\newcommand{\etabental}[1]{\textcolor{black}{\eta_{(#1)}}}
\newcommand{\counterbental}[0]{\textcolor{\paramcolor}{i}}
\newcommand{\countbental}[0]{\textcolor{\paramcolor}{\nu}}
\newcommand{\errorbental}[0]{\textcolor{\paramcolor}{\epsilon}}
\newcommand{\countpoly}[0]{\textcolor{\paramcolor}{\mu}}

\newcommand{\circumscribedopoly}[0]{\textcolor{\paramcolor}{\symSetting_{\circumss}}}
\newcommand{\rotatedopoly}[0]{\textcolor{\paramcolor}{\symSetting_{\rotatedss}}}
\newcommand{\genericbinary}[1]{\textcolor{black}{b_{#1}}}


\newcommand{\pipediameterijk}[0]{ \textcolor{\paramcolor}{  \symFluidPipeDiameter_{\indexGridLines}     }}
\newcommand{\pipeareaijk}[0]{ \textcolor{\paramcolor}{  \symFluidPipeArea_{\indexGridLines}     }}

\newcommand{\fluidpressure}[1]{\textcolor{black}{\symPressure_{#1}}}
\newcommand{\fluidpressureik}[0]{\textcolor{black}{\symPressure_{\indexGridNode }    }}
\newcommand{\fluidpressurejk}[0]{\textcolor{black}{\symPressure_{\indexGridNodeTwo }    }}

\newcommand{\fluidpressureikmin}[0]{\textcolor{\paramcolor}{\symPressure_{\indexGridNode }^{\minss}    }}
\newcommand{\fluidpressureikmax}[0]{\textcolor{\paramcolor}{\symPressure_{\indexGridNode }^{\maxss}    }}
\newcommand{\fluidpressurejkmin}[0]{\textcolor{\paramcolor}{\symPressure_{\indexGridNodeTwo }^{\minss}    }}
\newcommand{\fluidpressurejkmax}[0]{\textcolor{\paramcolor}{\symPressure_{\indexGridNodeTwo }^{\maxss}    }}

\newcommand{\fluidpressuresqik}[0]{\textcolor{black}{\symPressureSquared_{\indexGridNode }    }}
\newcommand{\fluidpressuresqjk}[0]{\textcolor{black}{\symPressureSquared_{\indexGridNodeTwo }    }}

\newcommand{\fluidpressuresqikmin}[0]{\textcolor{\paramcolor}{\symPressureSquared_{\indexGridNode }^{\minss}    }}
\newcommand{\fluidpressuresqjkmin}[0]{\textcolor{\paramcolor}{\symPressureSquared_{\indexGridNodeTwo }^{\minss}    }}
\newcommand{\fluidpressuresqikmax}[0]{\textcolor{\paramcolor}{\symPressureSquared_{\indexGridNode }^{\maxss}    }}
\newcommand{\fluidpressuresqjkmax}[0]{\textcolor{\paramcolor}{\symPressureSquared_{\indexGridNodeTwo }^{\maxss}    }}

\newcommand{\fluidflow}[1]{\textcolor{black}{\symVolumeFlow_{#1}}}
\newcommand{\fluidflowijk}[0]{\textcolor{black}{\fluidflow{\indexGridNode \indexGridNodeTwo}   }}

\newcommand{\fluidmassflow}[1]{\textcolor{black}{\symMassFlow_{#1}}}
\newcommand{\fluidmassflowijk}[0]{\textcolor{black}{\fluidmassflow{\indexGridNode \indexGridNodeTwo} }}

\newcommand{\fluidmassflowdirectionone}[0]{\textcolor{black}{y_{ij}^{+}}}
\newcommand{\fluidmassflowdirectiontwo}[0]{\textcolor{black}{y_{ij}^{-}}}

\newcommand{\fluiddensity}[1]{{\symDensity_{#1}  }}
\newcommand{\fluiddensityijk}[0]{\textcolor{\paramcolor}{\fluiddensity{\indexGridLines}   }}

\newcommand{\fluidviscosity}[1]{{\symViscosity_{#1}  }}
\newcommand{\fluidviscosityijk}[0]{\textcolor{\paramcolor}{\fluidviscosity{\indexGridLines}   }}

\newcommand{\heatflow}[1]{\textcolor{black}{\symHeatFlow_{#1}}}
\newcommand{\heatflowijk}[0]{\textcolor{black}{\symHeatFlow_{\indexGridNode \indexGridNodeTwo} }}

\newcommand{\heatresistance}[1]{\textcolor{\paramcolor}{\symHeatResistance_{#1}}}
\newcommand{\heatresistanceijk}[0]{ \textcolor{\paramcolor}{ \heatresistance{\indexGridLines}  }}   

\newcommand{\heattempijk}[0]{\textcolor{black}{\symTemperature_{\indexGridNode \indexGridNodeTwo}   }}
\newcommand{\heattempik}[0]{\textcolor{black}{\symTemperature_{\indexGridNode }   }}
\newcommand{\heattempjk}[0]{\textcolor{black}{\symTemperature_{ \indexGridNodeTwo}   }}

\newcommand{\fluidhagenresistance}[0]{ \textcolor{\paramcolor}{ w^{\hagenposeuilless}_{\indexGridLines}  }    }
\newcommand{\fluiddarcyresistance}[0]{ \textcolor{\paramcolor}{ w^{\darcyweisbachss}_{\indexGridLines}   }    }

\newcommand{\heatcapacity}[0]{ \textcolor{\paramcolor}{ \symHeatCapacity^{\text{\symPressure}} }}   


\newcommand{\symFlowVariable}[0]{   F }
\newcommand{\symFlowVariableSquare}[0]{   f }
\newcommand{\symPotentialVariable}[0]{ U }
\newcommand{\symPotentialAngleVariable}[0]{ \theta }
\newcommand{\symPotentialVariableSquare}[0]{ u }
\newcommand{\setFlowVariable}[0]{   \mathcal{F} }
\newcommand{\setPotentialVariable}[0]{ \mathcal{U} }
\newcommand{\indexFlowVariable}[0]{  f }
\newcommand{\indexPotentialVariable}[0]{ u }
\newcommand{\indexDomains}[0]{ d }
\newcommand{\setDomains}[0]{\mathcal{D}  }
\newcommand{\setDomainsExt}[0]{\mathcal{D}_0  }
\newcommand{\potentiali}[0]{ \symPotentialVariable_{\indexGridNode} }
\newcommand{\potentialz}[0]{ \symPotentialVariable_{\indexZIP} }
\newcommand{\potentialzlin}[0]{ \symPotentialVariable_{\indexZIP}^{\linss} }
\newcommand{\potentialilin}[0]{ \symPotentialVariable_{\indexGridNode}^{\linss} }
\newcommand{\potentialj}[0]{ \symPotentialVariable_{\indexGridNodeTwo} }
\newcommand{\potentialanglei}[0]{ \symPotentialAngleVariable_{\indexGridNode} }
\newcommand{\potentialanglej}[0]{ \symPotentialAngleVariable_{\indexGridNodeTwo} }
\newcommand{\potentialref}[0]{ \textcolor{\paramcolor}{ \symPotentialVariable^{\refss} }}
\newcommand{\potentialiref}[0]{ \textcolor{\paramcolor}{ \symPotentialVariable_{\indexGridNode}^{\refss} }}
\newcommand{\potentialzref}[0]{ \textcolor{\paramcolor}{ \symPotentialVariable_{\indexZIP}^{\refss} }}
\newcommand{\flowij}[0]{ \symFlowVariable_{\indexGridNode \indexGridNodeTwo} }
\newcommand{\flowi}[0]{ \symFlowVariable_{\indexGridNode } }
\newcommand{\flowz}[0]{ \symFlowVariable_{\indexZIP } }
\newcommand{\flowcomplexij}[0]{ \textcolor{\complexcolor}{\symFlowVariable_{\indexGridNode \indexGridNodeTwo} }}
\newcommand{\flowcomplexji}[0]{ \textcolor{\complexcolor}{\symFlowVariable_{\indexGridNodeTwo\indexGridNode } }}
\newcommand{\flowangleij}[0]{ {}{\symFlowVariable_{\indexGridNode \indexGridNodeTwo}^{\imagangle} }}
\newcommand{\flowji}[0]{ \symFlowVariable_{\indexGridNodeTwo \indexGridNode } }
\newcommand{\flowresistancetij}[0]{\textcolor{\paramcolor}{ r_{\indexGridLines } }}
\newcommand{\flowimpedancetij}[0]{\textcolor{\complexcolor}{ z_{\indexGridLines } }}
\newcommand{\flowreactancetij}[0]{\textcolor{\paramcolor}{ x_{\indexGridLines } }}
\newcommand{\potentialsqi}[0]{ \symPotentialVariableSquare_{\indexGridNode} }
\newcommand{\potentialsqz}[0]{ \symPotentialVariableSquare_{\indexZIP} }
\newcommand{\potentialsqj}[0]{ \symPotentialVariableSquare_{\indexGridNodeTwo} }
\newcommand{\flowsqij}[0]{ \symFlowVariableSquare_{\indexGridNode \indexGridNodeTwo} }
\newcommand{\potentialimin}[0]{ \textcolor{\paramcolor}{\symPotentialVariable_{\indexGridNode}^{\minss} }}
\newcommand{\potentialimax}[0]{ \textcolor{\paramcolor}{\symPotentialVariable_{\indexGridNode}^{\maxss} }}
\newcommand{\potentialirated}[0]{ \textcolor{\sizingcolor}{\symPotentialVariable_{\indexGridNode}^{\ratedss} }}
\newcommand{\potentialjmin}[0]{ \textcolor{\paramcolor}{\symPotentialVariable_{\indexGridNodeTwo}^{\minss} }}
\newcommand{\potentialjmax}[0]{ \textcolor{\paramcolor}{\symPotentialVariable_{\indexGridNodeTwo}^{\maxss} }}
\newcommand{\potentialjrated}[0]{ \textcolor{\sizingcolor}{\symPotentialVariable_{\indexGridNodeTwo}^{\ratedss} }}
\newcommand{\potentialijrated}[0]{ \textcolor{\sizingcolor}{\symPotentialVariable_{\indexGridNode\indexGridNodeTwo}^{\ratedss} }}
\newcommand{\potentialjirated}[0]{ \textcolor{\sizingcolor}{\symPotentialVariable_{\indexGridNodeTwo \indexGridNode}^{\ratedss} }}
\newcommand{\flowijmin}[0]{ \textcolor{\paramcolor}{\symFlowVariable_{\indexGridNode \indexGridNodeTwo}^{\minss} }}
\newcommand{\flowijmax}[0]{\textcolor{\paramcolor}{ \symFlowVariable_{\indexGridNode \indexGridNodeTwo}^{\maxss} }}
\newcommand{\flowijrated}[0]{\textcolor{\sizingcolor}{ \symFlowVariable_{\indexGridNode \indexGridNodeTwo}^{\ratedss} }}
\newcommand{\flowjirated}[0]{\textcolor{\sizingcolor}{ \symFlowVariable_{ \indexGridNodeTwo \indexGridNode }^{\ratedss} }}

\newcommand{\powerbalancefactorij}[0]{\textcolor{\paramcolor}{c_{\indexGridLines } }}

\newcommand{\enthalpycombustion}{\textcolor{\paramcolor}{\Delta \symEnthalpy^{\circ}_{\combustionss}}}
\newcommand{\enthalpyformation}{\textcolor{\paramcolor}{\Delta \symEnthalpy^{\circ}_{\formationss}}}

\newcommand{\slackconvexificationij}[0]{{\symSlack_{\indexGridNode \indexGridNodeTwo }^{\PFconvexss} }}
\newcommand{\slackconvexification}[0]{{\symSlack^{\PFconvexss} }}
\newcommand{\slackmaxconvexification}[0]{{\symSlack_{  }^{\maxss} }}
\newcommand{\slackminconvexification}[0]{{\symSlack_{  }^{\minss} }}
\newcommand{\slackabsconvexification}[0]{{\symSlack_{  }^{\absss} }}
\newcommand{\slackcomplexmagnitudeconvexification}[0]{{\symSlack_{  }^{\absss} }}
\newcommand{\slackbental}[0]{{\symSlack_{  }^{\bentallinearization} }}


\newcommand{\mapping}[0]{   \to }
\newcommand{\injection}[0]{   \rightarrowtail }
\newcommand{\surjection}[0]{   \twoheadrightarrow }
\newcommand{\bijection}[0]{   \leftrightarrow }

\newcommand{\symbollistlabelwith}{$\PPNeleclossk$}
\newcommand{\imagnumber}[0]{\textcolor{\paramcolor}{j}}
\newcommand{\imagangle}[0]{\textcolor{\paramcolor}{\angle}}

\newcommand{\timebra}{}

\newcommand{\constant}{}
\newcommand{\constantcomplex}{}
\newcommand{\linear}{}
\newcommand{\quadratic}{}
\newcommand{\convex}{ }
\newcommand{\nonconvex}{ }
\newcommand{\binary}{}
\newcommand{\convexquadratic}{}
\newcommand{\conic}{}
\newcommand{\circular}{}
\newcommand{\nonconvexquadratic}{}
\newcommand{\nonconvexnonlinear}{}
\newcommand{\linearcomplex}{}
\newcommand{\removed}{}
\newcommand{\MILP}{}
\newcommand{\SDP}{}
\newcommand{\SOCP}{}
\newcommand{\reformulated}{}
\newcommand{\substitution}{}
\newcommand{\convexified}{}
\newcommand{\penalization}{}
\newcommand{\indexZIP}{}

\clearpage 
 \section*{Nomenclature} \label{sec_nomenclature}
\begin{figure}
  \centering
    \includegraphics[width=0.99\columnwidth]{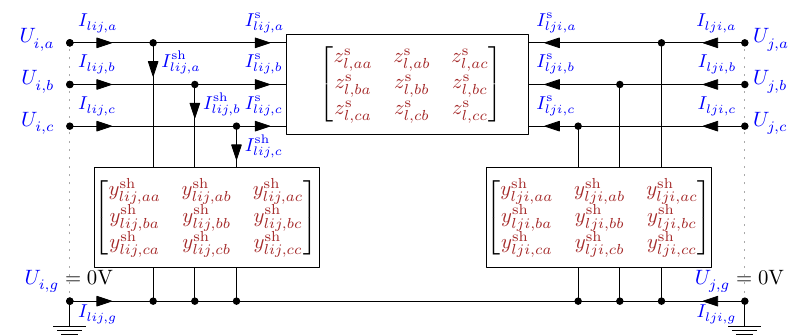}
  \caption{Unbalanced $3\times3$ $\Pi$-model branch in voltage and current variables. }  \label{fig_linemodel_iv_3x3}
\end{figure}

This article depends on the definition of a  variety of scalar, vector and matrix parameters and variables related to grid buses and branches (see Tables 1-5). The core variables are current, voltage and power, whereas parameters are mainly impedance and variable bounds. 
Fig.~\ref{fig_linemodel_iv_3x3} summarizes the variables and parameters defined in the fundamental $3\times3$ branch model for which the well known The $\Pi$-equivalent model used. With respect to balance networks, both series and shunt elements are represented by full complex-valued matrices including the mutual impedance coupling between the conductors. Using the $\Pi$-equivalent model, the branch current $\IijkSDP$ can be split into a series component $\IijskSDP$ and a shunt component $\IijshkSDP$, respectively.  All circuit element voltages are defined w.r.t. (local) ground voltage $\VikG=0\,V$. 

Table~\ref{tab_typography} illustrates typography and mathematical notation used throughout;
Table~\ref{tab_sets} defines sets and indices;
Table~\ref{tab_parameterss} defines parameters;
Table~\ref{tab_variables} defines typical engineering variables;
Table~\ref{tab_lifted_variables} defines lifted variables. Black and red colors indicate \emph{real-valued} variables and parameters, respectively. Blue and brown colors are used for \emph{complex-valued} variables and parameters, respectively. 

\begin{table}
 \centering 
  \caption{Typography and mathematical notation}\label{tab_typography}
\begin{IEEEeqnarraybox}[\IEEEeqnarraystrutmode]{L s}
\hline
x& real scalar variable \\
\mathbf{X}& real vector or matrix variable \\
\textcolor{\complexcolor}{x}& complex scalar variable \\
\textcolor{\complexcolor}{\mathbf{X}}& complex vector or matrix variable \\
\textcolor{\paramcolor}{x}& real scalar parameter \\
\textcolor{\paramcolor}{\mathbf{X}}& real vector or matrix parameter \\
\textcolor{\complexparamcolor}{x}& complex scalar parameter \\
\textcolor{\complexparamcolor}{\mathbf{X}}& complex vector or matrix parameter \\
\mathcal{X}& set \\
\mathbf{X}^{\transpose}&  transpose of   $\mathbf{X}$\\
\textcolor{\complexcolor}{\mathbf{X}}^{*}& conjugate  of   $\textcolor{\complexcolor}{\mathbf{X}}$\\
\textcolor{\complexcolor}{\mathbf{X}}^{\hermitiantranspose}= \textcolor{\complexcolor}{(\mathbf{X}^*)}^{\transpose}& conjugate transpose of  $\textcolor{\complexcolor}{\mathbf{X}}$\\
\circ & element-wise multiplication \\
\oslash & element-wise division \\
\imagnumber& imaginary unit, satisfies $ \imagnumber^2 = -1$ \\
a\imagangle b & polar notation of complex number $ a \cdot e^{\imagnumber b}$ \\
\setRealMatrix{n}{m} & set of real $n\times m$ matrices \\
 \setComplexMatrix{n}{m}  &  set of complex $n\times m$ matrices \\
 \setHermitianMatrix{n} \subset \setComplexMatrix{n}{n}  &  set of Hermitian $n\times n$  matrices \\
\diag(\textcolor{\complexcolor}{\mathbf{X}}) & extract diagonal of $\textcolor{\complexcolor}{\mathbf{X}}$, $\diag: \setComplexMatrix{n}{n} \rightarrow \setComplexMatrix{n}{1} $\\
 \hline
\end{IEEEeqnarraybox}
\end{table}

\begin{table}
 \centering 
  \caption{Sets and indices}\label{tab_sets}
    \begin{tabular}{l l  }
\hline
Phases & $\indexPhases,\indexPhasesTwo \in \setPhases = \{\symPhaseA, \symPhaseB, \symPhaseC  \}  $ \\
Branches &$\indexGridLines\in  \omegaGridLines $ \\		
Buses & $\indexGridNode \in  \omegaGridNodes $ \\
Topology (forward)& $\indexGridLines \indexGridNode\indexGridNodeTwo \in  \setTopology \subseteq  \omegaGridLines \times \omegaGridNodes \times \omegaGridNodes $\\		
Topology (reverse)  & $\indexGridLines  \indexGridNode \indexGridNodeTwo \in  \setTopologyReverse = \{\indexGridLines \indexGridNodeTwo \indexGridNode  \, | \,  \indexGridLines \indexGridNode\indexGridNodeTwo \in  \setTopology \} $ \\		
Topology & $ \indexGridLines \indexGridNode\indexGridNodeTwo \in \setTopologyBoth = \setTopology \cup \setTopologyReverse  $ \\		
Bus pairs & $ \indexGridNode\indexGridNodeTwo \in  \setBuspairs  = \{ ij \, | \,  \indexGridLines \indexGridNode\indexGridNodeTwo \in  \setTopology   \} \subseteq \omegaGridNodes \times \omegaGridNodes $\\
Units & $\indexUnit \in  \setUnits $\\		
Unit connectivity & $\indexUnit \indexGridNode \in  \setUnitTopology \subseteq  \setUnits \times \omegaGridNodes $\\		
Shunt connectivity & $\indexShunt \indexGridNode \in  \setShuntTopology \subseteq  \setUnits \times \omegaGridNodes $\\		
\hline
\end{tabular}
\end{table}

\begin{table}
\setlength{\tabcolsep}{1pt}
 \centering 
  \caption{Parameters}\label{tab_parameterss}
    \begin{tabular}{l l  }
\hline
Bus voltage magnitude min./max. (V) & $\VikSDPmin, \VikSDPmax     \in \setRealMatrix{|\setPhases|}{1}$ \\
Bus phase angle diff. min./max. (rad)& $\thetaiSDPmin, \thetaiSDPmax  \in \setRealMatrix{|\setPhases|}{1}$ \\
 Branch current rating (A) &$ \IijkSDPrated     \in \setRealMatrix{|\setPhases|}{1}$ \\
 Branch apparent power rating (VA) &$ \SijkSDPrated     \in \setRealMatrix{|\setPhases|}{1}$ \\
Branch series  impedance  ($\Omega$)    &$  \zijksSDP   \in \setComplexMatrix{|\setPhases|}{|\setPhases|} $  \\ 	
 Branch series  admittance   (S)   &$  \yijksSDP     \in \setComplexMatrix{|\setPhases|}{|\setPhases|} $ \\ 	
Branch from/to shunt  admittance  (S)    &$  \yijkshSDP,\yjikshSDP   \in \setComplexMatrix{|\setPhases|}{|\setPhases|}$   \\ 	
Bus pair angle diff. min./max. (rad)& $\thetaijSDPmin, \thetaijSDPmax  \in \setRealMatrix{|\setPhases|}{1}$ \\
Bus shunt  admittance (S)    &$  \ybSDP   \in \setComplexMatrix{|\setPhases|}{|\setPhases|}$   \\ 	
Unit current rating (A) &$ \IuunitSDPrated     \in \setRealMatrix{|\setPhases|}{1}$ \\
Unit active power bounds (W) &$ \PuunitSDPmin,\PuunitSDPmax      \in \setRealMatrix{|\setPhases|}{1}$ \\
Unit reactive power bounds (var) &$ \QuunitSDPmin,\QuunitSDPmax       \in \setRealMatrix{|\setPhases|}{1}$ \\
\hline
\end{tabular}
\end{table}


\begin{table}
 \centering 
  \caption{Optimization variables  }\label{tab_variables}
    \begin{tabular}{l l  }
\hline
 Bus voltage  (V)   & $  \VikSDP, \VjkSDP    \in \setComplexMatrix{|\setPhases|}{1}$  \\ 	
Branch current  (A)   &   $ \IijkSDP, \IjikSDP     \in \setComplexMatrix{|\setPhases|}{1} $\\ 	
Branch series current    (A) &$  \IijskSDP,\IjiskSDP       \in \setComplexMatrix{|\setPhases|}{1} $\\ 	
Branch shunt  current   (A)  &$ \IijshkSDP, \IjishkSDP   \in \setComplexMatrix{|\setPhases|}{1} $\\
Branch power flow   (W)  &  $   \SijkSDP, \SjikSDP   \in \setComplexMatrix{|\setPhases|}{|\setPhases|}  $  \\ 	
Branch series power flow     (W) & $    \SijksSDP, \SjiksSDP   \in \setComplexMatrix{|\setPhases|}{|\setPhases|} $   \\ 	
 Unit current   (A)&$  \IuunitSDP \in \setComplexMatrix{|\setPhases|}{1}$  \\
 Unit  power     (W) &  $   \SuunitSDP   \in \setComplexMatrix{|\setPhases|}{|\setPhases|}$  \\ 	
\hline
\end{tabular}
\end{table}

\begin{table}
 \centering 
  \caption{Lifted optimization variables  }\label{tab_lifted_variables}
\begin{equation*}
\begin{IEEEeqnarraybox}[
\IEEEeqnarraystrutmode
]{s L}
\hline
Bus voltage product  (V$^2$)  &    \VisqkSDP, \VjsqkSDP \in   \setHermitianMatrix{|\setPhases|} \\ 	
Bus pair voltage product  (V$^2$)  &    \VijsqkSDP, \VjisqkSDP \in  \setComplexMatrix{|\setPhases|}{|\setPhases|} \\ 	
Branch   current product    (A$^2$)   & \IijsqkSDP, \IjisqkSDP \in   \setHermitianMatrix{|\setPhases|} \\ 	
Branch series  current product    (A$^2$)   & \IijsqskSDP \in   \setHermitianMatrix{|\setPhases|} \\ 	
System voltage product  (V$^2$)  & \VsqkSDP\in  \setHermitianMatrix{|\setBuspairs||\setPhases|} \\ 	
Bus pair  voltage product  (V$^2$) &  \VbpsqkSDP\in  \setHermitianMatrix{2|\setPhases|} \\ 	
Branch voltage-current product   (mix) &  \VlinesqkSDP\in  \setHermitianMatrix{2|\setPhases|}\\ 	
\hline
\end{IEEEeqnarraybox}
\end{equation*}
\end{table}

 \section{Introduction }
Driven by the increased rollout of distributed energy resources (DERs) such as PV, battery storage as well as electric vehicles, power distribution grids are facing a number of challenges associated with the large scale integration of these technologies.
In low-voltage grids specifically, one can observe phase unbalance readily, due to the presence of a significant amount of single-phase loads and DERs. Many of the current electric vehicles use single phase charging and small rooftop solar systems are connected via single phase inverters to the low voltage grid, further increasing the  unbalance. Furthermore, phase  unbalance can stem from insufficient conductor transposition in radial distribution networks. 
Unbalance implies an underutilization of the grid's transfer capacity, as it leads to higher losses and to faster-than-expected congestion. Simulation techniques, i.e. deriving solutions to the power flow  (PF) equations in unbalanced networks has long been a topic of interest~\cite{Berg1967,Chen1991,Cheng1995}, and is used in deriving hosting capacity by means of scenario analysis.  The \emph{accuracy} of modeling of low-voltage grids has  been studied in-depth by Urquhart \cite{Urquhart2016}. 

Unbalanced \emph{optimal} power flow (OPF) refers to the mathematical optimization of problems subject to the physics of unbalanced grids, and serves as the core for a variety of problem classes such as benchmarking centralized or distributed optimal control solutions \cite{Dkhili2020}, determination of  expansion options and  hosting capacity analysis considering control actions.


\subsection{State of the Art on Unbalanced (O)PF }
Similar to the balanced (positive sequence) modeling, one can distinguish between the bus injection model (BIM) and the branch flow model (BFM) formulations of the unbalanced PF equations. BIM forms eliminate all current variables, which leads to active and reactive power flows being expressed purely as a function of the voltage differences between connected buses. The series impedance is consequently used in admittance form, which makes it impossible to represent zero-impedance branches. Conversely, the BFM forms keep (a representation of) the current variable through the series impedance. In this case, series impedance is represented in impedance form, therefore the edge case of zero series impedance remains representable. 

Most OPF problems are developed in the complex power-voltage variable space, instead of the current-voltage variable space common in PF solvers. 
Table \ref{tab_formulations} maps a number of published formulations to BIM/BFM categories and the variable spaces in which they are defined. 
We refer to \cite{Molzahn2017a,Bienstock2020}  for recent in-depth reviews of  mathematical formulations for the OPF problem. 

A nonlinear programming (NLP) unbalanced OPF formulation for  branches is presented in  \cite{Mahdad2006}, in which the shunt impedances have been neglected. The unbalanced current-voltage form, with generation power dispatch constraints, is derived in \cite{Geth2020b}. The nonconvex power-lifted voltage form is presented in \cite{Geth2020d}.
The first rank-constrained semi-definite programming (SDP) of unbalanced OPF, and its SDP relaxation, are proposed by Dall'Anese et al.  in \cite{Anese}. 
Gan et al. developed the BFM variant, and  defined both BFM and BIM in a consistent notation in  \cite{Gan2014}, without the provision of the real-value form. 
Extensions to these SDP relaxations have been proposed:
Zhao et al. develop models for delta-connected loads \cite{Zhao2017a};
Bazrafshan et al. propose extensions for voltage regulators \cite{Bazrafshan2018a}; 
Usman et al. discuss neutral conductor modeling \cite{Usman2020}; 
Claeys et al. present detailed transformer models for unbalanced OPF \cite{Claeys2020b};  
Claeys et al. also present convex relaxation of delta/wye ZIP unbalanced loads \cite{Claeys2021a};  
Vanin et al. explore further relaxation to SOC problems \cite{Vanin2020}.

It is noted that the formulations so far published are not easily compared due to the inherent complexity of the notation, the variety of notations used, and lacking details of approximations applied during implementation in modeling software. 
For instance, the branch shunt impedance is often neglected, assumed to be diagonal, or simply modeled as bus shunts. 
This leads to inaccuracies in the power and current flow values and in the enforcement of the proper branch flow limits, e.g. when only the series current is bounded, not the total including the shunt current. 
This in turn makes it hard to validate the \emph{feasibility} and correctness of results on published data sets which include such components without modification. Examples include the IEEE PES distribution test feeders \cite{Mather2017}; the IEEE123 bus system specifically contains non-diagonal branch shunt matrices.

\subsection{Bound definitions and data}
In an OPF problem, authors typically assume \cite{Molzahn2017a,Bienstock2020,Coffrin2015f} the power flow is subject to a subset of:
\begin{itemize}
    \item apparent power, upper bound;
    \item current magnitude, upper bound;
    \item voltage magnitude, lower and upper bound;
    \item voltage angle difference between adjacent buses, lower and upper bound.
\end{itemize}
These bounds can be generalized to the case with phase unbalance. Furthermore, in the context of phase unbalance, additional limits have been discussed \cite{Girigoudar}. The most obvious missing bound is that of phase angle differences at a certain bus.
Note that bounds represent a key aspect of optimization problems, however currently there are few unbalanced power flow data sets that include all of these parameters. This makes benchmarking \emph{optimality} of different unbalanced OPF engines challenging.  

\subsection{Implementation of Complex-Value  Optimization Problems}

\begin{figure}
  \centering
    \includegraphics[width=0.35\columnwidth]{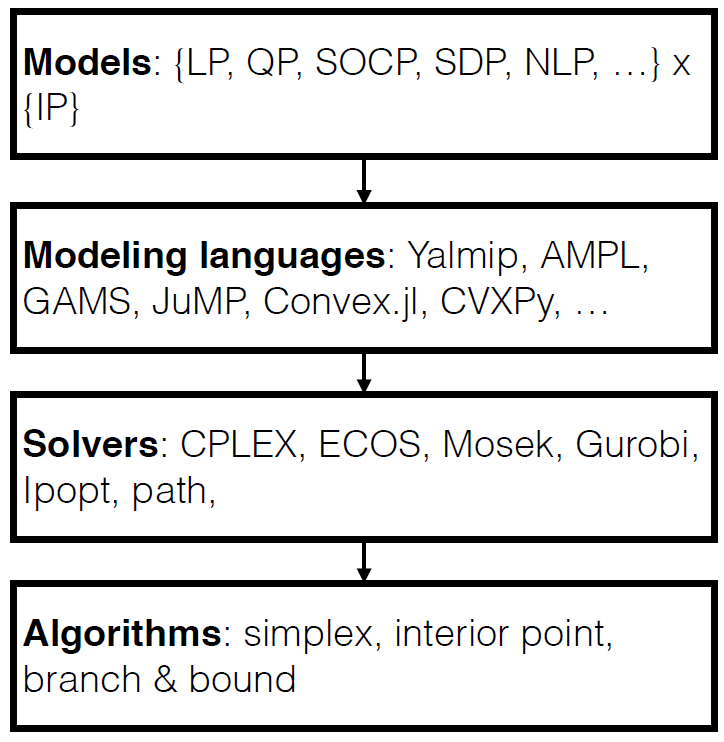}
  \caption{Formulation, implementation and solution of optimization problems.}  \label{fig_opf_flowchart}
\end{figure}

Figure~\ref{fig_opf_flowchart} illustrates the typical sequence of formulating, implementing and solving optimization problems as commonly applied to engineering problems. 
First, a mathematical model of the optimization problem is formulated. 
The problem  is formulated so that it can be handled by the chosen modeling language and optimization solver. 
Eventually, the modeling toolbox  translates the mathematical model into the appropriate form for the solver interface,   fills out parameter values, and dispatches the  optimization solver. 
Finally, the mathematical solver uses one or a variety of different algorithms to solve the optimization problem, and returns a solution.

In the literature, real-value formulations of the convex relaxations of unbalanced OPF are rarely discussed. 
CVX has been used as an optimization  modeling tool in \cite{Anese}. 
It is one of the very few optimization packages with support for complex optimization variable primitives, as well as matrix variables. 
CVX automates the reformulation process from a complex-value matrix to a set of real-value scalar variables and constraints. 
Functionality  to automate the real-value reformulation process was initially in a  JuMP  extension \cite{Sliwak2019}, but has recently become a built-in JuMP feature \cite{lubin_jump_2023}. 

Gilbert and Josz \cite{Gilbert2016a} illustrate with a limited-scope solver that there are potential advantages of solving SDP problems directly in the complex domain. 
Nevertheless, a mature implementation of such a solver has not been developed so far. 
Overall, the interest in real-value forms remains high.

\newgeometry{margin=1cm} 
\begin{landscape}


\begin{table*}
  \centering 
   \caption{Unbalanced power flow formulations, key references and variable spaces}\label{tab_formulations}
    \begin{tabular}{l  l  l  l  l  l  }
\hline
&Formulation & Ref. & Variable space  & Exact? & Section \\\hline
NLP& Kirchhoff unbalanced 	& \cite{Berg1967,Geth2020b} & $\IijkSDP, \VikSDP $		 & by definition, discussion: \cite{Urquhart} & \S\ref{sec_unbalanced_kirchhoff} \\
NLP& BIM unbalanced  scalar polar	&  \cite{Mahdad2006}& $\diag(\SijkSDP),\VikSDP $		 & if $\VikSDPmin \neq \mathbf{0}$ & \S \ref{sec_unbalanced_bim}\\
NLP& BIM unbalanced  scalar rect.	&  & $\diag(\SijkSDP),\VikSDP $		 & if $\VikSDPmin \neq \mathbf{0}$ &  \S \ref{sec_unbalanced_bim}\\
SDP& BIM unbalanced   	&\cite{Anese,Liu2018c,Gan2014b} & $\SijkSDP, \VisqkSDP  \succeq 0, \VijsqkSDP ,\VsqkSDP \succeq 0 $		  & $\rank(\VsqkSDP)=1$ &\S\ref{sec_bim_relax}  \\
SDP& BIM unbalanced  radial 	& \cite{Gan2014} & $\SijkSDP, \VisqkSDP \succeq 0, \VijsqkSDP, \VbpsqkSDP \succeq 0 $		  &   $\rank(\VbpsqkSDP) = 1$, only  if radial &\S\ref{sec_bim_relax} \\
SDP& BFM unbalanced     	& \cite{Gan2014}  & $ \SijkSDP, \VisqkSDP \succeq 0, \IijsqskSDP \succeq 0, \VlinesqkSDP \succeq 0 $		 &$\rank(\VlinesqkSDP)=1$, only  if radial & \S\ref{sec_unbal_bfm}   \\
\hline
\end{tabular}
\end{table*}

\end{landscape}
\restoregeometry

\subsection{Scope, Contributions and Report Structure}

In summary, we observe  the following gaps in the literature:
\begin{itemize}
    \item polar and rectangular real forms in the power-voltage variable space for  branches \emph{with shunts} ($\Pi$-model);
    \item real form of the power - lifted voltage BIM and BFM SDP relaxations; 
    \item definitions of current, power and voltage angle difference limits in both the power-voltage and power-lifted voltage variable spaces.
\end{itemize}

This work therefore provides consistent derivations of real-value formulations for variants of   unbalanced optimal power flow  in different variable spaces, to enable straight-forward  implementation in  optimization modeling toolboxes and  efficient comparison on published data sets.

We limit ourselves here to the three-wire case, and note that recent works \cite{CLAEYS2022108522,Geth2023} have explored OPF models \emph{with explicit representation of the neutral}.


First, the basic relationship between the different variables and parameters is defined in  \S \ref{sec_preliminaries}, which builds the mathematical foundation of this report; it also contains a generalized model to represent loads and generators in the system. 
Section \S \ref{sec_unbalanced_kirchhoff} provides  Kirchhoff's and Ohm's laws in the multiconductor form.
Next, \S \ref{sec_nlp_bfm} develops a formulation of nonlinear BFM unbalanced power flow and 
\S \ref{sec_unbalanced_bim} the equivalent BIM.
Furthermore, \S \ref{sec_bim_relax} derives the BIM and \S \ref{sec_unbal_bfm} the BFM in the lifted variable space which is the basis various convex relaxations. 
Moreover \S \ref{sec_feasible_sets} presents the feasible sets of the derived real-value formulations and discusses implementation aspects.
Finally, \S \ref{sec_conclusions} presents the conclusions.

 \section{Notation and Basic Relationships} \label{sec_preliminaries}

\subsection{Scalar and Matrix Variables}
  The  voltage $\VikSDP$ of a bus $i$ is a complex value, vector variable, encapsulating variables for each conductor:
    \begin{IEEEeqnarray}{RL?s}
  \VikSDP= 
  \begin{bmatrix} 
\VikA    \\
\VikB     \\
\VikC   
 \end {bmatrix} 
 =
 \begin{bmatrix} 
 \VikANmag    \imagangle \VikANangle \\
 \VikBNmag    \imagangle \VikBNangle \\
 \VikCNmag    \imagangle\VikCNangle
 \end {bmatrix} \label{eq_var_voltage_polar} 
 = \VikSDPmag \imagangle \VikSDPang
 =
 \begin{bmatrix} 
\VikANreal + \imagnumber \VikANimag    \\
\VikBNreal+ \imagnumber\VikBNimag     \\
\VikCNreal+ \imagnumber\VikCNimag  
 \end {bmatrix}
 = \VikSDPreal + \imagnumber \VikSDPimag
 \label{eq_var_voltage_rectangular}
   \end{IEEEeqnarray}   
  Similarly,  the total current of a branch $l$ connecting a pair of buses $i$ and $j$ is a complex value vector variable,
 \begin{IEEEeqnarray}{C?s}
  \IijkSDP = 
 \begin{bmatrix} 
\IijkA    \\
\IijkB     \\
\IijkC  
 \end {bmatrix}     .\nonumber
    \end{IEEEeqnarray}  
The  sending-side current ($i \rightarrow j$) through the series element in the $\Pi$-model is defined as $\IijskSDP$, the sending-side shunt  current as $ \IijshkSDP$, respectively. 
  The complex power flow in branch $l$ from bus $\indexGridNode$ to $\indexGridNodeTwo$  depends on the bus voltage of the sending-side $\VikSDP$ and the conjugate transpose (indicated with superscript $\hermitiantranspose$) of the current  $\IijkSDP$,
 \begin{IEEEeqnarray}{L?s}
\SijkSDP = 
 \VikSDP (\IijkSDP)^{\hermitiantranspose}
= \left[\begin{array}{ccc} 
\SijkAA\! \!& \SijkAB\!\! & \SijkAC  \\
\SijkBA\! \!& \SijkBB\!\! & \SijkBC \\
\SijkCA\! \!& \SijkCB\! \!& \SijkCC 
\end{array} \right].
  \end{IEEEeqnarray}
    We observe that $ \diag(\SijkSDP  )=\VikSDP \circ (\IijkSDP)^*$, 
    where $\circ$ is the element-wise (Hadamard) product. 
    It is noted that $\rank(\SijkSDP) = 1$, as it is defined as the outer product of two vectors.
The off-diagonals relate to the diagonal elements according to,
\begin{IEEEeqnarray}{C?s}
\frac{\SijkAA}{\VikA} = \frac{\SijkBA}{\VikB} = \frac{\SijkCA}{\VikC} = (\IijkA)^*, \IEEEyesnumber \IEEEyessubnumber  \\
\frac{\SijkBB}{\VikB} = \frac{\SijkAB}{\VikA} = \frac{\SijkCB}{\VikC} = (\IijkB)^*,  \IEEEyessubnumber \\
\frac{\SijkCC}{\VikC} = \frac{\SijkAC}{\VikA} = \frac{\SijkBC}{\VikB} = (\IijkC)^* ,  \IEEEyessubnumber
    \end{IEEEeqnarray}
    which means that the off-diagonals of $\SijkSDP$ are scaled and rotated versions of the more easily interpretable diagonal elements. 
    \subsection{Variable Bounds}
Voltage magnitudes have minimum and maximum operational limits, which are specific to each bus and phase,
 \begin{IEEEeqnarray}{C?s}
 \mathbf{0} \le   
 \underbrace{\begin{bmatrix} 
   \VikAmin \\
   \VikBmin \\
   \VikCmin\\
    \end{bmatrix}}_{\VikSDPmin}  \leq 
  \begin{bmatrix} 
|\VikA|    \\
|\VikB|     \\
|\VikC|  \\
 \end {bmatrix} 
 \leq
\underbrace{  \begin{bmatrix} 
   \VikAmax \\
   \VikBmax \\
   \VikCmax\\
 \end{bmatrix}}_{\VikSDPmax} . \label{eq_voltage_bounds_polar}
    \end{IEEEeqnarray}  
        In this report, we overload `$\leq$' for `$\geq$' for vectors \emph{and} matrices to indicate \emph{element-wise inequality} (conversely, `$\succeq$' is used as the symbol for matrix positive semidefiniteness).
Recognizing we obtain the magnitude squared by multiplying complex numbers with their own conjugates, the nodal voltage bounds can also be presented quadratically as:
  \begin{IEEEeqnarray}{C?s}
\mathbf{0} \le \VikSDPmin \circ \VikSDPmin  \leq \VikSDP \circ (\VikSDP)^* \leq \VikSDPmax \circ \VikSDPmax .  \label{eq_voltage_bounds_matrix}
    \end{IEEEeqnarray}  
    Note that the upper bound constraints are convex (space inside a circle), but the lower bound constraints are nonconvex (space outside of a circle). 
Next,  apparent power limits  
 are defined using the absolute value of the diagonals of the branch flow matrix,
 \begin{IEEEeqnarray}{C?s}
\mathbf{0} \le \begin{bmatrix} 
|\SijkAA|  \\
 |\SijkBB|\\
|\SijkCC| \\
 \end{bmatrix} \leq
 \begin{bmatrix} 
  \SijratedA \\
 \SijratedB \\
   \SijratedC\\
 \end{bmatrix} = \SijkSDPrated, 
     \end{IEEEeqnarray}  
which is equivalent to second order cone (SOC) constraints,
  \begin{IEEEeqnarray}{C?s}
  \begin{bmatrix} 
   |\SijkAA|^2  \\
   |\SijkBB|^2 \\
   |\SijkCC|^2  \\
 \end{bmatrix}
 =
\begin{bmatrix} 
   (\PijkAA)^2 + (\QijkAA)^2  \\
   (\PijkBB)^2 + (\QijkBB)^2\\
   (\PijkCC)^2 + (\QijkCC)^2 \\
 \end{bmatrix} \leq
 \begin{bmatrix} 
 (  \SijratedA )^2\\
(   \SijratedB )^2\\
(   \SijratedC)^2\\
 \end{bmatrix}, 
    \end{IEEEeqnarray}  
and can succinctly be written as,
  \begin{IEEEeqnarray}{C?s}
\mathbf{0}\le \diag(\SijkSDP) \circ \diag(\SijkSDP)^* \leq \SijkSDPrated \circ \SijkSDPrated . \label{eq_power_bounds}
    \end{IEEEeqnarray}  
    The magnitudes of the branch current should stay below rated values and are bounded for the diagonals of the branch current matrix,
 \begin{IEEEeqnarray}{C?s}
\mathbf{0}\le    \begin{bmatrix} 
|\IijkA |   \\
|\IijkB|     \\
|\IijkC  | \\
 \end {bmatrix}  \leq  
 \begin{bmatrix} 
   \IijratedA \\
   \IijratedB \\
   \IijratedC\\
 \end{bmatrix} =  \IijkSDPrated \label{ineq_current_a},
       \end{IEEEeqnarray}  
which again can be written as a set of SOC constraints,
  \begin{IEEEeqnarray}{C?s}
\mathbf{0} \le \IijkSDP \circ (\IijkSDP)^* \leq \IijkSDPrated \circ \IijkSDPrated . \label{eq_current_ratings_matrix}
    \end{IEEEeqnarray}  
A valid SOC representation of current magnitude limits  \eqref{ineq_current_a} that does not require explicit current variables $\IijkSDP$ is useful for BIM forms. 
Using the nodal voltage magnitudes, the branch current limits can directly be enforced on the power flows, 
  \begin{IEEEeqnarray}{C?s}
  \begin{bmatrix} 
   |\SijkAA|^2  \\
   |\SijkBB|^2 \\
   |\SijkCC|^2  \\
 \end{bmatrix} \leq
 \begin{bmatrix} 
 (  \IijratedA )^2\\
(   \IijratedB )^2\\
(   \IijratedC)^2\\
 \end{bmatrix}\circ
   \begin{bmatrix} 
|\VikA|^2    \\
|\VikB|^2     \\
|\VikC|^2  \\
 \end {bmatrix}, 
    \end{IEEEeqnarray}  which can also be developed in matrix notation  as, 
  \begin{IEEEeqnarray}{C?s}
 \diag(\SijkSDP) \circ \diag(\SijkSDP)^* \leq \IijkSDPrated \circ \IijkSDPrated \circ \VikSDP \circ (\VikSDP)^*. \label{eq_branch_currents_WS}
  \end{IEEEeqnarray}
  Valid bounds on all elements of the power flow matrix $\SijkSDP = \PijkSDP + \imagnumber \QijkSDP$ are,
        \begin{IEEEeqnarray}{C?s}
- \VikSDPmax (\IijkSDPrated)^{\transpose} \leq \PijkSDP, \QijkSDP \leq \VikSDPmax (\IijkSDPrated)^{\transpose} , \label{eq_complex_s_bounds}
   \end{IEEEeqnarray}  
   both for active and reactive power components, respectively.
    The voltage angle differences between connected buses $i$ and $j$ are bounded,
\begin{IEEEeqnarray}{C?s}
\!     \begin{bmatrix} 
-\frac{\pi}{2} \\
-\frac{\pi}{2} \\
-\frac{\pi}{2} 
 \end {bmatrix}  \leq
\underbrace{     \begin{bmatrix} 
\thetaijAAmin \\
\thetaijBBmin \\
\thetaijCCmin
 \end {bmatrix}}_{\thetaijSDPmin}
 \leq    \begin{bmatrix} 
 \VikANangle -\VjkANangle     \\
 \VikBNangle    -\VjkBNangle     \\
 \VikCNangle   -\VjkCNangle   
 \end {bmatrix}  \leq 
\underbrace{    \begin{bmatrix} 
\thetaijAAmax \\
\thetaijBBmax \\
\thetaijCCmax
 \end {bmatrix} }_{\thetaijSDPmax}
 \leq
   \begin{bmatrix} 
\frac{\pi}{2} \\
\frac{\pi}{2} \\
\frac{\pi}{2} 
 \end {bmatrix} \label{eq_angle_diff_eq}.
\end{IEEEeqnarray}

For a relatively balanced voltage phasor, we expect $ \VikANangle -\VikBNangle \approx  \VikBNangle    -\VikCNangle  \approx  \VikCNangle   -\VikANangle  \approx \frac{2\pi}{3}$. 
The voltage angle differences between phases on buses $i$ can be bounded to enforce angle balance relative to the expected 120 degrees,
\begin{IEEEeqnarray}{C?s}
\underbrace{     \begin{bmatrix} 
\thetaiAAmin \\
\thetaiBBmin \\
\thetaiCCmin
 \end {bmatrix}}_{\thetaiSDPmin}
 \leq    \begin{bmatrix} 
 \VikANangle -\VikBNangle -\frac{2\pi}{3}    \\
 \VikBNangle    -\VikCNangle   -\frac{2\pi}{3}  \\
 \VikCNangle   -\VikANangle  -\frac{2\pi}{3} 
 \end {bmatrix}   \leq 
\underbrace{    \begin{bmatrix} 
\thetaiAAmax \\
\thetaiBBmax \\
\thetaiCCmax
 \end {bmatrix} }_{\thetaiSDPmax} 
 \label{eq_angle_diff_eq_bus}.
\end{IEEEeqnarray}
  The voltage phasor in  reference buses $i \in \setReferenceNodes \subset \omegaGridNodes$ is assumed fixed, e.g.,
    \begin{IEEEeqnarray}{C?s}
\VikSDP=  \VikrefSDP=  
 \begin{bmatrix} 
\VikANref \imagangle   \thetaikANref \\
\VikBNref \imagangle \thetaikBNref   \\
\VikCNref  \imagangle\thetaikCNref  
\end {bmatrix} .\label{eq_ref_bus_voltage}
  \end{IEEEeqnarray}

   \subsection{Branch Impedance}
The  circuit series impedance matrix  is defined as a full matrix with no assumption on the particular structure:
\begin{IEEEeqnarray}{C?s}
\zijksSDP = \rijksSDP + \imagnumber  \xijksSDP
=
\begin{bmatrix} 
\zijksAA  &\zijksAB  &\zijksAC    \\
\zijksBA  &\zijksBB  &\zijksBC   \\
\zijksCA  & \zijksCB  & \zijksCC   \\
 \end {bmatrix} , \label{eq_4x4_impedance_def}
  \end{IEEEeqnarray}  
  where each element consists of a series resistive and reactive impedance. It is noted that in physical systems, we expect $\rijksSDP \succeq 0, \xijksSDP \succeq 0$.
  The impedance matrix can be rewritten in the corresponding admittance form,
  \begin{IEEEeqnarray}{C?s}
  \yijksSDP = {(\zijksSDP)^{-1}} =  \gijksSDP + \imagnumber  \bijksSDP . \nonumber
  \end{IEEEeqnarray}  
  where $ {(\zijksSDP)^{-1}}$ is the matrix inverse of $\zijksSDP$. In case of missing conductors, e.g. single, or two-conductor connections,  $\zijksSDP$ is not invertible but it is valid to use the Moore-Penrose inverse instead. 
The shunt admittances at the sending and receiving sides, respectively $\yijkshSDP, \yjikshSDP$, are defined, 
 \begin{IEEEeqnarray}{C?s}
 \yijkshSDP  = \gijkshSDP+ \imagnumber  \bijkshSDP 
= \begin{bmatrix} 
\yijkshAA   & \yijkshAB & \yijkshAC \\
\yijkshBA &  \yijkshBB    & \yijkshBC\\
\yijkshCA&  \yijkshCB  &  \yijkshCC  \\
 \end {bmatrix} .
  \end{IEEEeqnarray}  
 Although  for typical distribution lines and cables the shunt admittances are diagonal and have equal values for both sides, in this generalized derivation they can be considered as a two different and full matrices. This also allows for re-use of the representation for other elements, such as  transformers. 

\subsection{Loads and Generators as Units}
We define units $\indexUnit \in \setUnits$ to generalize loads, generators and storage elements. The current flowing from the connected bus into the unit is
%
      \begin{IEEEeqnarray}{C?s}
 \IuunitSDP=  
 \begin{bmatrix} 
 \IuunitAN  \\
\IuunitBN \\
\IuunitCN  \\
 \end {bmatrix}.
     \end{IEEEeqnarray}  
     The unit current is bounded by the corresponding current rating $\IuunitSDPrated$ similarly to \eqref{ineq_current_a}. 
Tuples of units and the buses they are connected to are defined in the connectivity set $\indexUnit \indexGridNode \in  \setUnitTopology $. 
The active and reactive power consumed by a unit is defined $\SuunitSDP = \PuunitSDP + \imagnumber \QuunitSDP = \VikSDP (\IuunitSDP)^\hermitiantranspose$. We  define bounds on active/reactive power dispatch separately,
\begin{IEEEeqnarray}{C?s}
\PuunitSDPmin \leq \diag(\PuunitSDP) \leq \PuunitSDPmax \IEEEyesnumber \IEEEyessubnumber \label{eq_unit_active_power},\\
\QuunitSDPmin \leq \diag( \QuunitSDP) \leq \QuunitSDPmax  \IEEEyessubnumber  \label{eq_unit_reactive_power}.
\end{IEEEeqnarray}

We do not provide any further detail on the specific modeling of delta or wye-connected units and refer to \cite{Zhao2017a, Claeys2021a}. 
\subsection{Shunts}
A shunt element, e.g. shunt capacitance, inductance or resistance, with index $\indexShunt$ has an admittance  $\ybSDP = \gbSDP + \imagnumber \bbSDP$. 
Tuples of shunts and the bus they are connected to are defined in the connectivity set $\indexShunt\indexGridNode \in \setShuntTopology$ .
   The current  from the bus to the shunt is, 
         \begin{IEEEeqnarray}{C?s}
 \IbSDP = 
 \begin{bmatrix} 
 \IbA  \\
\IbB \\
\IbC  \\
 \end {bmatrix}.
     \end{IEEEeqnarray}  
     and is bounded by the rated current  $\IhshuntSDPrated$.

   \section{Application of Circuit Laws} \label{sec_unbalanced_kirchhoff}
This section illustrates how Kirchhoff's and Ohm's laws are used in the context of multi-conductor branches with matrix impedances for the representation of the circuit physics, providing the relationship between the voltage, current and power variables.  
  \subsection{Branch Model}

\begin{figure}
  \centering
    \includegraphics[width=0.65\columnwidth]{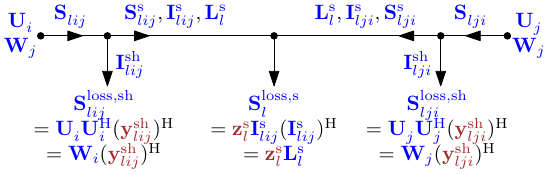}
  \caption{Power flow and losses in unbalanced  $\Pi$-model branch}  \label{fig_linemodel}
\end{figure}

Ohm's law for branch $lij$ representing the voltage drop along the branch is formulated in matrix form using nodal voltages, the branch current and the impedance matrix,  
    \begin{IEEEeqnarray}{C?s}
\VjkSDP     =\VikSDP  - \zijksSDP \IijskSDP     .\label{eq_ohms_matrix}
          \end{IEEEeqnarray}
  Kirchhoff's current law (KCL) is used to split up the series and shunt (to ground) currents in the $\Pi$-section,
          \begin{IEEEeqnarray}{C?s}
\IijkSDP =   \IijskSDP + \IijshkSDP \label{eq_total_current_from_def}, 
\IjikSDP =   \IjiskSDP + \IjishkSDP, \\
\IijskSDP + \IjiskSDP = 0. \label{eq_series_current_variables}
      \end{IEEEeqnarray}
Fig.~\ref{fig_linemodel} defines the harmonized single-wire equivalent  \emph{vector/matrix} variables in both natural and lifted (defined in upcoming sections)  variable spaces for a clear presentation of the basic relationships. The scalar representation of all variables and parameters is defined in Fig.~\ref{fig_linemodel_iv_3x3}.

\subsection{Shunts and Units at Nodes}

The unit model can be adapted for a variety of optimization problems, e.g. optimal dispatch of distributed generation, optimal scheduling of storage and electric vehicle charging and optimal demand management. Such extensions essentially define feasible sets dependent on $\diag(\PuunitSDP), \diag(\QuunitSDP)$. Without loss of generality, in this work, we focus on the feasible sets defined by \eqref{eq_unit_active_power}-\eqref{eq_unit_reactive_power}. 

Similarly, the shunt power $\SbSDP$ is defined as 
\begin{IEEEeqnarray}{C?s}
\SbSDP = \PbSDP + \imagnumber \QbSDP=   \VikSDP (\IbSDP)^{\hermitiantranspose}. \label{eq_bus_shunt_power}
\end{IEEEeqnarray}
The shunt current relates to the voltage through Ohm's law,
\begin{IEEEeqnarray}{C?s}
\IbSDP =  \ybSDP  \VikSDP . \label{eq_bus_shunt_current}
\end{IEEEeqnarray}
Substituting \eqref{eq_bus_shunt_current} into \eqref{eq_bus_shunt_power}, the shunt power consumption can be expressed  depending on the nodal voltage only,
\begin{IEEEeqnarray}{C?s}
\SbSDP =  \VikSDP (\VikSDP)^{\hermitiantranspose}  (\ybSDP)^{\hermitiantranspose} . \label{eq_shunt_power_def}
\end{IEEEeqnarray}

\subsection{Bus Model: Kirchhoff's Current Law}
Kirchhoff's Current Law (KCL) is conventionally expressed in current variables, but can also be lifted to the complex power variable space. 
KCL in current variables at each bus $i$ is,
    \begin{IEEEeqnarray}{C?s}
\sum_{lij \in \setTopologyBoth} \IijkSDP +  \sum_{\indexUnit \indexGridNode \in  \setUnitTopology } \IuunitSDP   +  \sum_{ bi \in \setShuntTopology } \IbSDP   = 0.
\label{eqnodalcurrentbalance} 
   \end{IEEEeqnarray}  
KCL in the complex power variable space is obtained by taking the conjugate transpose of  \eqref{eqnodalcurrentbalance} and element-wise multiplying with $\VikSDP\neq 0$ on the left,
    \begin{IEEEeqnarray}{C?s}
 \!\!\!\sum_{lij \in \setTopologyBoth} \!\!\diag(\SijkSDP) + \!\! \!\!\sum_{\indexUnit \indexGridNode \in  \setUnitTopology }\!\! \!\!\diag( \SuunitSDP ) + \!\!\!\! \sum_{ bi \in \setShuntTopology } \!\!\!\!\diag( \SbSDP)=0, \label{eq_power_balance_nodes}
   \end{IEEEeqnarray}  
   which means that the diagonal elements of the apparent power matrices of the connected branch flows, units and shunt need to sum to zero for each bus. Commonly, the shunt power expression \eqref{eq_shunt_power_def} is substituted into this equation. The real-value equivalent forms are obtained,
          \begin{IEEEeqnarray}{C?s}
 \!\!\!\sum_{lij \in \setTopologyBoth}\!\!\diag( \PijkSDP ) +\!\!\!\!  \sum_{\indexUnit \indexGridNode \in  \setUnitTopology }\!\!\!\!\diag( \PuunitSDP)   + \!\!\!\! \sum_{ bi \in \setShuntTopology }\!\!\!\! \diag( \PbSDP) =0, \label{eq_kcl_active_node} \IEEEyesnumber \IEEEyessubnumber\\
 \!\!\!\sum_{lij \in \setTopologyBoth} \!\!\diag(\QijkSDP) + \!\!\!\! \sum_{\indexUnit \indexGridNode \in  \setUnitTopology } \!\!\!\!\diag(\QuunitSDP)  +\!\!\!\!  \sum_{ bi \in \setShuntTopology }\!\!\!\! \diag( \QbSDP) =0.  \label{eq_kcl_reactive_node} \IEEEyessubnumber 
   \end{IEEEeqnarray}

    \section{Unbalanced Branch Flow Model} \label{sec_nlp_bfm}
    This section illustrates how the unbalanced BFM is derived. 
First, we define a variable for the complex power flow $\SijksSDP$ through the \emph{series} element of the $\Pi$-section, 
\begin{IEEEeqnarray}{C?s}
  \SijksSDP = \PijksSDP + \imagnumber  \QijksSDP = \VikSDP( \IijskSDP)^{\hermitiantranspose} .
  \label{eq_complex_power_flow_unbalanced} 
 \end{IEEEeqnarray}
              \subsection{Power Flow Model}
 Passive components in electrical circuits cause losses. In the BFM we distinguish between the losses associated with the from-side shunt, i.e. $\SijkshlossSDP$,  the series impedance, i.e. $\SijkslossSDP$, and the to-side shunt, i.e. $\SjikshlossSDP$. 
   The series losses  depend on the voltage drop over the impedance and the current flow through it, 
\begin{IEEEeqnarray}{C?s}
\SijkslossSDP  = (\VikSDP - \VjkSDP ) (\IijskSDP)^{\hermitiantranspose}. \label{eq_series_power_losses}
\end{IEEEeqnarray}
The voltage drop itself is can be derived from Ohm's law, 
\begin{IEEEeqnarray}{C?s}
(\VikSDP - \VjkSDP )  =  \zijksSDP   \IijskSDP. \label{eq_ohms_series}
\end{IEEEeqnarray}
Substituting \eqref{eq_ohms_series} into \eqref{eq_series_power_losses} we obtain the series loss
\begin{IEEEeqnarray}{C?s}
\SijkslossSDP  =  \zijksSDP   \IijskSDP (\IijskSDP)^{\hermitiantranspose} = \zijksSDP   \IjiskSDP (\IjiskSDP)^{\hermitiantranspose} ,   \label{eq_loss_series}
\end{IEEEeqnarray}
which is symmetric because of \eqref{eq_series_current_variables}. 
The shunt losses $\SijkshlossSDP$ are derived from the current through the shunt, i.e. $\IijshkSDP$, and the voltage at the shunt, i.e. $\VikSDP$,
\begin{IEEEeqnarray}{C?s}
\SijkshlossSDP =   \VikSDP (\IijshkSDP)^{\hermitiantranspose}, \label{eq_shunt_loss_power}
\end{IEEEeqnarray}
The shunt current relates to the voltage through Ohm's law,
\begin{IEEEeqnarray}{C?s}
\IijshkSDP =  \yijkshSDP  \VikSDP . \label{eq_shunt_current}
\end{IEEEeqnarray}
Substituting \eqref{eq_shunt_current} into \eqref{eq_shunt_loss_power}, the shunt losses are,
\begin{IEEEeqnarray}{C?s}
\SijkshlossSDP =  \VikSDP (\VikSDP)^{\hermitiantranspose}  (\yijkshSDP)^{\hermitiantranspose} , \,
\SjikshlossSDP =   \VjkSDP ( \VjkSDP)^{\hermitiantranspose}  (\yjikshSDP)^{\hermitiantranspose} . \label{eq_loss_shunt}
\end{IEEEeqnarray}
%
As the the sum of the sending and receiving side power flows need to equal the branch losses, the branch loss balance (see Fig.~\ref{fig_linemodel}) can be written using the different loss components \eqref{eq_loss_series} and \eqref{eq_loss_shunt},
\begin{IEEEeqnarray}{L?s}
      \SijkSDP + \SjikSDP = \SijkshlossSDP + \SijkslossSDP + \SjikshlossSDP  \nonumber \\ 
      = \VikSDP ( \VikSDP)^{\hermitiantranspose} (\yijkshSDP)^{\hermitiantranspose}  + \zijksSDP   \IijskSDP (\IijskSDP)^{\hermitiantranspose} +  \VjkSDP (\VjkSDP)^{\hermitiantranspose} (\yjikshSDP)^{\hermitiantranspose} \label{eq_bfm_line_power_balance}. 
\end{IEEEeqnarray}

     \section{Unbalanced Bus Injection Model} \label{sec_unbalanced_bim}
         This section illustrates how the unbalanced BIM is obtained for which the branch flow current is substituted by nodal voltage varibales. 
             \subsection{Power Flow Model}
Therefore, we write \eqref{eq_ohms_series} in admittance form,
    \begin{IEEEeqnarray}{C?s}
\IijskSDP = \yijksSDP(\VikSDP - \VjkSDP). \label{eq_volt_diff}
   \end{IEEEeqnarray}  
The sending end apparent power flow can be calculated using the sending end voltage and the sum of the series current and the sending end shunt current, respectively,
\begin{equation}
\SijkSDP = \VikSDP(   \IijskSDP + \IijshkSDP)^{\hermitiantranspose}. \label{eq_flow_from}
\end{equation}

By substituting \eqref{eq_shunt_current} and \eqref{eq_volt_diff} into \eqref{eq_flow_from} we derive,
\begin{equation}
\SijkSDP = \VikSDP(        \VikSDP)^{\hermitiantranspose} ( \yijkshSDP       )^{\hermitiantranspose} +  \VikSDP  (\VikSDP - \VjkSDP)^{\hermitiantranspose} (\yijksSDP)^{\hermitiantranspose} \label{eq_bim_complex_power},
\end{equation}  
   which is the nonlinear complex matrix form of the BIM. Note that this constraint is defined $\forall \indexGridLines \indexGridNode\indexGridNodeTwo \in \setTopologyBoth$, i.e. for  both the sending end $\SijkSDP$ and the receiving end $\SjikSDP$, respectively. 
We now choose rectangular coordinates for the voltage and power variables and obtain,
\begin{IEEEeqnarray}{C?s}
\PijkSDP  = 
\left(\VikSDPreal(\VikSDPreal)^{\transpose} + \VikSDPimag(\VikSDPimag)^{\transpose}\right) ( \gijkshSDP       )^{\transpose} 
+ \left(\VikSDPimag(\VikSDPreal)^{\transpose} - \VikSDPimag(\VikSDPimag)^{\transpose}\right) ( \bijkshSDP       )^{\transpose} \nonumber \\
+ \left( 
\VikSDPreal  (\VikSDPreal - \VjkSDPreal)^{\transpose} 
+  \VikSDPimag  (\VikSDPimag - \VjkSDPimag)^{\transpose}
\right) (\gijksSDP)^{\transpose} \nonumber \\
+ \left( 
\VikSDPimag  (\VikSDPreal - \VjkSDPreal)^{\transpose} 
-  \VikSDPreal  (\VikSDPimag - \VjkSDPimag)^{\transpose}
\right) (\bijksSDP)^{\transpose} 
\label{eq_bim_complex_power_real}\\
\QijkSDP  = 
- \left(\VikSDPreal(\VikSDPreal)^{\transpose} + \VikSDPimag(\VikSDPimag)^{\transpose}\right) ( \bijkshSDP       )^{\transpose} 
+ \left(\VikSDPimag(\VikSDPreal)^{\transpose} - \VikSDPimag(\VikSDPimag)^{\transpose}\right) ( \gijkshSDP       )^{\transpose} \nonumber \\
- \left( 
\VikSDPreal  (\VikSDPreal - \VjkSDPreal)^{\transpose} 
+  \VikSDPimag  (\VikSDPimag - \VjkSDPimag)^{\transpose}
\right) (\bijksSDP)^{\transpose} \nonumber \\
+ \left( 
\VikSDPimag  (\VikSDPreal - \VjkSDPreal)^{\transpose} 
-  \VikSDPreal  (\VikSDPimag - \VjkSDPimag)^{\transpose}
\right) (\gijksSDP)^{\transpose} 
\label{eq_bim_complex_power_imag},
\end{IEEEeqnarray}  

In the coming subsections we derive the \emph{diagonalized} and scalarized real-valued formulation of the BIM in both the polar and rectangular voltage coordinate systems.  

\subsection{Unbalanced BIM Rectangular Scalar Form}



   Using   $\textcolor{\complexcolor}{U_{i,\indexPhases}} =  \VikpNreal + \imagnumber\VikpNimag $, the  expressions for active power of the diagonal elements of $\SijkSDP$ can be parameterized, $\indexPhases,\indexPhasesTwo \in \setPhases$, as,
       \begin{multline}{}
\Pijkpp =  
\sum_{\indexPhasesTwo \in \setPhases} ( \VikpNreal   \VikppNreal + \VikpNimag   \VikppNimag)  \left( \gijksph + \gijkshph \right)  \\
+  \sum_{\indexPhasesTwo \in \setPhases} ( \VikpNimag   \VikppNreal - \VikpNreal   \VikppNimag)  \left( \bijksph + \bijkshph \right)   \\
-     \sum_{\indexPhasesTwo \in \setPhases}  ( \VikpNreal   \VjkppNreal + \VikpNimag   \VjkppNimag)\gijksph    \\
-     \sum_{\indexPhasesTwo \in \setPhases}  ( \VikpNimag   \VjkppNreal - \VikpNreal   \VjkppNimag) \bijksph   ,  \label{eq_active_rectangular_pf}
   \end{multline}  
   and
          \begin{multline}{}
\Qijkpp =  
- \sum_{\indexPhasesTwo \in \setPhases} ( \VikpNreal   \VikppNreal + \VikpNimag   \VikppNimag)  \left( \bijksph + \bijkshph \right)  \\
+  \sum_{\indexPhasesTwo \in \setPhases} ( \VikpNimag   \VikppNreal - \VikpNreal   \VikppNimag)  \left( \gijksph + \gijkshph \right)   \\
+    \sum_{\indexPhasesTwo \in \setPhases}  ( \VikpNreal   \VjkppNreal + \VikpNimag   \VjkppNimag)\bijksph    \\
-   \sum_{\indexPhasesTwo \in \setPhases} ( \VikpNimag   \VjkppNreal - \VikpNreal   \VjkppNimag) \gijksph   .\label{eq_reactive_rectangular_pf}
   \end{multline}  
   The bus shunt expressions are,
    \begin{IEEEeqnarray}{C?s}
\Pkpp =  
\sum_{\indexPhasesTwo \in \setPhases} ( \VikpNreal   \VikppNreal + \VikpNimag   \VikppNimag)  \gkshph \nonumber \\
+  \sum_{\indexPhasesTwo \in \setPhases} ( \VikpNimag   \VikppNreal - \VikpNreal   \VikppNimag)   \bkshph      ,  \label{eq_active_rectangular_pf_shunt}\\
\Qkpp =  
- \sum_{\indexPhasesTwo \in \setPhases} ( \VikpNreal   \VikppNreal + \VikpNimag   \VikppNimag)   \bkshph  \nonumber\\
+  \sum_{\indexPhasesTwo \in \setPhases} ( \VikpNimag   \VikppNreal - \VikpNreal   \VikppNimag)   \gkshph     .\label{eq_reactive_rectangular_pf_shunt}
    \end{IEEEeqnarray}
   In rectangular coordinates, the nodal voltage bounds as defined in \eqref{eq_voltage_bounds_matrix} become
 \begin{IEEEeqnarray}{C?s}
 \begin{bmatrix} 
  ( \VikAmin)^2 \\
   (\VikBmin)^2\\
   (\VikCmin)^2\\
    \end{bmatrix}  \leq 
  \begin{bmatrix} 
(\VikANreal)^2 + (\VikANimag)^2    \\
(\VikBNreal)^2 + (\VikBNimag)^2     \\
(\VikCNreal)^2 + (\VikCNimag)^2   \\
 \end {bmatrix} 
 \leq
  \begin{bmatrix} 
  ( \VikAmax)^2 \\
 (  \VikBmax)^2 \\
(   \VikCmax)^2\\
 \end{bmatrix} , \label{eq_voltage_bounds_rectangular}
    \end{IEEEeqnarray}  
    which are nonconvex for strictly positive voltage lower bounds.

The voltage angle difference constraint between buses \eqref{eq_angle_diff_eq} needs to be reformulated. 
We derive the unbalanced tangent identity \eqref{eq_tan_rect_su},
\begin{IEEEeqnarray}{C?s}
\tan\left(    \begin{bmatrix} 
\VikANangle \\
\VikBNangle \\
\VikCNangle
\end {bmatrix}  -
 \begin{bmatrix} 
\VjkANangle \\
\VjkBNangle \\
\VjkCNangle
\end {bmatrix}   \right) 
= 
 \begin{bmatrix} 
\VikANimag \VjkANreal - \VikANreal \VjkANimag    \\
\VikBNimag \VjkBNreal - \VikBNreal \VjkBNimag   \\
\VikCNimag \VjkCNreal - \VikCNreal \VjkCNimag  
 \end {bmatrix}
 \oslash
 \begin{bmatrix} 
\VikANreal \VjkANreal + \VikANimag \VjkANimag    \\
\VikBNreal \VjkBNreal + \VikBNimag \VjkBNimag     \\
\VikCNreal \VjkCNreal + \VikCNimag  \VjkCNimag
 \end {bmatrix} \label{eq_tan_rect_su} .
 \end{IEEEeqnarray}
Because the angle difference is bounded by $[-\pi/2, \pi/2]$ for voltage stability, we can rewrite this as quadratic nonconvex constraint \eqref{eq_su_bim_vad},  where the tangent function is applied element-wise to the components of the vector containing the bounds on the angle difference, 
\begin{IEEEeqnarray}{C?s}
  \tan\circ(\thetaijSDPmin) \circ \begin{bmatrix} 
\VikANreal \VjkANreal + \VikANimag \VjkANimag    \\
\VikBNreal \VjkBNreal + \VikBNimag \VjkBNimag     \\
\VikCNreal \VjkCNreal + \VikCNimag  \VjkCNimag
 \end {bmatrix}    \leq
 \begin{bmatrix} 
\VikANimag \VjkANreal - \VikANreal \VjkANimag    \\
\VikBNimag \VjkBNreal - \VikBNreal \VjkBNimag   \\
\VikCNimag \VjkCNreal - \VikCNreal \VjkCNimag  
 \end {bmatrix} \nonumber \\
 \leq  
   \tan\circ(\thetaijSDPmax) \circ \begin{bmatrix} 
\VikANreal \VjkANreal + \VikANimag \VjkANimag    \\
\VikBNreal \VjkBNreal + \VikBNimag \VjkBNimag     \\
\VikCNreal \VjkCNreal + \VikCNimag  \VjkCNimag
 \end {bmatrix}. \label{eq_su_bim_vad}
 \end{IEEEeqnarray}
 We derive the equivalent expression for the angle differences between phases \emph{on the same bus}, 
The phase angle difference constraint \eqref{eq_angle_diff_eq_bus} is equivalent to,
\begin{IEEEeqnarray}{C?s}
\frac{2\pi}{3} + \thetaiSDPmin
 \leq    \begin{bmatrix} 
 \VikANangle -\VikBNangle   \\
 \VikBNangle    -\VikCNangle     \\
 \VikCNangle   -\VikANangle  
 \end {bmatrix}   \leq 
 \frac{2\pi}{3} +   \thetaiSDPmax. \label{eq_pad_shifted}
\end{IEEEeqnarray}
For the tangent function to be increasing and invertible, we want to further restrict ourselves to,
\begin{IEEEeqnarray}{C?s}
\frac{\pi}{2} \leq  \frac{2\pi}{3} + \thetaiSDPmin
 \leq 
 \frac{2\pi}{3} +   \thetaiSDPmax \leq \pi,
\end{IEEEeqnarray}
which implies, 
\begin{IEEEeqnarray}{C?s}
 - \frac{\pi}{6} \leq \thetaiSDPmin \leq
  \thetaiSDPmax \leq \frac{\pi}{3}. \label{eq_pad_bounds_limit}
\end{IEEEeqnarray}
With this restriction we derive the quadratic nonconvex constraint \eqref{eq_su_bim_vad_bus} from \eqref{eq_pad_shifted}, using an identity similar to \eqref{eq_tan_rect_su},
\begin{IEEEeqnarray}{C?s}
   \tan\circ\left( \frac{2\pi}{3} +  \thetaiSDPmin \right) \circ \begin{bmatrix} 
\VikANreal \VikBNreal + \VikANimag \VikBNimag    \\
\VikBNreal \VikCNreal + \VikBNimag \VikCNimag     \\
\VikCNreal \VikANreal + \VikCNimag  \VikANimag
 \end {bmatrix}    \leq
 \begin{bmatrix} 
\VikANimag \VikBNreal - \VikANreal \VikBNimag    \\
\VikBNimag \VikCNreal - \VikBNreal \VikCNimag   \\
\VikCNimag \VikANreal - \VikCNreal \VikANimag  
 \end {bmatrix} \nonumber \\
 \leq  
   \tan\circ\left( \frac{2\pi}{3} +  \thetaiSDPmax \right) \circ \begin{bmatrix} 
\VikANreal \VikBNreal + \VikANimag \VikBNimag    \\
\VikBNreal \VikCNreal + \VikBNimag \VikCNimag     \\
\VikCNreal \VikANreal + \VikCNimag  \VikANimag
 \end {bmatrix}. \label{eq_su_bim_vad_bus}
\end{IEEEeqnarray}

\subsection{Unbalanced BIM Polar Scalar Form} \label{sec_unbal_bim_polar}
The active and reactive power flow in each conductor can be written in a paramterized way, using nodal voltages per conductor $\textcolor{\complexcolor}{U_{i,\indexPhases}}  = \VikpNmag \imagangle \VikpNangle$.
The expressions of the diagonal elements of $\SijkSDP$ can be parameterized, $\indexPhases,\indexPhasesTwo \in \setPhases = \{\symPhaseA, \symPhaseB, \symPhaseC  \} $ and the real-value formulation of active and reactive power is obtained using trigonometric functions, 
       \begin{multline}{}
\Pijkpp  =  
\sum_{\indexPhasesTwo \in \setPhases} \VikpNmag  
 \VikppNmag \cos(\VikpNangle \!- \VikhNangle)  \left( \gijksph + \gijkshph \right)  \\
+  \sum_{\indexPhasesTwo \in \setPhases} \VikpNmag  \VikppNmag \sin(\VikpNangle \!- \VikhNangle) \left( \bijksph + \bijkshph \right)   \\
-     \sum_{\indexPhasesTwo \in\setPhases}  \VikpNmag \VjkppNmag\!   \cos(\VikpNangle \!- \VjkhNangle)\gijksph    \\
-    \sum_{\indexPhasesTwo \in \setPhases}  \VikpNmag \VjkppNmag\!   \sin(\VikpNangle \!- \VjkhNangle) \bijksph   , \label{eq_active_polar_pf} 
   \end{multline}  
   and
          \begin{multline}{}
\Qijkpp =  
-\sum_{\indexPhasesTwo \in \setPhases} \VikpNmag  
 \VikppNmag \cos(\VikpNangle \!- \VikhNangle)  \left( \bijksph + \bijkshph \right)  \\
+  \sum_{\indexPhasesTwo \in\setPhases} \VikpNmag  \VikppNmag \sin(\VikpNangle \!- \VikhNangle) \left( \gijksph + \gijkshph \right)  \\
+   \sum_{\indexPhasesTwo \in\setPhases}  \VikpNmag \VjkppNmag\!   \cos(\VikpNangle \!- \VjkhNangle)\bijksph    \\
-   \sum_{\indexPhasesTwo \in \setPhases}  \VikpNmag \VjkppNmag\!    \sin(\VikpNangle \!- \VjkhNangle) \gijksph   . \label{eq_reactive_polar_pf}
   \end{multline}  
   Thus, the active and reactive power flow through each conductor is obtained using the mutual coupling of the nodal voltages and branch impedances. 
The bus shunt expressions are,
 \begin{IEEEeqnarray}{C?s}
\Pkpp  =  
\sum_{\indexPhasesTwo \in \setPhases} \VikpNmag  
 \VikppNmag \cos(\VikpNangle \!- \VikhNangle)  \gkshph  \nonumber \\
+  \sum_{\indexPhasesTwo \in \setPhases} \VikpNmag  \VikppNmag \sin(\VikpNangle \!- \VikhNangle)  \bkshph, \label{eq_active_polar_pf_shunt} \\
\Qkpp =  
-\sum_{\indexPhasesTwo \in \setPhases} \VikpNmag  
 \VikppNmag \cos(\VikpNangle \!- \VikhNangle)   \bkshph  \nonumber \\
+  \sum_{\indexPhasesTwo \in\setPhases} \VikpNmag  \VikppNmag \sin(\VikpNangle \!- \VikhNangle) \gkshph  . \label{eq_reactive_polar_pf_shunt}
 \end{IEEEeqnarray}
The bus pair voltage angle difference constraint \eqref{eq_angle_diff_eq} and phase angle difference constraint \eqref{eq_angle_diff_eq_bus} directly apply in this variable space.

        \section{Lifting of Bus Injection Model} \label{sec_bim_relax}
        This section defines new variables to represent products of voltages and currents for the BIM, to enable a lift-and-project approach. This approach is commonly used for SOC and SDP relaxations of the nonlinear power-voltage formulation.
        \subsection{Lifted Variables}
        \subsubsection{Bus Voltage}
We define an auxiliary variable for the  voltage products, $\VisqkSDP$, satisfying    
 \begin{IEEEeqnarray}{C?s}
\VisqkSDP =  \VisqkSDPreal + \imagnumber \VisqkSDPimag =     \VikSDP(\VikSDP)^{\hermitiantranspose} , \IEEEyesnumber \IEEEyessubnumber  \label{eq_w_definition} \\
  \VisqkSDP \succeq 0, \rank(\VisqkSDP) = 1.  \IEEEyessubnumber \label{eq_w_definition_sdp}
 \end{IEEEeqnarray}
 Note that \eqref{eq_w_definition} is a quadratic nonconvex constraint, and \eqref{eq_w_definition_sdp} is the well-known rank-constrained SDP equivalent form. 
We illustrate the structure of $\VisqkSDP =\VikSDP(\VikSDP)^{\hermitiantranspose} $ as a real-valued matrix in rectangular coordinates,
 \begin{IEEEeqnarray}{Cs}
\VisqkSDP =  \begin{bmatrix} 
\underline{\VisqksocpAAreal}	& \underline{\VisqksocpABreal}   & \underline{\VisqksocpACreal}  \\
\VisqksocpABreal  	&\underline{\VisqksocpBBreal}   & \underline{\VisqksocpBCreal }  \\
\VisqksocpACreal 	& \VisqksocpBCreal 	  & \underline{ \VisqksocpCCreal}
\end{bmatrix}   \nonumber
 +  \imagnumber 
  \begin{bmatrix} 
				0	\!\!&\!\!  \underline{\VisqksocpABimag  }\!\! &\!\! \underline{\VisqksocpACimag } \\
- \VisqksocpABimag \!\!	&  	\!\!		0		\!\! &\!\! \underline{\VisqksocpBCimag } \\
- \VisqksocpACimag		\!\!&\!\! - \VisqksocpBCimag	 \!\! & \!\!0
\end{bmatrix}  . \\
\label{eq_lifted_voltage_hermitian}
 \end{IEEEeqnarray}
 Note that this representation requires 9 unique scalar variables (underlined) and that the diagonal is real-valued,
  \begin{IEEEeqnarray}{C?s}
\diag ( \VisqkSDP )  =  
  \VikSDP \circ \VikSDP^* 
=  
\begin{bmatrix} 
|\VikA|^2  \\
|\VikB|^2	   \\
|\VikC |^2	\\
\end{bmatrix} 
 =
\begin{bmatrix} 
{\VisqksocpAAreal}	  \\
{\VisqksocpBBreal}    \\
 { \VisqksocpCCreal} \\
\end{bmatrix}
.\label{eq_lifted_voltage_diagonal}
 \end{IEEEeqnarray}
One can define bounds on the matrix entries as,
 \begin{IEEEeqnarray}{C?s}
\VikSDPmin \circ \VikSDPmin \leq   \diag ( \VisqkSDPreal ) \leq \VikSDPmax \circ \VikSDPmax \IEEEyesnumber \IEEEyessubnumber, \label{eq_voltage_bounds_lifted}\\
- \VikSDPmax (\VikSDPmax)^{\transpose} \leq   \VisqkSDPreal, \VisqkSDPimag \leq \VikSDPmax (\VikSDPmax)^{\transpose} .\IEEEyessubnumber  \label{eq_voltage_bounds_lifted2}\
  \end{IEEEeqnarray}
%
%
In this variable space, the phase voltage angle difference constraint \eqref{eq_angle_diff_eq_bus} becomes,
 \begin{IEEEeqnarray}{Cs}
  \!\!\tan\left( \frac{2\pi}{3} + \thetaiSDPmin\right)  \circ
  \begin{bmatrix} 
 {\VisqksocpABreal}   \\
 {\VisqksocpACreal}  \\
{\VisqksocpBCreal }  
\end{bmatrix}   \leq 
  \begin{bmatrix} 
  {\VisqksocpABimag  }\\
  {\VisqksocpACimag } \\
  {\VisqksocpBCimag } 
\end{bmatrix} \nonumber \\
\leq 
  \tan\left( \frac{2\pi}{3} + \thetaiSDPmax\right)  \circ
  \begin{bmatrix} 
 {\VisqksocpABreal}   \\
 {\VisqksocpACreal}  \\
{\VisqksocpBCreal }  
\end{bmatrix}  \label{eq_angle_diff_eq_bus_lifted}.
 \end{IEEEeqnarray}
 The reference bus phasor is fixed,
     \begin{IEEEeqnarray}{C?s}
\VisqkSDP=  \VikrefSDP (\VikrefSDP)^{\hermitiantranspose}.\label{eq_ref_bus_voltage_lifted}
  \end{IEEEeqnarray}

        \subsubsection{Bus Voltage Cross Product}

We define a variable $\VijsqkSDP$ for the cross-product of the voltages $\VikSDP$ and $\VjkSDP$ of the buses associated with a branch (i.e. all bus-pairs) as used in equations (\ref{eq_active_polar_pf}-\ref{eq_reactive_polar_pf}) and (\ref{eq_active_rectangular_pf}-\ref{eq_reactive_rectangular_pf}), 
 \begin{IEEEeqnarray}{C?s}
\VijsqkSDP =   \VikSDP (\VjkSDP)^{\hermitiantranspose}  \label{eq_wij_definition}.
 \end{IEEEeqnarray}
 Although $\VijsqkSDP$ is rank-1 by construction, it is not Hermitian. 
Note that this definition implies,
 \begin{IEEEeqnarray}{C?s}
   \VijsqkSDP = (\VjisqkSDP)^{\hermitiantranspose}  .
 \end{IEEEeqnarray}
One can define bounds on the matrix entries of $\VijsqkSDP$ as,
 \begin{IEEEeqnarray}{C?s}
- \VikSDPmax (\VjkSDPmax)^{\transpose} \leq   \VijsqkSDPreal, \VijsqkSDPimag \leq \VikSDPmax (\VjkSDPmax)^{\transpose} . \label{eq_w_cross_bounds}
  \end{IEEEeqnarray}

        \subsection{Power Flow Model}

 Now power flow equation \eqref{eq_bim_complex_power} can be written using the lifted variables,  
 \begin{IEEEeqnarray}{C?s}
  \SijkSDP = \VisqkSDP ( \yijkshSDP       )^{\hermitiantranspose} +  (\VisqkSDP -  \VijsqkSDP) (\yijksSDP)^{\hermitiantranspose} .
 \end{IEEEeqnarray}
     The real-value equivalents are,
        \begin{IEEEeqnarray}{C?s}
  \PijkSDP = \VisqkSDPreal ( \gijkshSDP    + \gijksSDP   )^{\transpose} + \VisqkSDPimag ( \bijkshSDP  +\bijksSDP     )^{\transpose} \nonumber \\
    -  \VijsqkSDPreal (\gijksSDP)^{\transpose}   -  \VijsqkSDPimag (\bijksSDP)^{\transpose} ,\IEEEyesnumber   \IEEEyessubnumber \label{eq_active_power_bim_lifted}\\
  \QijkSDP = \VisqkSDPimag ( \gijkshSDP   +\gijksSDP    )^{\transpose} - \VisqkSDPreal ( \bijkshSDP +\bijksSDP      )^{\transpose} \nonumber \\
   -  \VijsqkSDPimag (\gijksSDP)^{\transpose} +  \VijsqkSDPreal (\bijksSDP)^{\transpose} .\IEEEyessubnumber\label{eq_reactive_power_bim_lifted}
   \end{IEEEeqnarray}  
   Similarly, the bus shunt power with lifted variables is,
    \begin{IEEEeqnarray}{C?s}
\SbSDP =  \VisqkSDP  (\ybSDP)^{\hermitiantranspose},
 \end{IEEEeqnarray}
     and its real-value equivalents are,
        \begin{IEEEeqnarray}{C?s}
  \PbSDP = \VisqkSDPreal ( \gbSDP     )^{\transpose} + \VisqkSDPimag ( \bbSDP       )^{\transpose},  \IEEEyesnumber   \IEEEyessubnumber \label{eq_shunt_lifted_active}\\
    \QbSDP = \VisqkSDPimag ( \gbSDP   )^{\transpose} - \VisqkSDPreal ( \bbSDP       )^{\transpose}   .  \IEEEyessubnumber \label{eq_shunt_lifted_reactive}
       \end{IEEEeqnarray}  
       
    Using \eqref{eq_lifted_voltage_diagonal}, we lift the branch current limit \eqref{eq_branch_currents_WS} to the $\VisqkSDP$ variable space as
    \begin{IEEEeqnarray}{C?s}
 \diag(\SijkSDP) \circ \diag(\SijkSDP)^* \leq \IijkSDPrated \circ \IijkSDPrated \circ \diag(\VisqkSDP), \label{eq_lifted_current_limit_WS}
  \end{IEEEeqnarray}
  which is a convex SOC constraint.

\subsection{Voltage Angle Difference Bound} \label{sec_tangent_voltage_angle}
 We note the `tangent inequality' \cite{Jabr2007} can be extended to the three-phase case in the following way:
\begin{IEEEeqnarray}{C?s}
\tan\left(    \begin{bmatrix} 
\VikANangle \\
\VikBNangle \\
\VikCNangle
\end {bmatrix}  -
 \begin{bmatrix} 
\VjkANangle \\
\VjkBNangle \\
\VjkCNangle
\end {bmatrix}   \right) = {\diag(\VijsqkSDPimag)}\,{\oslash}\,{ \diag(\VijsqkSDPreal)}
\end{IEEEeqnarray}
The voltage angle difference bounds therefore are, 
\begin{IEEEeqnarray}{C?s}
\tan\left(  \thetaijSDPmin   \right) \circ \diag(\VijsqkSDPreal)
\leq  \diag(\VijsqkSDPimag) \nonumber \\
\leq
\tan\left(  \thetaijSDPmax   \right) \circ \diag(\VijsqkSDPreal). \label{eq_bim_vad_bounds}
\end{IEEEeqnarray}

         \subsection{Rank-Constrained SDP Model} \label{sec_mesh_grids_bim_sdp}
\subsubsection{Meshed Grids}

Note that \eqref{eq_wij_definition} can be generalized (to support meshed grids), by shaping the matrices into a block matrix $\VsqkSDP$,
\begin{IEEEeqnarray}{C?s}
\VsqkSDP =   
      \begin{bmatrix} 
 \VikSDP  \\
 \VjkSDP \\
 \vdots \\
 \VzkSDP
\end{bmatrix}
      \begin{bmatrix} 
 \VikSDP  \\
 \VjkSDP \\
  \vdots \\
 \VzkSDP
\end{bmatrix}^{\hermitiantranspose}
= \begin{bmatrix} 
 \VisqkSDP & \VijsqkSDP & \cdots & \VizsqkSDP \\
 \VjisqkSDP & \VjsqkSDP & \cdots &\VjzsqkSDP \\
\vdots & \vdots & \ddots &  \vdots \\
\VzisqkSDP & \VzjsqkSDP & \cdots & \VzsqkSDP  \\
\end{bmatrix}   \IEEEyesnumber \IEEEyessubnumber  ,\\
\VsqkSDP \succeq 0, \rank\left( \VsqkSDP \right) =1. \label{bim_rank_constraint_meshed} \IEEEyessubnumber
 \end{IEEEeqnarray}
 Note that $\VsqkSDP  \in \setHermitianMatrix{|\setPhases||\setBuspairs|}$ with $|\setBuspairs|$  the number of  unique bus pairs in the topology.
If there is no branch between buses $i$ and $j$, the corresponding $\VijsqkSDP = \mathbf{0} $. 
The real-value equivalent SDP constraint \cite{Fazel2001}, for $\VsqkSDP = \VsqkSDPreal + \imagnumber \VsqkSDPimag$, is, 
 \begin{IEEEeqnarray}{C?s}
 \begin{bmatrix} 
 \phantom{-}\VsqkSDPreal & \VsqkSDPimag \\
- \VsqkSDPimag & \VsqkSDPreal 
\end{bmatrix}  \succeq 0, 
\rank\left( \begin{bmatrix} 
 \phantom{-}\VsqkSDPreal & \VsqkSDPimag \\
- \VsqkSDPimag & \VsqkSDPreal 
\end{bmatrix}  \right) =1. \label{eq_psd_system_real}
 \end{IEEEeqnarray}
Chordal relaxation   can be employed to replace $\VsqkSDP$  with a set of smaller matrices \cite{Liu2018c,Gan2014b}.

\subsubsection{Radial Grids}
In case of radial grids, a decomposition of $\VsqkSDP$ is easily derived. 
Only constraints of type \eqref{eq_buspairs_W} are retained \cite{Gan2014b}:
 \begin{IEEEeqnarray}{C?s}
\forall ij \in \setBuspairs:  \VbpsqkSDP =
      \begin{bmatrix} 
 \VikSDP  \\
 \VjkSDP 
\end{bmatrix}
      \begin{bmatrix} 
 \VikSDP  \\
 \VjkSDP 
\end{bmatrix}^{\hermitiantranspose}
=
    \begin{bmatrix} 
 \VisqkSDP & \VijsqkSDP \\
 \VjisqkSDP & \VjsqkSDP 
\end{bmatrix} \nonumber  ,\\
\VbpsqkSDP \succeq 0, \rank\left(  \VbpsqkSDP \right) =1.  \label{eq_buspairs_W}
 \end{IEEEeqnarray}
 Note that $\VbpsqkSDP  \in \setHermitianMatrix{2|\setPhases|}$.
Note that $\VbpsqkSDP \succeq 0$ implies both  $\VisqkSDP \succeq 0$ and $\VjsqkSDP \succeq 0$.
The real-value equivalent form  \cite{Fazel2001} is, 
 \begin{IEEEeqnarray}{C?s}
\VbpsqkSDPtworeal =
\begin{bmatrix} 
 \VisqkSDPreal & \VijsqkSDPreal &  \VisqkSDPimag & \VijsqkSDPimag \\
 (\VijsqkSDPreal)^{\transpose} & \VjsqkSDPreal &   (-\VijsqkSDPimag)^{\transpose}  & \VjsqkSDPimag  \\
-  \VisqkSDPimag & -\VijsqkSDPimag &\VisqkSDPreal   & \VijsqkSDPreal \\
 (\VijsqkSDPimag)^{\transpose} &- \VjsqkSDPimag &   (\VijsqkSDPreal)^{\transpose}  & \VjsqkSDPreal 
\end{bmatrix}  \succeq0 .  \IEEEyesnumber \IEEEyessubnumber \label{eq_sdp_bim_radial_real} \\
\rank(\VbpsqkSDPtworeal) = 1   \IEEEyessubnumber \label{eq_rank_bim_radial_real}
 \end{IEEEeqnarray}

  \section{Lifting  of Branch Flow Model} \label{sec_unbal_bfm}
This section defines new variables to represent products of voltages and currents specific to the BFM and then details the lift-and-project approach taken.


\subsection{Lifted Variables}
\subsubsection{Series Current}

 The auxiliary variable for the  current products, $\IijsqskSDP$ satisfies,  
 \begin{IEEEeqnarray}{C?s}
 \IijsqskSDP
 = \IijskSDP (\IijskSDP)^{\hermitiantranspose} = \IjiskSDP (\IjiskSDP)^{\hermitiantranspose}  \IEEEyesnumber \IEEEyessubnumber = \IijsqskSDPreal + \imagnumber \IijsqskSDPimag \\
  \IijsqskSDP \succeq 0, \rank(\IijsqskSDP) = 1 \IEEEyessubnumber
  \end{IEEEeqnarray} 
 Which has the following  representation in scalar variables:
  \begin{IEEEeqnarray}{Cs}
\!\IijsqskSDP  \!= \! \underbrace{\begin{bmatrix} 
\underline{\IijsqsksocpAAreal	}	\!\!& \underline{ \IijsqsksocpABreal   }\!\!	& \underline{\IijsqsksocpACreal } \\
\IijsqsksocpABreal  		\!\!& \underline{\IijsqsksocpBBreal   } \!\!& \underline{ \IijsqsksocpBCreal  } \\
\IijsqsksocpACreal 		\!\!& \IijsqsksocpBCreal 		 \!\!& \underline{ \IijsqsksocpCCreal}
\end{bmatrix}  }_{\IijsqskSDPreal } \nonumber 
\!+ \imagnumber
\underbrace{  \begin{bmatrix} 
		0		\!\!	& \!\!\underline{ \IijsqsksocpABimag } 	\!\!& \!\!\underline{ \IijsqsksocpACimag } \\
- \IijsqsksocpABimag \!\!	& \!\!  			0		\!\!&\!\! \underline{\IijsqsksocpBCimag } \\
- \IijsqsksocpACimag	\!\!	& \!\!- \IijsqsksocpBCimag	\!\!  & \!\!0
\end{bmatrix}  }_{\IijsqskSDPimag  }\!. \\ \label{eq_lifted_series_current_variable}
 \end{IEEEeqnarray}
 This representation requires 9 unique variables (underlined).
  \subsubsection{Total Current}
The total current $\IijkSDP$ is lifted as,
            \begin{IEEEeqnarray}{C?s}
    \IijsqkSDP = \IijsqkSDPreal + \imagnumber \IijsqkSDPimag=        \IijkSDP (\IijkSDP)^{\hermitiantranspose} \IEEEyesnumber \IEEEyessubnumber \label{eq_total_current_variable_lifted} \\
    \IijsqkSDP \succeq 0 ,  \rank(\IijsqkSDP) = 1  \IEEEyessubnumber
          \end{IEEEeqnarray}
%
%
We start from the definition
$  \IijkSDP =   \IijskSDP + (\yijkshSDP)  \VikSDP $ and  multiply both sides with their conjugate transpose,
            \begin{IEEEeqnarray}{C?s}
          \IijkSDP (\IijkSDP)^{\hermitiantranspose} = (\IijskSDP + (\yijkshSDP)  \VikSDP)((\IijskSDP)^{\hermitiantranspose} +  (\VikSDP)^{\hermitiantranspose}(\yijkshSDP)^{\hermitiantranspose} ). \nonumber 
                    \end{IEEEeqnarray}
Substituting in the lifted variables we obtain,
            \begin{IEEEeqnarray}{C?s}
\IijsqkSDP          =\IijsqskSDP + (\yijkshSDP)  \VisqkSDP (\yijkshSDP)^{\hermitiantranspose}+ (\yijkshSDP)  \SijksSDP + (\SijksSDP)^{\hermitiantranspose} (\yijkshSDP)^{\hermitiantranspose} ,
\end{IEEEeqnarray}
 proving the lifted total current variable is a linear combination of $\IijsqskSDP$, $\VisqkSDP$ and $\SijksSDP$. In the real domain this becomes,
 \begin{IEEEeqnarray}{C?s}
\IijsqkSDPreal          = \IijsqskSDPreal +  \gijkshSDP \VisqkSDPreal (\gijkshSDP)^{\transpose}  + \gijkshSDP  \VisqkSDPimag (\bijkshSDP)^\transpose \nonumber \\
+  \bijkshSDP \VisqkSDPreal (\bijkshSDP)^{\transpose} - \bijkshSDP  \VisqkSDPimag (\gijkshSDP)^{\transpose}  \nonumber\\
 +  \gijkshSDP\PijksSDP - \bijkshSDP \QijksSDP +(\gijkshSDP\PijksSDP)^{\transpose}  - (\bijkshSDP\QijksSDP)^{\transpose}  , \IEEEyesnumber \IEEEyessubnumber \label{eq_tot_current_real}\\
\IijsqkSDPimag        = \IijsqskSDPimag +  \bijkshSDP \VisqkSDPreal (\gijkshSDP)^{\transpose}  + \bijkshSDP  \VisqkSDPimag (\bijkshSDP)^\transpose \nonumber\\ 
\gijkshSDP \VisqkSDPimag (\gijkshSDP)^{\transpose} - \gijkshSDP   \VisqkSDPreal (\bijkshSDP)^{\transpose} \nonumber\\ 
+ \gijkshSDP\QijksSDP + \bijkshSDP\PijksSDP   -  (\bijkshSDP\PijksSDP)^{\transpose}  -  (\gijkshSDP\QijksSDP)^{\transpose} . \IEEEyessubnumber \label{eq_tot_current_imag}
\end{IEEEeqnarray}
One can define current bounds \eqref{eq_current_ratings_matrix} on the matrix entries of $\IijsqkSDP$ as,
            \begin{IEEEeqnarray}{C?s}
0\le \diag(\IijsqkSDPreal) \leq \IijkSDPrated \circ \IijkSDPrated  \IEEEyesnumber \IEEEyessubnumber   \label{eq_tot_current_bound}, \\
- \IijkSDPrated  (\IijkSDPrated)^{\transpose} \leq \IijsqkSDPreal, \IijsqkSDPimag \leq \IijkSDPrated  (\IijkSDPrated)^{\transpose} \IEEEyessubnumber \label{eq_tot_current_bound2}.
\end{IEEEeqnarray}
Note that \eqref{eq_tot_current_bound}-\eqref{eq_tot_current_bound2} is equivalent to \eqref{eq_lifted_current_limit_WS}, however the former is linear due to \eqref{eq_tot_current_real} -\eqref{eq_tot_current_imag} and the latter is quadratic-convex (SOC).
There are no direct bounds on the auxiliary variable for the series current $\IijsqskSDP$, however we can derive valid bounds through the known bounds total current and voltage,
                \begin{IEEEeqnarray}{C?s}
 -( \IijkSDPrated (\IijkSDPrated)^{\transpose} + |\yijkshSDP|  \VikSDPmax (\VikSDPmax)^{\transpose} |\yijkshSDP|^{\transpose} )
 \leq 
 \IijsqskSDPreal, \IijsqskSDPimag       \nonumber \\
 \leq 
  \IijkSDPrated (\IijkSDPrated)^{\transpose} + |\yijkshSDP|  \VikSDPmax (\VikSDPmax)^{\transpose} |\yijkshSDP|^{\transpose}  \label{eq_implied_bounds_lifted_series_current},
\end{IEEEeqnarray}
where $|\yijkshSDP|$ indicates the element-wise application of the absolute value operation to obtain the magnitude. 
        \subsection{Power Flow Model}
Given $\SlkslossSDP  =  \zijksSDP  \IijsqskSDP$, 
 $\SijkshlossSDP = \VisqkSDP (\yijkshSDP)^{\hermitiantranspose}$, and 
$\SjikshlossSDP =  \VjsqkSDP (\yjikshSDP)^{\hermitiantranspose}$, the power balance \eqref{eq_bfm_line_power_balance} of a branch becomes,
\begin{IEEEeqnarray}{C}
 \SijkSDP + \SjikSDP  
      = \VisqkSDP (\yijkshSDP)^{\hermitiantranspose}  + \zijksSDP  \IijsqskSDP + \VjsqkSDP (\yjikshSDP)^{\hermitiantranspose}.
\end{IEEEeqnarray}
 The equivalent active and reactive power expressions are,
   \begin{IEEEeqnarray}{C?s}
\PijkSDP + \PjikSDP =  \VisqkSDPreal (\gijkshSDP)^{\transpose}  +  \VisqkSDPimag (\bijkshSDP)^\transpose + \rijksSDP \IijsqskSDPreal - \xijksSDP \IijsqskSDPimag  \nonumber \\
+\VjsqkSDPreal (\gjikshSDP)^{\transpose}  +  \VjsqkSDPimag (\bjikshSDP)^{\transpose} , \IEEEyesnumber \IEEEyessubnumber \label{eq_bfm_active_flow} \\
\QijkSDP + \QjikSDP =  \VisqkSDPimag (\gijkshSDP)^{\transpose} -    \VisqkSDPreal (\bijkshSDP)^{\transpose} +  \xijksSDP \IijsqskSDPreal + \rijksSDP \IijsqskSDPimag  \nonumber \\
+ \VjsqkSDPimag (\gjikshSDP)^{\transpose} -    \VjsqkSDPreal (\bjikshSDP)^{\transpose} .  \IEEEyessubnumber \label{eq_bfm_reactive_flow}
 \end{IEEEeqnarray}
Either $\SijkSDP$ or $\SijksSDP$ can be substituted out through,
    \begin{IEEEeqnarray}{C?s}
\PijkSDP = \PijksSDP +    \VisqkSDPreal (\gijkshSDP)^{\transpose}, \IEEEyesnumber  \label{eq_subst_series_active_power} 
\QijkSDP = \QijksSDP -    \VisqkSDPreal (\bijkshSDP)^{\transpose} .    \label{eq_subst_series_reactive_power}
 \end{IEEEeqnarray}

 \subsection{Ohm's Law}
       Ohm's law \eqref{eq_ohms_matrix} is reformulated  by multiplying both sides with their Hermitian adjoint,
\begin{IEEEeqnarray}{L}
  \! \VjkSDP (\VjkSDP)^{\hermitiantranspose} = (\VikSDP - \zijksSDP \IijskSDP) (\VikSDP - \zijksSDP \IijskSDP)^{\hermitiantranspose}  ,  \nonumber\\
 \! \!= \!\!  \VikSDP( \VikSDP)^{\hermitiantranspose} \!-\!  \SijksSDP (\zijksSDP)^{\hermitiantranspose} -  \zijksSDP  (\SijksSDP)^{\hermitiantranspose}  
\!+\! \zijksSDP  \IijskSDP (\IijskSDP)^{\hermitiantranspose}  \!(\zijksSDP)^{\hermitiantranspose}   \!      . \label{eq_ohms_squared}
 \end{IEEEeqnarray} Substituting in the lifted variables into \eqref{eq_ohms_squared}, we obtain,
\begin{IEEEeqnarray}{C?s}
 \VjsqkSDP = \VisqkSDP  -  \SijksSDP (\zijksSDP)^{\hermitiantranspose} -  \zijksSDP   (\SijksSDP)^{\hermitiantranspose}    + \zijksSDP \IijsqskSDP  (\zijksSDP)^{\hermitiantranspose}      .
 \end{IEEEeqnarray}
  The equivalent real expressions are,
 \begin{IEEEeqnarray}{RL}
 \VjsqkSDPreal   =  \VisqkSDPreal   - \PijksSDP(\rijksSDP)^{\transpose} - \QijksSDP(\xijksSDP)^{\transpose}   
  -\rijksSDP(\PijksSDP)^{\transpose} -\xijksSDP(\QijksSDP)^{\transpose}  \nonumber\\
  + \rijksSDP\IijsqskSDPreal(\rijksSDP)^{\transpose}-\xijksSDP\IijsqskSDPimag(\rijksSDP)^{\transpose}
  +  \xijksSDP\IijsqskSDPreal (\xijksSDP)^{\transpose} + \rijksSDP\IijsqskSDPimag(\xijksSDP)^{\transpose}, \nonumber \\     \IEEEyesnumber \IEEEyessubnumber \label{eq_ohms_real}\\
  \VjsqkSDPimag  =  \VisqkSDPimag   -\QijksSDP(\rijksSDP)^{\transpose} + \PijksSDP(\xijksSDP)^{\transpose}   
 - \xijksSDP(\PijksSDP)^{\transpose}  + \rijksSDP(\QijksSDP)^{\transpose}  \nonumber\\
  + \xijksSDP\IijsqskSDPreal(\rijksSDP)^{\transpose} + \rijksSDP\IijsqskSDPimag(\rijksSDP)^{\transpose} 
  - \rijksSDP\IijsqskSDPreal(\xijksSDP)^{\transpose} + \xijksSDP\IijsqskSDPimag(\xijksSDP)^{\transpose}     .\nonumber \\ \IEEEyessubnumber \label{eq_ohms_imag}
 \end{IEEEeqnarray}
These equations are symmetrical, so it is sufficient to generate the scalar constraints only on the upper triangle for $\VjsqkSDPimag$; upper triangle and diagonal for $\VjsqkSDPreal$ (generating all scalar constraints is redundant). 

\subsection{Voltage Angle Difference Bound}
We can derive $\VijsqkSDP$ as a function of $\VisqkSDP$ and $\SijksSDP$, 
 \begin{IEEEeqnarray}{C?s}
 \VijsqkSDP = \VisqkSDP -   \SijksSDP (\zijksSDP)^{\hermitiantranspose} , 
  \end{IEEEeqnarray}
which in the reals becomes,
  \begin{IEEEeqnarray}{C?s}
  \VijsqkSDPreal = \VisqkSDPreal -   \PijksSDP (\rijksSDP)^{\transpose} -   \QijksSDP (\xijksSDP)^{\transpose} , \IEEEyesnumber \IEEEyessubnumber \label{eq_wijreal_subst}\\
   \VijsqkSDPimag = \VisqkSDPimag -   \QijksSDP (\rijksSDP)^{\transpose} +   \PijksSDP (\xijksSDP)^{\transpose} .  \IEEEyessubnumber\label{eq_wijimag_subst}
 \end{IEEEeqnarray}
Therefore, to avoid introducing the variable $\VijsqkSDP$ in the BFM, we substitute this into \eqref{eq_bim_vad_bounds}, 
\begin{IEEEeqnarray}{C?s}
\tan\left(  \thetaijSDPmin   \right) \circ \diag(\VisqkSDPreal -   \PijksSDP (\rijksSDP)^{\transpose} -   \QijksSDP (\xijksSDP)^{\transpose}) \nonumber \\
\leq  \diag(\VisqkSDPimag -   \QijksSDP (\rijksSDP)^{\transpose} +   \PijksSDP (\xijksSDP)^{\transpose}) \nonumber \\
\leq
\tan\left(  \thetaijSDPmax   \right) \circ \diag(\VisqkSDPreal -   \PijksSDP (\rijksSDP)^{\transpose} -   \QijksSDP (\xijksSDP)^{\transpose}). \label{eq_bfm_vad_lifted}
\end{IEEEeqnarray}

         \subsection{Rank-Constrained SDP Model}
The product \eqref{eq_complex_power_flow_unbalanced}  formulated in the lifted variable space is,
 \begin{IEEEeqnarray}{C?s}
\forall lij \in \setTopology:  \VlinesqkSDP= \begin{bmatrix} 
 \VisqkSDP & \SijksSDP \\
(\SijksSDP)^{\hermitiantranspose} & \IijsqskSDP 
\end{bmatrix} ,  \IEEEyesnumber \IEEEyessubnumber \\
\VlinesqkSDP \succeq 0, \rank\left(  \VlinesqkSDP \right) =1.  \label{eq_bfm_rank_constraint} \IEEEyessubnumber
 \end{IEEEeqnarray}
Note that $\VlinesqkSDP  \in \setHermitianMatrix{2|\setPhases|}$.
Furthermore $\VlinesqkSDP\succeq 0$ necessitates $\VisqkSDP \succeq 0$ and $\IusqkSDP \succeq0 $ but not $\VjsqkSDP \succeq 0$. So for leaf buses, an explicit $\VjsqkSDP \succeq 0$ is needed. 
The real-value equivalent form is,
 \begin{IEEEeqnarray}{C?s}
\VlinesqkSDPrealtwo=  \begin{bmatrix} 
 \VisqkSDPreal & \PijksSDP &  \VisqkSDPimag & \QijksSDP \\
 (\PijksSDP)^{\transpose} &  \IijsqskSDPreal &  -(\QijksSDP)^{\transpose} &  \IijsqskSDPimag  \\
- \VisqkSDPimag & -\QijksSDP  &  \VisqkSDPreal & \PijksSDP \\
 (\QijksSDP)^{\transpose} & - \IijsqskSDPimag &  (\PijksSDP)^{\transpose} &  \IijsqskSDPreal 
\end{bmatrix}\succeq 0  , \IEEEyesnumber \IEEEyessubnumber \label{eq_complex_power_product_sdp} \\
\rank(\VlinesqkSDPrealtwo) = 1 . \IEEEyessubnumber  \label{eq_complex_power_product_rank} 
  \end{IEEEeqnarray}

  \section{Implementation of Real-Value Feasible Sets} \label{sec_feasible_sets}
  Table \ref{tab_var_constraints} section summarizes the feasible sets in the real domain. 
  The  shunt    $\PbSDP, \QbSDP$  feasible sets can be substituted into KCL and is therefore not listed explicitly. 
\newgeometry{margin=2cm} 
\begin{landscape}


\begin{table*}
\scriptsize
  \centering 
   \caption{Formulation feasible sets: variable spaces, bounds, constraints and substitution options (NCNL = nonconvex nonlinear, NCQ = nonconvex quadratic, QCP = quadratically constrained programming) }\label{tab_var_constraints}
    \begin{tabular}{l  l  l  l  l  l l l    }
\hline
& & \multicolumn{2}{c}{$\SijkSDP-\VikSDP$} 	 	&   	& \multicolumn{2}{c}{$\SijkSDP-\VisqkSDP$}  	  \\ 
 \cline {3-4}  \cline {6-7}   
& &AC polar & AC rectangular & &BIM lifted relax.&  BFM lifted relax.  \\
Overall complexity & & NLP & QCP & & SDP &  SDP \\\hline
Voltage variables & $\VikSDP$& $[ \VikpNmag \imagangle \VikpNangle]_{ \indexPhases \in \setPhases}$   \eqref{eq_var_voltage_polar}	& $[ \VikpNreal + \imagnumber\VikpNimag]_{ \indexPhases \in \setPhases}$    \eqref{eq_var_voltage_rectangular}    & 	 &  $\VisqkSDPreal, \VisqkSDPimag$ \eqref{eq_lifted_voltage_hermitian}   & $\VisqkSDPreal, \VisqkSDPimag$ \eqref{eq_lifted_voltage_hermitian}  , \\
Unit power variables &$\SuunitSDP$ &   $\PuunitSDP, \QuunitSDP $ \eqref{eq_unit_active_power}\eqref{eq_unit_reactive_power}	&   $\PuunitSDP, \QuunitSDP $ \eqref{eq_unit_active_power}\eqref{eq_unit_reactive_power} & 	 &    $\PuunitSDP, \QuunitSDP $ \eqref{eq_unit_active_power}\eqref{eq_unit_reactive_power} &     $\PuunitSDP, \QuunitSDP $ \eqref{eq_unit_active_power}\eqref{eq_unit_reactive_power} \\
Power flow variables & $\SijkSDP$ & $\PijkSDP, \QijkSDP$  \eqref{eq_complex_s_bounds}	&   $\PijkSDP, \QijkSDP$\eqref{eq_complex_s_bounds} & 	 & $\PijkSDP, \QijkSDP$   \eqref{eq_complex_s_bounds} &  $\PijkSDP, \QijkSDP$  \eqref{eq_complex_s_bounds} \\
Current variables & $\IijkSDP$ & eliminated 	&   eliminated  & 	 &  eliminated    & $\IijsqkSDPreal, \IijsqkSDPimag  $ \\
Additional variables & 	& - &  -&  & $\VijsqkSDPreal, \VijsqkSDPimag$  \eqref{eq_w_cross_bounds}   &$\IijsqskSDPreal, \IijsqskSDPimag$  \eqref{eq_lifted_series_current_variable} \eqref{eq_implied_bounds_lifted_series_current},  \\ \hline
Voltage bounds & $\VikSDPmin, \VikSDPmax$ 	&  LP \eqref{eq_voltage_bounds_polar} & NCQ\eqref{eq_voltage_bounds_rectangular}  & 	 &  LP\eqref{eq_voltage_bounds_lifted}, \eqref{eq_voltage_bounds_lifted2} & LP \eqref{eq_voltage_bounds_lifted}, \eqref{eq_voltage_bounds_lifted2},  \\
Phase angle difference bounds &  $\thetaiSDPmin, \thetaiSDPmax$	& LP \eqref{eq_angle_diff_eq_bus} &  NCQ \eqref{eq_su_bim_vad_bus}  & 	 & LP \eqref{eq_angle_diff_eq_bus_lifted}  & LP \eqref{eq_angle_diff_eq_bus_lifted} \\ 
Power flow bounds & $\SijkSDPrated$& SOC \eqref{eq_power_bounds} &  SOC \eqref{eq_power_bounds} & 	 &  SOC \eqref{eq_power_bounds}  &  SOC \eqref{eq_power_bounds} \\
Current bounds &  $\IijkSDPrated$	& SOC \eqref{eq_current_ratings_matrix} &  SOC\eqref{eq_current_ratings_matrix}  & 	 & SOC \eqref{eq_lifted_current_limit_WS}  & LP \eqref{eq_tot_current_bound}, \eqref{eq_tot_current_bound2} \\ 
Voltage angle difference bounds &  $\thetaijSDPmin, \thetaijSDPmax$	&  LP \eqref{eq_angle_diff_eq} &  NCQ \eqref{eq_su_bim_vad}  & 	 &  LP \eqref{eq_bim_vad_bounds}  & LP \eqref{eq_bfm_vad_lifted} \\ 
\hline
Kirchhoff's current law & 	& LP \eqref{eq_kcl_active_node}, \eqref{eq_kcl_reactive_node}  & LP \eqref{eq_kcl_active_node}, \eqref{eq_kcl_reactive_node}  & 	 & LP \eqref{eq_kcl_active_node}, \eqref{eq_kcl_reactive_node}   & LP \eqref{eq_kcl_active_node}, \eqref{eq_kcl_reactive_node} \\
Bus shunt equations (admittance) & 	& NCNL \eqref{eq_active_polar_pf_shunt}, \eqref{eq_reactive_polar_pf_shunt} & NCQ \eqref{eq_active_rectangular_pf_shunt}, \eqref{eq_reactive_rectangular_pf_shunt}  & 	 & LP \eqref{eq_shunt_lifted_active}, \eqref{eq_shunt_lifted_reactive} & LP \eqref{eq_shunt_lifted_active}, \eqref{eq_shunt_lifted_reactive} \\
Power flow equations  ($\Pi$-model)& 	& NCNL \eqref{eq_active_polar_pf}, \eqref{eq_reactive_polar_pf} & NCQ \eqref{eq_active_rectangular_pf}, \eqref{eq_reactive_rectangular_pf}  & 	 & LP \eqref{eq_active_power_bim_lifted}, \eqref{eq_reactive_power_bim_lifted} & LP \eqref{eq_bfm_active_flow}, \eqref{eq_bfm_reactive_flow} \\
Kirchhoff's voltage law & 	& implicit  & implicit  & 	 & SDP \eqref{eq_psd_system_real} or \eqref{eq_sdp_bim_radial_real}, \eqref{eq_rank_bim_radial_real}  & implied by Ohm's law when radial \\
Ohm's law & 	& implied   & implicit  & 	 &  implicit   & LP \eqref{eq_ohms_real}, \eqref{eq_ohms_imag} \\ 
Complex power definition &  & implicit   & implicit  & 	 &  implicit   & SDP \eqref{eq_complex_power_product_sdp} \\ 
Reference bus fixed phasor &  	& LP \eqref{eq_ref_bus_voltage}   & LP \eqref{eq_ref_bus_voltage}    & 	 & LP \eqref{eq_ref_bus_voltage_lifted}   & LP \eqref{eq_ref_bus_voltage_lifted} \\ 
\hline
Substitution & 	& - & -  & 	 & -  &  $\IijsqkSDPreal, \IijsqkSDPimag \rightarrow \IijsqskSDPreal, \IijsqskSDPimag$ \eqref{eq_tot_current_real}, \eqref{eq_tot_current_imag}\\
 & 	&  & & 	 & &$\PijkSDP, \QijkSDP \leftrightarrow \PijksSDP,   \QijksSDP$ \eqref{eq_subst_series_active_power}, \eqref{eq_subst_series_reactive_power}  \\
 & 	&  &  & 	 &   &$ \VijsqkSDPreal, \VijsqkSDPimag \rightarrow  \VisqkSDPreal, \VisqkSDPimag, \PijksSDP, \QijksSDP$ \eqref{eq_wijreal_subst}, \eqref{eq_wijimag_subst}  \\
\hline
\end{tabular}
\end{table*}
\end{landscape}
\restoregeometry

 The AC polar and rectangular forms can be implemented in optimization toolboxes with support for nonlinear programming. The core power flow equations \eqref{eq_active_polar_pf}, \eqref{eq_reactive_polar_pf} or \eqref{eq_active_rectangular_pf}, \eqref{eq_reactive_rectangular_pf} are nonlinear, however only the rectangular form is quadratically representable.
 
Solving the lifted BIM and BFM forms requires dropping the rank constraints \eqref{eq_rank_bim_radial_real} or \eqref{eq_complex_power_product_rank} to obtain a SDP formulation.   
It is noted that further SOC relaxation can be performed \cite{Kim2003,Vanin2020}. 
Furthermore,  for the tightness of the SOC relaxation of the SDP constraints, it is best to work in the complex domain as long as possible. 
The process for obtaining the original (nonlifted) current and voltage variables is detailed in \cite{Geth2019}. 
 
 Implementations of the real-value formulations are available in the package \textsc{PowerModelsDistribution.jl}~\cite{Fobes2020}. 
 This package extends \textsc{PowerModels.jl}~\cite{Dunning2015} and is built on top of JuMP ~\cite{Dunning2015}, a Julia package for mathematical programming. 
 In the implementation, where possible, the matrix forms are used, and the scalarization of the equations is performed by JuMP.  
 The discussed substitution and elimination of $\IijsqkSDPreal, \IijsqkSDPimag$, $\PijksSDP, \QijksSDP$ and $\VijsqkSDPreal, \VijsqkSDPimag$ in the BFM  are performed through JuMP expressions, which improves the readability of the mathematical model in the code. We note that nondimensionalization (per unit scaling) is performed in the implementation, in an effort to improve numerical conditioning.

 Using \textsc{PowerModelsDistribution}~\cite{Fobes2020}, one can read \textsc{OpenDSS} power flow case files and run power flow and optimal power flow case studies. 
 \textsc{OpenDSS} is used to validate the accuracy of the power flow formulations. 
The worst relative voltage error for IEEE 13, 34, 123 and LV test case  \cite{Mather2017} is  1.4E-7~\cite{Fobes2020}. 

  \section{Conclusions} \label{sec_conclusions}
  In this report, the electrical physics of three-phase grids are derived in complex-power matrix variables and different voltage variable spaces, i.e. polar, rectangular and lifted. 
  Using the different variable spaces, the unbalanced power flow formulations are derived for  generic $\Pi$-model branches with asymmetric shunt impedance (full matrices, no structure assumed). 
  The real-value matrix and scalar equivalent formulations are derived, and implemented in an open-source software package. 
  The derived mathematical framework can easily be extended and applied to a variety of continuous and discrete optimization problems in distribution networks.  
  Expressions are presented to enforce voltage magnitude, angle, power and current bounds across the formulations in an exact manner, without needing to resort to auxiliary variables. 
  Nevertheless, including auxiliary variables instead of eliminating them does not by default increase computational effort. 
  Comparing these variants numerically requires the development of an extensive unbalanced OPF test case library, which includes values for the different bounds discussed, not just voltage magnitudes. 
  Future work includes adding more bounds, e.g. for sequence voltage magnitudes and unbalance metrics \cite{Girigoudar}.

\section*{Acknowledgement}
Part of this work has been performed in the framework of the ADriaN project supported by Fluvius. Special thanks to  Reinhilde D'hulst, Carleton Coffrin, David Fobes, Sander Claeys, Thomas Brinsmead and Rahmat Heidari.

\bibliographystyle{IEEEtran}

\end{document}